\documentclass[graybox,envcountchap,sectrefs]{svmono}

\usepackage{setspace}
\usepackage{mathtools}
\usepackage{amsmath}
\usepackage{amsfonts}
\usepackage{amssymb}
\usepackage{helvet}
\usepackage{courier}
\usepackage{mdwlist}
\usepackage{longtable}
\usepackage{multirow}
\usepackage{booktabs}
\usepackage{tabularx}
\usepackage{algorithm}
\usepackage{algorithmic}
\usepackage{hyperref}
\usepackage{type1cm}

\usepackage{makeidx}         
\usepackage{graphicx}        
\usepackage{makecell}
\usepackage{multicol}        
\usepackage[bottom]{footmisc}

\makeindex             
\newtheorem{assumption}{Assumption}[chapter]
\numberwithin{algorithm}{chapter}

\begin{document}

\author{Wei Wei, Tsinghua University}
\title{Tutorials on Advanced Optimization Methods}
\subtitle{}
\maketitle

\frontmatter

\mainmatter

\appendix

\chapter*{}

\begin{center} 
{\Large \bf Tutorials on Advanced Optimization Methods}
\end{center}

\

\begin{center} 
{\bf Wei Wei, Tsinghua University}
\end{center}

\

This material is the appendix part of my book collaborated with Professor Jianhui Wang at Southern Methodist University:

{\color{blue} Wei Wei, Jianhui Wang. Modeling and Optimization of Interdependent Energy Infrastructures. Springer Nature Switzerland, 2020.}

{\small \url{https://link.springer.com/book/10.1007%2F978-3-030-25958-7}}

This material provides thorough tutorials on some optimization techniques frequently used in various engineering disciplines, including
\begin{enumerate}
\item[$\clubsuit$] Convex optimization
\item[$\clubsuit$] Linearization technique and mixed-integer linear programming
\item[$\clubsuit$] Robust optimization
\item[$\clubsuit$] Equilibrium/game problems
\end{enumerate}

It discusses how to reformulate a difficult (non-convex, multi-agent, min-max) problem to a solver-compatible form (semidefinite program, mixed-integer linear program) via convexification, linearization, and decomposition, so the original problem can be reliably solved by commercial/open-source software. Fundamental algorithms (simplex algorithm, interior-point algorithm) are not the main focus.
  
This material is a good reference for self-learners who have basic knowledge in linear algebra and linear programming. It is one of the main references for an optimization course taught at Tsinghua University. If you need teaching slides, please contact {\color{blue}wei-wei04@mails.tsinghua.edu.cn} or find the up-to-date contact information at {\small \url{https://sites.google.com/view/weiweipes/}}

\pdfbookmark[2]{Bookmarktitle}{internal_label}
\tableofcontents

%
%
%

\motto{The great watershed in optimization isn't between linearity and nonlinearity, but convexity and non-convexity. \\  \rightline{ $-$Ralph Tyrrell Rockafellar}}
\chapter{Basics of Linear and Conic Programs}
\label{App-A} 

The mathematical programming theory has been thoroughly developed in width and depth since its birth in 1940s, when George Dantzig invented simplex algorithm for linear programming. The most influential findings in the field of optimization theory can be summarized as \cite{App-A-CVX-Book-Ben}: 

1) Recognition of the fact that under mild conditions, a convex optimization program is computationally tractable: the computational effort under a given accuracy grows moderately with the problem size even in the worst case. In contrast, a non-convex program is generally computationally intractable: the computational effort of the best known methods grows prohibitively fast with respect to the problem size, and it is reasonable to believe that this is an intrinsic feature of such problems rather than a limitation of existing optimization techniques. 

2) The discovery of interior-point methods, which was originally developed in 1980s to solve LPs and could be generalized to solve convex optimization problems as well. Moreover, between these two extremes (LPs and general convex programs), there are many important and useful convex programs. Although nonlinear, they still possess nice structured properties, which can be utilized to develop more dedicated algorithms. These polynomial-time interior-point algorithms turn out to be considerably more efficient than those exploiting only the convex property.

The superiority of formulating a problem as a convex optimization problem is apparent. The most appealing advantage is that the problem can be solved reliably and efficiently. It is also convenient to build the associated dual problem, which gives insights on sensitivity information and may help develop distributed algorithm for solving the problem. Convex optimization has been applied in a number of energy system operational issues, and well acknowledged for its computational superiority. We believe that it is imperative for researchers and engineers to develop certain understanding on this important topic. 

As we have already learnt in previous chapters, many optimization problems in energy system engineering can be formulated as or converted to convex programs. The goal of this chapter is to help readers develop necessary background knowledge and skills to apply several well-structured convex optimization models, including LPs, SOCPs, and SDPs, i.e., to formulate or transform their problems as these specific convex programs. Certainly, convex transformation (or convexification) may be rather tricky and require special knowledge and skills. Nevertheless, the attempt often turns out to be worthwhile. We also pay special attention to nonconvex QCQPs, which can model various decision-making problems in engineering, such as optimal power flow and optimal gas flow. We discuss convex relaxation technique based on SDP,which is shown to be very useful to get a high-quality objective lower bound. We also present MILP formulations for some special QPs; because of the special problem structure, these MILP models can tackle practically sized problems in reasonable time.

Most materials regarding convex sets and functions come from \cite{App-A-CVX-Book-Boyd} and its solution manual \cite{App-A-CVX-Book-Solution}; extensions of duality theory from linear programming to conic programming follows from \cite{App-A-CVX-Book-Ben}. We consolidate  necessary contents in a convenient way to make this book self-contained and easy to follow.

\section{Basic Notations}
\label{App-A-Sect01}

\subsection{Convex Sets}
\label{App-A-Sect01-01}

\begin{figure}[!htp]
\centering
\includegraphics[scale=0.60]{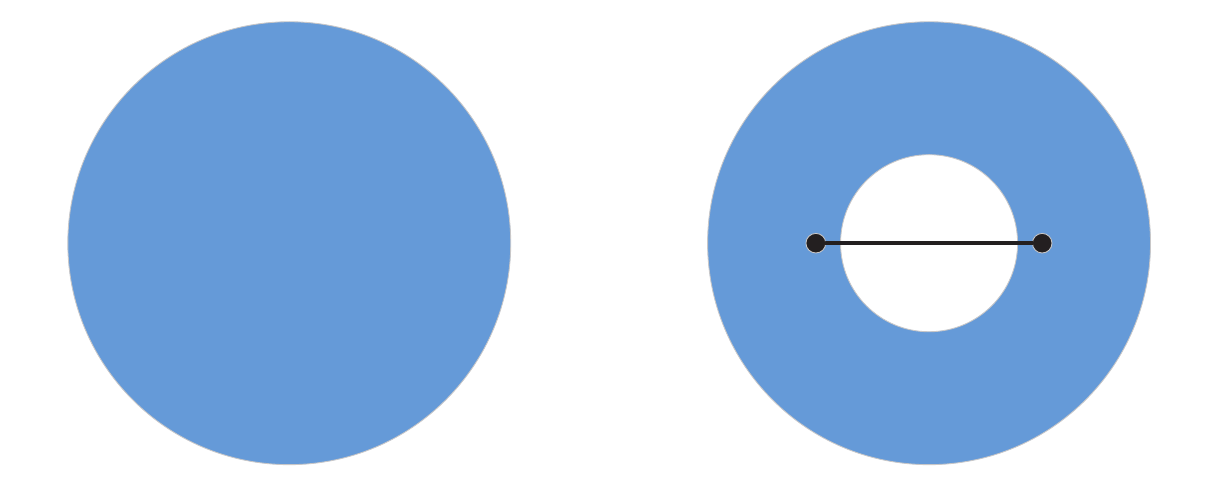}
\caption{Left: the circle is convex; right: the ring is non-convex.}
\label{fig:App-01-01}
\end{figure}

A set $C \in \mathbb R^n$ is convex if the line segment connecting any two points in $C$ is contained in $C$, i.e., for any $x_1,x_2 \in C$, we have $\theta x_1 + (1-\theta)x_2 \in C$, $\forall \theta \in [0,1]$. Roughly speaking, standing at anywhere in a convex set, you can see every other point in the set. Fig. \ref{fig:App-01-01} illustrates a simple convex set and a non-convex set in $\mathbb R^2$.

The convex combination of $k$ points $x_1,\cdots,x_k$ is defined as $\theta_1 x_1 + \cdots + \theta_k x_k$, where $\theta_1,\cdots,\theta_k \ge 0$, and $\theta_1 + \cdots +\theta_k = 1$. A convex combination of points can be regarded as a weighted average of the points, with $\theta_i$ the weight of $x_i$ in the mixture.  

\begin{figure}[!htp]
\centering
\includegraphics[scale=0.35]{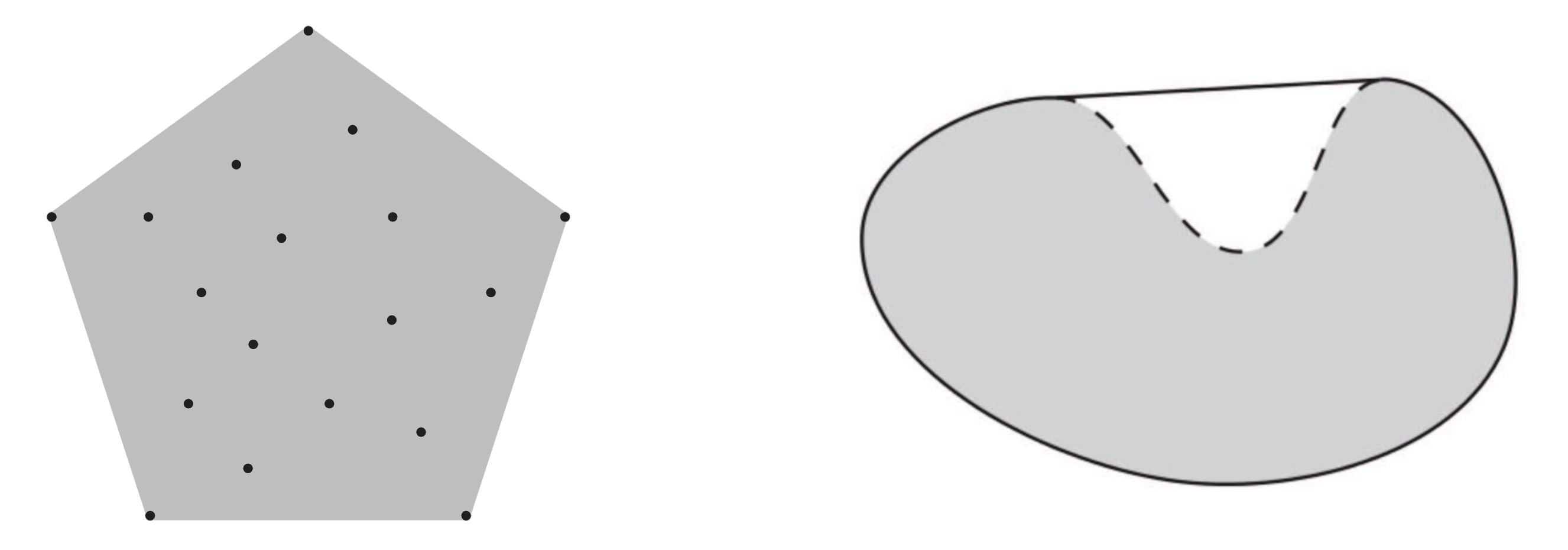}
\caption{Left: The convex hull of eighteen points. Right: The convex hull of a kidney shaped set.}
\label{fig:App-01-02}
\end{figure}

The convex hull of set $C$, denoted conv$(C)$, is the smallest convex set that contains $C$. Particularly, if $C$ has finite elements, then
\begin{equation}
\mbox{conv}(C) = \{ \theta_1 x_1 + \cdots + \theta_k x_k ~|~ x_i \in C,~ \theta_i \ge 0,~ i=1,~ \cdots,k,~ \theta_1 + \cdots + \theta_k = 1\} \notag
\end{equation}
Fig. \ref{fig:App-01-02} illustrates the convex hulls of two sets in $\mathbb R^2$.

Some useful convex sets are briefly introduced.

\vspace{12pt}
{\noindent \bf 1. Cones}

A set $C$ is called a cone, or nonnegative homogeneous, if for any $x \in C$, we have $\theta x \in C$, $\forall \theta \ge 0$. A set $C$ is a convex cone if it is convex and a cone: for any $x_1, x_2 \in C$ and $\theta_1,\theta_2 \ge 0$, we have $\theta_1 x_2 + \theta_2 x_2 \in C$.

\begin{figure}[!htp]
\centering
\includegraphics[scale=0.80]{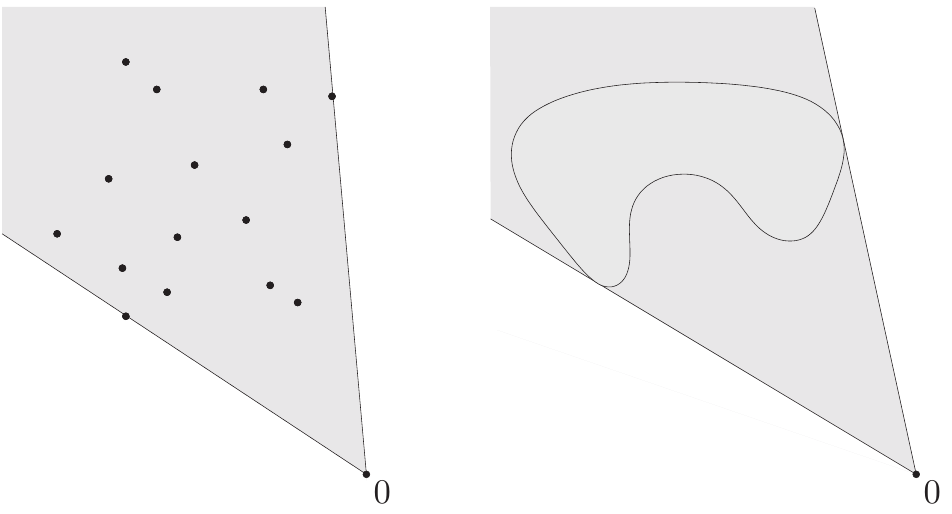}
\caption{The conic hulls of the two sets \cite{App-A-CVX-Book-Boyd}.}
\label{fig:App-01-03}
\end{figure}

The conic combination (or nonnegative linear combination) of $k$ points $x_1,\cdots,x_k$ is defined as $\theta_1 x_1 + \cdots + \theta_k x_k$, where $\theta_1,\cdots,\theta_k \ge 0$. If a set of finite points $\{x_i\}$, $i=1,2\cdots$ resides in a convex cone $C$, then every conic combination of $\{x_i\}$ remains in $C$. Conversely, a set $C$ is a convex cone if and only if it contains all conic combinations of its elements. The conic hull of set $C$ is the smallest convex cone that contains $C$. Fig. \ref{fig:App-01-03} illustrates the conic hulls of two sets in $\mathbb R^2$.

Some widely used cones are introduced.

\vspace{12pt}
{\noindent \bf a. The nonnegative orthant}

The nonnegative orthant is defined as
\begin{equation}
\label{eq:App-01-Orthant-P}
\mathbb R^n_+ = \{ x \in \mathbb R^n ~|~ x \ge 0 \}
\end{equation}
It is the set of vectors composed of non-negative entries. It is clearly a convex cone.

\vspace{12pt}
{\noindent \bf b. Second-order cone}

\begin{figure}[!htp]
\centering
\includegraphics[scale=0.70]{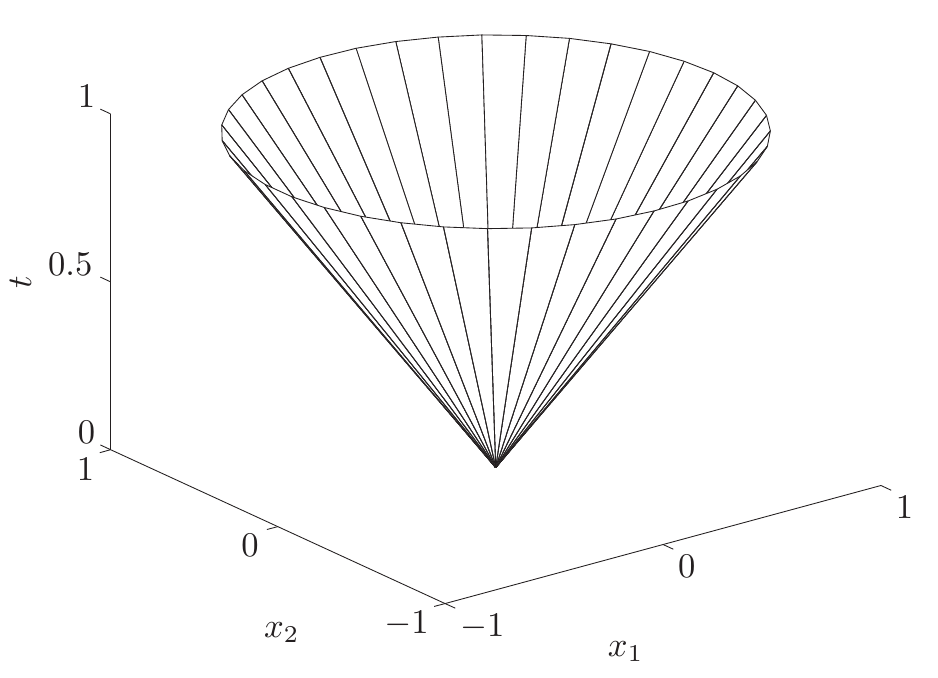}
\caption{$\mathbb L^3_C = \left\{ (x_1,x_2,t) ~\middle|~ \sqrt{x^2_1 + x^2_2} \le t \right\}$ in $\mathbb R^3$ \cite{App-A-CVX-Book-Boyd}.}
\label{fig:App-01-04}
\end{figure}

The unit second-order cone is defined as
\begin{equation}
\label{eq:App-01-Unit-SOC}
\mathbb L^{n+1}_C = \{(x,t) \in \mathbb R^{n+1} ~|~ \| x \|_2 \le t\}
\end{equation}
It is also called the Lorentz cone or ice-cream cone. Fig. \ref{fig:App-01-04} exhibits $\mathbb L^3_C$.

For any $(x,t) \in \mathbb L^{n+1}_C$ and $(y,z) \in \mathbb L^{n+1}_C$, we have
\begin{equation}
\| \theta_1 x + \theta_2 y \|_2 \le \theta_1 \|x\|_2 + \theta_2 \|y\|_2 \le \theta_1 t + \theta_2 z \Rightarrow 
\theta_1 \begin{bmatrix} x \\ t \end{bmatrix} + 
\theta_2 \begin{bmatrix} y \\ z \end{bmatrix} \in \mathbb L^{n+1}_C \notag
\end{equation}
which means that the unit second-order cone is a convex cone.

Sometimes, it is convenient to use the following inequality to represent a second-order cone in optimization problems 
\begin{equation}
\label{eq:App-01-General-SOC}
\|A x + b \|_2 \le c^T x + d
\end{equation}
where $A \in \mathbb R^{m \times n}$, $b \in \mathbb R^m$, $c \in \mathbb R^n$, $d \in \mathbb R$. It is the inverse image of the unit second-order cone under the affine mapping $f(x) =(A x + b, c^T x + d)$, and hence is convex. Second-order cones in forms of (\ref{eq:App-01-Unit-SOC}) and (\ref{eq:App-01-General-SOC}) are interchangeable.
\begin{equation}
\|A x + b \|_2 \le c^T x + d \Leftrightarrow 
\begin{bmatrix} A \\ c^T \end{bmatrix} x +
\begin{bmatrix} b \\ d \end{bmatrix} \in \mathbb L^{m+1}_C
\notag
\end{equation}
and hence is convex.

\vspace{12pt}
{\noindent \bf c. Positive semidefinite cone}

The set of symmetric $m \times m$ matrices is denoted by
\begin{equation}
\mathbb S^m = \{X \in \mathbb R^{m \times m}~|~ X = X^T\}   \notag
\end{equation} 
which is a vector space with dimension $m(m + 1)/2$. 

The set of symmetric positive semidefinite matrices is denoted by 
\begin{equation}
\mathbb S^m_+ = \{X \in \mathbb S^m ~|~ X \succeq 0 \}   \notag
\end{equation} 

The set of symmetric positive definite matrices is denoted by
\begin{equation}
\mathbb S^m_{++} = \{X \in \mathbb S^m ~|~ X \succ 0 \}   \notag
\end{equation} 

Clearly, $\mathbb S^m_+$ is a convex cone: if $A, B \in \mathbb S^m_+$, then for any $x \in \mathbb R^m$ and positive scalars $\theta_1,\theta_2 \ge 0$, we have
\begin{equation}
x^T(\theta_1 A + \theta_2 B) x = \theta_1 x^T A x + \theta_2 x^T B x \ge 0 \notag
\end{equation}
implying $\theta_1 A + \theta_2 B \in \mathbb S^m_+$.

\begin{figure}[!htp]
\centering
\includegraphics[scale=0.70]{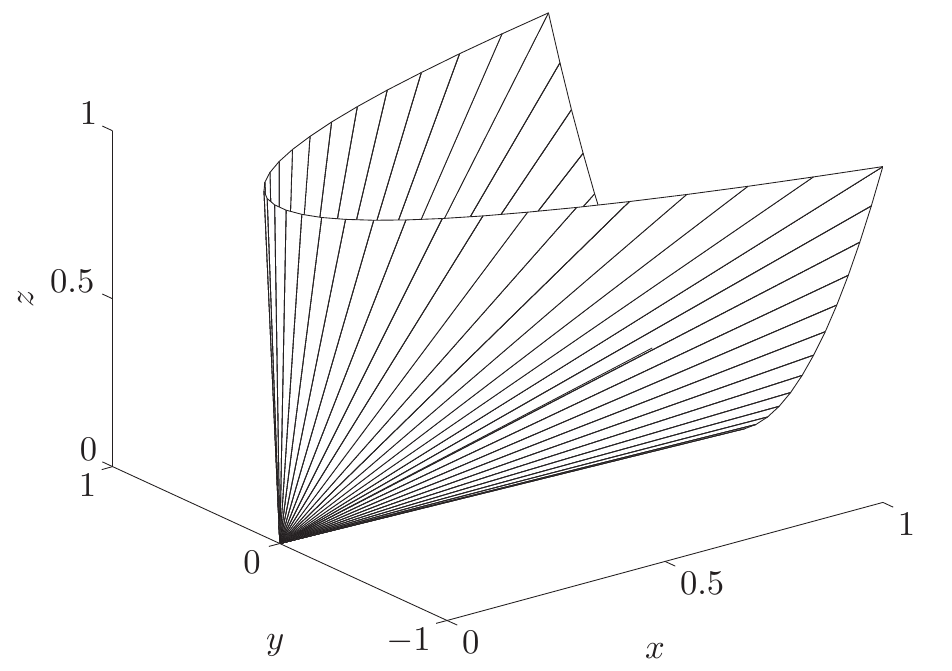}
\caption{Positive semidefinite cone in $\mathbb S^2$ (or in $\mathbb R^3$) \cite{App-A-CVX-Book-Boyd}.}
\label{fig:App-01-05}
\end{figure}

A positive semidefinite cone in $\mathbb R^2$ can be expressed via three variables $x,y,z$ as
\begin{equation}
\begin{bmatrix} 
x & y \\ y & z 
\end{bmatrix} \succeq 0 
\Leftrightarrow x \ge 0,~ xz \ge y^2 \notag
\end{equation}
which is plotted in Fig. \ref{fig:App-01-05}. In fact, $\mathbb L^3_C$ and $\mathbb S^2_+$ are equivalent to each other. To see this, the hyperbolic inequality $xz \ge y^2$ with $x \ge 0, z \ge 0$ defines the same feasible region in $\mathbb R^3$ as the following second-order cone
\begin{equation}
\left\| \begin{gathered}
2y \\ x-z 
\end{gathered} \right\|_2  \le  x + z,~
x \ge 0,~ z \ge 0  \notag
\end{equation}
In higher-order dimensions, every second-order cone can be written as an LMI via Schur complement as 
\begin{equation}
\label{eq:App-01-SOC-LMI}
\|A x + b \|_2 \le c^T x + d \Rightarrow 
\left[ \begin{gathered} (c^T x + d)I \\  (A x + b)^T \end{gathered}~~
\begin{gathered} A x + b \\ c^T x + d \end{gathered} \right] \succeq 0 
\end{equation}

In this sense of representability, positive semidefinite cones are more general than second-order cones. However, the transformation in (\ref{eq:App-01-SOC-LMI}) may not be superior from the computational perspective, because SOCPs are more tractable than SDPs.

\vspace{12pt}
{\noindent \bf d. Copositive cone}

A copositive cone $\mathbb C^n_+$ consists of symmetric matrices whose quadratic form is nonnegative over the nonnegative orthant $\mathbb R^n_+$:
\begin{equation}
\label{eq:App-01-COP-Cone}
\mathbb C^n_+ = \{ A ~|~ A \in \mathbb S^n,~ x^T A x \ge 0,~ 
\forall x \in \mathbb R^n_+ \} 
\end{equation}

The copositive cone $\mathbb C^n_+$ is closed, pointed, and convex \cite{App-A-COP-Cone-Property}. Clearly, $\mathbb S^n_+ \subseteq \mathbb C^n_+$, and every entry-wise nonnegative symmetric matrix $A$ belongs to $\mathbb C^n_+$. Actually, $\mathbb C^n_+$ is significantly larger than the positive semidefinite cone and the nonnegative symmetric matrix cone.

\vspace{12pt}
{\noindent \bf 2. Polyhedra}

A polyhedron is defined as the solution set of a finite number of linear inequalities:
\begin{equation}
\label{eq:App-01-Poly-H}
P = \{ x~|~ Ax \le b \}   
\end{equation}
(\ref{eq:App-01-Poly-H}) is also called a hyperplane representation for a polyhedron. It is easy to show that polyhedra are convex sets. Sometimes, a polyhedron is also called a polytope. The two concepts are often used interchangeably in this book. Because of physical bounds of decision variables, the polyhedral feasible regions in practical energy system optimization problems are usually bounded, which means that there is no extreme ray.

Polyhedra can be expressed via the convex combination as well. The convex hull of a finite number of points 
\begin{equation}
\label{eq:App-01-Poly-V-1}
\mbox{conv} \{v_1,\cdots,v_k\} = \{ \theta_1 v_1 + \cdots + \theta_k v_k ~|~ \theta \ge 0,~ 1^T \theta = 1 \}
\end{equation}
defines a polyhedron. (\ref{eq:App-01-Poly-V-1}) is called a convex hull representation. If the polyhedron is unbounded, a generalization of this convex hull representation is
\begin{equation}
\label{eq:App-01-Poly-V-2}
\{ \theta_1 v_1 + \cdots + \theta_k v_k ~|~ \theta \ge 0,~ \theta_1 + \cdots + \theta_m = 1, m \le k \} 
\end{equation}
which considers nonnegative linear combinations of $v_i$, but only the first $m$ coefficients whose summation is 1  are bounded, and the remaining ones can take arbitrarily large values. In view of this, the convex hull of points $v_1,\cdots,v_m$ plus the conic hull of points $v_{m+1}, \cdots,v_k$ is a polyhedron. The reverse is also correct: any
polyhedron can be represented by convex hull and conic hull.

How to represent a polyhedron depends on what information is available: if its boundaries are expressed via linear inequalities, the hyperplane representation is straightforward; if its extreme points and extreme rays are known in advance, the convex-conic hull representation is more convenient. With the growth in dimension, it is becoming more difficult to switch (derive one from the other) between the hyperplane representation and the hull representation.

\subsection{Generalized Inequalities}
\label{App-A-Sect01-02}

A cone $K \subseteq \mathbb R^n$ is called a proper cone if it satisfies:

1) $K$ is convex and closed.

2) $K$ is solid, i.e., it has non-empty interior.

3) $K$ is pointed, i.e., $x \in K$, $-x \in K$ $\Rightarrow x = 0$.

A proper cone $K$ can be used to define a generalized inequality, a partial ordering on $\mathbb R^n$, as follows
\begin{equation}
\label{eq:App-01-GNI}
x \preceq_K y  \Longleftrightarrow y - x \in K
\end{equation}
We denote $x \succeq_K y$ for $y \preceq_K x$. Similarly, a strict partial
ordering can be defined by 
\begin{equation}
\label{eq:App-01-GNI}
x \prec_K y  \Longleftrightarrow y - x \in \mbox{int}(K)
\end{equation}
where int$(K)$ stands for the interior of $K$, and write $x \succ_K y$ for $y \prec_K x$.

The nonnegative orthant $\mathbb R^n_+$ is a proper cone. When $K = \mathbb R^n_+$, the partial ordering $\preceq_K$ comes down to the element-wise comparison between vectors: for $x,y \in \mathbb R^n$, $x \preceq_{\mathbb R^n_+} y$ means $ x_i \le y_i$, $i = 1,\cdots n$, or the traditional notation $x \le y$. 

The positive semidefinite cone $\mathbb S^n_+$ is a proper cone in $\mathbb S^n$. When $K = \mathbb S^n_+$, the partial ordering $\preceq_K$ comes down to a linear matrix inequality  between symmetric matrices:  for $X,Y \in \mathbb S^n$, $X \preceq_{\mathbb S^n_+} Y$ means $ Y - X$ is positive semidefinite. Because it arises so frequently, we can drop the subscript $\mathbb S^n_+$ when we write a linear matrix inequality $Y \succeq X$ or $X \preceq Y$. It is understood that such a generalized inequality corresponds to the positive semidefinite cone without particular mention.

A generalized inequality is equivalent to linear constraints with $K = \mathbb R^n_+$; for other cones, such as the second-order cone $\mathbb L^{n+1}_C$ or  the positive semidefinite cone $\mathbb S^n_+$, the feasible region is nonlinear but remains convex.

\subsection{Dual Cones and Dual Generalized Inequalities}
\label{App-A-Sect01-03}

Let $K$ be a cone in $\mathbb R^n$. Its dual is defined as the following set
\begin{equation}
\label{eq:App-01-Dual-Cone}
K^* = \{y ~|~ x^T y \ge 0,~ \forall x \in K\} 
\end{equation}
Because $K^*$ is the intersection of homogeneous half spaces (half spaces passing through the origin). It is a closed convex cone.

The interior of $K^*$ is given by 
\begin{equation}
\label{eq:App-01-Interior-DCone}
\mbox{int}(K^*) = \{y ~|~ x^T y > 0,~ \forall x \in K,~ x \ne 0 \} 
\end{equation}
To see this, if $y^T x > 0$, $\forall x \in K$, then $(y+u)^T x > 0$, $\forall x \in K$ holds for all $u$ that is sufficiently small; hence $y \in \mbox{int}(K^*)$. Conversely, if $y \in K^*$ and $\exists x \in K : y^T x = 0$, $x \ne 0$, then $(y - tx)^T x < 0$, $\forall t > 0$, indicating $y \notin \mbox{int}(K^*)$.

If int$(K) \ne \emptyset$, then $K^*$ is pointed. If this is not true, suppose $\exists y \ne 0$: $y \in K^*$, $-y \in K^*$, i.e., $y^T x \ge 0$, $\forall x \in K$ and $-y^T x \ge 0$, $\forall x \in K$, so we have $x^T y = 0$, $\forall x \in K$, which is in contradiction with int$(K) \ne \emptyset$.

\begin{figure}[!htp]
\centering
\includegraphics[scale=0.50]{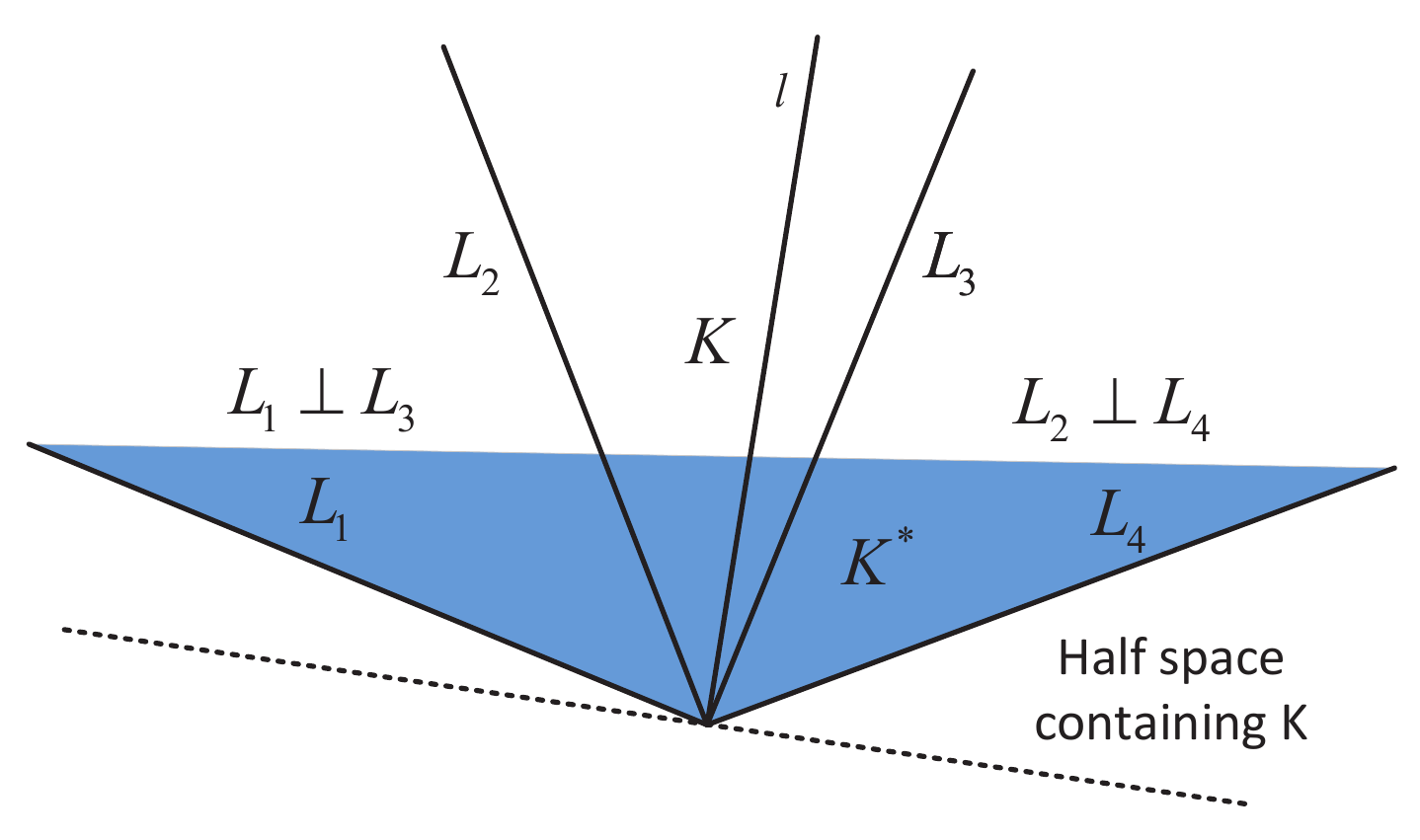}
\caption{Illustration of a cone and its dual cone in $\mathbb R^2$.}
\label{fig:App-01-06}
\end{figure}

In conclusion, $K^*$ is a proper cone, if the original cone $K$ is so; $K^*$ is closed and convex, regardless of the original cone $K$. Fig. \ref{fig:App-01-06} shows a cone $K$ (the region between $L_2$ and $L_3$) and its dual cone $K^*$ (the region between $L_1$ and $L_4$) in $\mathbb R^2$.

In light of the definition of $K^*$, a non-zero vector $y$ is the normal of a homogeneous half space which contains $K$ if and only if  $y \in K^*$.  The intersection of all such half spaces containing $K$ constitutes the cone $K$ (if $K$ is closed), in view of this
\begin{equation}
\label{eq:App-01-Double-Dual}
K = \bigcap_{y \in K^*} \left\{x ~|~ y^T x \ge 0 \right\} =
\{x ~|~ y^T x \ge 0,~ \forall y \in K^* \} = K^{**}
\end{equation}

This fact can be also understood in $\mathbb R^2$ from Fig. \ref{fig:App-01-06}. The extreme cases for the normal vector $y$ such that the corresponding half space contains $K$ are $L_1$ and $L_4$, and the intersection of these half spaces for all $y \in K^*$ turns out to be the original cone $K$.

Next, we investigate the dual cones of three special proper cones, i.e.,  $\mathbb R^n_+$, $\mathbb L^{n+1}_C$, and $\mathbb S^n_+$, respectively.

\vspace{12pt}
{\noindent \bf 1. The nonnegative orthant}

By observing the fact
\begin{equation}
x^T y \ge 0,~ \forall x \ge 0 \Longleftrightarrow y \ge 0 \notag
\end{equation}
we naturally have $(\mathbb R^n_+)^* = \mathbb R^n_+$; in other words, the nonnegative orthant is self-dual. 

\vspace{12pt}
{\noindent \bf 2. The second-order cone}

Now, we show that the second-order cone is also self-dual: $(\mathbb L^{n+1}_C)^* = \mathbb L^{n+1}_C$. To this end, we need to demonstrate  
\begin{equation}
x^T u + t v \ge 0,~ \forall (x,t) \in L^{n+1}_C 
\Longleftrightarrow \|u\|_2 \le v  \notag 
\end{equation}

$\Rightarrow$: Suppose the right-hand condition is false, and $\exists (u,v):\|u\|_2 > v$, by recalling Cauchy-Schwarz inequality $|a^T b| \le \|a\|_2 \|b\|_2$, we have
\begin{equation}
\min_{x} \left\{ x^T u ~\middle|~ \mbox{s.t.} \left\| x \right\|_2 \le t \right\} = - t \| u\|_2  \notag
\end{equation}
In such circumstance, $x^T u + t v = t(v - \|u\|_2) < 0$, $\forall t > 0$, which is in contradiction with the left-hand condition.

$\Leftarrow$: Again, according to Cauchy-Schwarz inequality, we have 
\begin{equation*}
x^T u + t v \ge -\|x\|_2 \|u\|_2 + t v \ge -\|x\|_2 \|u\|_2 + \|x\|_2 v = \|x\|_2 \left( v - \|u\|_2 \right) \ge 0
\end{equation*}

\vspace{12pt}
{\noindent \bf 3. The positive semidefinite cone}

We investigate the dual cone of $\mathbb S^n_+$. The inner product of $X,Y \in \mathbb S^n$ is defined by the element-wise summation 
\begin{equation}
\langle X,Y \rangle = \sum_{i=1}^n \sum_{j=1}^n X_{ij} Y_{ij}
= \mbox{tr}(X Y^T) \notag
\end{equation}

We establish this fact: $(\mathbb S^n_+)^* = \mathbb S^n_+$, which boils down to
\begin{equation}
\mbox{tr}(X Y^T) \ge 0,~ \forall X \succeq 0 
\Longleftrightarrow Y \succeq 0  \notag
\end{equation}

$\Rightarrow$: Suppose $Y \notin \mathbb S^n_+$, then $\exists q \in \mathbb R^n$ such that 
\begin{equation}
q^T Y q = \mbox{tr}(q q^T Y^T) < 0 \notag
\end{equation}
which is in contradiction with the left-hand condition because $X = q q^T \in \mathbb S^n_+$. 

$\Leftarrow$: Now suppose $X, Y \in \mathbb S^n_+$. $X$ can be expressed via its eigenvalues $\lambda_i \ge 0$ and eigenvectors $q_i$ as $X = \sum_{i=1}^n \lambda_i q_i q^T_i$, then we arrive at
\begin{equation}
\mbox{tr}(X Y^T) = \mbox{tr} \left( Y\sum_{i=1}^n \lambda_i q_i q^T_i \right)= \sum_{i=1}^n \lambda_i q^T_i Y q_i \ge 0 \notag
\end{equation}
In summary, it follows that the positive semidefinite cone is self-dual.

\vspace{12pt}
{\noindent \bf 4. The completely positive cone}

Following the same concept of matrix inner product, it is shown that $(\mathbb C^n_+)^*$ is the cone of so-called completely positive matrices and can be expressed as \cite{App-A-COP-Cone-Dual}
\begin{equation}
\label{eq:App-01-COMPL-Cone}
\mathbb (\mathbb C^n_+)^* = \mbox{conv} \{ x x^T ~|~ x \in \mathbb R^n_+\}
\end{equation}
In contrast to previous three cones, the copositive cone $\mathbb C^n_+$ is not self-dual.

\vspace{12pt}
When the dual cone $K^*$ is proper, it induces a generalized inequality $\preceq_{K^*}$, which is called the dual generalized inequality of the one induced by cone $K$ (if $K$ is proper). According to the definition of dual cone, an important fact relating a generalized inequality and its dual is

1) $x \preceq_K y$ if and only if $\lambda^T x \le \lambda^T y$, $\forall \lambda \in K^*$.

2) $x \prec_K y$ if and only if $\lambda^T x < \lambda^T y$, $\forall \lambda \in K^*$, $\lambda \ne 0$.

When $K = K^{**}$, the dual generalized inequality of $\preceq_{K^*}$ is $\preceq_K$, and the above property holds if the positions of $K$ and $K^*$ are swapped.

\subsection{Convex Function and Epigraph}
\label{App-A-Sect01-04}

A function $f:\mathbb R^n \to \mathbb R$ is convex if its feasible region $X$ is a convex set, and for all $x_1,x_2 \in X$, the following condition holds
\begin{equation}
\label{eq:App-01-Convex-Function}
f(\theta x_1 + (1-\theta)x_2) \le \theta f(x_1) + (1-\theta) f(x_2),~ 
\forall \theta \in [0,1]
\end{equation}

\begin{figure}[!htp]
\centering
\includegraphics[scale=0.60]{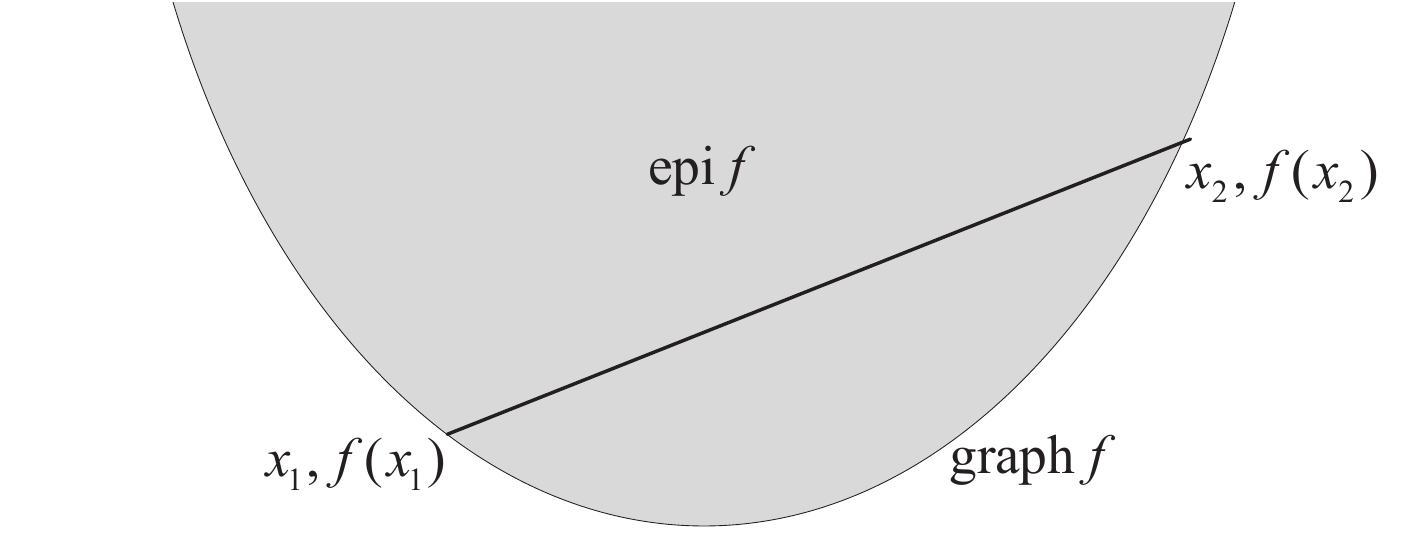}
\caption{Illustration of the graph of a convex function $f(x)$ (the solid line) and its epigraph (the shaded area) in $\mathbb R^2$.}
\label{fig:App-01-07}
\end{figure}

The geometrical interpretation of inequality (\ref{eq:App-01-Convex-Function}) is that the chord connecting points $(x_1, f(x_1))$ and $(x_2, f(x_2))$ always lies above the curve of $f$ between $x_1$ and $x_2$ (see Fig. \ref{fig:App-01-07}). Function $f$ is strictly convex if strict inequality holds in (\ref{eq:App-01-Convex-Function}) when $x_1 \ne x_2$ and $0 < \theta < 1$. Function $f$ is called (strictly) concave if $-f$ is (strictly) convex. An affine function is both convex and concave.

The graph of a function $f:\mathbb R^n \to \mathbb R$ is defined as
\begin{equation}
\label{eq:App-01-Graph-f}
\mbox{graph } f = \{(x,f(x))~|~ x \in X\}
\end{equation}
which is a subset of $\mathbb R^{n+1}$.

The epigraph of a function $f:\mathbb R^n \to \mathbb R$ is defined as
\begin{equation}
\label{eq:App-01-Epigraph-f}
\mbox{epi } f = \{(x,t)~|~ x \in X,~ f(x) \le t \}
\end{equation}
which is a subset of $\mathbb R^{n+1}$. These definitions are illustrated through Fig. \ref{fig:App-01-07}.

Epigraph bridges the concepts of convex sets and convex functions: A function is convex if and only if its epigraph is a convex set. Epigraph is frequently used in formulating optimization problems. A nonlinear objective function can be replaced by a linear objective and an additional constraint in epigraph form. In this sense, we can assume that any optimization problem has a linear objective function. Nonetheless, this does not facilitate solving the problem, as non-convexity moves to the constraints, if the objective function is not convex. Nonetheless, the solution to an optimization problem with a linear objective can always be found at the boundary of the convex hull of its feasible region, implying that if we can characterize the convex hull, a problem in epigraph form admits an exact convex hull relaxation. However, in general, it is difficult to express convex hull in an analytical form. 

Analyzing convex functions is a well developed field. Broadening the knowledge in convex analysis could be mathematically demanding, especially for readers who are primarily  interested in applications. We will not pursue in sophisticated theories in depth any more. Readers are referred to the literature suggested at the end of this chapter for further information.

\section{From Linear to Conic Program}
\label{App-A-Sect02}

Linear programming is one of the most mature and tractable mathematical programming problems. In this section, we first investigate and explain the motivation of linear programming duality theory, then provide a unified model for conic programming problems. LPs, SOCPs, and SDPs are special cases of conic  programs associated with generalized inequalities $\preceq_K$ where $K=\mathbb R^n$, $\mathbb L^{n+1}_C$, and $\mathbb S^n_+$, respectively. Our aim is to help readers who are not familiar with conic programs build their decision-making problems in these formats with structured convexity, and  write out their dual problems more conveniently. The presentation logic is consistent with \cite{App-A-CVX-Book-Ben}, and most of the presented materials in this section also come from \cite{App-A-CVX-Book-Ben}. 

\subsection{Linear Program and its Duality Theory}
\label{App-A-Sect02-01}

A linear program is an optimization program with the form
\begin{equation}
\label{eq:App-01-LP-Compact}
\min \{ c^T x ~|~ A x \ge b \}
\end{equation}
where $x$ is the vector of decision variables, $A$, $b$, $c$ are constant coefficient matrices with compatible dimensions. We assume LP (\ref{eq:App-01-LP-Compact}) is feasible, i.e., its feasible set $X = \{x~|~ Ax \ge b\}$ is a non-empty polyhedron; moreover, because of the limited ranges of decision variables representing physical quantities, we assume $X$ is bounded. In such circumstance, LP (\ref{eq:App-01-LP-Compact}) always has a finite optimum. LPs can be solved by mature algorithms, such as the simplex algorithm and the interior-point algorithm, which are not the main focus of this book.     

A question which is important both in theory and practice is: how to find a systematic way to bound the optimal value of (\ref{eq:App-01-LP-Compact})? Clearly, if $x$ is a feasible solution, an instant upper bound is given by $c^T x$. Lower bounding is to find a value $a$, such that $c^T x \ge a$ holds for all $x \in X$. 

A trivial answer is to solve the problem and retrieve its optimal value, which is the tightest lower bound. However, there may be a smarter way to retrieve a valid lower bound with much cheaper computational expense. To outline the basic motivation, let us consider the following example 
\begin{equation}
\label{eq:App-01-LP-Example}
\min \left\{ \sum_{i=1}^6 x_i ~\middle|~ \begin{gathered}
2 x_1 + 1 x_2 + 3 x_3 + 8 x_4 + 5 x_5 + 3 x_6 \ge 5  \\
6 x_1 + 2 x_2 + 6 x_3 + 1 x_4 + 1 x_5 + 4 x_6 \ge 2  \\
2 x_1 + 7 x_2 + 1 x_3 + 1 x_4 + 4 x_5 + 3 x_6 \ge 1  \\
\end{gathered} \right\}  
\end{equation}
Although LP (\ref{eq:App-01-LP-Example}) is merely a toy case for modern solvers and computers, one may guess it is still a little bit complicated for mental arithmetic. In fact, we can claim the optimal value is 0.8 at a glance without any sophisticated calculation: summing up the three constraints yields an inequality
\begin{equation}
\label{eq:App-01-LP-Example-Weigh-Sum}
10(x_1 + x_2 + x_3 + x_4 + x_5 + x_6) \ge 8 
\end{equation}
which immediately gives the optimal value is 0.8. To understand why such a value is indeed the optimum, by adding the constraints together and dividing both sides by 10, inequality (\ref{eq:App-01-LP-Example-Weigh-Sum}) implies that the objective function must get a value which is greater than or equal to 0.8 at any feasible point; moreover, to demonstrate that 0.8 is attainable, we can find a point $x^*$ which activates the three constraints simultaneously, so (\ref{eq:App-01-LP-Example-Weigh-Sum}) becomes an equality. LP duality is merely a formal generalization of this simple trick.

Multiplying each constraint in $Ax \ge b$ with a non-negative weight $\lambda_i$, and adding all constraints together, we will see  

\begin{equation}
\lambda^T A x \ge \lambda^T b \notag
\end{equation}

If we choose $\lambda$ elaborately such that $\lambda^T A = c^T$, then $\lambda^T b$ will be a valid lower bound of the optimal value of (\ref{eq:App-01-LP-Compact}). To improve the lower bound estimation, one may optimize the weighting vector $\lambda$, giving rise to the following problem  
\begin{equation}
\label{eq:App-01-LP-Dual}
\max_{\lambda} \{ \lambda^T b ~|~ A^T \lambda = c,~ \lambda \ge 0 \}
\end{equation}
where $\lambda$ is the vector of decision variables or dual variables, and the feasible region $D=\{\lambda ~|~ A^T \lambda = c,~ \lambda \ge 0\}$ is a polyhedron. Clearly, (\ref{eq:App-01-LP-Dual}) is also an LP, and is called the dual problem of LP (\ref{eq:App-01-LP-Compact}). Correspondingly, (\ref{eq:App-01-LP-Compact}) is called the primal problem. From above construction, we immediately conclude $c^T x \ge \lambda^T b$.

\begin{proposition}
\label{pr:App-01-LP-Weak-Duality}
(Weak duality): The optimal value of (\ref{eq:App-01-LP-Dual}) is less than or equal to the optimal value of (\ref{eq:App-01-LP-Compact}).
\end{proposition}

In fact, the optimal bound offered by (\ref{eq:App-01-LP-Dual}) is tight.

\begin{proposition}
\label{pr:App-01-LP-Strong-Duality}
(Strong duality): Optimal values of (\ref{eq:App-01-LP-Dual}) and (\ref{eq:App-01-LP-Compact}) are equal.
\end{proposition}

To see this, an explanation is given in \cite{App-A-CVX-Book-Ben}. If a real number $a$ is the optimal value of the primal LP (\ref{eq:App-01-LP-Compact}), the system of linear inequalities

\begin{equation}
S^P: \left\{ \begin{gathered} 
-c^T x > - a : \lambda_0 \\
 A x \ge b: \lambda
\end{gathered} \right. \notag
\end{equation}
must have an empty solution set, indicating that at least one of the following two systems does have a solution (called separation property later)
\begin{equation}
S^D_1: \left\{ \begin{gathered} 
-\lambda_0 c + A^T \lambda  =  0 \\
-\lambda_0 a + b^T \lambda \ge 0 \\
~~ \lambda_0 > 0,~~ \lambda \ge 0
\end{gathered} \right. \notag
\end{equation}
\begin{equation}
S^D_2: \left\{ \begin{gathered} 
-\lambda_0 c + A^T \lambda  =  0 \\
-\lambda_0 a + b^T \lambda  > 0 \\
~~ \lambda_0 \ge 0,~~ \lambda \ge 0
\end{gathered} \right. \notag
\end{equation}
  
We can show that $S^P$ has no solutions if and only if $S^D_1$ has a solution.

$S^D_1$ has a solution $\Rightarrow$ $S^P$ has no solution is clear. Otherwise, suppose that $S^P$ has a solution $x$, because $\lambda_0$ is strictly positive, the weighted summation of inequalities in $S^P$ leads to 
\begin{equation*}
0 = 0^T x =(-\lambda_0 c + A^T \lambda)^T x =-\lambda_0 c^T x + \lambda^T A x > -\lambda_0 a + \lambda^T b   
\end{equation*}
which is in contradiction with the second inequality in $S^D_1$. 

$S^P$ has no solution $\Rightarrow$ $S^D_1$ has a solution. Suppose $S^D_1$ has no solution, $S^D_2$ must have a solution owing to the separation property (Theorem 1.2.1 in \cite{App-A-CVX-Book-Ben}). Moreover, if $\lambda_0 > 0 $, the solution of system $S^D_2$ also solves system $S^D_1$, so there must be $\lambda_0 = 0$. As a result, the solution of $S^D_2$ is independent of the values of $a$ and $c$. Let $c = 0$ and $a = 0$, the solution $\lambda$ of $S^D_2$ satisfies $A^T \lambda = 0$, $b^T \lambda > 0$. Therefore, for any $x$ with a compatible dimension, $\lambda^T(Ax-b) = \lambda^T Ax - \lambda^T b <0$ holds. In addition, because $\lambda \ge 0$, we can conclude that $Ax \ge b$ has no solution, a contradiction to the assumption that (\ref{eq:App-01-LP-Compact}) is feasible. 

Now, consider the solution of $S^D_1$. Without loss of generality, we can assume $\lambda_0 = 1$; otherwise, if $\lambda_0 \ne 1$, ($1,\lambda/\lambda_0$) also solves $S^D_1$. In view of this, in normalized condition ($\lambda_0=1$), $S^D_1$  comes down to
\begin{equation}
S^D_3: \left\{
\begin{gathered} 
A^T \lambda  =  c \\
b^T \lambda \ge a \\
\lambda \ge 0
\end{gathered} \right. \notag
\end{equation}

Now we can see the strong duality: Let $a^*$ be the optimal solution of (\ref{eq:App-01-LP-Compact}). For any $a < a^*$, $S^P$ has no solution, so $S^D_1$ has a solution $(1,\lambda^*)$. According to $S^D_3$, the optimal value of (\ref{eq:App-01-LP-Dual}) is no smaller than $a$, i.e., $a \le b^T \lambda^* \le a^*$. When $a$ tends to $a^*$, we can conclude that the primal and dual optimal values are equal. Since the primal problem always has a finite optimum (as we assumed before), so does the dual problem, as they share the same optimal value. Nevertheless, even if the primal feasible region is bounded, the dual feasible set $D$ may be unbounded, and the dual problem is always bounded above. Please refer to \cite{App-A-LP-Book-Dantzig,App-A-LP-Book-Bertsimas,App-A-LP-Book-Vanderbei} for more information on duality theory in linear programming. 

\begin{proposition}
\label{pr:App-01-LP-OCD-Primal-Dual}
(Primal-dual optimality condition) If LP (\ref{eq:App-01-LP-Compact}) is feasible and $X$ is bounded, then any feasible solution to the following system 
\begin{equation}
\label{eq:App-01-LP-OCD-Primal-Dual}
\begin{lgathered}
A x \ge b \\
A^T \lambda = c,~ \lambda \ge 0 \\
c^T x = b^T \lambda 
\end{lgathered} 
\end{equation}
solves the original primal-dual pair of LPs: $x^*$ is the optimal solution of (\ref{eq:App-01-LP-Compact}), and $\lambda^*$ is the optimal solution of (\ref{eq:App-01-LP-Dual}).
\end{proposition}

\noindent  (\ref{eq:App-01-LP-OCD-Primal-Dual}) is also called the primal-dual optimality condition of LPs. It consists of linear inequalities and equalities, and there is no objective function to be optimized. 

Substituting $c = A^T \lambda$ into the last equation of (\ref{eq:App-01-LP-OCD-Primal-Dual}) gives $\lambda^T A x = \lambda^T b$, i.e.
\begin{equation}
\lambda^T (b-Ax) = 0 \notag 
\end{equation}
Since $\lambda \ge 0$ and $Ax \ge b$, above equation is equivalent to 
\begin{equation}
\lambda_i (b-Ax)_i = 0 \notag 
\end{equation}
where notation $(b-Ax)_i$ and $\lambda_i$ stand for the $i$-th components of vectors $b-Ax$ and $\lambda$, respectively. This condition means that at most one of $\lambda_i$ and $(b-Ax)_i$ can take a strictly positive value. In other words, if the $i$-th inequality constraint is inactive, then its dual multiplier $\lambda_i$ must be 0; otherwise, if $\lambda_i >0$, then the corresponding inequality constraint must be binding. This phenomenon is called the complementarity and slackness condition.  

Applying KKT optimality condition for general nonlinear programs to LP (\ref{eq:App-01-LP-Compact}) we have:
\begin{proposition}
\label{pr:App-01-LP-OCD-KKT}
(KKT optimality condition) If LP (\ref{eq:App-01-LP-Compact}) is feasible and $X$ is bounded, the following system 
\begin{equation}
\label{eq:App-01-LP-OCD-KKT}
\begin{gathered}
0 \le \lambda \bot A x - b \ge 0 \\
A^T \lambda = c
\end{gathered} 
\end{equation}
has a solution ($x^*,\lambda^*$) (may not be unique), where $a \bot b$ means $a^T b = 0$, $x^*$ solves (\ref{eq:App-01-LP-Compact}) and $\lambda^*$ solves (\ref{eq:App-01-LP-Dual}).
\end{proposition}

The question that which one of (\ref{eq:App-01-LP-OCD-Primal-Dual}) and (\ref{eq:App-01-LP-OCD-KKT}) is better can be subtle and has very different practical consequences. At the first look, the former one seems more tractable because (\ref{eq:App-01-LP-OCD-Primal-Dual}) is a linear system while (\ref{eq:App-01-LP-OCD-KKT}) contains complementarity and slackness conditions. However, the actual situation in practice is more complicated. For example, to solve a bilevel program with an LP lower level, the LP is often replaced by its optimality condition. In a bilevel optimization structure, some of the coefficients $A$, $b$, and $c$ are optimized by the upper-level agent, say, the coefficient vector $c$ representing the price is controlled by the upper level decision maker, while $A$ and $b$ are constants. If we use (\ref{eq:App-01-LP-OCD-Primal-Dual}), the term $c^T x$ in the single-level equivalence becomes non-convex, although $c$ is a constant in the lower level, preventing a global optimal solution from being found easily. In contrast to this, if we use (\ref{eq:App-01-LP-OCD-KKT}) and linearize the complementarity and slackness condition via auxiliary integer variables, the single-level equivalent problem can be formulated as an MILP, whose global optimal solution can be procured with reasonable computation effort.

The dual problem of LPs which maximize its objective can be derived in the same way. Consider the LP
\begin{equation}
\label{eq:App-01-LP-Max-Primal}
\max \{ c^T x ~|~ A x \le b \}
\end{equation}
For this problem, we need an upper bound on the objective function. To this end, associating a non-negative dual vector $\lambda$ with the constraint, and   adding the weighted inequalities together, we have  

\begin{equation}
\lambda^T A x \le \lambda^T b \notag
\end{equation}
If we intentionally choose $\lambda$ such that $\lambda^T A = c^T$, then $\lambda^T b$ will be a valid upper bound of the optimal value of (\ref{eq:App-01-LP-Max-Primal}). The dual problem   
\begin{equation}
\label{eq:App-01-LP-Max-Dual}
\min_{\lambda} \{ \lambda^T b ~|~ A^T \lambda = c,~ \lambda \ge 0 \}
\end{equation}
optimizes the weighting vector $\lambda$ to offer the tightest upper bound. 

Constraints in the form of equality  and $\ge$ inequality can be considered using the same paradigm. Bearing in mind that we are seeking an upper bound, so we need a certification for $c^T x \le a$, so the dual variables for equalities have no signs and those for $\ge$ inequalities should be negative.  

Sometimes it is useful to define the dual cone of a polyhedron, despite that a bounded polyhedron is not a cone. Recall its definition, the dual cone of a polyhedron $P$ can be defined as 
\begin{equation}
\label{eq:App-01-Dual-Polytope-1}
P^* = \{y ~|~ x^T y \ge 0,~ \forall x \in P \} 
\end{equation}
where $P = \{x ~|~ Ax \ge b \}$. As we have demonstrated in Sect. \ref{App-A-Sect01-03}, the dual cone is always closed and convex; however, for a general set, its dual cone does not have an analytical expression. 

For polyhedral sets, the condition in (\ref{eq:App-01-Dual-Polytope-1}) holds if and only if the minimal value of $x^T y$ over $P$ is non-negative. For a given vector $y$, let us investigate the minimum of $x^T y$ through an LP 
\begin{equation}
\min_x \{y^T x ~|~ Ax \ge b \} \notag
\end{equation}
It is known from Proposition \ref{pr:App-01-LP-Weak-Duality} that $y^T x \ge b^T \lambda$, $\forall \lambda \in D_P$, where $D_P =\{\lambda ~|~ A^T \lambda = y, \lambda \ge 0 \}$. Moreover, if $\exists \lambda \in D_P$ such that $b^T \lambda < 0$, Proposition \ref{pr:App-01-LP-Strong-Duality} certifies the existence of $x \in P$ such that $y^T x = b^T \lambda < 0$. In conclusion, the dual cone of polyhedron $P$ can be cast as 
\begin{equation}
\label{eq:App-01-Dual-Polytope-2}
P^* = \{ y ~|~ \exists \lambda: b^T \lambda \ge 0,~  
A^T \lambda = y,~ \lambda \ge 0 \} 
\end{equation}
which is also a polyhedron. It can be observed from (\ref{eq:App-01-Dual-Polytope-2}) that all constraints in $P^*$ are homogeneous, so $P^*$ is indeed a polyhedral cone.

\subsection{General Conic Linear Program}
\label{App-A-Sect02-02}

Linear programs cover vast topics in engineering optimization problems. Its duality program provides informative quantifications and valuable insights of the problem at hand, which help develop efficient algorithms for itself and facilitate building tractable reformulations for more complicated mathematical programming models, such as robust optimization, multi-level optimization and equilibrium problems. The algorithms of LPs, which are perfectly developed by now, can solve quite large instances (with up to hundreds of thousands of variables and constraints). Nevertheless, there are practical problems which cannot be modeled by LPs. To cope with these essentially nonlinear cases, one needs to explore new models and computational methods
beyond the reach of LPs.

The broadest class of optimization problems which the LP can be compared with is the class of convex optimization problems. Convexity marks whether a problem can be solved efficiently, and any local optimizer of a convex program must be a global optimizer. Efficiency is quantified by the number of arithmetic operations required to solve the problem. Suppose that all we know about the problem is its convexity: its objective and constraints are convex functions in decision variables $x \in \mathbb R^n$, and their values along with their derivatives at any given point can be evaluated within $M$ arithmetic operations. The best known complexity for finding an $\epsilon$-solution turns out to be \cite{App-A-CVX-Book-Ben} 
\begin{equation}
O(1) n(n^3+M) \ln \left( \frac{1}{\epsilon} \right) \notag
\end{equation}
Although this bound grows polynomially with $n$, the computation time may be still unacceptable for a large $n$ like $n = 1,000$, which is in contrast to LPs which are solvable with $n = 100,000$. The reason is: linearity are much stronger than convexity; the structure of an affine function $a^T x + b$ solely depends on its constant coefficients $a$  and $b$; function values and derivatives are never evaluated in a state-of-the-art LP solver. There are many classes of convex programs which are essentially nonlinear, but still possess nice analytical structure, which can be used to develop more dedicated algorithms. These algorithms may perform much more efficiently  than those exploiting only convexity. In what follows, we consider such a class of convex program, i.e., the conic program, which is a simple extension of LP. Its general form and mathematical model are briefly introduced, while the details about interior-point algorithms is beyond the scope of this book, which can be found in \cite{App-A-CVX-Book-Ben,App-A-CVX-Book-Boyd}. 

\vspace{12pt}
{\noindent \bf 1. Mathematical model}

When we consider to add some nonlinear factors in LP (\ref{eq:App-01-LP-Compact}), the most common way is to replace a linear function $a^T x$ with a nonlinear but convex function $f(x)$. As what has been explained, this may not be advantageous from a computational perspective. In contrast to this, we sustain all functions to be linear, but inject nonlinearity in the comparative operators $\ge$ or $\le$.  Recall the definition of generalized inequalities $\succeq_K$ with cone $K$, we consider the following problem in this section 
\begin{equation}
\label{eq:App-01-CP-Primal}
\min_x \{ c^T x ~|~ A x \succeq_K b \}
\end{equation}
which is called a conic programming problem. An LP is a special case of the conic program with $K = \mathbb R^n_+$. With this generalization, we are able to formulate a much wider spectrum of optimization problems which cannot be modeled as LPs, while enjoy nice properties of structured convexity. 

\vspace{12pt}
{\noindent \bf 2. Conic duality}

Aside from developing high-performance algorithms, the most important and elegant theoretical result in the area of LP is its duality theorem. In view of their similarities in mathematical appearances, how can the LP duality theorem be extended to conic programs? Similarly, the motivation of duality is the desire of a systematic way to  certify a lower bound on the optimal value of conic program (\ref{eq:App-01-CP-Primal}). Let us try the same trick: multiplying the dual vector $\lambda$ on both sides of $Ax \succeq_K b$, and adding them together, we obtain $\lambda^T A x$ and $b^T \lambda$; moreover, if we are lucky to get $A^T \lambda = c$, we guess $b^T \lambda$ can serve as a lower bound of the optimum of (\ref{eq:App-01-CP-Primal}) under some condition. The condition can be translated into: what is the admissible region of $\lambda$, such that the inequality $\lambda^T A x \ge b^T \lambda$ is a consequence of $Ax \succeq_K b$? A nice answer has been given at the end of Sect. \ref{App-A-Sect01-03}. Let us explain the problem from some simple cases.

Particularly, when $K = \mathbb R^n_+$, the admissible region of $\lambda$ is also $\mathbb R^n_+$, because we have already known the fact that the dual variable of $\ge$ inequalities in an LP which minimizes its objective should be non-negative. However, $\mathbb R^n_+$ is no longer a feasible region of $\lambda$ for conic programs with generalized inequality $\succeq_K$ if $K \ne \mathbb R^n_+$. To see this, consider $\mathbb L^3_C$ and the corresponding generalized inequality
\begin{equation}
\begin{bmatrix}  x \\  y  \\ z \end{bmatrix}  \succeq_{\mathbb L^3_C}
\begin{bmatrix}  0 \\  0  \\ 0 \end{bmatrix}  \Longleftrightarrow
z \ge \sqrt{x^2+y^2} \notag
\end{equation}
$(x,y,z)=(-1,-1,1.5)$ is a feasible solution. However, the weighted summation of both sides with $\lambda = [1,1,1]^T$ gives a false inequality $-0.5 \ge 0$.

To find the feasible region of $\lambda$, consider the condition 
\begin{equation}
\label{eq:App-01-CP-Lambda}
\forall a \succeq_K 0 ~\Rightarrow~  \lambda^T a \ge 0
\end{equation}
If (\ref{eq:App-01-CP-Lambda}) is true, we have the following logical inferences
\begin{equation}
\begin{gathered} 
\\ \Leftrightarrow \\ \Rightarrow \\ \Leftrightarrow 
\end{gathered}  \quad
\begin{aligned}
     Ax & \succeq_K b        \\
Ax - b & \succeq_K 0         \\
\lambda^T (Ax-b) & \ge 0     \\
\lambda^T A x  & \ge \lambda^T b
\end{aligned}   \notag
\end{equation}
Conversely, if $\lambda$ is an admissible vector for certifying
\begin{equation}
\forall (a,b: a \succeq_K b) ~\Rightarrow~  \lambda^T a \ge \lambda^T b \notag
\end{equation}
then, (\ref{eq:App-01-CP-Lambda}) is clearly true by letting $b=0$. Therefore, the admissible set of $\lambda$ for generalized inequality $\succeq_K$ with cone $K$ can be written as
\begin{equation}
\label{eq:App-01-CP-Dual-Cone}
K^* =\{ \lambda ~|~ \lambda^T a \ge 0,~ \forall a \in K \}
\end{equation}
which contains vectors whose inner products with all vectors belonging to $K$ are nonnegative. Recall the definition in (\ref{eq:App-01-Dual-Cone}), we can observe that the set $K^*$ is actually the dual cone of cone $K$.

Now we are ready to setup the dual problem of conic program (\ref{eq:App-01-CP-Primal}). As in the case of LP duality, we try to recover the objective function from the linear combination of constraints by choosing a proper dual variable $\lambda$, i.e., $\lambda^T A x= c^T x$, in addition, $\lambda \in K^*$ ensures $\lambda^T A x \ge \lambda^T b$, implying that $\lambda^T b$ is a valid lower bound of the objective function. The best bound one can expect is the optimum of the problem  
\begin{equation}
\label{eq:App-01-CP-Dual}
\max_\lambda \{ b^T \lambda ~|~ A^T \lambda = c,~ \lambda \succeq_{K^*} 0\}
\end{equation}
which is also a conic program, and called the dual problem of conic program (\ref{eq:App-01-CP-Primal}). From above construction, we have already known that $c^T x \ge b^T \lambda$ is satisfied  for all feasible $x$ and $\lambda$, which is the weak duality of conic programs.

In fact, the primal-dual pair of conic programs has following properties:
\begin{proposition}
\label{pr:App-01-CP-Conic-Duality}
(Conic Duality Theorem) \cite{App-A-CVX-Book-Ben} : The following conclusions hold true for conic program (\ref{eq:App-01-CP-Primal}) and its dual (\ref{eq:App-01-CP-Dual}).

1) Conic duality is symmetric: the dual problem is still a conic one, and the primal and dual problems are dual to each other.

2) Weak duality holds: the duality gap $c^T x - b^T \lambda$ is nonnegative over the primal and dual feasible sets.

2) If either of the primal problem or the dual problem is strictly feasible and has a finite optimum, then the other is solvable, and the duality gap is zero: $c^T x^* = b^T \lambda^*$ for some $x^*$ and $\lambda^*$.

3) If either of the primal problem or the dual problem is strictly feasible and has a finite optimum, then a pair of primal-dual feasible solutions ($x, \lambda$) solves the respective problems if and only if
\begin{equation}
\label{eq:App-01-CP-OCD-PD}
\begin{gathered}
Ax \succeq_K b \\
A^T \lambda = c \\
\lambda \succeq_{K^*} 0  \\
c^T x = b^T \lambda
\end{gathered} 
\end{equation}
or
\begin{equation}
\label{eq:App-01-CP-OCD-KKT}
\begin{gathered}
0 \preceq_{K^*} \lambda \bot Ax - b \succeq_K 0 \\
A^T \lambda = c
\end{gathered}
\end{equation}
where (\ref{eq:App-01-CP-OCD-PD}) is called the primal-dual optimality condition, and (\ref{eq:App-01-CP-OCD-KKT}) is called the KKT optimality condition.
\end{proposition}

The proof can be found in \cite{App-A-CVX-Book-Ben} and is omitted here. To highlight the role of strict feasibility in Proposition \ref{pr:App-01-CP-Conic-Duality}, consider the following example  
\begin{equation}
\min_x \left\{ x_2 ~\middle|~ 
\begin{bmatrix} x_1 \\ x_2 \\ x_1 \end{bmatrix}
\succeq_{\mathbb L^3_C} 0
 \right\} \notag
\end{equation}
The feasible region is 
\begin{equation*}
\sqrt{x^2_1 + x^2_2} \le x_1 \Leftrightarrow x_2 = 0,~ x_1 \ge 0
\end{equation*}
 So its optimal value is 0. As explained before, second-order cones are self-dual: $(\mathbb L^3_C)^* = \mathbb L^3_C$, it is easy to see the dual problem is 
\begin{equation}
\max_\lambda \left\{ 0 ~\middle|~ \lambda_1 + \lambda_3 = 0,~
\lambda_2 = 1,~ \lambda \succeq_{\mathbb L^3_C} 0
 \right\} \notag
\end{equation}
The feasible region is 
\begin{equation}
\left\{ \lambda ~\middle|~ \sqrt{\lambda_1^2 + \lambda^2_2} \le \lambda_3,~
\lambda_3 \ge 0,~ \lambda_2 = 1,~ \lambda_1 = - \lambda_3 \right\} \notag
\end{equation}
which is empty, because $ \sqrt{(-\lambda_3)^2 + 1} > \lambda_3$.

This example demonstrates that the existence of a strictly feasible point is indispensable for conic duality. But this condition is not necessary in LP duality, which means strong duality holds in conic programming with stronger assumptions.  

Several classes of conic programs with particular cones are of special interests. The cones in these problems are self-dual, so we can set up the dual program  directly, which allows to explore deeply into the original problem, or convert it into equivalent formulations which are more computationally friendly. The structure of these relatively simple cones also helps develop efficient algorithms for corresponding conic programs. In what follows, we will investigate two
extremely important classes of conic programs.

\subsection{Second-order Cone Program}
\label{App-A-Sect02-03}

{\bf 1. Mathematical models of the primal and dual problems}

Second-order cone program is a special class of conic problem with $K =\mathbb L^{n+1}_C$. It minimizes a linear function over the intersection
of a polytope and the Cartesian product of second-order cones, and can be formulated as 
\begin{equation}
\label{eq:App-01-SOCP-Conic-Primal}
\min_x \left\{ c^T x ~\middle|~ Ax-b \succeq_K 0 \right\}
\end{equation}
where $x \in \mathbb R^n$, and $K = \mathbb L^{m_1}_C \times \cdots \times \mathbb L^{m_k}_C \times \mathbb R^{m_p}_+$, in other words, the conic constraints in (\ref{eq:App-01-SOCP-Conic-Primal}) can be expressed as $k$ second-order cones $A_i x - b_i \succeq_{\mathbb L^{m_i}} 0$, $i = 1,\cdots,k$ plus one polyhedron $A_p x - b_p \ge 0$ with the following matrix partition
\begin{equation}
\begin{bmatrix} A; b\end{bmatrix} =
\begin{bmatrix}
[A_1;b_1]  \\  \vdots \\ [A_k;b_k] \\ [A_p;b_p]
\end{bmatrix}
\notag
\end{equation}

Recall the definition of second-order cone, we further partition the sub-matrices $A_i,b_i$ into 
\begin{equation}
\begin{bmatrix} A_i; b_i \end{bmatrix} = \left[ 
\begin{gathered} D_i \\ p^T_i  \end{gathered} ~~
\begin{gathered} d_i \\ q_i    \end{gathered} ~  \right],~
i = 1,\cdots,k  \notag
\end{equation}
where $D_i \in \mathbb R^{(m_i-1) \times n}$, $p_i \in \mathbb R^n$, $d_i \in \mathbb R^{m_i-1}$, $q_i \in \mathbb R$. Then we can write (\ref{eq:App-01-SOCP-Conic-Primal}) as
\begin{equation}
\label{eq:App-01-SOCP-MP-Primal}
\begin{aligned}
\min_x ~~ & c^T x   \\
\mbox{s.t.}~~ & A_p x \ge b_p  \\
& \|D_i x - d_i\|_2 \le p^T_i x -q_i,~ i=1,\cdots,k
\end{aligned}   
\end{equation}
(\ref{eq:App-01-SOCP-MP-Primal}) is often more convenient for model builders.

It is easy to see that the cone $K$ in (\ref{eq:App-01-SOCP-Conic-Primal}) is self-dual, as both second-order cone and non-negative orthant are self-dual. In this regard, the dual problem of SOCP (\ref{eq:App-01-SOCP-Conic-Primal}) can be expressed as
\begin{equation}
\label{eq:App-01-SOCP-Conic-Dual}
\max_\lambda \left\{ b^T \lambda ~\middle|~ A^T \lambda = c,~ 
\lambda \succeq_K 0 \right\}
\end{equation}

Partitioning the dual vector as
\begin{equation}
\lambda = \begin{bmatrix} \lambda_1 \\ \vdots \\ \lambda_k \\ \lambda_p \end{bmatrix},~
\lambda_i \in \mathbb L^{m_i}_C,~ i = 1, \cdots, k,~ \lambda_p \ge 0 \notag
\end{equation}
We can write the dual problem as
\begin{equation}
\label{eq:App-01-SOCP-Conic-Dual-Decomp}
\begin{aligned}
\max_\lambda~ & \sum_{i=1}^k b^T_i \lambda_i + b^T_p \lambda_p \\
\mbox{s.t.}~~ & \sum_{i=1}^k A^T_i \lambda_i + A^T_p \lambda_p = c \\
& \lambda_i \in \mathbb L^{m_i}_C,~ i = 1, \cdots, k  \\
& \lambda_p \ge 0
\end{aligned} 
\end{equation}

We further partition $\lambda_i$ according to the norm representation in (\ref{eq:App-01-SOCP-MP-Primal})
\begin{equation}
\lambda_i = \begin{bmatrix} \mu_i \\ \nu_i \end{bmatrix},~
\mu_i \in \mathbb R^{m_i-1},~ \nu_i \in \mathbb R \notag
\end{equation}
all second-order cone constraints are associated with dual variables as
\begin{equation}
\begin{bmatrix} D_i x \\ p^T_i x \end{bmatrix} - 
\begin{bmatrix} d_i   \\ q_i     \end{bmatrix}
\in \mathbb L^{m_i}_C : 
\begin{bmatrix} \mu_i \\ \nu_i \end{bmatrix},~ 
i = 1, \cdots, k  \notag
\end{equation}
so the admissible region of dual variables $(\mu_i,\nu_i)$ is 
\begin{equation}
\begin{bmatrix} \mu_i \\ \nu_i \end{bmatrix} \in (\mathbb L^{m_i}_C)^*  
\Rightarrow \|\mu_i\|_2 \le \nu_i \notag
\end{equation}

Finally, we arrive at the dual form of (\ref{eq:App-01-SOCP-MP-Primal})
\begin{equation}
\label{eq:App-01-SOCP-MP-Dual}
\begin{aligned}
\max_\lambda~ & \sum_{i=1}^k \left(\mu^T_i d_i + \nu_i q_i \right) + b^T_p \lambda_p \\
\mbox{s.t.}~~ & \sum_{i=1}^k \left(D^T_i \mu_i + \nu_i p_i \right) + A^T_p \lambda_p = c \\
& \|\mu_i\|_2 \le \nu_i,~ i = 1, \cdots, k  \\
& \lambda_p \ge 0
\end{aligned} 
\end{equation}

(\ref{eq:App-01-SOCP-MP-Primal}) and (\ref{eq:App-01-SOCP-MP-Dual}) are more convenient than (\ref{eq:App-01-SOCP-Conic-Primal}) and (\ref{eq:App-01-SOCP-Conic-Dual-Decomp}) respectively because norm constraints can be recognized by most commercial solvers, whereas generalized inequalities $\succeq_K$ and constraints with the form $\in \mathbb L^{m_i}_C$ are supported only in some dedicated packages. Strict feasibility can be expressed in a more straightforward manner via norm constraints: the primal problem is strictly feasible if $\exists x: \|D_i x - d_i\|_2 < p^T_i x -q_i,~ i=1,\cdots,k,~A_p x > b_p$; the dual problem is strictly feasible if $\|\mu_i\|_2 < \nu_i,~ i = 1, \cdots, k,~ \lambda_p > 0$. In view of this, (\ref{eq:App-01-SOCP-MP-Primal}) and (\ref{eq:App-01-SOCP-MP-Dual}) are treated as the standard forms of an SOCP and its dual by practitioners whose primary interests are applications.

\vspace{12pt}
{\noindent \bf 2. What can be expressed via SOCPs?}

Mathematical programs raised in engineering applications may not always appear in standard convex forms, and convexity may be hidden in seemingly non-convex expressions. Therefore, an important step is to recognize the potential existence of a convex form that is equivalent to the original formulation. This task can be rather tricky. We introduce some frequently used functions and constraints that can be represented by second-order cone constraints.

\vspace{6pt}
{\noindent \bf a. Convex quadratic constraints}

A convex quadratic constraint has the form
\begin{equation}
\label{eq:App-01-SOCP-Convex-QC-1}
x^T P x + q^T x + r \le 0
\end{equation}
where $P \in \mathbb S^n_+$, $q \in \mathbb R^n$, $r \in \mathbb R$ are constant coefficients. Let $t = q^Tx + r$, we have
\begin{equation}
t = \frac{(t+1)^2}{4} - \frac{(t-1)^2}{4} \notag
\end{equation}
Performing the Cholesky factorization $P = D^T D$, (\ref{eq:App-01-SOCP-Convex-QC-1}) can be represented by 
\begin{equation}
\|Dx\|^2_2 +  \frac{(t+1)^2}{4} \le \frac{(t-1)^2}{4}  \notag
\end{equation}
So (\ref{eq:App-01-SOCP-Convex-QC-1}) is equivalent to the following second-order cone constraint
\begin{equation}
\label{eq:App-01-SOCP-Convex-QC-2}
\left\| \begin{gathered}
2 D x \\  q^Tx + r + 1
\end{gathered} \right\|_2 \le
q^Tx + r - 1 
\end{equation}

However, not every second-order cone constraint can be expressed via a convex quadratic constraint. By squaring $\|D x - d\|_2 \le p^T x -q$ we get an equivalent quadratic inequality
\begin{equation}
x^T (D^T D - p p^T ) x + 2(q p^T - d^T D)x +d^T d -q^2 \le 0
\end{equation}
with $p^T x - q \ge 0$. The matrix $M = D^T D - p p^T$ is not always positive semidefinite. Indeed, $M \succeq 0$ if and only if $\exists u, \|u\|_2 \le 1: p=D^T u$.  On this account, SOCPs are more general than convex QCQPs. 

\vspace{6pt}
{\noindent \bf b. Hyperbolic constraints}

Hyperbolic constraints are frequently encountered in engineering optimization problems. They are non-convex in their original forms but can be represented by a second-order cone constraint. A hyperbolic constraint has the form
\begin{equation}
\label{eq:App-01-SOCP-Hyper-1}
x^T x \le yz,~ y > 0,~ z > 0
\end{equation}
where $x \in \mathbb R^n$, $y,z \in \mathbb R_{++}$. Noticing the fact that $4yz = (y+z)^2 - (y-z)^2$, (\ref{eq:App-01-SOCP-Hyper-1}) is equivalent to the following second-order cone constraint
\begin{equation}
\label{eq:App-01-SOCP-Hyper-2}
\left\| \begin{gathered}
 2 x \\  y-z
\end{gathered} \right\|_2 \le
y+z,~ y > 0,~z > 0 
\end{equation}
However, a hyperbolic constraint can not be expressed via a convex quadratic constraint, because the compact quadratic form of (\ref{eq:App-01-SOCP-Hyper-1}) is 
\begin{equation}
\begin{bmatrix} x \\ y \\ z \end{bmatrix}^T P
\begin{bmatrix} x \\ y \\ z \end{bmatrix} \le 0,~
P = \begin{bmatrix} 
2I &  0  &  0  \\
 0 &  0  & -1  \\
 0 & -1  &  0 
\end{bmatrix}  \notag
\end{equation}
where the matric $P$ is indefinite.   

Many instances can be regarded as special cases of hyperbolic constraints, such as the upper branch of hyperbola 
\begin{equation}
\{ (x,y) ~|~ xy \ge 1,~ x > 0 \}  \notag
\end{equation}
and the epigraph of a fractional-quadratic function $g(x,s)=x^T x/s$, $s>0$
\begin{equation}
\left\{ (x,s,t) ~\middle|~ t \ge \frac{x^T x}{s},~ s > 0 \right\}  \notag
\end{equation}

\vspace{6pt}
{\noindent \bf c. Composition of second-order cone representable functions}

A function is called second-order cone representable if its epigraph can be represented by second-order cone constraints. Second-order cone representable functions are closed under composition \cite{App-A-SOCP-Boyd}. Suppose two univariate convex functions $f_1(x)$ and $f_2(x)$ are second-order cone representable, and $f_1(x)$ is monotonically increasing, the composition $g(x) = f_1(f_2(x))$ is also second-order cone representable, because its epigraph $\{(x,t) ~|~ g(x) \le t\}$ can be expressed by
\begin{equation}
\{(x,t) ~|~ \exists s: f_1(s) \le t,~ f_2(x) \le s \} \notag
\end{equation}
where $f_1(s) \le t$ and $f_2(x) \le s$ essentially come down to second-order cone constraints.

\vspace{12pt}
{\noindent \bf d. Maximizing the production of concave functions}

Suppose two functions $f_1(x)$
and $f_2(x)$  are concave with $f_1(x) \ge 0,~ f_2(x) \ge 0$, and $-f_1(x)$ and $-f_2(x)$ are second-order cone representable [which means $f_1(x) \ge t_1$ and $f_2(x) \ge t_2$ are (equivalent to) second-order cone constraints]. Consider the maximum of their production 
\begin{equation}
\label{eq:App-01-SOCP-Bargain-1}
\max_x \{ f_1(x) f_2(x) ~|~ x \in X, f_1(x) \ge 0, f_2(x) \ge 0 \}
\end{equation}
where the feasible region $X$ is the intersection of a polyhedron and second-order cones. It is not instantly clear whether problem (\ref{eq:App-01-SOCP-Bargain-1}) is a convex optimization problem or not. This formulation frequently arises in engineering applications, such as the Nash Bargaining problem and multi-objective optimization problems.

By introducing auxiliary variables $t,t_1,t_2$, it is immediately seen that problem (\ref{eq:App-01-SOCP-Bargain-1}) is equivalent to the following SOCP
\begin{equation}
\label{eq:App-01-SOCP-Bargain-2}
\begin{aligned}
\max_{x,t,t_1,t_2}~~ & t \\
\mbox{s.t.}~~ & x \in X  \\
&   t_1  \ge 0,~ t_2 \ge 0,~ t_1t_2 \ge t^2  \\
& f_1(x) \ge t_1,~ f_2(x) \ge t_2 
\end{aligned}
\end{equation}
At the optimal solution, $f_1(x)f_2(x)=t^2$.

\vspace{12pt}
{\noindent \bf 3. Polyhedral approximation of second-order cones}

Although SOCPs can be solved very efficiently, the state-of-the-art in numerical computing of SOCPs is still incomparable to that in LPs. The salient computational superiority of LPs inspires a question: can we approximate an SOCP by an LP without dramatically increasing the problem size? There have been other reasons to explore LP approximations for SOCPs. For example, to solve a bilevel program with an SOCP lower level, the SOCP should be replaced by its optimality conditions. However, the primal-dual optimality condition (\ref{eq:App-01-CP-OCD-PD}) may introduce bilinear terms, while the second-order cone complementarity constraints in KKT optimality condition (\ref{eq:App-01-CP-OCD-KKT}) cannot be linearized easily. If the SOCP can be approximated by an LP, then the KKT optimality condition can be linearized and the original bilevel program can be reformulated as an MILP. Clearly, if we only work in original variables, the number of additional constraints would quickly grow unacceptable with the increasing problem dimension and required accuracy.  In this section, we introduce the technique developed in \cite{App-A-SOCP-LP}, which lifts the problem into higher dimensions with moderate numbers of auxiliary variables and constraints. 

We start with the basic question: find a polyhedral $\epsilon$-approximation ${\rm \Pi}$ for $\mathbb L^3_C$ such that: 

1) If $x \in \mathbb L^3_C$, then $\exists u: (x,u) \in {\rm \Pi}$.

2) If $(x,u) \in {\rm \Pi}$ for some $u$, then $\sqrt{x^2_1+x^2_2} \le (1+\epsilon) x_3$.

Geometrically, the polyhedral cone $\rm \Pi$ includes a system of homogeneous linear equalities and inequalities in variables $x,u$; its projection on $x$-space is an $\epsilon$-outer approximation of $\mathbb L^3_C$, and the error bound is quantified by $\epsilon x_3$. The answer to this question is given in \cite{App-A-SOCP-LP}. It is shown that ${\rm \Pi}$ can be expressed by
\begin{equation}
\label{eq:App-01-SOCP-Polyhedral}
\begin{aligned}
(a) \quad & \left\{ \begin{lgathered} 
\xi^0 \ge |x_1|  \\  \eta^0 \ge  |x_2|  \end{lgathered} \right. \\
(b) \quad & \left\{ \begin{lgathered}
\xi^j   =  \cos \left( \frac{\pi}{2^{j+1}} \right) \xi^{j-1} +
           \sin \left( \frac{\pi}{2^{j+1}} \right) \eta^{j-1} \\
\eta^j \ge \left| -\sin \left( \frac{\pi}{2^{j+1}} \right) \xi^{j-1} +
           \cos \left( \frac{\pi}{2^{j+1}} \right) \eta^{j-1} \right|   
    \end{lgathered} \right. ,~ j=1,\cdots,v \\
(c) \quad & \left\{ \begin{lgathered}  \xi^v  \le x_3  \\  
\eta^v  \le \tan \left( \frac{\pi}{2^{v+1}} \right) \xi^v 
\end{lgathered} \right.
\end{aligned}
\end{equation}

Formulation (\ref{eq:App-01-SOCP-Polyhedral}) can be understood from an geometric point of view:

1) Given $x \in \mathbb L^3_C$, set $\xi^0 = |x_1|$, $\eta^0 = |x_2|$, which satisfies (a) in (\ref{eq:App-01-SOCP-Polyhedral}), and point $P^0 = (\xi^0,\eta^0)$ belongs to the first quadrant. Let 
\begin{equation}
\left\{ \begin{lgathered}
\xi^j   =  \cos \left( \frac{\pi}{2^{j+1}} \right) \xi^{j-1} +
           \sin \left( \frac{\pi}{2^{j+1}} \right) \eta^{j-1} \\
\eta^j  = \left| -\sin \left( \frac{\pi}{2^{j+1}} \right) \xi^{j-1} +
           \cos \left( \frac{\pi}{2^{j+1}} \right) \eta^{j-1} \right|   
    \end{lgathered} \right.   \notag
\end{equation}
which ensures (b). Point  $P^j = (\xi^j,\eta^j)$ is obtained from $P^{j-1}$ according to following operation: rotate $P^{j-1}$ by angle $\phi_j = \pi/2^{j+1}$ clockwise and get a mediate point $Q^{j-1}$; if $Q^{j-1}$ resides in the upper half-plane, $P^j = Q^{j-1}$; otherwise $P^j$ is the reflection of $Q^{j-1}$ with respect to the $x$-axis. By this construction, it is clear that all vectors from the origin to $P^j$ have the same Euclidean norm, i.e., $\|[x_1,x_2]\|_2$.
Moreover, as $P^0$ belongs to the first quadrant, the angle of $Q^0$ must satisfy  $-\pi/4 \le \arg(Q^0) \le \pi/4$, and $0 \le \arg(P^1) \le \pi/4$. With the procedure going on, we have $|\arg(Q^j)| \le \pi/2^{j+1}$, and $0 \le \arg(P^{j+1}) \le \pi/2^{j+1}$, for $j=1,\cdots,v$. In the last step, $\xi^v \le \|P^v\|_2 = \|[x_1,x_2]\|_2 \le x_3$ and $ 0 \le \arg(P^v) \le \pi/2^{v+1}$ hold, ensuring condition (c). In this manner, a point in $\mathbb L^3_C$ has been extended to a solution of (\ref{eq:App-01-SOCP-Polyhedral}).

2) Given $(x,u) \in {\rm \Pi}$, where $u = \{\xi^j,\eta^j\},j=1,\cdots,v$. Define $P^j = [\xi^j, \eta^j]$, and it directly follows from (a) and (b) that all $P^j$ belongs to the first quadrant, and $\left\|P^0\right\|_2 \ge \sqrt{x^2_1+x^2_2}$. Moreover, recall the construction of $Q^j$ in previous analysis, it is seen $\|P^j\|_2 = \|Q^j\|_2$; the absolute value of the vertical coordinate of $P^{j+1}$ is no less than that of $Q^j$; therefore, $\|P^{j+1}\|_2 \ge \|Q^j\|_2 = \|P^j\|_2$. At last
\begin{equation}
\left\| P^v \right\|_2 \le \dfrac{x_3}{\cos \left( \dfrac{\pi}{2^{v+1}} \right)}\notag
\end{equation} 
so we arrive at $\sqrt{x^2_1+x^2_2} \le (1+\epsilon) x_3$, where 
\begin{equation}
\epsilon =  \dfrac{1}{\cos \left( \dfrac{\pi}{2^{v+1}} \right)} - 1
\end{equation}
In this way, a solution of (\ref{eq:App-01-SOCP-Polyhedral}) has been approximately extended to $\mathbb L^3_C$.

Now, let us consider the general case: approximating 
\begin{equation}
\mathbb L^{n+1}_C = \left\{(y,t)~\middle|~ \sqrt{y^2_1+\cdots+y^2_n} \le t \right\}  \notag
\end{equation} via a polyhedral cone. Without loss of generality, we assume $n=2^K$. To make use of the outcome in (\ref{eq:App-01-SOCP-Polyhedral}), $y$ is split into $2^{K-1}$ pairs $(y_1,y_2),\cdots,(y_{n-1},y_n)$, which are called variables of generation 0. A successor variable is associated with each pair, which is called variable of generation 1, and is further divided into $2^{K-2}$ pairs and associated with variable of generation 2, and so on. After $K-1$ steps of dichotomy, we complete variable splitting with two variables of generation $K-1$. The only variable of generation $K$ is $t$. For notation convenience, let $y^l_i$ be $i$-th variable of generation $l$, the original vector $y=[y^0_1,\cdots,y^0_n]$, and $t=y^K_1$. The ``parents'' of $y^l_i$ are variables $y^{l-1}_{2i-1},y^{l-1}_{2i}$. The total number of variables in the ``tower'' is $2n-1$.

Using the tower of variables $y^l$, $\forall l$, the system of constraints
\begin{equation}
\label{eq:App-01-SOCP-Tower-Cone}
\sqrt{(y^{l-1}_{2i-1})^2 + (y^{l-1}_{2i})^2} \le y^l_i,~
i = 1,\cdots,2^{K-l},~ l = 1,\cdots,K 
\end{equation}
gives the same feasible region on $y$ as $L^{n+1}_C$, and each second-order cone in $\mathbb L^3_C$ in (\ref{eq:App-01-SOCP-Tower-Cone}) can be approximated by a polyhedral cone given in (\ref{eq:App-01-SOCP-Polyhedral}). 

The size of this polyhedral approximation is unveiled in \cite{App-A-SOCP-LP}:

1) The dimension of the lifted variable is $p \le n+O(1)\sum_{l=1}^K 2^{K-l} v_l$.

2) The number of constraints is $q \le O(1) \sum_{l=1}^K 2^{K-l} v_l$.

The quality of the approximation is \cite{App-A-SOCP-LP}
\begin{equation}
\beta = \prod_{l=1}^K \dfrac{1}{\cos \left( \dfrac{\pi}{2^{v_l+1}} \right)} - 1 \notag
\end{equation}

Given a desired tolerance $\epsilon$, choose
\begin{equation}
v_l = \lfloor O(1) l \ln \frac{2}{\epsilon} \rfloor \notag
\end{equation}
with a proper constant $O(1)$, we can guarantee the following bounds:
\begin{equation}
\begin{aligned}
\beta  & \le \epsilon \\
    p  & \le O(1)n \ln \frac{2}{\epsilon}  \\
    q  & \le O(1)n \ln \frac{2}{\epsilon}  \\
\end{aligned} \notag
\end{equation}
which implies that the required numbers of variables and constraints grow
linearly in the dimension of the target second-order cone.

\subsection{Semidefinite Program}
\label{App-A-Sect02-04}

{\bf 1. Notation clarification}

In this section, variables appear in the form of symmetric matrices, some notations should be clarified first. 

The Frobenius inner product of two matrices $A,B \in \mathbb M^n$ is defined by
\begin{equation}
\label{eq:App-01-SDP-Frobenius}
\langle A, B \rangle = \mbox{tr}(AB^T) = 
\sum_{i=1}^n \sum_{j=1}^n A_{ij} B_{ij}
\end{equation}

The Euclidean norm of a matrix $X \in \mathbb M^n$ can be defined through the Frobenius inner product as follows
\begin{equation}
\label{eq:App-01-SDP-Matrix-Norm}
\left\| X \right\|_2 = \sqrt{\langle X,X \rangle} = \sqrt{\mbox{tr}(X^T X)}
\end{equation}

Equipped with the Frobenius inner product, the dual cone of a given cone $K \subset \mathbb S^n$ is defined by
\begin{equation}
\label{eq:App-01-SDP-Dual-Cone}
K^*=\{Y \in \mathbb S^n ~|~ \langle Y, X \rangle \ge 0,~ \forall X \in K\}
\end{equation}
Among the cones in $\mathbb S^n$, this section talks about the positive semidefinite cone $\mathbb S^n_+$. As what has been demonstrated in Sect. \ref{App-A-Sect01-03}, $S^n_+$ is self-dual, i.e., $(\mathbb S^n_+)^* = \mathbb S^n_+$. The interior of cone $\mathbb S^n_+$ consists of all $n \times n$ matrices that are positive definite, and is denoted by $\mathbb S^n_{++}$.

\vspace{12pt}
{\noindent \bf 2. Primal and dual formulations of SDPs}

When $K=\mathbb S^n_+$, conic program (\ref{eq:App-01-CP-Primal}) boils down to an SDP
\begin{equation}
\min_x \{ c^T x ~|~ A x - b \in \mathbb S^n_+ \}  \notag
\end{equation}
 which minimizes a linear objective over the intersection of affine plane $y = Ax -b$ and the positive semidefinite cone $\mathbb S^n_+$. However, the notation in such a form is a little confusing:  $Ax-b$ is a vector, which is not dimensionally compatible with the cone $\mathbb S^n_+$. In fact, we have met a similar difficulty at the very beginning: the vector inner product does not apply to matrices, which is consequently replaced with the Frobenius inner product.   There are two prevalent ways to resolve the confliction in dimension, leading to different formulations which will be discussed.

\vspace{12pt}
{\noindent \bf a. Formulation based on vector decision variables}

In this formulation, $b$ is replaced with a matrix $B \in \mathbb S^n$, and $Ax$ is replaced with a linear mapping $\mathcal A x: \mathbb R^n \to \mathbb S^n$. In this way, $\mathcal A x - B$ becomes an element of $\mathbb S^n$. A simple way to specify the linear mapping $\mathcal A x$ is  
\begin{equation}
\mathcal A x = \sum_{j=1}^n x_j A_j,~ x = [x_1,\cdots,x_n]^T, ~
A_1,\cdots,A_n \in \mathbb S^n \notag
\end{equation}

With all these input matrices, an SDP can be written as 
\begin{equation}
\label{eq:App-01-SDP-LMI-Primal}
\min_x \{ c^T x ~|~ x_1 A_1 + \cdots + x_n A_n - B \succeq 0 \} 
\end{equation}
where the cone $\mathbb S^n_+$ is omitted in the operator $\succeq$ without causing confusion. The constraint in (\ref{eq:App-01-SDP-LMI-Primal}) is an LMI.
This formulation is general enough to capture the situation in which multiple LMIs exist, because  
\begin{equation}
\mathcal A_i x - B_i \succeq 0,~ i=1,\cdots,k 
\Leftrightarrow  \mathcal A x - B \succeq 0  \notag
\end{equation}
with $\mathcal A x = \mbox{Diag} (\mathcal A_1x,\cdots,\mathcal A_k x)$ and  $B = \mbox{Diag} (B_1,\cdots,B_k)$.

The general form of conic duality can be specified in the case when the cone $K = \mathbb S^n_+$. Associating a matrix dual variable $\rm \Lambda$ with the LMI constraint, and recalling the fact that $(\mathbb S^n_+)^*  = \mathbb S^n_+$, the dual program of SDP (\ref{eq:App-01-SDP-LMI-Primal}) reads:
\begin{equation}
\label{eq:App-01-SDP-LMI-Dual}
\max_{\rm \Lambda} \{ \langle B, {\rm \Lambda} \rangle ~|~ 
\langle A_i, {\rm \Lambda} \rangle =c_i,~ i = 1, \cdots,n, ~ 
{\rm \Lambda} \succeq 0 \} 
\end{equation}
which remains an SDP.

Apply conic duality theorem given in Proposition \ref{pr:App-01-CP-Conic-Duality}
to SDPs (\ref{eq:App-01-SDP-LMI-Primal}) and (\ref{eq:App-01-SDP-LMI-Dual}).

1) Suppose $A_1,\cdots,A_n$ are linearly independent， i.e., no nontrivial linear combination of $A_1,\cdots,A_n$ gives an all zero matrix.

2) The primal SDP (\ref{eq:App-01-SDP-LMI-Primal}) is strict feasible, i.e., $\exists x: x_1 A_1 + \cdots + x_n A_n \succ B$, and is solvable (the minimum is attainable) 

3) The dual SDP (\ref{eq:App-01-SDP-LMI-Dual}) is strict feasible, i.e.,
$\exists {\rm \Lambda} \succ 0: \langle A_i, {\rm \Lambda} \rangle =c_i,i = 1,\cdots,n$, and is solvable (the maximum is attainable). 

The optimal values of (\ref{eq:App-01-SDP-LMI-Primal}) and (\ref{eq:App-01-SDP-LMI-Dual}) are equal, and the complementarity and slackness condition
\begin{equation}
\label{eq:App-01-SDP-Complement}
\langle {\rm \Lambda}, x_1 A_1 + \cdots + x_n A_n - B \rangle = 0
\end{equation}
is necessary and sufficient for a pair of primal and dual feasible solutions ($x,{\rm \Lambda})$ to be optimal for their respective problems. For a pair of positive semidefinite matrices, it can be shown that 
\begin{equation}
\langle X Y \rangle = 0 \Leftrightarrow  XY = YX = 0  \notag
\end{equation}
indicating that the
eigenvalues of these two matrices in some certain basis are “complementary”: for every common eigenvector, at most one of the two eigenvalues of $X$ and $Y$ can be strictly positive.

\vspace{12pt}
{\noindent \bf b. Formulation based on matrix decision variables}

This formulation directly incorporates a matrix decision variable $X \in \mathbb S^n_+$, and imposes other restrictions on $X$ through linear equations. In the objective function, the vector inner product $c^T x$ is replaced by a  Frobenius inner product $\langle C, X \rangle$. In this way, an SDP can be written as 
\begin{equation}
\label{eq:App-01-SDP-Hybrid-Primal}
\begin{aligned}
\min_X  ~~  &   \langle C, X \rangle   \\
\mbox{s.t.}~~ & \langle A_i, X \rangle = b_i : \lambda_i,~ i=1,\cdots,m \\
            &  X \succeq 0 : {\rm \Lambda}
\end{aligned} 
\end{equation}
By introducing dual variables (following the colon) for individual constraints, the dual program of (\ref{eq:App-01-SDP-Hybrid-Primal}) can be constructed as 
\begin{equation}
\begin{aligned}
\max_{\lambda,\rm \Lambda}  ~~  &   
b^T \lambda + \langle 0, {\rm \Lambda} \rangle   \\
\mbox{s.t.}~~ & 
{\rm \Lambda} + \lambda_1 A_1 + \cdots + \lambda_n A_n = C \\
              & {\rm \Lambda} \succeq 0 
\end{aligned}   \notag
\end{equation}

Eliminating $\rm \Lambda$, we obtain
\begin{equation}
\label{eq:App-01-SDP-Hybrid-Dual}
\begin{aligned}
\max_\lambda  ~~  &   b^T \lambda   \\
\mbox{s.t.}~~ & C - \lambda_1 A_1 - \cdots - \lambda_n A_n \succeq 0 
\end{aligned} 
\end{equation}

It is observed that (\ref{eq:App-01-SDP-Hybrid-Primal}) and (\ref{eq:App-01-SDP-Hybrid-Dual}) are in the same form compared with (\ref{eq:App-01-SDP-LMI-Dual}) and (\ref{eq:App-01-SDP-LMI-Primal}), respectively, except for the signs of some coefficients. 

SDP handles positive semidefinite matrices, so it is especially powerful in eigenvalue related problems, such as Lyapunov stability analysis and controller design, which are the main field of control theorists. Moreover, every SOCP can be formulated as an SDP because 
\begin{equation}
\| y \|_2 \le t  \Leftrightarrow  
\begin{bmatrix} tI & y \\ y^T & t\end{bmatrix} \succeq 0 \notag
\end{equation}
Nevertheless, solving SOCPs via SDP may not be a good idea. Interior-point algorithms for SOCPs have much better worst-case complexity than those for SDPs. In fact, SDPs are extremely popular in the convex relaxation technique for non-convex quadratic optimization problems, owing to its ability to offer a nearly global optimal solution in many practical applications, such as the OPF problem in power systems. The SDP based convex relaxation method for non-convex QCQPs will be discussed in the next section. Here we talk about some special cases involving homogeneous quadratic functions or at most two non-homogeneous quadratic functions.

\vspace{12pt}
{\noindent \bf 3. Homogeneous quadratic programs}

Consider the following quadratic program 
\begin{equation}
\label{eq:App-01-SDP-Homo-QP}
\begin{aligned}
\min~~ &  x^T  B  x   \\
\mbox{s.t.}~~ &  x^T A_i  x \ge 0,~ i = 1,\cdots,m  
\end{aligned}
\end{equation}
where $A_1,\cdots,A_m,B \in \mathbb S^n$ are constant coefficients. Suppose that problem (\ref{eq:App-01-SDP-Homo-QP}) is feasible. Due to its homogeneity, the optimal value is clear: $-\infty$ or 0, depending on whether there is a feasible solution $x$ such that $x^T B x < 0$ or not. But it is unclear which situation takes place, i.e., to judge $x^T B x \ge 0$ over the intersection of homogeneous inequalities $x^T A_i  x \ge 0$, $i=1,\cdots,m$, or whether the implication
\begin{equation}
\label{eq:App-01-SDP-H-SLemma-1}
x^T A_i  x \ge 0,~ i = 1,\cdots,m  \Rightarrow x^T B x \ge 0
\end{equation}
holds.

\begin{proposition}
\label{pr:App-01-SDP-H-SLemma-1}
If there exist $\lambda_i \ge 0$, $i=1,2,\cdots$ such that $B \succeq \sum_i \lambda_i A_i$, then the indication in (\ref{eq:App-01-SDP-H-SLemma-1}) is true.
\end{proposition}

To see this, $B \succeq \sum_i \lambda_i A_i \Leftrightarrow x^T (B - \sum_i \lambda_i A_i ) x \ge 0 \Leftrightarrow x^T B x \ge \sum_i \lambda_i x^T A_i x$; therefore, $x^T B x$ is a direct consequence of $x^T A_i  x \ge 0, i = 1,\cdots,m$, as the right-hand side of the last inequality is non-negative. Proposition \ref{pr:App-01-SDP-H-SLemma-1} provides a sufficient condition for (\ref{eq:App-01-SDP-H-SLemma-1}), and necessity is generally not guaranteed. Nevertheless, if $m=1$, the condition is both necessary and sufficient.

\begin{proposition}
\label{pr:App-01-SDP-H-SLemma-2}
(S-Lemma) Let $A,B \in \mathbb S^n$ and a homogeneous quadratic inequality
\begin{equation}
(a) \quad x^T A x \ge 0 \notag
\end{equation}
is strictly feasible. Then the homogeneous quadratic inequality
\begin{equation}
(b) \quad x^T B x \ge 0 \notag
\end{equation}
is a consequence of (a) if and only if $\exists \lambda \ge 0: B \succeq \lambda A$.
\end{proposition}

Proposition \ref{pr:App-01-SDP-H-SLemma-2} is called the S-Lemma or S-Procedure. It can be proved by many means. The most instructive one, in our tastes, is based on the semidefinite relaxation, which can be found in \cite{App-A-CVX-Book-Ben}.

\vspace{12pt}
{\noindent \bf 4. Non-homogeneous quadratic programs with a single constraint}

Consider the following quadratic program 
\begin{equation}
\label{eq:App-01-SDP-Heter-QP}
\begin{aligned}
\min~~ &  f_0 (x) = x^T A_0 x + 2 b^T_0 x + c_0 \\
\mbox{s.t.}~~ & f_1(x) = x^T A_1 x + 2 b^T_1 x + c_1 \le 0  
\end{aligned}
\end{equation}
Let $f^*$ denote the optimal solution, so $f_0(x) - f^* \ge 0$ is a consequence of $-f_1(x) \ge 0$. A sufficient condition for this implication is $\exists \lambda \ge 0: f_0(x) - f^* + \lambda f_1(x) \ge 0$. The left-hand side is a quadratic function with matrix form
\begin{equation}
\left[ \begin{gathered} x \\  1  \end{gathered}  \right]^T
\left[ \begin{gathered} A_0 + \lambda A_1 \\ (b_0 + \lambda b_1)^T \end{gathered}\quad
\begin{gathered} b_0 + \lambda b_1 \\ c_0 + \lambda c_1 - f^* \end{gathered}\right]
\left[ \begin{gathered} x \\  1  \end{gathered}  \right]  \notag
\end{equation}
Its non-negativeness is equivalent to 
\begin{equation}
\left[ \begin{gathered} 
A_0 + \lambda A_1 \\ (b_0 + \lambda b_1)^T 
\end{gathered} \quad
\begin{gathered} 
b_0 + \lambda b_1 \\ c_0 + \lambda c_1 - f^* 
\end{gathered} \right]  \succeq 0
\end{equation}

Similar to the homogeneous case, this condition is also sufficient.   In view of this, the optimal value $f^*$ of (\ref{eq:App-01-SDP-Heter-QP}) solves the following SDP 
\begin{equation}
\begin{aligned}
\min_{\lambda,f}~~ & f  \\
\mbox{s.t.} ~~ &
\left[ \begin{gathered} 
A_0 + \lambda A_1 \\ (b_0 + \lambda b_1)^T 
\end{gathered} \quad
\begin{gathered} 
b_0 + \lambda b_1 \\ c_0 + \lambda c_1 - f 
\end{gathered} \right]  \succeq 0
\end{aligned}
\label{eq:App-01-SDP-Heter-QP-LMI}
\end{equation}

This conclusion is known as the non-homogeneous S-Lemma:

\begin{proposition}
\label{pr:App-01-SDP-N-SLemma-2}
(Non-homogeneous S-Lemma) Let $A_i \in \mathbb S^n$, $b_i \in \mathbb R^n$, and $c_i \in \mathbb R$, $i=0,1$, if $\exists x: x^T A_1 x + 2 b^T_1 x + c_1 < 0$, the implication
\begin{equation}
x^T A_1 x + 2 b^T_1 x + c_1 \le 0 \Rightarrow 
x^T A_0 x + 2 b^T_0 x + c_0 \le 0 \notag
\end{equation}
holds if and only if 
\begin{equation}
\exists \lambda \ge 0 : 
\left[ \begin{gathered} A_0 \\ b_0^T  \end{gathered} ~~
\begin{gathered} b_0 \\ c_0 \end{gathered} \right] \preceq \lambda
\left[ \begin{gathered} A_1 \\ b_1^T  \end{gathered} ~~
\begin{gathered} b_1 \\ c_1 \end{gathered} \right]  
\end{equation}
\end{proposition}

Because the implication can boil down to the maximum of quadratic function $x^T A_0 x + 2 b^T_0 x + c_0$ being non-positive over set $\{x|x^T A_1 x + 2 b^T_1 x + c_1 \le 0 \}$, which is a special case of (\ref{eq:App-01-SDP-Heter-QP}) by letting $f_0(x) = -x^T A_0 x - 2 b^T_0 x - c_0$, Proposition \ref{pr:App-01-SDP-N-SLemma-2} is a particular case of (\ref{eq:App-01-SDP-Heter-QP-LMI}) with the optimum $f^*=0$.

A formal proof based on semidefinite relaxation is given in \cite{App-A-CVX-Book-Ben}. Since a quadratic inequality describes an ellipsoid, Proposition \ref{pr:App-01-SDP-N-SLemma-2} can be used to test whether an  ellipsoid is contained in another one.

\begin{figure}[!htp]
\centering
\includegraphics[scale=0.50]{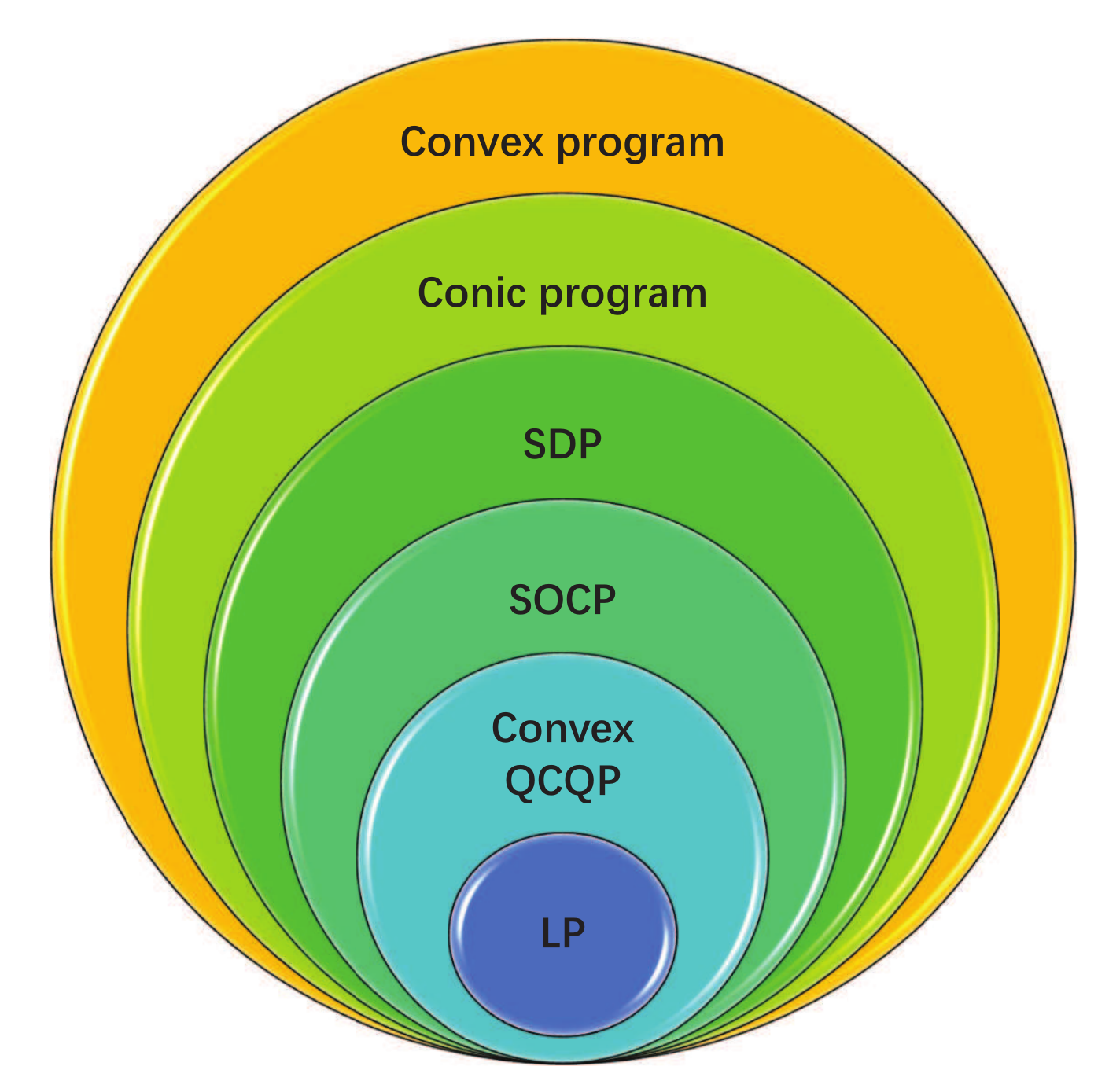}
\caption{Relations of the discussed convex programs.}
\label{fig:App-01-08}
\end{figure}

As a short conclusion, we summarize the relation of discussed convex programs in Fig. \ref{fig:App-01-08}. 

\section{Convex Relaxation Methods for Non-convex QCQPs}
\label{App-A-Sect03}

One of the most prevalent and promising applications of SDP is to build tractable approximations of computationally intractable optimization problems. One of the most quintessential appliances is the convex relaxation of
quadratically constrained  quadratic programs (QCQPs), which cover vast engineering optimization problems. QCQPs are generally non-convex and could have more than one locally optimal solution, and each of them may yield significant different objective values. However, gradient based algorithms can only find a local solution which largely depends on the initial point. One primary interest is to identify the global optimal solution or determine a high-quality bound for the optimum, which can be used to quantify the optimality gap of a given local optimal solution. The SDP relaxation technique for solving
non-convex QCQPs are briefly reviewed in this section.

\subsection{SDP Relaxation and Valid Inequalities}
\label{App-A-Sect03-01}

A standard fact of quadratic expression is
\begin{equation}
\label{eq:App-01-SDPr-Quad}
x^T Q x = \langle Q, x x^T   \rangle 
\end{equation}
where $\langle \cdot \rangle$ stands for the Frobenius inner product.

Following the logic in \cite{App-A-SDP-Relaxation-Tutor}, we focus our attention on QCQPs in the following form
\begin{equation}
\label{eq:App-01-SDPr-QCQP}
\min~~ \{x^T C x + c^T x ~|~ x \in F \}
\end{equation}
where
\begin{equation}
\label{eq:App-01-SDPr-QCQP-Fea}
F = \left\{ x \in \mathbb R^n ~\middle|~ x^T A_k x + a^T_k x \le b_k,~ 
k = 1, \cdots, m,~  l \le x \le u  \right\}
\end{equation}
All coefficient matrices and vectors have compatible dimensions. If $A_k = 0$ in all constraints, then the feasible set $F$ is a polyhedron, and (\ref{eq:App-01-SDPr-QCQP}) reduces to a quadratic program (QP); If $A_k \succeq 0$, $k=1,\cdots,m$ and $C \succeq 0$, (\ref{eq:App-01-SDPr-QCQP}) is a convex QCQP, which is easy to solve. Without loss of generality, we assume $A_k$, $k=1,\cdots,m$ and $C$ are indefinite, $F$ is a non-convex set, and the objective is a non-convex function. In fact, a number of hard optimization problems can be cast as non-convex QCQP (\ref{eq:App-01-SDPr-QCQP}). For example, a polynomial optimization problem can be reduced to a QCQP by introducing a tower of condensing variables, e.g., $x_1 x_2 x_3 x_4$ could be replaced by quadratic term $x_{12} x_{34}$ with $x_{12} = x_1 x_2$ and $x_{34}= x_3 x_4$. Moreover, a binary constraint $x \in \{0,1\}$ is equivalent to quadratic equality $x(x-1)=1$ where $x$ is continuous. 

A common idea to linearize non-convex terms $x^T A_k x$ is to define new variables $X_{ij}=x_i x_j$, $i=1,\cdots,n$, $j=1,\cdots,n$. In this way, $x^T A_k x = \sum_i \sum_j A_{ij} x_i x_j = \sum_{ij} A_{ij}X_{ij}$, and the last term is linear.  Recall (\ref{eq:App-01-SDPr-Quad}), this fact can be written in a compact form
\begin{equation}
x^T A_k x = \langle A_k, X\rangle, ~ X = x x^T  \notag
\end{equation}
With this transformation, QCQP (\ref{eq:App-01-SDPr-QCQP}) becomes
\begin{equation}
\label{eq:App-01-SDPr-QCQP-Ext}
\min~~ \{ \langle C, X \rangle + c^T x ~|~ (x,X) \in \hat F \}
\end{equation}
where
\begin{equation}
\label{eq:App-01-SDPr-QCQP-Fea-Ext}
\hat F = \left\{ (x,X) \in \mathbb R^n \times \mathbb S^n ~\middle|~ 
\begin{gathered} 
\langle A_k, X \rangle + a^T_k x \le b_k,~ k = 1, \cdots, m~ \\
 l \le x \le u,~~ X = x x^T  
\end{gathered} \right\}
\end{equation}
In problem (\ref{eq:App-01-SDPr-QCQP-Ext}), non-convexity are concentrated in the relation between the lifting variable $X$ and the original variable $x$, whereas all other constraints are linear. Moreover, if we replace $\hat F$ with its convex hull conv($\hat F$), the optimal solution of (\ref{eq:App-01-SDPr-QCQP-Ext}) will not change, because its objective function is linear. However, conv($\hat F$) does not have a closed form expression. Convex relaxation approaches can be interpreted as attempting to approximate conv($\hat F$) through structured convex constraints which can be recognized by existing solvers.

We define the following linear relaxation
\begin{equation}
\label{eq:App-01-SDPr-QCQP-Fea-LPr}
\hat L = \left\{ (x,X) \in \mathbb R^n \times \mathbb S^n ~\middle|~ 
\begin{gathered} 
\langle A_k, X \rangle + a^T_k x \le b_k \\
 k = 1, \cdots, m~ \\  l \le x \le u  
\end{gathered} \right\}
\end{equation}
which contains only linear constraints. Now, let us consider the lifting constraint
\begin{equation}
\label{eq:App-01-SDPr-X=xxT}
X = x x^T 
\end{equation}
which is called a rank-1 constraint. However, a rank constraint is non-convex and cannot be accepted by most solvers. Notice the fact that if (\ref{eq:App-01-SDPr-X=xxT}) holds, then
\begin{equation}
\begin{bmatrix} 1 & x^T \\ x  &   X   \end{bmatrix} =
\begin{bmatrix} 1 & x^T \\ x  & x x^T \end{bmatrix} =
\begin{bmatrix} 1 \\ x \end{bmatrix} 
\begin{bmatrix} 1 \\ x \end{bmatrix}^T  \succeq  0 \notag
\end{equation}

Define an LMI constraint
\begin{equation}
\label{eq:App-01-SDPr-QCQP-Fea-LMI}
\mbox{LMI} = \left\{ (x,X) ~\middle|~ Y=
\begin{bmatrix} 
1 & x^T \\ x  & X   
\end{bmatrix} \succeq  0 \right\}
\end{equation}
The positive semi-definiteness condition is true over conv($\hat F$). 

The basic SDP relaxation of (\ref{eq:App-01-SDPr-QCQP-Ext}) replaces the rank-1 constraint in $\hat F$ with a weaker but convex constraint (\ref{eq:App-01-SDPr-QCQP-Fea-LMI}), giving rise to the following SDP  
\begin{equation}
\label{eq:App-01-SDPr-QCQP-Basic}
\begin{aligned}
\min~~ &  \langle C, X \rangle + c^T x  \\ 
\mbox{s.t.}~~ & (x,X) \in \hat L \cap \mbox{LMI}
\end{aligned}
\end{equation}

Clearly, the LMI constraint enlarges the feasible region defined by (\ref{eq:App-01-SDPr-X=xxT}), so the optimal solution to (\ref{eq:App-01-SDPr-QCQP-Basic}) may not be feasible in the original QCQP, and the optimal value is a strict lower bound. In this situation, the SDP relaxation is inexact. Conversely, if matrix $Y$ is indeed rank-1 at the optimal solution, then the SDP relaxation is exact and $x$ solves the original QCQP (\ref{eq:App-01-SDPr-QCQP}). 

The basic SDP relaxation model  (\ref{eq:App-01-SDPr-QCQP-Basic}) can be further improved by enforcing additional linkages between $x$ and $X$, which are called valid inequalities. Suppose linear inequalities $\alpha^T x \le \alpha_0$ and $\beta^T x \le \beta_0$ are chosen from $\hat L$, then the quadratic inequality
\begin{equation}
(\alpha_0 - \alpha^T x) (\beta_0 - \beta^T x) = \alpha_0 \beta_0 - \alpha_0 \beta^T x -\beta_0 \alpha^T x + x^T \alpha \beta^T x \ge 0 \notag
\end{equation}
holds for all $x \in \hat L$. The last quadratic term can be linearized via the lifting variable $X$, resulting in the following linear inequality
\begin{equation}
\label{eq:App-01-SDPr-QCQP-VIN}
\alpha_0 \beta_0 - \alpha_0 \beta^T x - \beta_0 \alpha^T x + 
\langle \beta \alpha^T, X \rangle \ge 0 
\end{equation}

Any linear inequality in $\hat L$ (possibly the same) can be used to construct valid inequalities. Because additional constraints are imposed on $X$, the relaxation could be tightened, and the feasible region shrinks but may still be larger than conv($\hat F$). 

If we construct valid inequality (\ref{eq:App-01-SDPr-QCQP-VIN}) from side constraint $l \le x \le u$, we get
\begin{equation}
\label{eq:App-01-SDPr-QCQP-VIN-Quad}
\left. \begin{gathered} 
(x_i - l_i) (x_j - l_j)  \ge  0  \\
(x_i - l_i) (u_j - x_j)  \ge  0  \\
(u_i - x_i) (x_j - l_j)  \ge  0  \\
(u_i - x_i) (u_j - x_j)  \ge  0  \\
\end{gathered}  \right\},~~   
\forall i,j = 1,\cdots,n,~ i \le j  
\end{equation}
Expanding these quadratic inequalities,  the coefficients of quadratic terms $x_i x_j$ are equal to 1, and we obtain simple bounds on $X_{ij}$
\begin{equation}
\begin{gathered}
x_i l_j + x_j l_i - l_i l_j  \le  X_{ij}  \\
x_i u_j + x_j l_i - l_i u_j  \ge  X_{ij}  \\
u_i x_j - u_i l_j + x_i l_j  \ge  X_{ij}  \\
u_i x_j + x_i u_j - u_i u_j  \le X_{ij}
\end{gathered}  \notag
\end{equation}
or in a compact matrix form \cite{App-A-SDP-Relaxation-Tutor}
\begin{equation}
\label{eq:App-01-SDPr-QCQP-VIN-RLT}
\mbox {RLT} = \left\{ (x,X) ~\middle|~ 
\begin{gathered}
l x^T + x l^T - l l^T \le X  \\
u x^T + x u^T - u u^T \le X  \\
x u^T + l x^T - l u^T \ge X  
\end{gathered}  \right\}  
\end{equation}
(\ref{eq:App-01-SDPr-QCQP-VIN-RLT}) is known as the reformulation-linearization technique after the term appeared in \cite{App-A-APP-RLT}. These constraints have been extensively studied since it was proposed in \cite{App-A-Convex-Concave-Envelop}, due to the simple structure and satisfactory performance in various applications.
The improved SDP relaxation with valid inequalities can be written as
\begin{equation}
\label{eq:App-01-SDPr-QCQP-LMI+RLT}
\begin{aligned}
\min~~ &  \langle C, X \rangle + c^T x  \\ 
\mbox{s.t.}~~ & (x,X) \in \hat L \cap \mbox{LMI} \cap \mbox{RLT}
\end{aligned}
\end{equation}

From the construction of $\hat L$, LMI, and RLT, it is directly concluded that
\begin{equation}
\label{eq:App-01-SDPr-QCQP-Inclusion}
\mbox{conv}(\hat F) \subseteq \hat L \cap \mbox{LMI} \cap \mbox{RLT}
\end{equation}
The inclusion becomes tight only in some very special situations, such as those encountered in the homogeneous and non-homogeneous S-Lemma. Nevertheless, what we really need is the equivalence between the optimal solution of the relaxed problem (\ref{eq:App-01-SDPr-QCQP-LMI+RLT}) and that of the original problem (\ref{eq:App-01-SDPr-QCQP}): if the optimal matrix variable of (\ref{eq:App-01-SDPr-QCQP-LMI+RLT}) allows a rank-1 decomposition  
\begin{equation}
\begin{bmatrix} 1 & x^T \\ x  &   X   \end{bmatrix} =
\begin{bmatrix} 1 \\ x \end{bmatrix} 
\begin{bmatrix} 1 \\ x \end{bmatrix}^T   \notag
\end{equation}
which indicates $X$ has a rank-1 decomposition $X = x x^T$, then $x$ is optimal in (\ref{eq:App-01-SDPr-QCQP}), and the SDP relaxation is said to be exact, although $\mbox{conv}(\hat F)$ may be a strict subset of $\hat L \cap \mbox{LMI} \cap \mbox{RLT}$.

\subsection{Successively Tightening the Relaxation}
\label{App-A-Sect03-02}

If the matrix $X$ has a rank higher than 1, the corresponding optimal solution $x$ in (\ref{eq:App-01-SDPr-QCQP-LMI+RLT}) may be infeasible in (\ref{eq:App-01-SDPr-QCQP}).
The rank-1 constraint on $X$ can be exactly described by a pair of LMIs $X \succeq x x^T$ and $X \preceq x x^T$. The former one is redundant to (\ref{eq:App-01-SDPr-QCQP-Fea-LMI}) indicated by the Schur complement theorem; the latter one is non-convex, which is simply neglected in the SDP relaxation.

\vspace{12pt}
{\noindent \bf 1. A dynamical valid inequality generation approach}

An approach is proposed in \cite{App-A-SDP-Relaxation-Tutor} to generate valid inequalities dynamically by harnessing the constraint violations in $X \preceq x x^T$. The motivation comes from the fact that 
\begin{equation}
X -x x^T \preceq 0 \Leftrightarrow \langle X, v_i v^T_i \rangle \le (v^T_i x)^2,~ i=1,\cdots,n   \notag
\end{equation}
where  $\{v_1,\cdots,v_n\}$ is a set of orthogonal basis of $\mathbb R^n$. To see this, any vector $h \in \mathbb R^n$ can be expressed as the linear combination of the orthogonal basis as $h = \sum_{i=1}^n \lambda_i v_i$, therefore, $h^T (X-x x^T) h = \langle X, h h^T \rangle - (h^T x)^2 = \sum_{i=1}^n \lambda^2_i [\langle X, v_i v^T_i \rangle - (v^T_i x)^2] \le 0$. In view of this,
\begin{equation}
\hat F = \hat L \cap \mbox{LMI} \cap \mbox{NSD} \notag
\end{equation}
where
\begin{equation}
\begin{aligned}
\mbox{NSD} & = \{(x,X)~|~ X -x x^T \preceq 0 \}  \\
& = \{ (x,X) ~|~ \langle X, v_i v^T_i \rangle \le 
(v^T_i x)^2, i=1,\cdots,n \} 
\end{aligned} \notag
\end{equation}
If $\{v_1,\cdots,v_n\}$ is the standard orthogonal basis, 
\begin{equation}
\label{eq:App-01-SDPr-QCQP-NSD-1}
\mbox{NSD} = \{ (x,X) ~|~ X_{ii} \le x_i^2,~i=1,\cdots,n \} 
\end{equation}
It is proposed in \cite{App-A-SDP-Relaxation-Tutor} to construct NSD as 
\begin{equation}
\label{eq:App-01-SDPr-QCQP-NSD-2}
\mbox{NSD} = \{ (x,X) ~|~ \langle X, \eta_i \eta^T_i \rangle \le 
(\eta^T_i x)^2,~i=1,\cdots,n \} 
\end{equation}
where $\{\eta_1,\cdots,\eta_n\}$ are the eigenvectors of matrix $X - x x^T$, because they exclude infeasible points with respect to $X -x x^T \preceq 0$ most effectively.

Non-convex constraints in (\ref{eq:App-01-SDPr-QCQP-NSD-1}) and (\ref{eq:App-01-SDPr-QCQP-NSD-2}) can be handled by a special disjunctive programming derived in \cite{App-A-QCQP-Extended} and the convex-concave procedure investigated in \cite{App-A-CCP-Boyd}. The former one is an exact approach which requires binary variables to formulate disjunctive constraints; the latter is a heuristic approach which only solves convex optimization problems. We do not further detail these techniques here.

\vspace{12pt}
{\noindent \bf 2. A rank penalty method \cite{App-A-SDP-Rank-CCP}}

In view of the rank-1 exactness condition, another way to tighten SDP relaxation is to work on the rank of the optimal solution. A successive rank penalty approach is proposed in \cite{App-A-SDP-Rank-CCP}. We consider problem (\ref{eq:App-01-SDPr-QCQP-Ext}) as a rank-constrained SDP
\begin{equation}
\label{eq:App-01-SDP-Rank}
\min~~ \{ \langle {\rm \Omega}, Y \rangle ~|~ Y \in \hat L \cap \mbox{LMI} \cap \mbox{RLT},~ \mbox{rank}(Y) = 1 \}
\end{equation}
where
\begin{equation}
{\rm \Omega} = 
\begin{bmatrix}  
0 & 0.5 c^T \\ 0.5 c  &  C
\end{bmatrix},~~  Y =
\begin{bmatrix} 
1 & x^T \\ x  & X    
\end{bmatrix}   \notag
\end{equation}
constraints $\hat L$, LMI, and RLT (rearranged for variable $Y$) are defined in (\ref{eq:App-01-SDPr-QCQP-Fea-LPr}), (\ref{eq:App-01-SDPr-QCQP-Fea-LMI}), and (\ref{eq:App-01-SDPr-QCQP-VIN-RLT}), respectively.  The last constraint in (\ref{eq:App-01-SDP-Rank}) ensures that $Y$ has a rank-1 decomposition such that $X = x x^T$. Actually, LMI and RLT are redundant to the rank-1 constraint, but will give a high quality convex relaxation when the rank constraint is relaxed. 

To treat the rank-1 constraint in a soft manner, we introduce a dummy variable $Z$, and penalize the matrix rank in the objective function, giving rising to the following problem
\begin{equation}
\label{eq:App-01-SDP-Rank-Penalty}
\begin{aligned}
\min_{Y}~~ & \left\{ \langle {\rm \Omega}, Y \rangle + 
\min_{Z} \frac{\rho}{2} \| Y-Z \|^2_2 \right\}   \\
\mbox{s.t.} ~~ & Y \in \hat L \cap \mbox{LMI} \cap \mbox{RLT},~ 
\mbox{rank}(Z) = 1 
\end{aligned}
\end{equation}
If the penalty parameter $\rho$ is sufficiently large, the penalty term will be zero at the optimal solution, so $Y=Z$ and rank$(Z)=1$. One advantage of this treatment is that the constraints on $Y$ and $Z$ are decoupled, and the inner rank minimization problem has a closed-form solution. 

To see this, if rank$(Y) = k > 1$, the singular value decomposition of $Y$ has the form $Y = U {\rm \Sigma} V^T$, where
\begin{equation}
{\rm \Sigma} = \mbox{diag}(S,0),~ 
S = \mbox{diag}(\sigma_1,\cdots,\sigma_k),~
\sigma_1 \ge \cdots \ge \sigma_k > 0 \notag
\end{equation}
$U$ and $V$ are orthogonal matrices. Let matrix $D$ have the same dimension as $Y$, $D_{11}=\sigma_1$, and $D_{ij}=0$, $\forall (i,j) \ne (1,1)$, we have  
\begin{equation}
\begin{aligned}
& \quad  \min_{Z}~ \left\{ \frac{\rho}{2} \| Y - Z \|^2_2 ~\middle|~ 
\mbox{rank}(Z) = 1 \right\}   \\
& = \min_{Z}~ \left\{ \frac{\rho}{2} \| U(Y - Z)V^T \|^2_2 ~\middle|~ 
\mbox{rank}(Z) = 1 \right\}   \\
& = \min_{Z}~ \left\{ \frac{\rho}{2} \| {\rm \Sigma} - U Z V^T \|^2_2 ~\middle|~ \mbox{rank}(Z) = 1 \right\}   \\
& = \frac{\rho}{2} \| {\rm \Sigma} - D \|^2_2  
  = \frac{\rho}{2} \sum_{i=2}^k \sigma_i^2 \\
& =  \frac{\rho}{2} \| Y \|^2_2 - \frac{\rho}{2} \sigma_1^2 (Y)
\end{aligned}  \notag
\end{equation}

To represent the latter term via a convex function, let matrix $\rm \Theta$ have the same dimension as $Y$, ${\rm \Theta}_{11} = 1$, and ${\rm \Theta}_{ij}=0$, $\forall (i,j) \ne (1,1)$, we have 
\begin{equation}
\mbox{tr}(Y^T U {\rm \Theta} U^T Y) = \mbox{tr}(V {\rm \Sigma} U^T U {\rm \Theta} U^T U {\rm \Sigma} V^T) = \mbox{tr}(V {\rm \Sigma} {\rm \Theta} {\rm \Sigma} V^T) 
= \mbox{tr}({\rm \Sigma} {\rm \Theta} {\rm \Sigma}) = \sigma^2_1 (Y) \notag
\end{equation}

Define two functions $f(Y) = \langle {\rm \Omega}, Y \rangle + \frac{\rho}{2} \|Y\|^2_2$ and $g(Y) = \mbox{tr}(Y^T U {\rm \Theta} U^T Y)$. Because $\|Y\|_2$ is convex in $Y$ (Example 3.11, \cite{App-A-CVX-Book-Boyd}), so is $\|Y\|^2_2$ (composition rule, page 84, \cite{App-A-CVX-Book-Boyd}); clearly,  $f(Y)$ is a convex function in $Y$, as it is the sum of a linear function and a convex function. For the latter one, the Hessian matrix of $g(Y)$ is 
\begin{equation*}
\nabla^2_Y g(Y) = U {\rm \Theta} U^T = U {\rm \Theta}^T {\rm \Theta} U^T = ({\rm \Theta} U^T)^T {\rm \Theta} U^T \succeq 0
\end{equation*}
so $g(Y)$ is also convex in $Y$. Substituting above results into problem (\ref{eq:App-01-SDP-Rank-Penalty}), the rank constrained SDP (\ref{eq:App-01-SDP-Rank}) boils
down to
\begin{equation}
\label{eq:App-01-SDP-Rank-DCP}
\min_{Y} \left\{ \langle {\rm \Omega}, Y \rangle + \frac{\rho}{2} \|Y\|^2_2 - \frac{\rho}{2} \mbox{tr}(Y^T U {\rm \Theta} U^T Y) ~\middle|~ Y \in \hat L \cap \mbox{LMI} \cap \mbox{RLT} \right\}
\end{equation}

The objective function is a DC function, and the feasible region is convex, so (\ref{eq:App-01-SDP-Rank-DCP}) is a DC program. One can employ the convex-concave procedure discussed in \cite{App-A-CCP-Boyd} to solve this problem. The flowchart is summarized in Algorithm \ref{Ag:App-01-SDP-Rank-CCP}.

\begin{algorithm}[!htp]
\normalsize
\caption{\bf : Sequential SDP}
\begin{algorithmic}[1]
\STATE Choose an initial penalty parameter $\rho^0$, a penalty growth rate $\tau > 0$, and solve the following SDP relaxation model
\begin{equation}
\min~~ \{ \langle {\rm \Omega}, Y \rangle ~|~ 
Y \in \hat L \cap \mbox{LMI} \cap \mbox{RLT} \}  \notag
\end{equation}  
The optimal solution is $Y^*$.
\STATE Construct the linear approximation of $g(Y)$ as 
\begin{equation}
\label{eq:App-01-SDP-DCP-Concave-Lin}
g_L (Y,Y^*) = g(Y^*) + \langle \nabla g(Y^*),Y-Y^* \rangle
\end{equation}
 Solve the the following SDP
\begin{equation}
\label{eq:App-01-SDP-DCP-Master}
\min_{Y} \left\{ f(Y) - \frac{\rho}{2} g_L(Y,Y^*) ~\middle|~ Y \in \hat L \cap \mbox{LMI} \cap \mbox{RLT} \right\}
\end{equation}
The optimal solution is $Y^*$. 
\STATE If rank$(Y^*)=1$, terminate and report the optimal solution $Y^*$; otherwise, update $\rho \leftarrow (1+\tau) \rho$, and go to step 2.
\end{algorithmic}
\label{Ag:App-01-SDP-Rank-CCP}
\end{algorithm} 

For the convergence of Algorithm \ref{Ag:App-01-SDP-Rank-CCP}, we have the following properties.

\begin{proposition}
\label{pr:App-01-SDP-Rank-CCP-1}
\cite{App-A-SDP-Rank-CCP} The optimal value sequence $F(Y^i)$ generated by Algorithm \ref{Ag:App-01-SDP-Rank-CCP} is monotonically decreasing.
\end{proposition}

Denote by $F(Y) = f(Y) - \frac{\rho}{2} g(Y)$ the objective function of (\ref{eq:App-01-SDP-Rank-DCP}) in the DC form, and $H(Y,Y^i) = f(Y) - \frac{\rho}{2} g_L(Y,Y^i)$ the convexified objective function in (\ref{eq:App-01-SDP-DCP-Master}) by linearizing the concave term in $F(Y)$. Two basic facts help explain this proposition:

1) $g_L(Y^*,Y^*) = g(Y^*)$, $\forall Y^*$ which directly follows from the definition in (\ref{eq:App-01-SDP-DCP-Concave-Lin}).

2) For any given $Y^*$, $g(Y) \ge g_L(Y,Y^*)$, $\forall Y$, because the graph of a convex function must lie over its tangent plane at any fixed point.

First we can asset inequality $H(Y^{i+1},Y^i) \le H(Y^i,Y^i)$,  because $H(Y,Y^i)$ is optimized in problem (\ref{eq:App-01-SDP-DCP-Master}). The optimum $H(Y^{i+1},Y^i)$ deserves a value no greater than that at any feasible point.
Furthermore, with the definition of $H(Y^i,Y^i)$, we have 
\begin{equation}
H(Y^{i+1},Y^i) \le H(Y^i,Y^i) = f(Y^i) - \frac{\rho}{2} g_L(Y^i,Y^i) = 
f(Y^i) - \frac{\rho}{2} g(Y^i) = F(Y^i)  \notag
\end{equation}
On the other hand, 
\begin{equation}
H(Y^{i+1},Y^i)=f(Y^{i+1}) - \frac{\rho}{2} g_L(Y^{i+1},Y^i) \ge f(Y^{i+1}) - \frac{\rho}{2} g(Y^{i+1})=F(Y^{i+1})  \notag
\end{equation}
Consequently, we arrive at the monotonic property 
\begin{equation}
F(Y^{i+1}) \le F(Y^i) \notag
\end{equation}

\begin{proposition}
\label{pr:App-01-SDP-Rank-CCP-2}
\cite{App-A-SDP-Rank-CCP} The solution sequence $Y^i$ generated by Algorithm \ref{Ag:App-01-SDP-Rank-CCP} approaches to the optimal solution of problem (\ref{eq:App-01-SDP-Rank}) when $\rho \rightarrow \infty$.
\end{proposition}

It is easy to understand that whenever $\rho$ is sufficiently large, the penalty term will tend to 0, and the rank-1 constraint in (\ref{eq:App-01-SDP-Rank}) is met. A formal proof can be found in \cite{App-A-SDP-Rank-CCP}. A few more remarks are given below.

1) The convex-concave procedure in \cite{App-A-SDP-Boyd} is a local algorithm under mild conditions and needs a manually supplied initial point. Algorithm \ref{Ag:App-01-SDP-Rank-CCP}, however, is elaborately initiated at the solution offered by the SDP relaxation model, which usually appears to be close to the global optimal one for many engineering optimization problems. Therefore, Algorithm \ref{Ag:App-01-SDP-Rank-CCP} generally performs well and will identify the global optimal solution, although a provable guarantee is non-trivial.  

2) In practical applications, Algorithm \ref{Ag:App-01-SDP-Rank-CCP} could converge without the penalty parameter approaching infinity, because when some constraint quantification holds, there exists an exact penalty parameter $\rho^*$, such that the optimal solution leads to a zero penalty term for any $\rho \ge \rho^*$ \cite{App-A-Exact-Penalty-1,App-A-Exact-Penalty-2}, and Algorithm \ref{Ag:App-01-SDP-Rank-CCP} converges in a finite number of steps.
If the exact penalty parameter does not exist, Algorithm \ref{Ag:App-01-SDP-Rank-CCP} may fail to converge. In such circumstance, one can impose an upper bound on $\rho$, and use an alternative convergence criterion: the change of the objective value $F(Y)$ in two consecutive steps is less than a given threshold value. As a result, Algorithm \ref{Ag:App-01-SDP-Rank-CCP} will be able to find an approximate solution of problem (\ref{eq:App-01-SDP-Rank}), and the rank-1 constraint may not be enforced.

3) From the numeric computation perspective, a very large $\rho$ may cause ill-conditioned problem and lead to numerical instability, so it is useful to gradually increase $\rho$ from a small value. Another reason for the moderate growth of $\rho$ is that it does not cause dramatic change of optimal solutions in two successive iterations. As a result, $g_L(Y,Y^*)$ can provide relatively accurate approximation for $g(Y)$ in every iteration. 
 
4) The penalty term $\rho_i p(Y^i)/2 = \rho_i \left[ \|Y^i\|^2_2 - \mbox{tr}(Y^{iT} U {\rm \Theta} U^T Y^i) \right]/2$ gives an upper bound on the optimality gap induced by rank relaxation. To see this, let $\rho^*$ and $Y^*$ be the exact penalty parameter and corresponding optimal solution of (\ref{eq:App-01-SDP-DCP-Master}), i.e., $p(Y^*)=0$; $\rho_i$ and $Y^i$ be the penalty parameter and optimal solution in $i$-th iteration. According to Proposition \ref{pr:App-01-SDP-Rank-CCP-1}, we have $\langle {\rm \Omega}, Y^* \rangle \le \langle {\rm \Omega}, Y^i \rangle + \rho_i p(Y^i)/2$; moreover, since the rank-1 constraint is relaxed before Algorithm \ref{Ag:App-01-SDP-Rank-CCP} could converge, $\langle {\rm \Omega}, Y^i \rangle \le \langle {\rm \Omega}, Y^* \rangle$ holds. Therefore, $\langle {\rm \Omega}, Y^i \rangle$ and $\langle {\rm \Omega}, Y^i \rangle + \rho_i p(Y^i)/2$ are lower and upper bounds for the optimal value of problem (\ref{eq:App-01-SDP-Rank}). In this regard, $\rho_i p(Y^i)/2$ is an estimation on the optimality gap.

\subsection{Completely Positive Program Relaxation}
\label{App-A-Sect03-03}

Inspired by the convex hull expression in (\ref{eq:App-01-COMPL-Cone}), researchers have shown that most non-convex QCQPs can be modeled as linear programs over the intersection of a completely positive cone and a polyhedron \cite{App-A-COPr-1,App-A-COPr-2,App-A-COPr-3}. For example, consider minimizing a quadratic function over a standard simplex
\begin{equation}
\label{eq:App-01-COPr-Example-1}
\begin{aligned}
\min~~ &  x^T Q x  \\
\mbox{s.t.} ~~ &  e^T x = 1  \\
       &  x \ge 0
\end{aligned}
\end{equation}
where $Q \in \mathbb S^n$, and $e$ denotes the all-one vector with $n$ entries. Following the paradigm similar to (\ref{eq:App-01-SDPr-QCQP-Ext}), let $X = x x^T$, and then we can construct a valid inequality
\begin{equation*}
1 = x^T e e^T x = x^T E x = \langle E, X \rangle  
\end{equation*}
where $E=ee^T$ is the all-one matrix. According to (\ref{eq:App-01-COMPL-Cone}), conv$\{xx^T| x \in \mathbb R^n_+ \}$ is given by $(\mathbb C^n_+)^*$. Therefore, problem (\ref{eq:App-01-COPr-Example-1}) transforms to
\begin{equation}
\label{eq:App-01-COPr-Example-2}
\begin{aligned}
\min~~ &  \langle Q, X \rangle  \\
\mbox{s.t.} ~~ &  \langle E, X \rangle = 1 \\
       &  X \in  \mathbb (\mathbb C^n_+)^*
\end{aligned}
\end{equation}
Problem (\ref{eq:App-01-COPr-Example-2}) is a convex relaxation of (\ref{eq:App-01-COPr-Example-1}). Because the objective is linear, the optimal solution must be located at one extremal point of the convex hull of the feasible region. In view of the representation in (\ref{eq:App-01-COMPL-Cone}), the extremal points are exactly rank-1, so the convex relaxation (\ref{eq:App-01-COPr-Example-2}) is always exact.  

Much more general results are demonstrated in \cite{App-A-COPr-2} that every quadratic program with linear and binary constraints can be rewritten as a completely positive program. More precisely, a mixed-integer quadratic program 
\begin{equation}
\label{eq:App-01-COPr-Example-3}
\begin{aligned}
\min~~ &  x^T Q x + 2 c^T x \\
\mbox{s.t.} ~~ &  a^T_i x = b_i,~ i = 1,\cdots,m  \\
       &  x \ge 0,~ x_j \in \{0,1\}, j \in B
\end{aligned}
\end{equation}
and the following completely positive program
\begin{equation}
\label{eq:App-01-COPr-Example-4}
\begin{aligned}
\min~~ &  \langle Q, X \rangle + 2 c^T x \\
\mbox{s.t.} ~~ &  a^T_i x = b_i,~ i = 1,\cdots,m \\
       &  \langle a_i a^T_i, X \rangle = b^2_i,~ i = 1,\cdots,m \\ 
       &  x_j = X_{jj},~ j \in B  \\ 
       &  X \in  \mathbb (C^n_+)^*
\end{aligned}
\end{equation}
have the same optimal solution, as long as problem (\ref{eq:App-01-COPr-Example-3}) satisfies: $a^T_i x = b_i$, $\forall i$ and $x \ge 0$ implies $x_j \le 1$, $\forall j \in B$. Actually, this is a relatively mild condition \cite{App-A-COPr-2}. Complementarity constraints can be handled in the similar way. Whether problems with general quadratic constraints can be restated as completely positive programs in the similar way remains an open question.  

The NP-hardness of problem (\ref{eq:App-01-COPr-Example-3}) makes (\ref{eq:App-01-COPr-Example-4}) NP-hard itself. The complexity has been encapsulated into the last cone constraint. The relaxation model is still interesting due to its convexity. Furthermore, it can be approximated via a sequence of SDPs with growing sizes \cite{App-A-COPr-SOS} given an arbitrarily small error bound.

\subsection{MILP Approximation}
\label{App-A-Sect03-04}

SDP relaxation technique introduces a squared matrix variable that contains $n(n+1)/2$ independent variables. Although exploiting the sparse pattern of $X$ via graphic theory is helpful to expedite problem solution, the computational burden is still high especially when the initial relaxation is inexact and a sequence of SDPs should be solved. Inspired by difference-of-convex programming an alternative choice is to express the non-convexity of QCQP by univariate concave functions, and approximate these concave functions via PWL functions compatible with mixed-integer programming solvers. This approach has been expounded in \cite{App-A-QCQP-MILP}.

Consider nonconvex QCQP
\begin{equation}
\label{eq:App-01-QCQP-MILP-Appr-1}
\begin{aligned}
\min~~ & x^T A_0 x + a^T_0 x   \\
\mbox{s.t.} ~~ & x^T A_k x + a^T_k x \le b_k,~ k = 1, \cdots, m
\end{aligned}
\end{equation}
We can always find $\delta_0$, $\delta_1$, $\cdots$, $\delta_m \in \mathbb R^+$, such that $A_k + \delta_k I \succeq 0$, $k=0,\cdots,m$. For example, $\delta_k$ can take the absolute value of the most negative eigenvalue of $A_k$, and $\delta_k=0$ if $A_k \succeq 0$. Then, problem (\ref{eq:App-01-QCQP-MILP-Appr-1}) can be cast as  
\begin{equation}
\label{eq:App-01-QCQP-MILP-Appr-2}
\begin{aligned}
\min~~ & x^T (A_0+\delta_0I) x + a^T_0 x -\delta_0 1^T y  \\
\mbox{s.t.} ~~ & x^T (A_k + \delta_k I) x + a^T_k x - \delta_k 1^T y \le b_k,~ k = 1, \cdots, m \\
& y_i = x^2_i, i=1, \cdots, n
\end{aligned}
\end{equation}
Problem (\ref{eq:App-01-QCQP-MILP-Appr-2}) is actually a difference-of-convex program; however, the nonconvex terms are consolidated in much simpler parabolic equalities, which can be linearized via the SOS2 based PWL approximation technique discussed in Appendix \ref{App-B-Sect01}. Except for the last $n$ quadratic equalities, remaining constraints and objective function of problem (\ref{eq:App-01-QCQP-MILP-Appr-2}) are all convex, so the linearized problem gives rise to a mixed-integer convex quadratic program.  

Alternatively, we can first perform convex relaxation by replacing $y_i = x^2_i$ with $y_i \ge x^2_i$, $i=1, \cdots, n$; if strict inequality holds  at the optimal solution, a disjunctive cut is generated to remove this point from the feasible region. However, the initial convex relaxation can be very weak ($y=+\infty$ is usually an optimal solution). Predefined disjunctive cuts can be added \cite{App-A-QCQP-MILP}. 

Finally, nonconvex QCQP is a hard optimization problem. Developing an efficient algorithm should leverage specific problem structure. For example, SDP relaxation is suitable for OPF problems; MILP approximation can be used for small and dense problems. Unlike SDP relaxation works on a squared matrix variable, the number of auxiliary variables in (\ref{eq:App-01-QCQP-MILP-Appr-2}) and its mixed-integer convex quadratic program approximation is moderate. Therefore, this approach is promising to tackle practical problems whose coefficient matrices are usually sparse. Furthermore, no particular assumption is needed to guarantee the exactness of relaxation, so this method is general enough to tackle a wide spectrum of engineering optimization problems.

\section{MILP Formulation of Nonconvex QPs}
\label{App-A-Sect04}

In a non-convex QCQP, if the constraints are all linear, it is called a nonconvex QP. There is no doubt that convex relaxation methods presented in the previous section can be applied to nonconvex QPs. However, the relaxation is generally inexact. In this section, we introduce exact MILP formulations to globally solve such a nonconvex  optimization problem; unlike the mixed-integer programming approximation method in Sect. \ref{App-A-Sect03-04}, in which approximation error is inevitable, by using duality theory, the MILP models will be completely equivalent to the original QP. Thanks to the advent of powerful MILP solvers, this method is becoming increasingly competitive compared to existing global solution methods and is attracting more attentions from the research community.

\subsection{Nonconvex QPs over polyhedra}

The presented approach is devised in \cite{App-A-QP-MILP}. A nonconvex QP with linear constraints has the form of 
\begin{equation}
\label{eq:App-01-NC-QP}
\begin{aligned}
\min~~ & \frac{1}{2} x^T Q x  + c^T x  \\
\mbox{s.t.} ~~ & Ax \le b
\end{aligned}
\end{equation}
where $Q$ is a symmetric, but indefinite matrix; $A$, $b$, $c$ are constant coefficients with compatible dimensions. We assume that finite lower and upper limits of the decision variable $x$ have been included, and thus the feasible region is a bounded polyhedron. The KKT conditions of (\ref{eq:App-01-NC-QP}) can be written as:
\begin{equation}
\label{eq:App-01-NC-QP-KKT}
\begin{gathered}
c + Qx + A^T \xi = 0 \\
0 \le \xi \bot b - Ax \ge 0
\end{gathered}
\end{equation}
If there is a multiplier $\xi$ so that the pair $(x,\xi)$ of primal and dual variables satisfies KKT condition (\ref{eq:App-01-NC-QP-KKT}), then $x$ is said to be a KKT point or a stationary point. The complementarity and slackness condition in  (\ref{eq:App-01-NC-QP-KKT}) gives $b^T \xi = x^T A^ T \xi$. For any primal-dual pair $(x,\xi)$ that satisfies (\ref{eq:App-01-NC-QP-KKT}), the following relations hold
\begin{equation}
\label{eq:App-01-NC-QP-Lin}
\begin{aligned}
\frac{1}{2} x^T Qx + c^T x &= \frac{1}{2} c^T x +\frac{1}{2} x^T (c + Qx)\\
& = \frac{1}{2} c^T x - \frac{1}{2}  x^T A^ T \xi = \frac{1}{2} \left( c^T x - b^T \xi \right)
\end{aligned}
\end{equation}

As such, the non-convex quadratic objective function is equivalently stated as a linear function in the primal and dual variables without loss of accuracy. 
Thus, if problem (\ref{eq:App-01-NC-QP}) has an optimal solution, then the solution can be retrieved by solving an LPCC
\begin{equation}
\label{eq:App-01-NC-QP-LPCC}
\begin{aligned}
\min~~ & \frac{1}{2} \left( c^T x - b^T \xi \right) \\
\mbox{s.t.} ~~ & c + Qx + A^T \xi = 0 \\
&0 \le \xi \bot b - Ax \ge 0
\end{aligned}
\end{equation}
which is equivalent to the following MILP
\begin{equation}
\label{eq:App-01-NC-QP-MILP}
\begin{aligned}
\min~~ & \frac{1}{2} c^T x - b^T \xi  \\
\mbox{s.t.} ~~ & c + Qx + A^T \xi = 0 \\
&  0 \le \xi \le M(1-z)   \\
&  0 \le b - Ax \le Mz   \\
&  z ~\mbox{ binary}
\end{aligned}
\end{equation}
where $M$ is a sufficiently large  constant; $z$ is a vector of binary variables. Regardless of the value of $z_i$, at most one of $\xi_i$ and $(b-Ax)_i$ can take a strictly positive value. For more rigorous discussions on this method, please see \cite{App-A-QP-MILP}, in which an unbounded feasible region is considered. More tricks in MILP reformulation technique can be found in the next chapter.

It should be pointed out that the set of optimal solutions of (\ref{eq:App-01-NC-QP}) is a subset of stationary points described by (\ref{eq:App-01-NC-QP-KKT}), because (\ref{eq:App-01-NC-QP-KKT}) is only a necessary condition for optimality but not sufficient. Nevertheless, as we assumed that the feasible region is a bounded polytope (thus compact), QP (\ref{eq:App-01-NC-QP}) must have a finite optimum, then according to \cite{App-A-QP-LPCC-Opt-Eqv}, the optimal value is equal to the minimum of objective function values perceived at stationary points. Therefore, MILP (\ref{eq:App-01-NC-QP-MILP}) provides an exact solution to (\ref{eq:App-01-NC-QP}).

Finally, we shed some light on the selection of $M$, since it has notable impact on the computational efficiency of (\ref{eq:App-01-NC-QP-MILP}). An LP based bound preprocessing method is thoroughly discussed in \cite{App-A-QP-LPCC-Bounding}, which is used in a finite branch-and-bound method for solving LPCC (\ref{eq:App-01-NC-QP-LPCC}). Here we briefly introduce the bounding method.

For the primal variable $x$ which represents physical quantities or measures, its bounds depends on practical situations and security considerations, and we assume that the bound is $0 \le x \le U$. The bound can be tightened by solving
\begin{equation}
\label{eq:App-01-NC-QP-LPCC-Primal-Bounding}
\min (\max)~ \{x_j ~|~ Ax \le b,~ 0 \le x \le U \}
\end{equation}
In (\ref{eq:App-01-NC-QP-LPCC-Primal-Bounding}), we can incorporate individual bounds for the components of vector $x$, which never wrecks the optimal solution and can be supplemented  in (\ref{eq:App-01-NC-QP-MILP}).    

For the dual variables, we consider (\ref{eq:App-01-NC-QP}) again with explicit bounds on primal variable $x$
\begin{equation*}
\begin{aligned}
\min~~ & \frac{1}{2} x^T Q x  + c^T x  \\
\mbox{s.t.} ~~ & Ax \le b: \xi \\
& 0 \le x \le U: \lambda,\rho
\end{aligned}
\end{equation*}
where $\xi$, $\lambda$, $\rho$ following the colon are dual variables. Its KKT condition reads 
\begin{subequations}
\label{eq:App-01-NC-QP-KKT-Primal-Bound}
\begin{gather}
c + Qx + A^T \xi - \lambda + \rho = 0 \label{eq:App-01-NC-QP-KKT-Primal-Bound-1}\\
0 \le \xi \bot b - Ax \ge 0 \label{eq:App-01-NC-QP-KKT-Primal-Bound-2}\\
0 \le x \bot \lambda \ge 0 \label{eq:App-01-NC-QP-KKT-Primal-Bound-3}  \\
0 \le U - x \bot \rho \ge 0 \label{eq:App-01-NC-QP-KKT-Primal-Bound-4} 
\end{gather}
Multiplying both sides of (\ref{eq:App-01-NC-QP-KKT-Primal-Bound-1}) by a feasible solution $x^T$
\begin{equation}
\label{eq:App-01-NC-QP-KKT-Primal-Bound-5}
c^T x + x^T Q x + x^T A^T \xi - x^T \lambda + x^T \rho = 0
\end{equation}
Substituting $\xi^T A x = \xi^T b$, $x^T \lambda = 0$, and $x^T \rho = \rho^T U$ concluded from (\ref{eq:App-01-NC-QP-KKT-Primal-Bound-2})-(\ref{eq:App-01-NC-QP-KKT-Primal-Bound-4}) into (\ref{eq:App-01-NC-QP-KKT-Primal-Bound-5}) outcomes
\begin{equation}
\label{eq:App-01-NC-QP-KKT-Primal-Bound-6}
c^T x + x^T Q x + b^T \xi + U^T \rho = 0
\end{equation}
\end{subequations}

The upper bounds (lower bounds are 0) on the dual variables required for MILP (\ref{eq:App-01-NC-QP-MILP})  can be computed from the following LP:
\begin{subequations}
\label{eq:App-01-NC-QP-KKT-Dual-Bounding}
\begin{align}
\max~~ & \lambda_j \label{eq:App-01-NC-QP-KKT-Dual-Bounding-Obj} \\
\mbox{s.t.}~~ & c + Qx + A^T \xi - \lambda + \rho = 0 
\label{eq:App-01-NC-QP-KKT-Dual-Bounding-Cons-1} \\
& \mbox{tr}(Q^T X) + c^T x + b^T \xi + U^T \rho = 0
\label{eq:App-01-NC-QP-KKT-Dual-Bounding-Cons-2} \\
& \mbox{Cons-RLT}=\{ (x,X)~|~(\ref{eq:App-01-SDPr-QCQP-VIN-RLT})\}
\label{eq:App-01-NC-QP-KKT-Dual-Bounding-Cons-3} \\
& 0 \le x \le U,~ Ax \le b, \lambda, \xi, \rho \ge 0
\label{eq:App-01-NC-QP-KKT-Dual-Bounding-Cons-4}  
\end{align}
\end{subequations}
In (\ref{eq:App-01-NC-QP-KKT-Dual-Bounding-Cons-2}), quadratic equality (\ref{eq:App-01-NC-QP-KKT-Primal-Bound-6}) is linearized by letting $X = xx^T$, and (\ref{eq:App-01-NC-QP-KKT-Dual-Bounding-Cons-3}) is a linear relaxation for above rank-1 condition, as explained in Sect. \ref{App-A-Sect03-01}. 

By exploiting the relaxation revealed in (\ref{eq:App-01-NC-QP-KKT-Dual-Bounding-Cons-2}), it bas been proved that problem (\ref{eq:App-01-NC-QP-KKT-Dual-Bounding}) always has a finite optimum, because the recession cone of the set comprised of the primal and dual variables as well as their associated valid inequalities is empty, see the proof of Proposition 3.1 in \cite{App-A-QP-LPCC-Bounding}. This is a pivotal theoretical guarantee. Other bounding techniques which only utilize KKT conditions hardly ensure a finite optimum.

\subsection{Standard Nonconvex QPs}

The presented approach is devised in \cite{App-A-Standard-QP}. A standard nonconvex QP entails minimizing a nonconvex quadratic function over a unit probability simplex  
\begin{equation}
\label{eq:App-01-Stand-QP}
\begin{aligned}
v(Q) = \min~ &  x^T Q x  \\
\mbox{s.t.} ~ & x \in {\rm \Delta}_n  
\end{aligned}
\end{equation}
where $Q$ is a symmetric matrices, and unit simplex
\begin{equation*}
{\rm \Delta}_n = \{x \in \mathbb R^n_+~|~e^T x = 1 \}
\end{equation*}
where $e$ is all-one vector. A nonhomogeneous objective can always be transformed to a quadratic form given the simplex constraint ${\rm \Delta}_n$:
\begin{equation*}
x^T Q x + 2c^T x = x^T (Q + e c^T + c e^T ) x,~ \forall x \in {\rm \Delta}_n
\end{equation*}

Standard nonconvex QPs have wide applications in portfolio optimization, quadratic resource allocation, graphic theory and so on. In addition, for a
given symmetric matrix $Q$, a necessary and sufficient condition for $Q$ being copositive is $v(Q) \ge 0$. Copositive programming is a young and active research field, and can help the research in convex relaxation. A fundamental problem is copositivity test, which entails solving (\ref{eq:App-01-Stand-QP}) globally. 

Problem (\ref{eq:App-01-Stand-QP}) is a special case of nonconvex QP (\ref{eq:App-01-NC-QP}), so the methods in previous subsection also work for (\ref{eq:App-01-Stand-QP}). The core trick is to select a big-M parameter in linearizing complementarity and slackness conditions. Due to its specific structure, the valid big-M parameter for problem (\ref{eq:App-01-Stand-QP}) can be chosen in a much more convenient way. To see this, the KKT condition of (\ref{eq:App-01-Stand-QP}) reads as
\begin{subequations}
\label{eq:App-01-Stand-QP-KKT}
\begin{align}
Qx - \lambda e - \mu & = 0 \label{eq:App-01-Stand-QP-KKT-1} \\
e^T x & = 1 \label{eq:App-01-Stand-QP-KKT-2}  \\
x & \ge 0 \label{eq:App-01-Stand-QP-KKT-3}  \\
\mu & \ge 0 \label{eq:App-01-Stand-QP-KKT-4}  \\
x_j \mu_j & = 0,~ j=1,\cdots,n \label{eq:App-01-Stand-QP-KKT-5}
\end{align}
\end{subequations}
where $\lambda$ and $\mu$ are dual variables associated with equality constraint $e^T x = 1$ and inequality constraint $x \ge 0$. Because the feasible region is polyhedral, constraint quantification always holds, and any optimal solution of (\ref{eq:App-01-Stand-QP}) must solve KKT system (\ref{eq:App-01-Stand-QP-KKT}).

Multiplying both sides of (\ref{eq:App-01-Stand-QP-KKT-1}) by $x$ results in $x^T Q x = \lambda x^T e - x^T \mu$; substituting (\ref{eq:App-01-Stand-QP-KKT-2}) and (\ref{eq:App-01-Stand-QP-KKT-5}) into the right-hand side concludes $x^T Q x = \lambda$. Provided with eligible big-M parameter, problem (\ref{eq:App-01-Stand-QP}) is (exactly) equivalent to the following MILP
\begin{equation}
\label{eq:App-01-Stand-QP-MILP-1}
\begin{aligned}
\min~~ & \lambda  \\  
\mbox{s.t.} ~~ & Qx - \lambda e - \mu = 0  \\ 
 & e^T x = 1,~ 0  \le  x  \le y \\
 &0 \le \mu_j \le M_j(1-y_j),~ j = 1, \cdots, n 
\end{aligned}
\end{equation}
where $y \in \{0,1\}^n$, and $M_j$ is the big-M parameter. It is the upper bound of dual variable $\mu_j$. To estimate such a bound, according to (\ref{eq:App-01-Stand-QP-KKT-1})
\begin{equation*}
\mu_j = e^{T}_j Q x - \lambda,~ j=1,\cdots,n 
\end{equation*}
where $e_j$ is the $j$-th column of $n \times n$ identity matrix. For the first term, 
\begin{equation*}
x^T Q e_j \le \max_{i \in \{1,\cdots,n\}} Q_{ij},~ j=1,\cdots,n 
\end{equation*}
As for the second term, we know $\lambda \ge v(Q)$, so any known lower bound of $v(Q)$ can be used to obtain an upper bound of $M_j$. One possible lower bound of $v(Q)$ is suggested in  \cite{App-A-Standard-QP} as 
\begin{equation*}
l(Q) = \min_{1 \le i,j \le n} Q_{ij} + \dfrac{1}{\sum_{k=1}^n \left(Q_{kk}- \min\limits_{1 \le i,j \le n} Q_{ij} \right)^{-1}} 
\end{equation*}
If the minimal element of $Q$ locates on the main diagonal, the second term vanishes and $l(Q) = \min_{1 \le i,j \le n} Q_{ij}$.

In summary, a valid choice of $M_j$ would be 
\begin{equation}
\label{eq:App-01-Stand-QP-Big-M}
M_j = \max_{i \in \{1,\cdots,n\}} Q_{ij} - l(Q),~ j = 1,\cdots,n
\end{equation}

It is found in \cite{App-A-Standard-QP} that if we relax (\ref{eq:App-01-Stand-QP-KKT-1}) as an inequality and solve the following MILP
\begin{equation}
\label{eq:App-01-Stand-QP-MILP-2}
\begin{aligned}
\min~~ & \lambda  \\  
\mbox{s.t.} ~~ & Qx - \lambda e - \mu \le 0  \\ 
 & e^T x = 1,~ 0  \le  x  \le y \\
 &0 \le \mu_j \le M_j(1-y_j),~ j = 1, \cdots, n 
\end{aligned}
\end{equation}
which is an relaxed version of (\ref{eq:App-01-Stand-QP-MILP-1}), the optimal solution will not change. However, in some instances, solving (\ref{eq:App-01-Stand-QP-MILP-2}) is significantly faster than solving (\ref{eq:App-01-Stand-QP-MILP-1}). More thorough theoretical analysis can be found in \cite{App-A-Standard-QP}.

\section{Further Reading}
\label{App-A-Sect05}

Decades of wonderful research has resulted in elegant theoretical developments and sophisticated computational softwares, which have brought convex optimization to an unprecedented dominating stage where it serves as the baseline and reference model for optimization problems in almost every discipline. Only problems which can be formulated as convex programs are regarded as theoretically solvable. We suggest following materials for readers who want to build a solid mathematical background or know more about applications in the field of convex optimization.   

1. Convex analysis and convex optimization. Convex analysis is a classic topic in mathematics, and focuses on basic concepts and topological properties of convex sets and convex functions. We recommend monographs \cite{App-A-Convex-Analysis-1,App-A-Convex-Analysis-2,App-A-Convex-Analysis-3}. The last one sheds more light on optimization related topics, including DC programming, polynomial programming, and equilibrium constrained programming, which are originally non-convex. The most popular textbooks on convex optimization include \cite{App-A-CVX-Book-Ben,App-A-CVX-Book-Boyd}. They contain important materials that everyone who wants to apply this technique should know.

2. Special convex optimization problems. The most mature convex optimization problems are LPs, SOCPs, and SDPs. We recommend \cite{App-A-LP-Book-Dantzig,App-A-LP-Book-Bertsimas,App-A-LP-Book-Vanderbei} for the basic knowledge of duality theory, simplex algorithm, interior-point algorithm, and applications of LPs. The modeling abilities of SOCPs and SDPs have been well discussed in \cite{App-A-CVX-Book-Ben,App-A-CVX-Book-Boyd}. A geometric program is a type of optimization problem whose objective and constraints are characterized by special monomials and posynomial functions. Through a logarithmic variable transformation, a geometric program can be mechanically converted to a convex optimization problem. Geometric programming is relatively restrictive in structure, and it may not be apparent to see whether a given problem can be expressed by a geometric program. We recommend a tutorial paper \cite{App-A-GOP-Boyd} and references therein on this topic. Copositive program is a relatively young field in operational research. It is a special class of conic programming which is more general than SDP. Basic information on copositive/completely positive programs is introduced in \cite{App-A-Copositive-1,App-A-Copositive-2,App-A-Copositive-3}. They are particularly useful in combinatorial and quadratic optimization. Though very similar to SDPs in appearances, copositive programs are NP-hard. Algorithms and applications of copositive and completely positive programs have continued to be highly active research fields \cite{App-A-COP-New-1,App-A-COP-New-2,App-A-COP-New-3}.

3. General convex optimization problems. Beside above mature convex optimization models that can be specified without high level of expertise, recognizing the convexity of a general mathematical programming problem may be rather tricky. A deep understanding on convex analysis is unavoidable. Furthermore, to solve the problem using off-the-shelf solvers, a user must find a way to transform the problem into one of the standard forms (if a general purpose NLP solver fails to solve it). The so-called disciplined convex programming method is proposed in \cite{App-A-Disp-CVX} to lower this expertise barrier. The method consists of a set of rules and conventions that one must follow when setting up the problem such that the convexity is naturally sustained. This methodology has been implemented in cvx toolbox under Matlab environment.

4. Convex relaxation methods. One major application of convex optimization is to derive tractable approximations for non-convex programs, so as to facilitate problem resolution in terms of computational efficiency and robustness. A general QCQP is a quintessential non-convex optimization problem. Among various convex relaxation approaches, the SDP relaxation is shown to be able to offer high quality solutions for many QCQPs raised in signal process \cite{App-A-SDPr-Signal-1,App-A-SDPr-Signal-2} and power system energy management \cite{App-A-SDPr-Power-1,App-A-SDPr-Power-2}. Decades of excellent studies on SDP relaxation methods for QCQPs are comprehensively reviewed in \cite{App-A-SDP-Relaxation-Tutor,App-A-SDPr-QCQP-Rev-1,App-A-SDPr-QCQP-Rev-2}. Some recent advances are reported in \cite{App-A-SDPr-QCQP-1,App-A-SDPr-QCQP-2,App-A-SDPr-QCQP-3,App-A-SDPr-QCQP-4,App-A-SDPr-QCQP-5,App-A-SDPr-QCQP-6,App-A-SDPr-QCQP-7}. The rank of the matrix variable has a decisive impact on the exactness (or tightness) of the SDP relaxation. Low rank SDP method are attracting increasing attentions from researchers, and many approaches are proposed to recover a low-rank solution. More information can be found in \cite{App-A-SDP-Rank-1,App-A-SDP-Rank-2,App-A-SDP-Rank-3,App-A-SDP-Rank-4,App-A-SDP-Rank-5} and references therein.

5. Sum-of-squares (SOS) programming is originally devised in \cite{App-A-SOS-1} to decompose a polynomial $f(x)$ as the square of another polynomial $g(x)$ (if there exists), such that $f(x)=[g(x)]^2$ must be non-negative. Non-negativity of a polynomial over a semi-algebraic set can be certified in a similar way via Positivstellensatz refutations. This can be done by solving a structured SDP \cite{App-A-SOS-1}, and implemented in a Matlab based toolbox \cite{App-A-SOS-2}. Based on these outcomes, a promising methodology is quickly developed for polynomial programs, which cover a broader class of optimization problems than QCQPs. It is proved that the global solution of a polynomial program can be found by solving a hierarchy of SDPs under mild conditions. This is very inspiring since polynomial programs are generally non-convex while SDPs are convex. We recommend \cite{App-A-Poly-SDP-1} for a very detailed discussion on this approach, and \cite{App-A-Poly-SDP-2,App-A-Poly-SDP-3,App-A-Poly-SDP-4,App-A-Poly-SDP-5} for some recent advances. However, users should be aware that this approach may be unpractical because the size of the relaxed SDP quickly becomes unacceptable after a few steps. Nonetheless, the elegant theory still marks a milestone in the research field.


%
%
%

\motto{There is no problem in all mathematics that cannot be solved by direct counting. But with the present implements of mathematics many operations can be performed in a few minutes which without mathematical methods would take a lifetime. \\  \rightline{ $-$Ernst Mach}}

\chapter{Formulation Recipes in Integer Programming}
\label{App-B}

As stated in Appendix \ref{App-A}, generally speaking, convex optimization problems can be solved efficiently. However, the majority of optimization problems encountered in practical engineering are non-convex, and gradient based NLP solvers terminate at a local optimum, which may be far away from the global one. In fact, any nonlinear function can be approximated by a PWL function with adjustable errors by controlling the granularity of partitions. A PWL function can be expressed via a logic form or incorporating integer variables. Thanks to the latest progress in branch-and-cut algorithms and the development of state-of-the-art MILP solvers, a large-scale MILP can often be solved globally within reasonable computational efforts \cite{App-MILP-Solver-Perform}, although the MILP itself is proved to be NP-hard. In view of this fact, PWL/MILP approximation serves as a viable option to tackle real-world non-convex optimization problems, especially those with special structures.  

This chapter introduces PWL approximation methods for nonlinear functions and linear representations of special non-convex constraints via integer programming techniques. When the majority of a problem at hand is linear or convex, while non-convexity arises from nonlinear functions with only one or two variables, linear complementarity constraints, logical inferences and so on, it is worth trying the methods in this chapter, in view of the fact that MILP solvers are becoming increasingly efficient to retrieve a solution with a pre-defined optimality gap.

\section{Piecewise Linear Approximation of Nonlinear Functions}
\label{App-B-Sect01}

\subsection{Univariate Continuous Function}
\label{App-B-Sect01-01}

Considering a nonlinear continuous  function $f(x)$ in a single variable $x$, we can evaluate the function values $f(x_0)$, $f(x_1)$, $\cdots$, $f(x_n)$ at given breakpoints $x_0$, $x_1$, $\cdots$, $x_k$, and replace $f(x)$ with the following PWL function 
\begin{equation}
\label{eq:App-02-PWL-Logic}
f(x) = \begin{cases}
m_1 x + c_1,   &  x \in [x_0,x_1]  \\ 
m_2 x + c_2,   &  x \in [x_1,x_2]  \\ 
\qquad \vdots  &   \qquad  \vdots  \\
m_k x + c_k,   &  x \in [x_{k-1},x_k]   
\end{cases}  
\end{equation}

\begin{figure}[!t]
\centering
\includegraphics[scale=0.45]{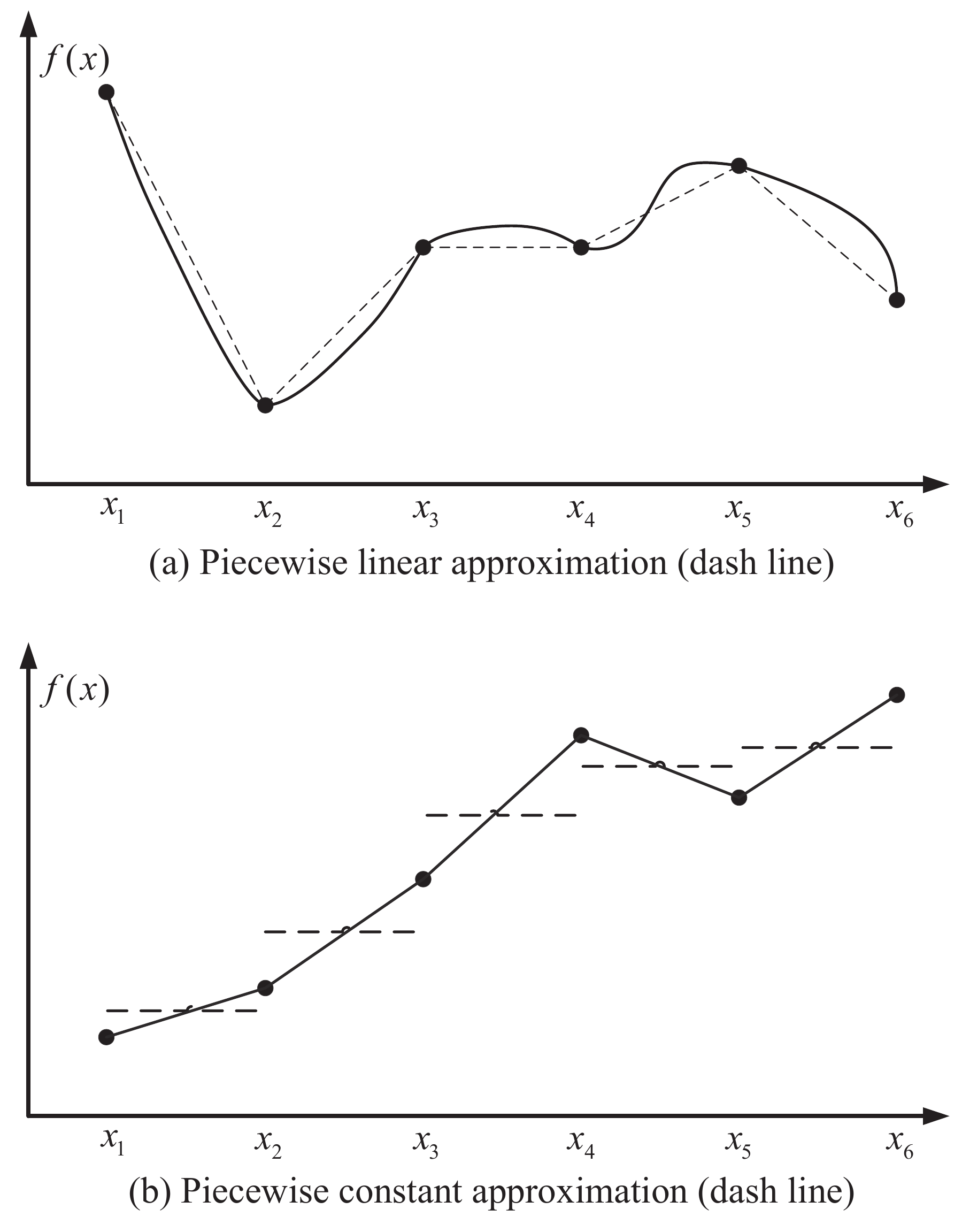}
\caption{Piecewise linear and piecewise constant approximations.}
\label{fig:App-02-01}
\end{figure}

As an illustrative example, two curves of the original nonlinear function and its PWL approximation are portrayed in part (a), Fig. \ref{fig:App-02-01}. The PWL function in (\ref{eq:App-02-PWL-Logic}) is a finite union of line segments, but still non-convex. Moreover, the logic representation in (\ref{eq:App-02-PWL-Logic}) is not compatible with commercial solvers. Given the fact that any point on a line segment can be expressed as a convex combination of two terminal points, (\ref{eq:App-02-PWL-Logic}) can be written as
\begin{equation}
\label{eq:App-02-PWL-CC}
\begin{aligned}
x & = \sum_i  \lambda_i x_i      \\
y & = \sum_i  \lambda_i f(x_i)   \\
\lambda & \ge 0,~ \sum_i \lambda_i =1 \\
\lambda & \in \mathbb{SOS}_2
\end{aligned}
\end{equation}
where $\mathbb{SOS}_2$ stands for the special ordered set of type 2, describing a vector of variables with at most two adjacent ones being able to take nonzero values. The $\mathbb{SOS}_2$ constraint on $\lambda$ can be declared via the build-in module of commercial solvers such as CPLEX or GUROBI. Please note that if $f(x)$ is convex and to be minimized, then the last $\mathbb{SOS}_2$ requirement is naturally met (thus can be relaxed), because the epigraph of $f(x)$ is a convex region. Otherwise, relaxing the last $\mathbb{SOS}_2$ constraint in (\ref{eq:App-02-PWL-CC}) gives rise to the convex hull of the sampled points $(x_0,y_0)$, $\cdots$, $(x_k,y_k)$. In general, the relaxation is inexact.

Branch-and-bound algorithms which directly working on SOS variables exhibit good performance \cite{App-MILP-SOS}, but it is desirable to explore equivalent MILP formulations to leverage the superiority of state-of-the-art solvers.
To this end, we first provide an explicit form using additional integer variables.
\begin{equation}
\label{eq:App-02-SOS2-MILP}
\begin{aligned}
\lambda_0 & \le z_1  \\
\lambda_1 & \le z_1 + z_2  \\
\lambda_2 & \le z_2 + z_3  \\
\vdots    &  \\
\lambda_{k-1} & \le z_{k-1} + z_k  \\   
\lambda_k & \le z_k    \\
 z_i  & \in \{ 0,1 \},~ \forall i,~ \sum\nolimits_{i=1}^k z_i = 1  \\
\lambda_i & \ge 0,~ \forall i,~ \sum\nolimits_{i=0}^k \lambda_i = 1   
\end{aligned}
\end{equation}
Formulation (\ref{eq:App-02-SOS2-MILP}) illustrates how integer variables can be used to enforce $\mathbb{SOS}_2$ requirements on the weighting coefficients. This formulation does not involve any manually supplied  parameter, and often gives stronger bounds when the integrality of binary variables are relaxed.

Sometimes, it is more convenient to use a piecewise constant approximation, especially when the original function $f(x)$ is not continuous. An example is exhibited in part (b), Fig. \ref{fig:App-02-01}. In this approach, the feasible interval of $x$ is partitioned into $S-1$ segments (associated with binary variables $\theta_s$, $s=1$, $\cdots$, $S-1$) by $S$ breakpoints $x_1$, $\cdots$, $x_S$ (associated with $S$ continuous weight variables $\lambda_s$, $s=1$, $\cdots$, $S$); In the $s$-th interval between $x_s$ and $x_{s+1}$, the function value $f(x)$ is approximated by the arithmetic mean $f_s = 0.5[f(x_s)+f(x_{s+1})]$, $s=1,\cdots,$ $S-1$, which is a constant as illustrated in Fig. \ref{fig:App-02-01}. With an appropriate number of partitions, an arbitrary function $f(x)$ can be approximated by a piecewise constant function as follows
\begin{subequations}
\label{eq:App-02-PWC-CC}
\begin{equation}
x = \sum_{s=1}^S \lambda_s x_s,~
y = \sum_{s=1}^{S-1} \theta_s f_s
\label{eq:App-02-PWC-CC-1}
\end{equation}
\begin{equation}
\lambda_1 \le \theta_1,~ \lambda_S \le \theta_{S-1}
\label{eq:App-02-PWC-CC-2}
\end{equation}
\begin{equation}
\lambda_s \le \theta_{s-1} + \theta_s,~ s = 2,\cdots, S-1 
\label{eq:App-02-PWC-CC-3}
\end{equation}
\begin{equation}
\lambda_s \ge 0,~ s = 1, \cdots, S,~ 
\sum\nolimits_{s=1}^S \lambda_s = 1
\label{eq:App-02-PWC-CC-4}
\end{equation}
\begin{equation}
\theta_s \in \{0,1\}, s=1,\cdots,S-1,~
\sum\nolimits_{s=1}^{S-1} \theta_s = 1
\label{eq:App-02-PWC-CC-5}
\end{equation}
\end{subequations}
In (\ref{eq:App-02-PWC-CC}), binary variable $\theta_s=1$ indicates interval $s$ is activated, and constraint (\ref{eq:App-02-PWC-CC-5}) ensures that only one interval will be activated; Furthermore, constraints (\ref{eq:App-02-PWC-CC-2})-(\ref{eq:App-02-PWC-CC-4}) enforce weigh coefficients $\alpha_s$, $s=1$, $\cdots$, $S$ to be $\mathbb{SOS}_2$; Finally, constraint (\ref{eq:App-02-PWC-CC-1}) expresses $y$ and $x$ via the linear combination of sampled values. The advantage of piecewise constant formulation (\ref{eq:App-02-PWC-CC}) lies in the binary expression of function value $y$, such that the product of $y$ and another continuous variable can be easily linearized via integer programming technique, which can be seen in Sect. \ref{App-B-Sect02-03}.

Clearly, the required number of binary variables introduced in formulation (\ref{eq:App-02-SOS2-MILP}) is $k$, which grows linearly with respect to the number of breakpoints, and the final MILP model may suffer from computational overheads due to the presence of a large number of binary variables when more breakpoints are involved for improving accuracy. In what follows, we present a useful formulation that only engages a logarithmic number of binary variables and constraints. This technique is proposed in \cite{App-MILP-SOS2-LogCC-1,App-MILP-SOS2-LogCC-2,App-MILP-SOS2-LogCC-3}. Consider the following constraints:
\begin{equation}
\label{eq:App-02-SOS2-Log}
\begin{aligned}
\sum_{i \in L_n} \lambda_i & \le z_n,~ \forall n \in N  \\
\sum_{i \in R_n} \lambda_i & \le 1 - z_n,~ \forall n \in N \\
                       z_n & \in \{0, 1\},~ \forall n \in N \\
   \lambda & \ge  0,~ \sum\nolimits^k_{i=0} \lambda_i = 1 
\end{aligned} 
\end{equation}
where $L_n$ and $R_n$ are index sets of weights $\lambda_i$, $N$ is an index set corresponding to the number of binary variables. The dichotomy sequences $\{L_n, R_n \}_{n \in N}$ constitute a branching scheme on the indices of weights, such that constraint (\ref{eq:App-02-SOS2-Log}) guarantees that at most two adjacent elements of $\lambda$ can take strictly positive values, so as to meet the $\mathbb{SOS}_2$ requirement.  The required number of binary variables $z_n$ is $\lceil \log_2 k \rceil$, which is significantly smaller than that involved in formulation (\ref{eq:App-02-SOS2-MILP}).

\begin{figure}[!t]
\centering
\includegraphics[scale=0.50]{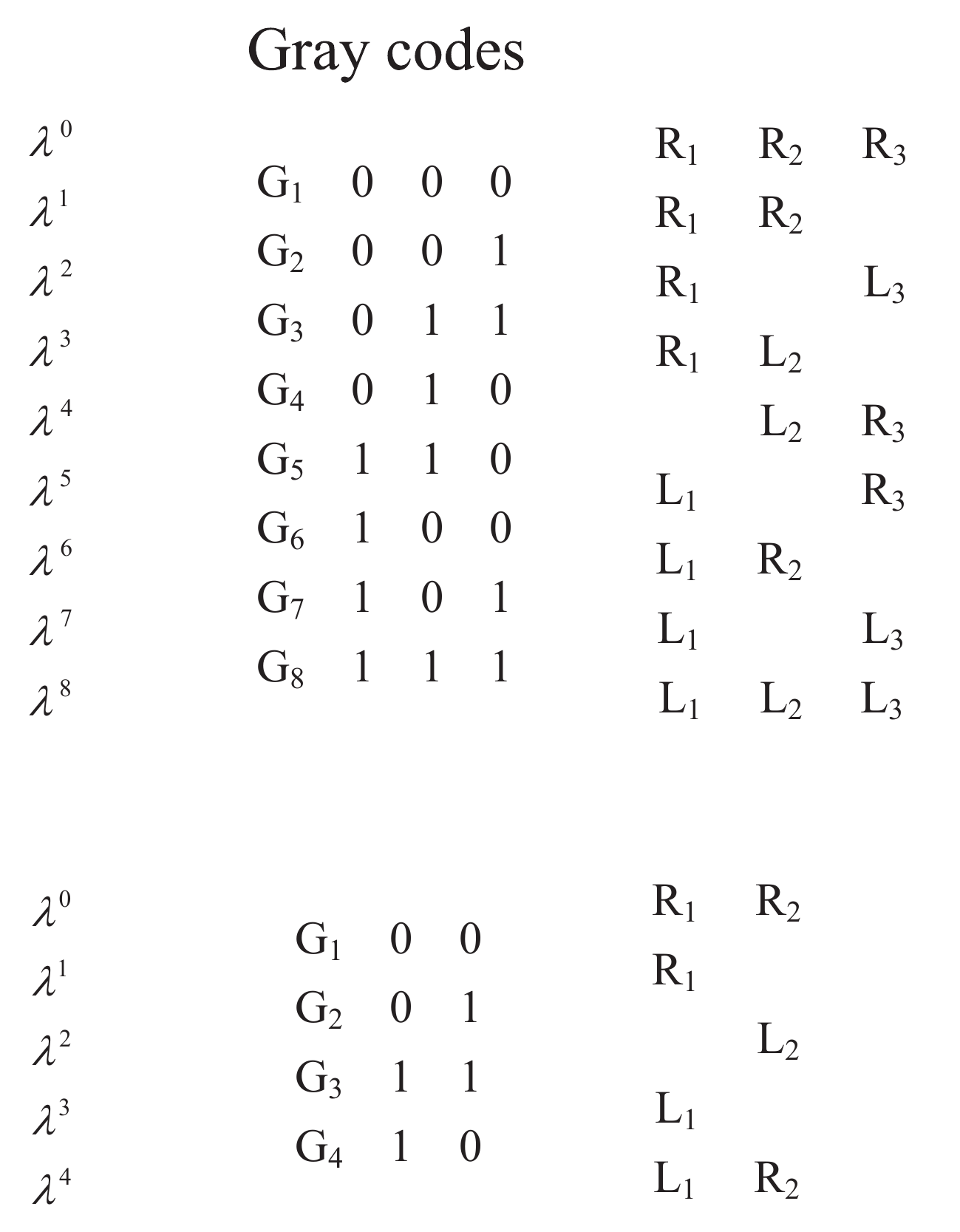}
\caption{Gray codes and sets $L_n$, $R_n$ for two and three binary variables.}
\label{fig:App-02-02}
\end{figure}

Next, we demonstrate how to design the sets $L_n$ and $R_n$ based on the concept of Gray codes. For notation brevity, we restrict the discussion to the instances with 2 and 3 binary variables (which are shown in Fig. \ref{fig:App-02-02}), indicating 5 and 9 breakpoints (or 4 and 8 intervals) in consequence. 

As shown in Fig. \ref{fig:App-02-02}, Gray codes G$_1$ - G$_8$ form a binary system where any two adjacent numbers only differ in one bit. For example, G$_4$ and G$_5$ differ in the first bit, and G$_5$ and G$_6$ differ in the second bit. Such Gray codes are used to describe which two adjacent weights are activated. In general, sets $R_n$ and $L_n$ are constructed as follows: the index $v \in L_n$ if the binary values of the $n$-th bit of two successive codes $G_n$ and $G_{n+1}$ are equal to 1, or $v \in R_n$ if they are equal to 0. This principle can be formally defined in a mathematical way as
\begin{equation}
\label{eq:App-02-SOS2-Log-L}
L_n = \left\{v~\middle|~ \begin{aligned}
      (G^n_v  &= 1 \mbox{ and } G^n_{v+1} = 1)   \\
\cup~ (v  & = 0 \mbox{ and } G^n_{1} =1)   \\
\cup~ (v  & = k \mbox{ and } G^n_{k} =1) 
\end{aligned}  \right\}
\end{equation}
\begin{equation}
\label{eq:App-02-SOS2-Log-R}
R_n = \left\{v~\middle|~ \begin{aligned}
      (G^n_v  &= 0 \mbox{ and } G^n_{v+1} = 0)   \\
\cup~ (v  & = 0 \mbox{ and } G^n_{1} = 0)   \\
\cup~ (v  & = k \mbox{ and } G^n_{k} = 0) 
\end{aligned}  \right\}
\end{equation}
where $G^n_v$ stands for the $n$-th bit of code $G_v$. 

For example, sets R$_1$, R$_2$, R$_3$ and L$_1$, L$_2$, L$_3$ for Gray codes G$_1$-G$_8$ are shown in Fig. \ref{fig:App-02-02}. In such a way, we can establish the rule that only two adjacent weights can be activated via (\ref{eq:App-02-SOS2-Log}). To see this, consider that if $\lambda_i > 0$ for $i=4,5$ and $\lambda_i = 0$ for other indices, we let $z_1=1$, $z_2=1$, $z_3$ = 0, which leads to the following constraint set:
\begin{equation*}
\left\{ \begin{lgathered}
\lambda_0 + \lambda_1 + \lambda_2 + \lambda_3 \le 1 - z_1 = 0 \\  
\lambda_5 + \lambda_6 + \lambda_7 + \lambda_8 \le z_1 = 1 \\ 
\lambda_0 + \lambda_1 + \lambda_6 \le 1 - z_2 = 0  \\
\lambda_3 + \lambda_4 + \lambda_8 \le z_2 = 1      \\
\lambda_0 + \lambda_4 + \lambda_5 \le 1 - z_3 = 1  \\ 
\lambda_2 + \lambda_7 + \lambda_8 \le z_3 = 0  \\ 
\lambda_i \ge 0, \forall i, \sum\nolimits^8_{i=0}  \lambda_i = 1
\end{lgathered}  \right.
\end{equation*}
Thus we can conclude that
\begin{gather*}
\lambda_4 + \lambda_5 = 1, \lambda_4 \ge 0, \lambda_5 \ge 0, \\  
\lambda_0=\lambda_1=\lambda_2=\lambda_3=\lambda_6=\lambda_7=\lambda_8=0
\end{gather*}
This mechanism can be interpreted as follows: $z_1=1$ enforces $\lambda_i = 0$, $i=0,1,2,3$
through set $R_1$; $z_2=1$ further enforces $\lambda_6 = 0$ through set $R_2$; finally, $z_3=0$ enforces $\lambda_7 = \lambda_8 = 0$ through set $L_3$. Then the remaining weights $\lambda_4$ and $\lambda_5$ constitute the positive coefficients. In this regard, only  $\log_2 8 =3$ binary variables and $2 \log_2 8 = 6$ additional constraints are involved. Compared with formulation (\ref{eq:App-02-SOS2-MILP}), the gray code can be regarded as extra branching operation enabled by problem structure, so the number of binary variables in expression (\ref{eq:App-02-SOS2-Log}) is greatly reduced in the case with a large value of $k$.

As a special case, consider the following problem
\begin{equation}
\label{eq:App-02-Example-1-NLP}
\min \left\{ \sum_i f_i(x_i) ~\middle|~ x \in X \right\}
\end{equation}
where $f_i(x_i),i=1,2,\cdots$ are convex univariate functions, and $X$ is a polytope. This problem is convex but nonlinear. The DCOPF problem, a fundamental issue in power market clearing, is given in this form, in which $f_i(x_i)$ is a convex quadratic function. Although any local NLP algorithm can find the global optimal solution of (\ref{eq:App-02-Example-1-NLP}), there are still reasons to seek approximated LP formulations. One is that problem (\ref{eq:App-02-Example-1-NLP}) may be embedded in another optimization problem and serve as its constraint. This is a pervasive modeling paradigm to study the strategic behaviors and market powers of energy providers, where the electricity market is cleared according to a DCOPF, and delivered energy of generation companies and nodal electricity prices are extracted from the optimal primal variables and dual variables associating with power balancing constraints, respectively. An LP representation allows to exploit the elegant LP duality theory for further analysis, and helps characterize optimal solution through primal-dual or KKT optimality conditions. To this end, we can opt to solve the following LP
\begin{equation}
\label{eq:App-02-Example-1-LP-Approx}
\begin{aligned}
\min_{x,y,\lambda} ~~ & \sum_i y_i  \\
\mbox{s.t.}~~ & y_i = \sum_{k}  \lambda_{ik} f_i(x_{ik}),~ \forall i   \\
& x_i = \sum_k  \lambda_{ik} x_{ik},~\forall i,~ x \in X   \\
& \lambda \ge 0,~ \sum_k \lambda_{ik} =1,~ \forall i \\
\end{aligned}
\end{equation}
where $x_{ik}, k=1,2,\cdots$ are break points (constants) for variable $x_i$, and the associated weights are $\lambda_{ik}$. Because $f_i(x_i)$ are convex functions, the $\mathbb{SOS}_2$ requirement on the weight variable $\lambda$ is naturally met, so it is relaxed from the constraints.   

\subsection{Bivariate Continuous Nonlinear Function}
\label{App-B-Sect01-02}

\begin{figure}[!t]
\centering
\includegraphics[scale=0.60]{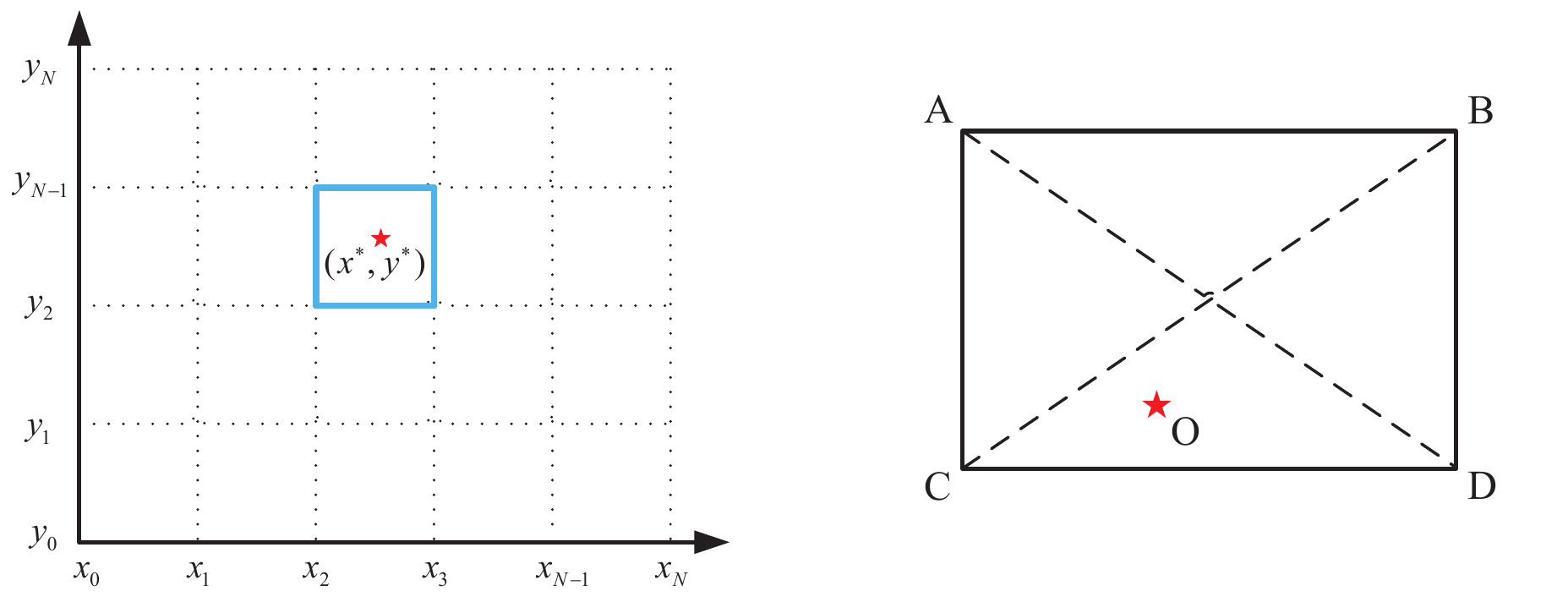}
\caption{Breakpoints and active rectangle for PWL approximation.}
\label{fig:App-02-03}
\end{figure}

Consider a continuous nonlinear function $f(x,y)$ in two variables $x$ and $y$. The entire feasible region is partitioned into $M \times N$ disjoint sub-rectangles by $M+N+2$ breakpoints $x_n$, $n=0,1,\cdots,N$ and $y_m$, $m=0,1,\cdots,M$, as illustrated in Fig. \ref{fig:App-02-03}, and the corresponding function values are $f_{mn} = f(x_m,y_n)$. By introducing a planar weighting coefficient matrix $\{ \lambda_{mn}\}$ for each grid point that satisfies 
\begin{subequations}
\begin{equation}
\label{eq:App-02-f(x,y)-Weight}
\begin{gathered}
\lambda_{mn} \ge 0,~ \forall m, \forall n   \\
\sum^M_{m=0} \sum^N_{n=0} \lambda_{mn} =1 
\end{gathered}   
\end{equation}
we can present any point $(x,y)$ in the feasible region by a convex combination of the extreme points of the sub-rectangle it resides in:
\begin{equation}
\label{eq:App-02-f(x,y)-Variable}
\begin{gathered}
x = \sum^M_{m=0} \sum^N_{n=0} \lambda_{mn} x_n = \sum^N_{n=0} 
\left( \sum^M_{m=0} \lambda_{mn} \right) x_n  \\
y = \sum^M_{m=0} \sum^N_{n=0} \lambda_{mn} y_m = \sum^M_{m=0} 
\left( \sum^N_{n=0} \lambda_{mn} \right) y_m
\end{gathered} 
\end{equation}
and its function value
\begin{equation}
\label{eq:App-02-f(x,y)-Value}
f(x,y) = \sum^M_{m=0} \sum^N_{n=0} \lambda_{mn} f_{mn}  
\end{equation}
\end{subequations}
is also a convex combination of the function values at the corner points.

As we can see from Fig. \ref{fig:App-02-03}, in a valid representation, if ($x^*, y^*$) belongs to a sub-rectangle, only the weight parameter associated with the four corner points can be non-negative, while others should be forced at 0. In such a pattern, the sum of columns/rows of matrix ${\rm \Lambda} = [\lambda_{mn}],\forall m,n$, which remains a vector, should constitute an $\mathbb{SOS}_2$, and $\rm \Lambda$ is called a planar $\mathbb{SOS}_2$, which can be implemented via two $\mathbb{SOS}_1$ constraints on the marginal weight vectors. In fact, at most three of the four corner points can be associated with uniquely determined non-negative weights. Consider point O and the active rectangle ABCD shown in Fig. \ref{fig:App-02-03}. The location of O can be expressed by a linear combination of the coordinates of corner points $x_A,x_B,x_C,x_D$ associating with non-negative weights $\lambda_A,\lambda_B,\lambda_C,\lambda_D$ as:
\begin{subequations}
\begin{equation}
\label{eq:App-02-Rectangle}
x_O = \lambda_A x_A + \lambda_B x_B + \lambda_C x_C + \lambda_D x_D \\
\end{equation}
In the first case 
\begin{equation}
\label{eq:App-02-Rectangle-Weight-Case1}
\lambda^1_A,~ \lambda^1_B,~ \lambda^1_C,~ \lambda^1_D \ge 0,~  
\lambda^1_A + \lambda^1_B + \lambda^1_C + \lambda^1_D = 1 
\end{equation}
In the second case
\begin{equation}
\label{eq:App-02-Rectangle-Weight-Case2}
\lambda^2_A, \lambda^2_C, \lambda^2_D \ge 0,~ \lambda^2_B = 0,~ 
\lambda^2_A + \lambda^2_C + \lambda^2_D = 1  
\end{equation}
In the third case
\begin{equation}
\label{eq:App-02-Rectangle-Weight-Case3}
\lambda^3_B, \lambda^3_C, \lambda^3_D \ge 0,~ \lambda^3_A = 0,~   
\lambda^3_B + \lambda^3_C + \lambda^3_D = 1 
\end{equation}
We use superscripts 1, 2, 3 to distinguish values of weights in different representations. According to Caratheodory theorem, the non-negative weights are uniquely determined in (\ref{eq:App-02-Rectangle-Weight-Case2}) and (\ref{eq:App-02-Rectangle-Weight-Case3}), and in the former (latter) case, we say $\rm \Delta$ACD ($\rm \Delta$BCD) is activated or selected. Denote function values in these three cases by  
\begin{gather}
f_1(x_O) = \lambda^1_A f(x_A) + \lambda^1_B f(x_B) + \lambda^1_C f(x_C) + \lambda^1_D f(x_D) \label{eq:App-02-Rectangle-Value-Case1} \\
f_2(x_O) = \lambda^2_A f(x_A) + \lambda^2_C f(x_C) + \lambda^2_D f(x_D) 
\label{eq:App-02-Rectangle-Weight-Case2} \\
f_3(x_O) = \lambda^3_B f(x_B) + \lambda^3_C f(x_C) + \lambda^3_D f(x_D)
\label{eq:App-02-Rectangle-Weight-Case3} 
\end{gather}
\end{subequations}
Suppose $f(x_A)< f(x_B)$, the plane defined by points B, C, D lies above that defined by points A, C, D, hence $f_2(x_O) < f_1(x_O) < f_3(x_O)$. If a smaller (larger) function value is in favor, then $\rm \Delta$ACD ($\rm \Delta$BCD) will be activated at the optimal solution. Please bear in mind that as long as A, B, C, D are not in the same plane, $f_1(x_O)$ will be strictly less (greater) than $f_3(x_O)$ ($f_2(x_O)$). Therefore, (\ref{eq:App-02-Rectangle-Value-Case1}) will not become binding at the optimal solution, and the weights for active corners are uniquely determined. If rectangle ABCD is small enough, such a discrepancy can be neglected. Nonetheless, non-uniqueness of the corner weights has little injury on its application, because the optimal solution $x_O$ and optimal value will be consistent with the original problem. The weights do not correspond to physical strategies that need to be deployed, and the linearization method can be considered as a black box to the decision maker, who provides function values at $x_A,x_B,x_C,x_D$ and receives a unique solution $x_O$. 

Detecting the active sub-rectangle that $(x^*, y^*)$ resides in requires  additional constraints on the weight parameter $\lambda_{mn}$. The aforementioned integer formulation is used to impose planar $\mathbb{SOS}_2$ constraints. Let $\lambda^n$ and $\lambda^m$ be the aggregated weights for $x$ and $y$, respectively, i.e.,
\begin{equation}
\label{eq:App-02-f(x,y)-fmn}
\begin{gathered}
\lambda^n = \sum^M_{m=0} \lambda_{mn},~ \forall n  \\  
\lambda^m = \sum^N_{n=0} \lambda_{mn},~ \forall m 
\end{gathered}
\end{equation}
which are also called the marginal weight vectors, and introduce the following constraints:
\begin{equation}
\label{eq:App-02-f(x,y)-MILP-x}
\mbox{For $x$: } \left\{
\begin{aligned}
\sum_{n \in L^1_k}  \lambda^n & \le z^1_k  \\
\sum_{n \in R^1_k}  \lambda^n & \le 1 - z^1_k  \\
                      z^1_k   & \in \{0, 1\} 
\end{aligned}   \right\},~ \forall k \in K_1
\end{equation}
\begin{equation}
\label{eq:App-02-f(x,y)-MILP-y}
\mbox{For $y$: } \left\{
\begin{aligned}
\sum_{m \in L^2_k}  \lambda^m & \le z^2_k  \\
\sum_{m \in R^2_k}  \lambda^m & \le 1 - z^2_k  \\
                      z^2_k   & \in \{0, 1\} 
\end{aligned}   \right\},~ \forall k \in K_2
\end{equation}
where $L^1_k$, $L^2_k$ and $R^2_k$, $R^2_k$ are index sets of weights $\lambda^n$ and $\lambda^m$, $K_1$ and $K_2$ are index sets of binary variables. The dichotomy sequences $\{L^1_k,R^1_k \}_{k \in K_1}$ and $\{L^2_k,R^2_k\}_{k \in K_2}$ constitute a branching scheme on the indices of weights, such that constraints (\ref{eq:App-02-f(x,y)-MILP-x}) and (\ref{eq:App-02-f(x,y)-MILP-y}) would guarantee that at most two adjacent elements of $\lambda^n$  and $\lambda^m$ can take strictly positive values, so as to detect the active sub-rectangle. In this approach, the required number of binary variables is $\lceil \log_2 M \rceil + \lceil \log_2 N \rceil$. The construction of these index sets has been explained in the univariate case.

Likewise, the discussions for problems (\ref{eq:App-02-Example-1-NLP}) and (\ref{eq:App-02-Example-1-LP-Approx})  can be easily extended if the objective function is the sum of bi-variate convex functions, implying that the planar $\mathbb{SOS}_2$ condition is naturally met.

\subsection{Approximation Error}

This section answers a basic question: For a given function, how many intervals (break points) are needed to achieve certain error bound $\varepsilon$? For the ease of understanding, we restrict our attention to univariate function, including the quadratic function $f(x)= ax^2$, and more generally, the continuous function $f(x)$  that is three times continuously differentiable. Let $\phi_f(x)$ be the PWL approximation for function $f(x)$ on $X=\{x|x_l \le x \le x_m \}$, the absolute maximum approximation error is defined by ${\rm \Delta} = \max_{x \in X} |f(x)-\phi_f (x)|$.

First let us consider the quadratic function $f(x) = a x^2$, which has been thoroughly studied in \cite{App-MILP-PWL-Error-1}. The analysis is briefly introduced here. Choose an arbitrary interval $[x_{i-1},x_i] \subset X$, the PWL approximation
can be parameterize in a single variable $t \in [0,1]$ as 
\begin{gather*}
x(t) = x_{i-1} + t (x_i-x_{i-1})  \\
\phi_f(x(t)) = ax_{i-1}^2 + at (x_i^2 - x_{i-1}^2) 
\end{gather*}
Clearly, $f(x(0))=\phi_f(x(0))=x^2_{i-1}$, $f(x(1))=\phi_f(x(1))=x^2_i$, and $\phi_f(x(t)) > f(x(t))$, $\forall t \in (0,1)$. The maximal approximation error in the interval must be found at a critical point which satisfies 
\begin{equation*}
\begin{aligned}
   & \dfrac{d}{dt}\left( \phi_f(x(t)) - f(x(t)) \right) \\
=~ & \dfrac{d}{dt} a \left[ x_{i-1}^2 + t (x_i^2 - x_{i-1}^2) - 
    ( x_{i-1} + t (x_i-x_{i-1}))^2 \right]  \\
=~ & \dfrac{d}{dt} a (x_i - x_{i-1})^2 (t-t^2) \\    
=~ & a (x_i - x_{i-1})^2 (1-2t) \\
=~ & 0 \Rightarrow t = \frac{1}{2}
\end{aligned}
\end{equation*}
implying that $x(1/2)$ is always a critical point where the approximation error reaches maximum, regardless of the partition of intervals, and the error is given by
\begin{equation*}
\begin{aligned}
{\rm \Delta} & = \phi \left( x \left( \frac{1}{2} \right) \right) 
                 - f\left( x \left( \frac{1}{2} \right) \right)  \\
 & = a \left[ x_{i-1}^2 + \frac{1}{2} (x_i^2 - x_{i-1}^2) 
     - \frac{1}{4}  ( x_{i-1} + x_i)^2 \right]   \\
 & = \frac{a}{4} (x_i - x_{i-1})^2      
\end{aligned}
\end{equation*}
which is quadratic in the length of the interval and independent of its location. In this regard, the intervals must be evenly distributed with equal length in order to get the best performance. If $X$ is divided into $n$ intervals, the absolute maximum approximation error is
\begin{equation*}
{\rm \Delta} =  \frac{a}{4n^2} (x_m - x_l)^2
\end{equation*}
Therefore, for a given tolerance $\epsilon$, the number of intervals should satisfy 
\begin{equation*}
n \ge \sqrt{\frac{a}{\varepsilon}}  \frac{x_m-x_l}{2}
\end{equation*}

For quadratic function $f(x)=a x^2$, coefficient $a$ determines its second-order derivative. For more general situations, above discussion implies that the number of intervals needed to perform a PWL approximation for function $f(x)$ may depend on its second-order derivative. This problem has been thoroughly studied in \cite{App-MILP-PWL-Error-2}. The conclusion is: for a three times continuously differentiable $f(x)$ in the interval $[x_l, x_m]$, the optimal number of segments $s(\varepsilon)$ under  given error tolerance $\varepsilon$ can be selected as
\begin{equation*}
s(\varepsilon) \propto \dfrac{c}{\sqrt{\varepsilon}},~ \varepsilon \to 0^+
\end{equation*}
where
\begin{equation*}
c = \dfrac{1}{4} \int^{x_m}_{x_l} \sqrt{|f^{\prime \prime}(x)|}
\end{equation*}
The conclusion still holds if $\sqrt{|f^{\prime \prime}(x)|}$ has integrable singularities at the endpoints.

\section{Linear Formulation of Product Terms}
\label{App-B-Sect02}

Product of two variables, or a bilinear term, naturally arises in optimization models from various disciplines. For one example, in economic studies, if the price $c$ and the quantity $q$ of a commodity are variables, then the cost $cq$ would be a bilinear term. For another, in circuit analysis, if both of the voltage $v$ and the current $i$ are variables, then the electric power $vi$ would be a bilinear term. Bilinear terms are non-convex. Throughout history, linearizing bilinear terms using linear constraints and integer variables is a frequently used technique in optimization community. This section presents several techniques for the central question of product linearization: how to enforce constraint $z=xy$, depending on the types of $x$ and $y$.

\subsection{Product of Two Binary Variables}
\label{App-B-Sect02-01}

If $x \in \mathbb {B}$ and $y \in \mathbb {B}$, then $z = x y$ is equivalent to the following linear inequalities  
\begin{equation}
\label{eq:App-02-xy-BB}
\begin{gathered}
0 \le z  \le  y    \\
0 \le x-z \le 1-y  \\
x \in \mathbb{B},~y \in \mathbb{B},~ z \in \mathbb{B} 
\end{gathered}
\end{equation}
It can be verified that if $x=1$, $y=1$, then $z = 1$ is achieved; if $x = 0$ or $y = 0$, then $z = 0$ is enforced, regardless of the value of $y$ or $x$. This is equivalent to the requirement $z = x y$.   

If $x \in \mathbb {Z}^+$ belongs to interval $[x^L,x^U]$, and $y \in \mathbb {Z}^+$ belongs to interval $[y^L,y^U]$, given the following binary expansion
\begin{equation}
\label{eq:App-02-xy-ZZ-BE}
\begin{aligned}
x = x^L + \sum_{k=1}^{K_1} 2^{k-1} u_k  \\
y = y^L + \sum_{k=1}^{K_2} 2^{k-1} v_k  
\end{aligned}
\end{equation}
where $K_1 = \lceil \log_2 (x^U-x^L) \rceil$, $K_2 = \lceil \log_2 (y^U-y^L) \rceil$. To develop a vector expression, define vectors $b^1 = [2^0, 2^1,\cdots,2^{K_1-1}]$, $u = [u_1,u_2,\cdots,u_{K_1}]^T$,
$b^2 = [2^0, 2^1,\cdots,2^{K_2}]$, $v = [v_1,v_2,\cdots,v_{K_2}]^T$, and matrices $B = (b^1)^T b^2$, $z = u v^T$, then
\begin{equation}
x y = x^L y^L + x^L b^2 v + y^L b^1 u + \langle B, z \rangle \notag
\end{equation}
where product matrix $z = u v^T$, and $\langle B, z \rangle = \sum_i \sum_j B_{ij} z_{ij}$. Relation among $u$, $v$, and $z$ can be  linearized via equation (\ref{eq:App-02-xy-BB}) element-wise. Its compact form is given by
\begin{equation}
\label{eq:App-02-xy-ZZ-Comp}
\begin{gathered}
{\bf 0}^{K_1 \times K_2} \le z  \le {\bf 1}^{K_1 \times 1} v^T  \\
{\bf 0}^{K_1 \times K_2} \le u {\bf 1}^{1 \times K_2} - z \le {\bf 1}^{K_1 \times K_2} - {\bf 1}^{K_1 \times 1} v^T  \\
u \in \mathbb{B}^{K_1 \times 1},~ v \in \mathbb{B}^{K_2 \times 1},~ 
z \in \mathbb{B}^{K_1 \times K_2} 
\end{gathered}
\end{equation}

\subsection{Product of Integer and Continuous Variables}
\label{App-B-Sect02-02}

We consider the binary-continuous case. If $x \in \mathbb {R}$ belongs to interval $[x^L,x^U]$, and $y \in \mathbb {B}$, and then $z = x y$ is equivalent to the following linear inequalities
\begin{equation}
\label{eq:App-02-xy-BC}
\begin{gathered}
x^L y \le z  \le  x^U y   \\
x^L (1-y) \le x-z \le x^U (1-y)  \\
x \in \mathbb{R},~y \in \mathbb{B},~ z \in \mathbb{R} 
\end{gathered}
\end{equation}
It can be verified that if $y = 0$, then $z$ is enforced to be 0 and $x^L \le x \le x^U$ is naturally met; if $y = 1$, then $z = x$ and $x^L \le z \le x^U$ must be satisfied, indicating the same relationship on $x$, $y$, and $z$. As for the integer-continuous case, the integer variable can be represented as (\ref{eq:App-02-xy-ZZ-BE}) using binary variables, yielding a linear combination of binary-continuous  products.

It should be mentioned that the upper bound $x^U$ and the lower bound $x^L$ are crucial for creating linearization inequalities. If explicit bounds are not available at hand, one can incorporate a constant $M$ that is big enough. The value of $M$ will have a notable impact on the computation time. To enhance efficiency, a desired value should be the minimal $M$ that ensures that inequality $-M \le x \le M$ never becomes binding at optimum, as it leads to the strongest bound if integrality of binary variables is neglected, expediting the converge of the branch-and-bound procedure. However, such a value is generally unclear before we solve the problem. Nevertheless, we do not actually need to find the smallest value $M^{*}$. Any $M \ge M^*$ produces the same optimal solution and is valid for linearization. Please bear in mind that an over-large $M$ not only deteriorates the computation time, but also cause numeric instability due to a large conditional number. So a proper tradeoff must be made between efficiency and accuracy.  A proper $M$ can be determined from estimating the bound of $x$ from certain heuristics, which is problem-dependent.

\subsection{Product of Two Continuous Variables}
\label{App-B-Sect02-03}

If $x \in \mathbb {R}$ belongs to interval $[x^L,x^U]$, and $y \in \mathbb {R}$ belongs to interval $[y^L,y^U]$, there are three options for linearizing their product $xy$. The first one considers $z = xy$ as a bivariate function $f(x,y)$, and applies the planar $\mathbb{SOS}_2$ method in Sect. \ref{App-B-Sect01-02}. The second one discretizes $y$, for example, as follows
\begin{equation}
\begin{gathered}
y = y^L + \sum_{k=1}^{K} 2^{k-1}  u_k {\rm \Delta} y \\
{\rm \Delta} y = \dfrac{y^U - y^L}{2^K},~ u_k \in \mathbb {B},~ \forall k
\end{gathered}
\label{eq:App-02-xy-CC-Discrete-y}
\end{equation}
and 
\begin{equation}
xy = x y^L + \sum_{k=1}^{K} 2^{k-1} v_k {\rm \Delta} y \\
\end{equation}
where $v_k = u_k x$ can be linearized through equation (\ref{eq:App-02-xy-BC})
as
\begin{equation}
\label{eq:App-02-xy-CC-BE}
\begin{gathered}
x^L u_k \le v_k  \le  x^U u_k,~ \forall k   \\
x^L (1 - u_k) \le x - v_k \le x^U (1 - u_k),~ \forall k   \\
 x  \in \mathbb{R},~
u_k \in \mathbb{B},~ \forall k,~
v_k \in \mathbb{R},~ \forall k 
\end{gathered}
\end{equation}

In practical problems, bilinear terms often appear as the inner production of two vectors. For convenience, we present the compact linearization of $x^T y$ via binary expansion. Let $y$ be the candidate vector variable to be discretized; perform (\ref{eq:App-02-xy-CC-Discrete-y}) on each element of $y$ 
\begin{equation*}
y_j = y^L_j + \sum_{k=1}^{K} 2^{k-1}  u_{jk} {\rm \Delta} y_j,~ \forall j
\end{equation*}
and thus
\begin{equation*}
x_j y_j = x_j y^L_j + \sum_{k=1}^{K} 2^{k-1} v_{jk} {\rm \Delta} y_j,~
\forall j,~ v_{jk} = u_{jk} x_j,~ \forall j, \forall k   
\end{equation*}
Relation $v_{jk} = u_{jk} x_j$ can be expressed via linear constraints 
\begin{equation*}
x^L_j u_{jk} \le v_{jk} \le x^U_j u_{jk},~ 
x^L_j(1-u_{jk}) \le x_j - v_{jk} \le x^U_j(1-u_{jk}),~ 
\forall j, \forall k
\end{equation*}

Denote by $V$ and $U$ are matrix variables consisting of $v_{jk}$ and $u_{jk}$, respectively; $1_K$ stands for all-one column vector with a dimension of $K$; ${\rm \Delta}_Y$ is a diagonal matrix with ${\rm \Delta} y_j$ being non-zero entries; vector $\zeta = [2^0,2^1,\cdots,2^{K-1}]$. Combining all above element-wise expressions together, we have the linear formulation of $x^T y$ in a compact matrix form
\begin{equation*}
x^T y = x^T y^L + \zeta V^T {\rm \Delta} y
\end{equation*}
in conjunction with
\begin{equation}
\begin{gathered}
y = y^L + {\rm \Delta}_Y U \zeta^T \\
(x^L \cdot 1_K^T) \otimes U \le V \le (x^U \cdot 1_K^T) \otimes U \\
(x^L \cdot 1_K^T) \otimes (1-U) \le x \cdot 1^T_K - V \le  (x^U \cdot 1_K^T) \otimes (1-U)  \\
x \in \mathbb{R}^J,~ y \in \mathbb R^J,~
U \in \mathbb{B}^{J \times K},~ V \in \mathbb{R}^{J \times K}
\end{gathered}
\label{eq:App-02-xy-CC-Vector-BE}
\end{equation}
where $\otimes$ represents element-wise product of two matrices with the same dimension. 

One possible drawback of this formulation is that the discretized variable is no longer continuous. The approximation accuracy can be improved by increasing the number of breakpoints without introducing too many binary variables, whose number is given by $\lceil \log_2 (y^U-y^L)/{\rm \Delta} y \rceil$. Furthermore, the variable to be discretized must have clear upper and lower bounds.  This is not restrictive because decision variables of engineering problems are subject to physical operating limitations, such as the maximum and minimum output of a generator. Nevertheless, if $x$, for example, is unbounded in formulation, but the problem has a finite optimum, we can replace $x^U(x^L)$  in (\ref{eq:App-02-xy-CC-BE}) with a large enough big-M parameter $M(-M)$, so that the true optimal solution remains feasible. It should be pointed out that the value of $M$ may influence the computational efficiency of the equivalent MILP, as mentioned previously. The optimal choice of $M$ in general cases remains an open problem, but there could be heuristic methods for specific instances. For example, if $x$ stands for the marginal production cost, which is a dual variable whose bounds are unclear, one can alternatively determine a suitable bound from historical data or price forecast. 

An alternative formulation for the second option deals with product term $x f(y)$, and $xy$ is a special case when $f(y)=y$. By performing the piecewise constant approximation  (\ref{eq:App-02-PWC-CC}) on function $f(y)$, the product becomes $xy = \sum_{s=1}^{S-1} x \theta_s y_s$, where $y_s$ is constant, $x$ is continuous, and $\theta_s$ is binary. The products $x \theta_s$, $s=1,\cdots,S-1$ can be readily linearized via the method in \ref{App-B-Sect02-02}. In this approach, the continuity of $x$ and $y$ are retained. However, the number of binary variables in the piecewise constant approximation for $f(y)$ grows linearly in the number of samples on $y$.

The third one converts the product into a separable form, and then performs PWL approximation for univariate nonlinear functions. To see this, consider a bilinear term $x y$. Introduce two continuous variables $u$ and $v$ defined as follows
\begin{equation}
\label{eq:App-02-xy-CC-Decomp-1}
\begin{gathered}
u = \frac{1}{2} (x+y) \\
v = \frac{1}{2} (x-y) 
\end{gathered}
\end{equation}
Now we have 
\begin{equation}
\label{eq:App-02-xy-CC-Decomp-2}
x y = u^2 - v^2
\end{equation}
In (\ref{eq:App-02-xy-CC-Decomp-2}), $u^2$ and $v^2$ are univariate nonlinear functions, and can be approximated by the PWL method  presented in Sect. \ref{App-B-Sect01-01}. Furthermore, if $x_l \le x \le x_u$, $y_l \le y \le y_u$, then the lower and upper bounds of $u$ and $v$ are given by
\begin{equation*}
\begin{gathered}
\frac{1}{2} (x_l + y_l) \le u \le \frac{1}{2} (x_u + y_u) \\
\frac{1}{2} (x_l - y_u) \le v \le \frac{1}{2} (x_u - y_l) 
\end{gathered}
\end{equation*}

Formulation (\ref{eq:App-02-xy-CC-Decomp-2}) has a connotative advantage. If $xy$ appears in the objective function which is to be minimized and is not involved in constraints, we only need to approximate $v^2$ because $u^2$ is convex and $-v^2$ is concave. The minimum amount of binary variables in this method is a logarithmic function in the number of break points, as explained in Sect. \ref{App-B-Sect01-01}.

The bilinear term $xy$ can be replaced by a single variable $z$ in the following situation: 1) if the lower bounds $x_l$ and $y_l$ are nonnegative; 2) either $x$ or $y$ is not referenced anywhere else except in $xy$. For instance, $y$ is such a variable, then the bilinear term $xy$ can be replaced by variable $z$ and constraint $x y_l \le z \le x y_u$. Once the problem is solved, $y$ can be recovered by $y = z/x$ if $x > 0$, and the inequality constraint on $z$ guarantees $y \in [y_l,y_u]$; otherwise if $x = 0$, then $y$ is undetermined and has no impact on the optimum.

\subsection{Monomial of Binary Variables}
\label{App-B-Sect02-04}

Previous cases discuss linearizing the product of two variables. Now we consider a binary monomial with $n$ variables
\begin{equation}
\label{eq:App-02-Monomial-Binary}
z = x_1 x_2 \cdots x_n,~ x_i \in \{0,1\},~ i=1,2,\cdots,n 
\end{equation}
Clearly, this monomial takes a binary value. Since the product of two binary can be expressed by a single one in light of (\ref{eq:App-02-xy-BB}), the monomial can be linearized recursively. Nevertheless, by making full use of the binary property of $z$, a smarter and concise way to represent (\ref{eq:App-02-Monomial-Binary}) is given by  
\begin{align}
z & \in \{0,1\}   \label{eq:App-02-Monomial-BL-1} \\
z & \le \dfrac{x_1 + x_2 + \cdots + x_n}{n} \label{eq:App-02-Monomial-BL-2}\\
z & \ge \dfrac{x_1 + x_2 + \cdots + x_n -n +1}{n} \label{eq:App-02-Monomial-BL-3}
\end{align}

If at least one of $x_i$ is equal to 0, because $\sum_{i=1}^n x_i - n + 1 \le 0$, (\ref{eq:App-02-Monomial-BL-3}) becomes redundant; moreover, $\sum^n_{i=1} x_i /n \le 1 - 1/n$, which removes $z=1$ from the feasible region, so $z$ will take a value of 0; otherwise, if all $x_i$ are equal to 1, $\sum^n_{i=1} x_i /n = 1$, and the right-hand side of (\ref{eq:App-02-Monomial-BL-3}) is $1/n$, which removes $z=0$ from the feasible region. Hence $z$ is forced to be 1. In conclusion, linear constraints (\ref{eq:App-02-Monomial-BL-1})-(\ref{eq:App-02-Monomial-BL-3}) have the same effect as (\ref{eq:App-02-Monomial-Binary}). 

In view of the above transformation technique, a binary polynomial program can always be reformulated as a binary linear program. Moreover, if a single continuous variable appears in the monomial,  the problem can be reformulated as an MILP.

\subsection{Product of Functions in Integer Variables}
\label{App-B-Sect02-05}

First, let us consider $z = f_1(x_1)f_2(x_2)$, where decision variables are positive integers, i.e., $x_i \in \{d_{i,1},d_{i,2},\cdots,d_{i,r_i}\},i=1,2$. Without particular tricks, $f_1$ and $f_2$ can be expressed as  
\begin{gather*}
f_1 = \sum_{j=1}^{r_1} f_1(d_{1,j})u_{1,j},~ u_{1,j} \in \{0,1\},~
\sum_{j=1}^{r_1} u_{1,j} = 1 \\
f_2 = \sum_{j=1}^{r_2} f_2(d_{2,j})u_{2,j},~ u_{2,j} \in \{0,1\},~
\sum_{j=1}^{r_2} u_{2,j} = 1 
\end{gather*}
and the product of two binary variables can be linearized via (\ref{eq:App-02-xy-BB}). Above formulation introduces a lot of intermediary binary variables, and is not propitious to represent a product with more functions recursively.   

Ref. \cite{App-MILP-Fun-Prod} suggests another choice 
\begin{equation}
z = \sum_{i=1}^{r_2} f_2(d_{2,j}) \sigma_{2,j},~ 
\sum_{i=1}^{r_2} \sigma_{2,j} = f_1(x_1),~ 
\sigma_{2,j} = f_1(x_1) u_{2,j}
\label{eq:App-02-Fun-Prod}
\end{equation}
where $u_{2,j} \in \{0,1\}$ and $\sum_{j=1}^{r_2} u_{2,j} = 1$. Although $f_1(x_1) u_{2,j}$ remains nonlinear because of decision variable $x_1$, (\ref{eq:App-02-Fun-Prod}) can be used to linearize a product with more than two nonlinear functions. 

To see this, Denote by $z_1 = f_1(x_1)$, $z_i = z_{i-1}f_i(x_i)$, $i = 2,\cdots,n$; integer variable $x_i \in \{d_{i,1},d_{i,2},\cdots,d_{i,r_i}\}$, $f_i(x_i)>0$, $i=1,\cdots,n$. By using (\ref{eq:App-02-Fun-Prod}), $z_i,i=1,\cdots,n$ have the following expressions \cite{App-MILP-Fun-Prod} 
\begin{equation}
\begin{aligned}
& z_1 = \sum_{j=1}^{r_1} f_1(d_{1,j})u_{1,j}  \\
& z_2 = \sum_{j=1}^{r_2} f_2(d_{2,j})\sigma_{2,j},~
\sum_{i=1}^{r_2} \sigma_{2,j} = z_1,~ \cdots \\
& z_n = \sum_{j=1}^{r_n} f_n(d_{n,j})\sigma_{n,j},~
\sum_{i=1}^{r_n} \sigma_{n,j} = z_{n-1}  \\ 
& \left. \begin{lgathered}
0 \le z_{i-1} - \sigma_{i,j} \le \bar z_{i-1}(1-u_{i,j}) \\
0 \le \sigma_{i,j} \le \bar z_{i-1} u_{i,j},~ u_{i,j} \in \{0,1\}
\end{lgathered}  \right\},~ 
\begin{lgathered}
j = 1, \cdots, r_i, \\
i = 2, \cdots, n 
\end{lgathered}  \\
& x_i = \sum_{j=1}^{r_i} d_{i,j} u_{i,j}, ~ \sum_{j=1}^{r_i} u_{i,j} = 1,~ i = 1,2, \cdots,n
\end{aligned}
\label{eq:App-02-Funs-Prod}
\end{equation}

In (\ref{eq:App-02-Funs-Prod}), the number of binary variables is $\sum_{i=1}^n r_i$, and grow linearly in the dimension of $x$ and the interval length of each $x_i$. To reduce the number of auxiliary binary variable $u_{i,j}$, the dichotomy procedure in Sect. \ref{App-B-Sect01-01} for SOS2 can be applied, which is discussed in \cite{App-MILP-Fun-Prod}.

\subsection{Log-sum Functions}
\label{App-B-Sect02-06}

We consider log-sum function $\log(x_1+x_2+\cdots+x_n)$, which arises from solving a signomial geometric programming problem. The basic element in such a problem has a form of
\begin{equation}
\label{eq:App-02-Signomial}
c_k \prod_{j=1}^l y_j^{a_{jk}}
\end{equation}
where $y_j > 0$, $c_k$ is a constant, and $a_{jk} \in \mathbb R$. Non-integer value of $a_{jk}$ makes signomial geometric programming problem even harder than polynomial programs. Under some variable transformation, the non-convexity of a signomial geometric program can be concentrated in some log-sum functions \cite{App-MILP-Signomial}. In view of the form in (\ref{eq:App-02-Signomial}), we discuss log-sum function in Sect. \ref{App-B-Sect02}. 

We aim to represent function $\log(x_1+x_2+\cdots+x_n)$ in terms of $\log x_1$, $\log x_2$, $\cdots$, $\log x_n$. Following the method in \cite{App-MILP-Signomial}, define a univariate function $F(X) = \log(1+e^X)$ and let $X_i=\log x_i$, ${\rm \Gamma}_i=\log(x_1+\cdots+x_i)$, $i = 1,\cdots,n$. The relation between $X_i$ and ${\rm \Gamma}_i$ can be revealed. Because
\begin{equation*}
\begin{aligned}
F(X_{i+1}-{\rm \Gamma}_i) & = \log\left( 1+ e^{\log x_{i+1} - \log(x_1+\cdots+x_i)} \right) \\
& = \log \left( 1+ \dfrac{x_{i+1}}{x_1+\cdots+x_i}\right) 
  = {\rm \Gamma}_{i+1} - {\rm \Gamma}_i
\end{aligned}
\end{equation*}
By stipulating $W_i = X_{i+1}-{\rm \Gamma}_i $, we have the following recursive equations
\begin{equation}
\label{eq:App-02-Signomial-Recursive}
\begin{aligned}
{\rm \Gamma}_{i+1} & = {\rm \Gamma}_i + F(W_i),~ i = 1,\cdots,n-1 \\
W_i & = X_{i+1}-{\rm \Gamma}_i,~ i = 1,\cdots,n-1 
\end{aligned}
\end{equation}

Function $F(W_i)$ can be linearized using the method in Sect. \ref{App-B-Sect01-01}. Based on this technique, an outer-approximation approach is proposed in \cite{App-MILP-Signomial} to solve signomial geometric programming problem via MILP.

\section{Other Frequently used Formulations}
\label{App-B-Sect03}

\subsection{Minimum Values}
\label{App-B-Sect03-01}

Let $x_1$, $x_2$, $\cdots$, $x_n$ be continuous variables with known lower bound $x^L_i$ and upper bound $x^U_i$, and $L = \min \{x^L_1, x^L_2, \cdots, x^L_n\}$, then their minimum $y = \min \{x_1, x_2, \cdots, x_n\}$ can be expressed via linear constraints
\begin{equation}
\label{eq:App-02-Minimum-PWL}
\begin{gathered}
x^L_i \le x_i \le x^U_i,~ \forall i  \\
y \le x_i, \forall i \\
x_i-(x^U_i - L)(1-z_i) \le y, \forall i \\
z_i \in \mathbb {B},~ \forall i,~  \sum_i z_i = 1 
\end{gathered}
\end{equation}
The second inequality guarantees $y \le \min \{x_1, x_2, \cdots, x_n\}$; in addition, if $z_i =1$, then $y \ge x_i$, hence $y$ achieves the minimal value of $\{x_i\}$. According to the definition of $L$, $x_i - y \le x^U_i - L, \forall i$ holds, thus the third inequality is inactive for the remaining $n-1$ variables with $z_i = 0$.

\subsection{Maximum Values}
\label{App-B-Sect03-02}

Let $x_1$, $x_2$, $\cdots$, $x_n$ be continuous variables with known lower bound $x^L_i$ and upper bound $x^U_i$,  and $U = \max \{x^U_1, x^U_2, \cdots, x^U_n\}$, then their maximum $y = \max \{x_1, x_2, \cdots, x_n\}$ can be expressed via linear constraints
\begin{equation}
\label{eq:App-02-Maximum-PWL}
\begin{gathered}
x^L_i \le x_i \le x^U_i,~ \forall i  \\
y \ge x_i, \forall i \\
x_i + (U - x^L_i)(1-z_i) \ge y, \forall i \\
z_i \in \mathbb {B},~ \forall i,~  \sum_i z_i = 1 
\end{gathered}
\end{equation}
The second inequality guarantees $y \ge \max \{x_1, x_2, \cdots, x_n\}$; in addition, if $z_i =1$, then $y \le x_i$, hence $y$ achieves the maximal value of $\{x_i\}$. According to the definition of $U$, $y - x_i \le U - x^L_i, \forall i$ holds, thus the third inequality is inactive for the remaining $n-1$ variables with $z_i = 0$.

\subsection{Absolute Values}
\label{App-B-Sect03-03}

Suppose $x \in \mathbb {R}$ and $|x| \le U$, the absolute value function $y = |x|$, which is nonlinear, can be expressed via PWL function as
\begin{equation}
\label{eq:App-02-Abs-PWL}
\begin{gathered}
 0 \le y - x \le 2 U z,~  U (1-z) \ge x  \\
 0 \le y + x \le 2 U (1-z),~ -U z  \le x  \\
 -U \le x \le U,~ z \in \mathbb{B}   
\end{gathered}
\end{equation}
When $x > 0$, the first line yields $z = 0$ and $y = x$, while the second line is inactive. When $x < 0$, the second line yields $z = 1$ and $y = -x$, while the first line is inactive. When $x = 0$, either $z = 0$ or $z = 1$ gives $y = 0$. In conclusion, (\ref{eq:App-02-Abs-PWL}) has the same effect as $y=|x|$.

\subsection{Linear Fractional of Binary Variables}
\label{App-B-Sect03-04}

A linear fractional of binary variables takes the form of 
\begin{equation}
\label{eq:App-02-BLF-1}
\dfrac{a_0+\sum_{i=1}^n a_i x_i}{b_0+\sum_{i=1}^n b_i x_i}
\end{equation}
We assume $b_0+\sum_{i=1}^n b_i x_i \ne 0$ for all $x \in \{0,1\}^n$. Define a new continuous variable
\begin{equation}
\label{eq:App-02-BLF-2}
y = \dfrac{1}{b_0+\sum_{i=1}^n b_i x_i}
\end{equation}
The lower bound and upper bound of $y$ can be easily computed. Then the linear fractional shown in (\ref{eq:App-02-BLF-1}) can be replaced with a linear expression
\begin{equation}
\label{eq:App-02-BLF-3}
a_0 y + \sum_{i=1}^n a_i z_i \\
\end{equation}
with constraints 
\begin{equation}
\label{eq:App-02-BLF-4}
b_0 y + \sum_{i=1}^n b_i z_i = 1
\end{equation}
\begin{equation}
\label{eq:App-02-BLF-5}
z_i = x_i y,~ \forall i
\end{equation}
where (\ref{eq:App-02-BLF-5}) describes a product of a binary variable and a continuous variable, which can be linearized through equation (\ref{eq:App-02-xy-BC}).

\subsection{Disjunctive Inequalities}
\label{App-B-Sect03-05}

Let $\{P^i\}, i = 1, 2, \cdots, m$ be a finite set of bounded polyhedra. Disjunctive inequalities usually arise when the solution space is characterized by the union $S = \cup_{i=1}^m P^i$ of these polyhedra. Unlike intersection operator which preserves convexity, disjunctive inequalities form a non-convex region. It can be represented by MILP model using binary variables. We introduce three emblematic methods.

\vspace{12pt}
{\noindent \bf 1. Big-M formulation}

\begin{figure}[!t]
\centering
\includegraphics[scale=0.52]{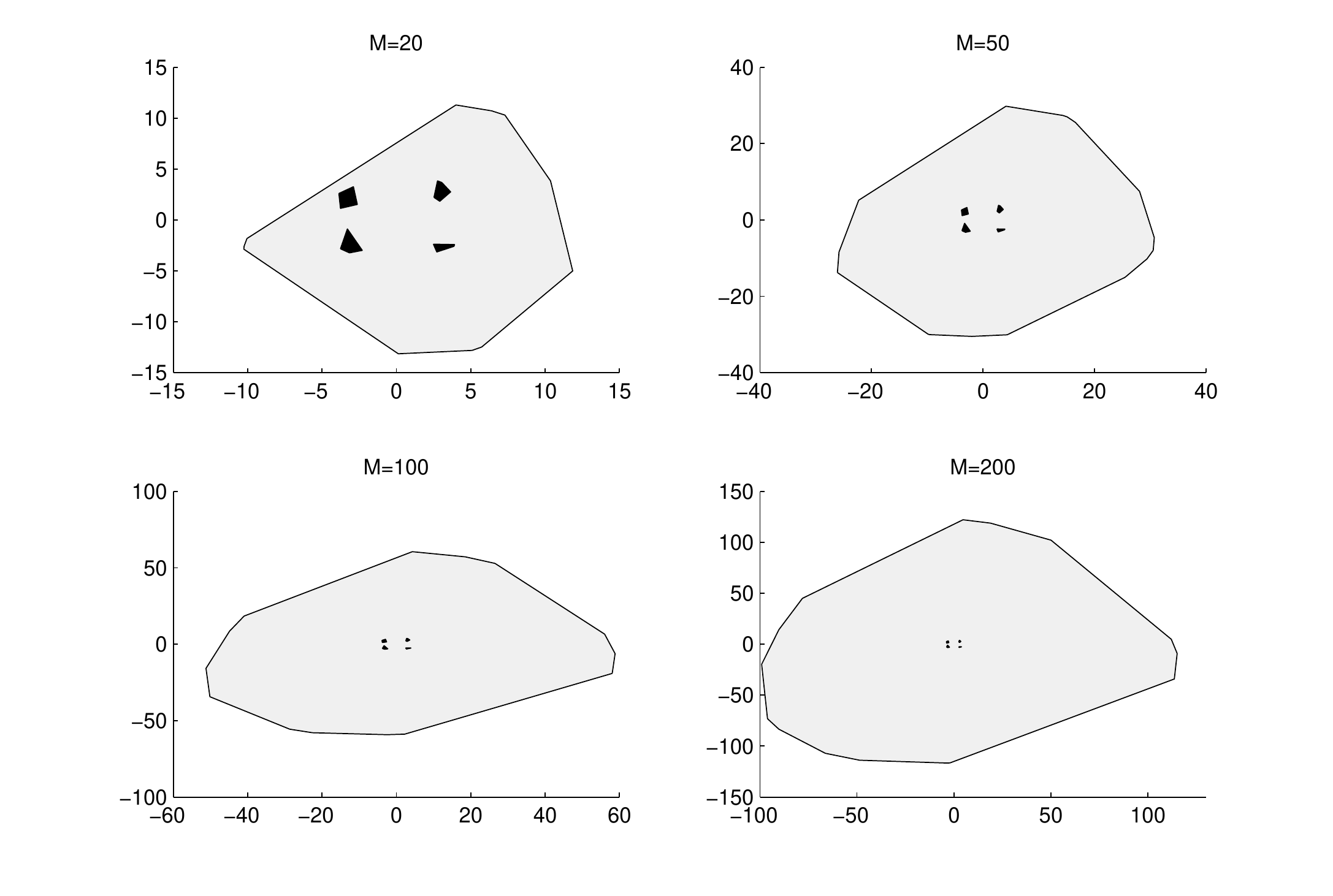}
\caption{Big-M formulation and their relaxed regions.}
\label{fig:App-02-04}
\end{figure}

The hyperplane representations of polyhedra are given by $P^i = \{ x \in \mathbb{R}^n | A^i x \le b^i \}, i = 1, 2, \cdots, m$. By introducing binary variables $z_i, i = 1, 2, \cdots, m $, an MILP formulation for $S$ can be written as
\begin{equation}
\label{eq:App-02-Disj-Big-M}
\begin{gathered}
A^i x \le b^i + M^i(1-z_i) ,~ \forall i  \\
z_i \in \mathbb{B},~\forall i,~ \sum^m_{i=1} z_i = 1
\end{gathered}
\end{equation}
where $M^i$ is a vector such that when $z_i = 0$, $A^i x \le b^i + M^i$ holds. To show the impact of the value of $M$ on the tightness of formulations (\ref{eq:App-02-Disj-Big-M}) when integrality constraints $z_i \in \mathbb{B}, \forall i$ are relaxed as $z_i \in [0,1], \forall i$, we contrivedly construct 4 polyhedra in $\mathbb {R}^2$, which are depicted in Fig. \ref{fig:App-02-04}. The continuous relaxations of (\ref{eq:App-02-Disj-Big-M}) with different values of $M$ are illustrated in the same graph, showing that the smaller the value of $M$, the tighter the relaxation of (\ref{eq:App-02-Disj-Big-M}). 

From a computational perspective, the element in $M$ should be as small as possible, because a huge constant without any insights about problem data will feature a bad conditional number. Furthermore, the continuous relaxation of MILP model will be very weak, resulting in poor objective value bounds and excessive branch-and-bound computation. The goal of big-M parameter selection is to create a model whose continuous relaxation is close to the convex hull of the original constraint, i.e. the smallest convex set that contains the original feasible region. A possible selection of the big-M parameter is 
\begin{equation}
\label{eq:App-02-Value-Big-M}
\begin{gathered}
M^i_l = \left(\max_{j \ne i} M^{ij}_l \right) - b^i_l \\
M^{ij}_l = \max_x \left\{ [A^i x]_l: A^j x \le b^j \right\}
\end{gathered}
\end{equation}
where subscript $l$ stands for the $l$-th element of a vector or $l$-th row of a matrix. As polyhedron $P^i$ are bounded, all bound parameters in (\ref{eq:App-02-Value-Big-M}) are well defined. 

However, even the tightest big-M parameter will yield a relaxed solution space that is generally larger than the convex hull of the original feasible set. In many applications, good variable bounds can be estimated from certain heuristic methods which explore specific problem data and structure.

\vspace{12pt}
{\noindent \bf 2. Convex hull formulation}

Let $\mbox{vert}(P^i) = \{ v^i_l \}, l = 1,2, \cdots,L^i$ denote sets of vertices of polyhedra $\{P^i\}, i = 1,2,\cdots,m$, where $L^i$ is the number of vertices of $P^i$. The set of extreme rays is empty since $P^i$ is bounded. By introducing binary variables $z_i, i = 1, 2, \cdots, m $, an MILP formulation for $S$ is given by
\begin{equation}
\label{eq:App-02-Disj-Conv}
\begin{gathered}
\sum_{i=1}^m  \sum_{l=1}^{L^i} \lambda^i_l v^i_l  = x \\
\sum_{l=1}^{L^i} \lambda^i_l = z_i,~ \forall i  \\
\lambda^i_l \ge 0,~ \forall i,~ \forall l  \\ 
z_i \in \mathbb{B},~\forall i,~ \sum^m_{i=1} z_i = 1
\end{gathered}
\end{equation}

Formulation (\ref{eq:App-02-Disj-Conv}) does not rely on manually supplied parameter. Instead, it requires enumerating all extreme points of polyhedra $P^i$. Although the vertex representation and hyperplane representation of a polyhedron are interchangeable, given the fact that vertex enumeration is time consuming for high-dimensional polyhedra, (\ref{eq:App-02-Disj-Conv}) is useful only if $P^i$ are originally represented by extreme points.

\vspace{12pt}
{\noindent \bf 3. Lifted formulation}

A smarter formulation exploits the fact that bounded polyhedra $P^i$ share the same recession cone $\{ 0 \}$, i.e., equation $A^i x = 0$ has no non-zero solutions. Otherwise, suppose $A^i x^* = 0$, $x^* \ne 0$, and $y \in P^i$, then $y + \lambda x^* \in P^i$, $\forall \lambda > 0$, because $A^i(y + \lambda x^*) = A^i y \le b^i$. As a result, $P^i$ is unbounded. Bearing this in mind, an MILP formulation for $S$ is given by
\begin{equation}
\label{eq:App-02-Disj-Lift}
\begin{gathered}
A^i x^i \le b^i z^i,~ \forall i  \\
\sum_{i=1}^m x^i = x  \\
z_i \in \mathbb{B},~\forall i \\
\sum^m_{i=1} z_i = 1
\end{gathered}
\end{equation}

Formulation (\ref{eq:App-02-Disj-Lift}) is also parameter-free. Since it incorporates additional continuous variable for each polytope, we call it a lifted formulation. It is easy to see that the feasible region of $x$ is the union of $P^i$: if $z_i = 0$, $x^i = 0$ as analyzed before; otherwise, if $z_i = 1$, $x = x^i \in P^i$.  

\vspace{12pt}
{\noindent \bf 4. Complementarity and slackness condition}

Complementarity and slackness condition naturally arises in the KKT optimality condition of a mathematical programming problem, an equilibrium problem, a hierarchical optimization problem, and so on. It is a quintessential law to characterize the logic condition under which a rational decision-making progress must obey. Here we pay attention to the linear case and equivalent MILP formulation, because nonlinear cases give rise to MINLPs, which are challenging to solve and not superior from the computational point of view.

A linear complementarity and slackness condition can be written as
\begin{equation}
\label{eq:App-B-LCP-MILP-1}
0 \le y \bot Ax - b \ge 0
\end{equation}
where vectors $x$ and $y$ are decision variables; $A$ and $b$ are constant coefficients with compatible dimensions; notation $\bot$ stands for the orthogonality of two vectors. In fact, (\ref{eq:App-B-LCP-MILP-1}) encompasses the following nonlinear constraints in traditional form
\begin{equation}
\label{eq:App-B-LCP-MILP-2}
y \ge 0,~ Ax - b \ge 0,~ y^T(Ax-b) = 0
\end{equation}
In view of the non-negativeness of $y$ and $Ax-b$, the orthogonality condition is equivalent to the element-wise logic form $y_i = 0$ or $a_i x-b_i=0$, $\forall i$, where $a_i$ is the $i$-th row of $A$; in other words, at most one of $y_i$ and $a_i x-b_i$ can take a strictly positive value,  implying that the feasible region is either the slice $y_i = 0$ or the slice $a_i x-b_i=0$. Therefore, (\ref{eq:App-B-LCP-MILP-1}) can be regarded as a special case of the disjunctive constraints.

In practical application, (\ref{eq:App-B-LCP-MILP-1}) usually serves as constraints in an optimization problem. For example, in a sequential decision making or a linear bilevel program, the KKT condition of the lower-level LP appears in the form of  (\ref{eq:App-B-LCP-MILP-1}), which is the constraint of the upper-level optimization problem. The main computation challenge arises from the  orthogonality condition, which is nonlinear and non-convex, and violates the linear independent constraint qualification, see Appendix \ref{App-D-Sect03} for an example. Nonetheless, in view of the switching logic between $y_i$ and $a_i x-b_i$, we can introduce a binary variable $z_i$ to select which slice is active \cite{App-MILP-Fortuny-Amat}
\begin{equation}
\label{eq:App-B-LCP-MILP-3}
\begin{gathered}
0 \le a_i x-b_i \le M z_i,~ \forall i   \\
0 \le y_i \le M(1-z_i),~ \forall i
\end{gathered}
\end{equation}
where $M$ is a large enough constant. According to (\ref{eq:App-B-LCP-MILP-3}), if $z_i=0$, then $(Ax - b)_i=0$ must hold, and the second inequality is redundant; otherwise, if $z_i=1$, then we have $y_i=0$, and the first inequality becomes redundant. (\ref{eq:App-B-LCP-MILP-3}) can be written in a compact form as
\begin{equation}
\label{eq:App-B-LCP-MILP-4}
\begin{gathered}
0 \le A x - b \le M z   \\
0 \le y \le M(1-z)
\end{gathered}
\end{equation}

It is worth mentioning that the big-M parameter $M$ has a notable impact on the feasible region of the relaxed problem as well as the computational efficiency of the MILP model, as illustrated in Fig. \ref{fig:App-02-04}. One should make sure that (\ref{eq:App-B-LCP-MILP-4}) would not remove the optimal solution from the feasible set. If both $x$ and $y$ have clear bounds, then $M$ can be easily estimated; otherwise, we may prudently employ a large $M$, at the cost of sacrificing the computational efficiency.

Furthermore, if we are aiming to solve (\ref{eq:App-B-LCP-MILP-1}) without an objective function and other constraints, such a problem is called a linear complementarity problem (under some proper transformation), for which we can build parameter-free MILP models. More details can be found in Appendix \ref{App-D-Sect04-02}.

\subsection{Logical Conditions}

Logical conditions are associated with indicator constraints with a statement like ``if event A then event B''. An event can be described in many ways. For example, a binary variable $a=1$ can stand for event A happens, and otherwise $a=0$; a point $x$ belongs to a set $X$ can denote a system is under secure operating condition, and otherwise $x \notin X$. In view of this, the disjunctive constraints discussed above is a special case of logical condition. In this section, we expatiate on how some usual logical conditions can be expressed via linear constraints. Let A, B, C, $\cdots$ associated with binary variables $a$, $b$, $c$, $\cdots$ represent events. Main results for linearizing typical logical conditions are summarized in \ref{tab:App-B-Logic-MILP} \cite{App-MILP-Logic-Cons}.

\begin{table}[htp]
\small
\renewcommand{\arraystretch}{1.3}
\renewcommand{\tabcolsep}{1em}
\caption{Linear form of some typical logic conditions}
\centering
\begin{tabular}{ll}
\hline
If A then B           &  $b \ge a$      \\
Not B                 &  $1-b$          \\ 
If A then not B       &  $a+b \le 1$    \\
If not A then B       &  $a+b \ge 1$    \\
A if and only if B    &  $a=b$          \\
If A then B and C     &  $b+c \ge 2a$   \\
If A then B or C      &  $b+c \ge a$    \\
If B or C then A      &  $2a \ge b+c$   \\
If B and C then A     &  $a \ge b+c-1$  \\
If M or more of N events then A  & $(N-M+1)a \ge b+c+\cdots-M+1$ \\  
\hline
\end{tabular}
\label{tab:App-B-Logic-MILP}
\end{table}

Logical AND is formulated as a function of two binary inputs. Specifically, $c = a$ AND $b$ can be expressed as $c=\min\{a,b\}$ or $c=ab$. The former one can be linearized via (\ref{eq:App-02-Minimum-PWL}) and the letter one through (\ref{eq:App-02-xy-BB}), and both of them renders 
\begin{equation}
\label{eq:App-B-Logic-AND-2}
c \le a,~ c \le b,~ c \ge a+b-1,~ c \ge 0  \\
\end{equation}

For the case with multiple binary inputs, i.e., $c=\min\{c_1,\cdots,c_n\}$, or $c = \prod_{i=1}^n c_i$, (\ref{eq:App-B-Logic-AND-2}) can be generalized as
\begin{equation}
\label{eq:App-B-Logic-AND-N}
c \le c_i,~ \forall i,~ c \ge \sum\nolimits_{i=1}^n c_i -n +1,~ c \ge 0
\end{equation}

Logical OR is formulated as a function of two binary inputs, i.e., $c=\max \{a,b\}$, which can be linearized via (\ref{eq:App-02-Maximum-PWL}), yielding 
\begin{equation}
\label{eq:App-B-Logic-OR-2}
c \ge a,~ c \ge b,~ c \le a+b,~ c \le 1  \\
\end{equation}

For the case with multiple binary inputs, i.e., $c=\max \{c_1,\cdots,c_n\}$, (\ref{eq:App-B-Logic-OR-2}) can be generalized as
\begin{equation}
\label{eq:App-B-Logic-OR-N}
c \ge c_i,~ \forall i,~ c \le \sum\nolimits_{i=1}^n c_i,~ c \le 1
\end{equation}

\section{Further Reading}
\label{App-B-Sect04}
 
Throughout the half-century long research and development, MILP has become an indispensable and unprecedentedly powerful modeling tool in mathematics and engineering, thanks to the advent of efficient solvers that encapsulate many state-of-the-art techniques  \cite{{App-MILP-History}}. This chapter aims to provide an overview on formulation recipes that transform complicated conditions into MILPs, so as to take full advantages of off-the-shelf solvers. The paradigm is able to deal with a fairly broad class of hard optimization problems. 

Readers who are interested in the  strength of MILP model, may find in-depth discussions in \cite{App-MILP-Strength} and references therein. For those interested in the PWL approximation of nonlinear functions, we refer to \cite{App-MILP-PWL-Function-1,App-MILP-PWL-Function-2,App-MILP-PWL-Function-3} and references therein, for various models and methods. The most promising one may be the convex combination model with a logarithmic number of binary variables, whose implementation has been thoroughly discussed in \cite{App-MILP-SOS2-LogCC-1,App-MILP-SOS2-LogCC-2,App-MILP-SOS2-LogCC-3}. For those who are interested in the polyhedral study of single-term bilinear sets and MILP based methods for bilinear programs may find extensive information in \cite{App-MILP-MIBLP-1,App-MILP-MIBLP-2} and references therein.
For those who need more knowledge about mathematical program with disjunctive constraints, in which  constraint activity is controlled by logical conditions, we recommend \cite{App-MILP-Disj-Review}; specifically, the choice of big-M parameter is discussed in \cite{App-MILP-Disj-Big-M}. For those who wish to learn more about integer programming techniques, we refer to \cite{App-MILP-Union} for the formulation of union of polyhedra, \cite{App-MILP-Representability} for the representability of MILP, and \cite{App-MILP-MICQP,App-MILP-Duality} for the more general mixed-integer conic programming as well as its duality theory. To the best of our knowledge, dissertation \cite{App-MILP-Dissertation-MIT} launches the most comprehensive and in-depth study on MILP approximation of non-convex optimization problems. State-of-the-art MILP formulations which balance problem size, strength, and branching behavior are developed and compared, including those mentioned above. The discussions in \cite{App-MILP-Dissertation-MIT} offer insights on designing efficient MILP models that perform extremely well in practice, despite of their theoretically non-polynomial complexity in the worst case.

%
%
%

\motto{To be uncertain is to be uncomfortable,
but to be certain is to be ridiculous.}

\chapter{Basics of Robust Optimization}
\label{App-C} 

Real-world decision-making models often involve unknown data. Reasons for data uncertainty could come from inexact measurements or forecast errors. For example, in power system operation, the wind power generation and system loads are barely known exactly at the time when the generation schedule should be made; in inventory management, market price and demand volatility is the main source of financial risks. In fact, optimal solutions to mathematical programming problems can be highly sensitive to parameter perturbations \cite{RO-Detail-1}. The optimal solution to the nominal problem may be highly suboptimal or even infeasible in reality due to parameter inaccuracy. Consequently, there is a great need of a systematic methodology that is capable of quantifying the impact of data inexactness on the solution quality, and is able to produce robust solutions that are insensitive to data uncertainty. 

Optimization under uncertainty has been a focus of the operational research community for a long time. Two approaches are prevalent to deal with uncertain data in optimization, namely stochastic optimization (SO) and robust optimization (RO). They differ in the ways of modeling uncertainty. The former one assumes that the true probability distribution of uncertain data is known or can be estimated from available information, and minimizes the expected cost in its objective function. SO provides strategies that are optimal in the sense of statistics. However, the probability distribution itself may be inexact owing to the lack of enough data, and the performance of the optimal solution could be sensitive to the probability distribution chosen in the SO model. The latter one considers uncertain data resides in a pre-defined uncertainty set, and minimizes the cost in the worst-case scenario in its objective function. Constraint violation is not allowed for all possible data realizations in the uncertainty set. RO is popular because it relies on simple data and distribution-free. From the computational perspective, it is equivalent to convex optimization problems for a variety of uncertainty sets and problem types; for the intractable cases, it can be solved via systematic iteration algorithms. For more technical details about RO, we refer to \cite{RO-Detail-1, RO-Detail-2, RO-Guide,RO-Convex}, survey articles \cite{RO-Survey,RO-Survey-2018}, and many references therein. Recently, distributionally robust optimization (DRO), an emerging methodology that inherits the advantages of SO and RO, has attracted wide attention. In DRO, uncertain data are described by probability distribution functions which are not known exactly and restricted in a functional ambiguity set constructed from available information and structured properties. The expected cost associated with the worst-case distribution is minimized, and the probability of constraint violations can be controlled via robust chance constraints. In many cases, the DRO can be reformulated as a convex optimization problem, or solved iteratively via convex optimization. RO and DRO approaches are young and active research fields, and the challenge is to explore tractable reformulations with various kinds of uncertainties. SO is a relatively mature technique, and the current research is focusing on probabilistic modeling of uncertainty, chance constrained programming, multi-stage SO such as stochastic dual dynamic programming, as well as more efficient computational methods.

There are several ways to categorize robust optimization methods. According to how uncertainty is dealt with, they can be classified into static (single-stage) RO and dynamic (multi-stage) RO. According to how uncertainty is modeled, they can be divided into RO and DRO. In the latter category, the ambiguity set for probability distribution can be further classified into the moment based one and the divergence based one. We will shed light on each of them in this chapter. Specifically, RO will be discussed in Sect. \ref{App-C-Sect01} and Sect. \ref{App-C-Sect02}, moment-based DRO will be presented in Sect. \ref{App-C-Sect03}, and divergence-based DRO, also called  robust SO will be illuminated in Sect. \ref{App-C-Sect04}. In the operations research community, DRO and robust SO refer to the same thing: optimization problem with distributional uncertainty, and can be used interchangeably, although DRO is preferred by the majority of researchers. In this book, we intentionally distinguish them because the moment ambiguity set can be set up with little information and is more likely a RO;  the divergence based set relies on an empirical distribution (may be inexact), so is more similar to an SO. In fact, the gap between SO and RO has been significantly narrowed by recent research progress in the sense of data-driven optimization.

\section{Static Robust Optimization}
\label{App-C-Sect01}

For the purpose of clarity, we begin to explain the paradigm of static RO from LPs, the best known and most frequently used mathematical programming problem in engineering applications. It is relatively easy to derive tractable robust counterparts with various uncertainty sets. Nevertheless, most results can be readily generalized to robust conic programs. The general form of an LP with uncertain parameters can be written as follows:
\begin{equation}
\label{eq:App-03-SRO-ULP}
\min_x \left\{ c^T x ~\middle|~ Ax \le b \right\}:(A,b,c) \in W
\end{equation}
where $x$ is the decision variable, $A$, $b$, $c$ are coefficient matrices with compatible dimensions, and $W$ denotes the set of all possible data realizations constructed from available information or historical data, or merely a rough estimation. 

Without loss of generality, we can assume that the objective function and the constraint right-hand side in  (\ref{eq:App-03-SRO-ULP}) are certain, and uncertainty only exists in coefficient matrix $A$. To see this,  it is not difficult to observe that  problem (\ref{eq:App-03-SRO-ULP}) can be written as an epigraph form 
\begin{equation*}
\min_{t,x,y} \{t~|~c^T x - t \le 0,~ Ax-by \le 0,~ y=1 \}: (A,b,c) \in W
\end{equation*}
By introducing additional scalar variables $t$ and $y$, coefficients appearing in the objective function and constraint right-hand side are constants.
With this transformation, it will be more convenient to define the feasible solution and the optimal solution to (\ref{eq:App-03-SRO-ULP}). Hereinafter, we neglect the uncertainty in cost coefficient vector $c$ and constraint right-hand vector $b$ without particular mention, and consider problem
\begin{equation}
\label{eq:App-03-SRO-LP-UA}
\min_x \left\{ c^T x ~\middle|~ Ax \le b \right\}: A \in W
\end{equation}

 Next we present solution concepts of static RO under uncertain data.

\subsection{Basic Assumptions and Formulations}
\label{App-C-Sect01-01}
Basic assumptions and definitions in static RO \cite{RO-Detail-1} are summarized as follows.

\begin{assumption}
\label{ap:App-03-SRO-1}
Vector $x$ represents ``here-and-now'' decisions: they should be determined without knowing exact values of uncertain parameters.
\end{assumption}

\begin{assumption}
\label{ap:App-03-SRO-2}
Once the decisions are made, constraints must be feasible when the actual data is within the uncertainty set $W$, and may be either feasible or not when the actual data step outside the uncertainty set $W$. 
\end{assumption}

These assumptions bring about the definition for a feasible solution of (\ref{eq:App-03-SRO-LP-UA}).

\begin{definition}
\label{df:App-03-SRO-Feasibility}
A vector $x$ is called a robust feasible solution to (\ref{eq:App-03-SRO-LP-UA}) if the following condition holds: 
\begin{equation}
\label{eq:App-03-SRO-Robust-Fea}
A x \le b,~ \forall A \in W 
\end{equation}
\end{definition}

To prescribe an optimal solution, the worst-case criterion is widely accepted in RO studies, leading to the following definition:
\begin{definition}
\label{df:App-03-SRO-Optimality}
The robust optimal value of (\ref{eq:App-03-SRO-LP-UA}) is the minimum value of the objective function over all possible $x$ that satisfies (\ref{eq:App-03-SRO-Robust-Fea}).
\end{definition}

After we have agreed on the meanings of feasibility and optimality of  (\ref{eq:App-03-SRO-LP-UA}), we can seek the optimal solution among all robust feasible solutions to the problem. Now, the robust counterpart (RC) of the uncertain LP (\ref{eq:App-03-SRO-LP-UA}) can be described as:
\begin{equation}
\label{eq:App-03-SRO-RC}
\begin{aligned}
\min_x ~~ & c^T x \\
\mbox{s.t.} ~~ & a^T_i x \le b_i,~ \forall i,~ \forall A \in W
\end{aligned} 
\end{equation}
where $a^T_i$ is the $i$-th row of matrix $A$, and $b_i$ is the $i$-th element of vector $b$. We have two observations on the formulation of robust constraints in (\ref{eq:App-03-SRO-RC}).

\begin{proposition}
\label{pr:App-03-Invariance-1}
Robust feasible solutions of (\ref{eq:App-03-SRO-RC}) remain the same if we replace $W$ with the Cartesian product $\hat W = W_1 \times \cdots \times W_n$, where $W_i = \{a_i | \exists A \in W \}$ is the projection of $W$ on the coefficient space of $i$-th row of $A$.  
\end{proposition}

This is called the constraint-wise property in static RO \cite{RO-Detail-1}. The reason is 
\begin{equation}
a^T_i x \le b_i,~ \forall A \in W~ 
\Leftrightarrow~  \max_{A \in W} a^T_i x \le b_i  
\Leftrightarrow~  \max_{a_i \in W_i} a^T_i x \le b_i
\notag  
\end{equation}
As a result, problem (\ref{eq:App-03-SRO-RC}) comes down to
\begin{equation}
\label{eq:App-03-SRO-RC-Hull}
\begin{aligned}
\min_x ~~ & c^T x \\
\mbox{s.t.} ~~ & a^T_i x \le b_i,\forall a_i \in W_i,~ \forall i
\end{aligned} 
\end{equation}

Proposition \ref{pr:App-03-Invariance-1} seems rather counter-intuitive. One may perceive that (\ref{eq:App-03-SRO-RC}) will be less conservative with uncertainty set $W$ since it is a subset of $\hat W$. In fact, later we will see that this intuition is true for adjustable robustness.

\begin{proposition}
\label{pr:App-03-Invariance-2}
Robust feasible solutions of (\ref{eq:App-03-SRO-RC-Hull}) remain the same if we replace $W_i$ with its convex hull ${\rm conv}(W_i)$. 
\end{proposition}

To see this, let vector $a^j_i$, $j=1,2,\cdots$ be the extreme points of $W_i$, then any point $\bar a_i \in$ conv$(W_i)$ can be expressed by $\bar a_i = \sum_j \lambda_j a^j_i$, where $\lambda_j \ge 0$, $\sum_j \lambda_j = 1$ are weight coefficients. If $x$ is feasible for all extreme points $a^j_i$, i.e., $a^j_i x \le b_i$, $\forall j$, then  
\begin{equation}
\bar a^T_i x = \sum_j \lambda_j a^j_i x \le \sum_j \lambda_j b_i = b_i \notag
\end{equation}
which indicates that the constraint remains intact for all uncertain parameters reside in conv$(W_i)$.

Combining Propositions \ref{pr:App-03-Invariance-1} and \ref{pr:App-03-Invariance-2}, we can conclude that the robust counterpart of an uncertain LP with a certain objective remains intact even if sets $W_i$ of uncertain data are extended to their closed convex hulls, and $W$ to the Cartesian product of the resulting sets. In other words, we can make a further assumption on the uncertainty set without loss of generality.

\begin{assumption}
\label{ap:App-03-SRO-3}
The uncertainty set $W$ is the Cartesian product of closed and convex sets.
\end{assumption}

\subsection{Tractable Reformulations}
\label{App-C-Sect01-02}

The constraint-wise property enables us to analyze the robustness of each constraint $a^T_i x \le b_i$, $\forall a_i \in W_i$ separately. Without particular mention, we will omit the subscript $i$ for brevity. To facilitate discussion, it is convenient to parameterize the uncertain vector as $a = \bar a + P \zeta$, where $\bar a$ is the nominal value of $a$, $P$ is a constant matrix, $\zeta$ is a new variable that is uncertain. This section will focus on how to derive tractable reformulation for robust constraints in the form of
\begin{equation}
\label{eq:App-03-SRO-RC-Single}
(\bar a + P \zeta)^T x \le b,~ \forall \zeta \in Z
\end{equation}
where $Z$ is the uncertainty set of variable $\zeta$. For same reasons, we can assume that $Z$ is closed and convex. A “computationally tractable” problem means that there are known solution algorithms which can solve the problem with polynomial running time in its input size even in the worst case. It has been shown in \cite{RO-Detail-1} that problem (\ref{eq:App-03-SRO-RC-Hull}) is generally intractable even if each $W_i$ is closed and convex. Nevertheless, tractability can be preserved for some special classes of uncertainty sets. Some well-known results are summarized in the following. 

Condition (\ref{eq:App-03-SRO-RC-Single}) contains an infinite number of constraints due to the enumeration over set $Z$. Later we will see that for some particular uncertainty sets, the $\forall$ quantifier as well as the uncertain parameter $\zeta$ can be eliminated by using duality theory, and the resulting constraint in variable $x$ is still convex.

\vspace{12pt}
{\noindent \bf 1. Polyhedral uncertainty set} 

We start with a commonly used uncertainty set: a polyhedron
\begin{equation}
\label{eq:App-03-SRO-US-Polyhedron}
Z = \{ \zeta ~|~ D \zeta + q \ge 0\}
\end{equation}
where $D$ and $q$ are constant matrices with compatible dimensions. 

To exclude the $\forall$ quantifier for variable $\zeta$, we investigate the worst case of the left-hand side and require
\begin{equation}
\label{eq:App-03-SRO-Poly-1}
\bar a^T x + \max_{\zeta \in Z} (P^T x)^T \zeta \le b
\end{equation}

For a fixed $x$, the second term is the optimum of an LP in variable $\zeta$. Duality theory of LP says that the following relation holds
\begin{equation}
\label{eq:App-03-SRO-Poly-2}
(P^T x)^T \zeta \le q^T u,~ \forall \zeta \in Z,~ \forall u \in U
\end{equation}
where $u$ is the dual variable, and $U = \{ u ~|~ D^T u + P^T x =0,~ u \ge 0\}$ is the feasible region of the dual problem. Please be cautious on the sign of $u$. We actually replace $u$ with $-u$ in the original dual LP. Therefore, a necessary condition to validate (\ref{eq:App-03-SRO-RC-Single}) is 
\begin{equation}
\label{eq:App-03-SRO-Poly-3}
\exists u \in U: \bar a^T x + q^T u \le b
\end{equation}
It is also sufficient if the second term takes its minimum value over $U$, because strong duality always holds for LPs, i.e. $(P^T x)^T \zeta = q^T u$ is satisfied at the optimal solution. In this regard, (\ref{eq:App-03-SRO-Poly-1}) is equivalent to
\begin{equation}
\label{eq:App-03-SRO-Poly-4}
\bar a^T x + \min_{u \in U}~ q^T u \le b
\end{equation}
In fact, the ``min'' operator in (\ref{eq:App-03-SRO-Poly-4}) can be omitted in a RC optimization problem that minimizes the objective function, and thus renders polyhedral constraints, although (\ref{eq:App-03-SRO-Poly-1}) is not given in a closed form and seems non-convex. 

In summary, the RC problem of an uncertain LP with polyhedral uncertainty
\begin{equation}
\label{eq:App-03-SRO-RC-LP-1}
\begin{aligned}
\min_x ~~ & c^T x \\
\mbox{s.t.} ~~ & 
(\bar a_i + P_i \zeta_i)^T x \le b_i,~ \forall \zeta_i \in Z_i,~\forall i\\
& Z_i = \{ \zeta_i ~|~ D_i \zeta_i + q_i \ge 0\},~ \forall i
\end{aligned} 
\end{equation}
can be equivalently formulated as
\begin{equation}
\label{eq:App-03-SRO-RC-LP-2}
\begin{aligned}
\min_x ~~ & c^T x \\
\mbox{s.t.} ~~ & 
\bar a^T_i x +  q^T_i u_i \le b_i,~\forall i \\
& D^T_i u_i + P^T_i x = 0,~ u_i \ge 0,~\forall i
\end{aligned} 
\end{equation}
which is still an LP.

\vspace{12pt}
{\noindent \bf 2. Cardinality constrained uncertainty set} 

Cardinality constrained uncertainty set is a special class of
polyhedral uncertainty set which incorporates a budget constraint and defined as follows
\begin{equation}
\label{eq:App-03-SRO-US-Card}
Z({\rm \Gamma}) = \left\{ \zeta ~\middle|~ -1 \le \zeta_j \le 1,~ \forall j,~ \sum_j |\zeta_j| \le \rm \Gamma \right \}
\end{equation}
where $\rm \Gamma$ is called the budget of uncertainty \cite{RO-Price-Robust}. Motivated by the fact that each entry $\zeta_j$ is unlikely to reach 1 or $-1$ at the same time, the budget constraint controls the total data deviation from their forecast values. In other words, the decision maker can achieve a compromise between the level of solution robustness and the optimal cost by adjusting the value of $\rm \Gamma$, which should be less than the dimension of $\zeta$, otherwise the the budget constraint will be redundant.

Although the cardinality constrained uncertainty set $Z({\rm \Gamma})$ is essentially a polyhedron, the number of its facets, or the number of linear constraints in (\ref{eq:App-03-SRO-US-Polyhedron}), grows exponentially in the dimension of $\zeta$, leading to a huge and dense coefficient matrix for the uncertainty set. To circumvent this difficulty, we can lift it into a higher dimensional space as follows by introducing auxiliary variables
\begin{equation}
\label{eq:App-03-SRO-US-Card-Lift}
Z({\rm \Gamma}) = \left\{ \zeta,\sigma ~\middle|~ - \sigma_j \le \zeta_j \le \sigma_j,~ \sigma_j \le 1,~\forall j,~ \sum_j \sigma_j \le \rm \Gamma \right \}
\end{equation}
The first inequality naturally suggests $\sigma_j \ge 0$, $\forall j$.
It is easy to see the equivalence of (\ref{eq:App-03-SRO-US-Card}) and (\ref{eq:App-03-SRO-US-Card-Lift}), and the numbers of variables and constraints in the latter one grows linearly in the dimension of $\zeta$. 

Following a similar paradigm, certifying constraint robustness with a cardinality constrained uncertainty set requires the optimal value function of the following LP in variables $\zeta$ and $\sigma$ representing the uncertainty   
\begin{equation}
\label{eq:App-03-SRO-Card-1}
\begin{aligned}
\max_{\zeta,\sigma}~~  (P^T & x)^T \zeta \\
\mbox{s.t.}~~ 
- \zeta_j - \sigma_j & \le 0,~ \forall j : u^n_j \\
  \zeta_j - \sigma_j & \le 0,~ \forall j : u^m_j \\
            \sigma_j & \le 1,~ \forall j : u^b_j \\
     \sum_j \sigma_j & \le {\rm \Gamma}  : u_r
\end{aligned} 
\end{equation}
where $u^n_j$, $u^m_j$, $u^b_j$, $\forall j$, and $u_r$ following a colon are the dual variables associated with each constraint. The dual problem of (\ref{eq:App-03-SRO-Card-1}) is given by
\begin{equation}
\label{eq:App-03-SRO-Card-2}
\begin{aligned}
\min_{u^n,u^m,u^b,u_r}~~ & u_r {\rm \Gamma} + \sum_j u^b_j \\
\mbox{s.t.}~~ & u^m_j - u^n_j = (P^T x)_j,~ \forall j \\
       & -u^m_j - u^n_j + u^b_j + u_r = 0,~ \forall j \\
       & u^m_j,~ u^n_j,~ u^b_j \ge 0,~ \forall j,~ u_r \ge 0 
\end{aligned} 
\end{equation}

In summary, the RC problem of an uncertain LP with cardinality constrained uncertainty
\begin{equation}
\label{eq:App-03-SRO-RC-Card-1}
\begin{aligned}
\min_x ~~ & c^T x \\
\mbox{s.t.} ~~ & 
(\bar a_i + P_i \zeta_i)^T x \le b_i,~ \forall \zeta_i \in Z_i({\rm \Gamma}_i),~\forall i
\end{aligned} 
\end{equation}
can be equivalently formulated as
\begin{equation}
\label{eq:App-03-SRO-RC-Card-2}
\begin{aligned}
\min_x ~~ & c^T x \\
\mbox{s.t.} ~~ & \bar a^T_i x +  u_{ri} {\rm \Gamma}_i + \sum_j u^b_{ij} \le b_i,~\forall i \\
 & u^m_{ij} - u^n_{ij} = (P^T_i x)_j,~ \forall i,~ \forall j  \\
 &-u^m_{ij} - u^n_{ij} + u^b_{ij} + u_{ir} = 0,~ \forall i,~ \forall j \\
 & u^m_{ij},~ u^n_{ij},~ u^b_{ij} \ge 0,~ \forall i,~ \forall j,~
   u_{ir} \ge 0,~ \forall i 
\end{aligned} 
\end{equation}
which is still an LP.

\vspace{12pt}
{\noindent \bf 3. Several other uncertainty sets} 

Equivalent convex formulations of the uncertain constraint (\ref{eq:App-03-SRO-RC-Single}) with some other uncertainty sets are summarized in Table \ref{tab:App-03-SRO-RCs} \cite{RO-Guide}. These outcomes are derived using the similar method described previously.

\begin{table}[!t]
\scriptsize
\renewcommand{\arraystretch}{2.0}
\renewcommand{\tabcolsep}{1em}
\caption{Equivalent convex formulations with different uncertainty sets}
\centering
\begin{tabular}{|c|c|c|c|}
\hline
Uncertainty  &   $Z$    &    Robust reformulation &  Tractability  \\ 
\hline
Box & $\|\zeta\|_\infty \le 1$ & $\bar a^T x + \| P^T x \|_1 \le b$ &  LP\\
Ellipsoidal & $\|\zeta\|_2 \le 1$ & $\bar a^T x+\| P^T x \|_2 \le b$ & LP\\
$p$-norm & $\|\zeta\|_p \le 1$ & $\bar a^T x + \| P^T x \|_q \le b$ &  Convex program \\
Proper cone & $D \zeta + q \in K$ & $ \left\{ \begin{lgathered}
\bar a^T x + q^T u \le b \\ D^T u + P^T x = 0 \\ u \in K^* \end{lgathered} \right.$ & Conic LP \\
Convex constraints & $h_k(\zeta) \le 0, \forall k$ & $ \left\{ \begin{lgathered}
\bar a^T x + \sum_k \lambda_k h^*_k \left( \frac{u^k}{\lambda_k} \right) \le b \\ \sum_k u^k = P^T x \\ \lambda_k \ge 0, \forall k \end{lgathered} \right.$ & Convex program \\
\hline
\end{tabular}
\label{tab:App-03-SRO-RCs}
\end{table}

Table \ref{tab:App-03-SRO-RCs} includes three cases: the p-norm uncertainty, the conic uncertainty, and general convex uncertainty. 
In the $p$-norm case, the H{\"o}lder's inequality is used, i.e.: 
\begin{equation}
(P^T x)^T \zeta \le \|P^T x\|_p \| \zeta \|_q 
\end{equation}
where $\|\cdot\|_p$ and $\|\cdot\|_q$ with $p^{-1} + q^{-1} = 1$ are a pair of dual norms. Since norm function of any order is convex \cite{CVX-Book-Boyd}, the resulting RC is a convex program. Moreover, if $q$ is a  positive rational number, the $q$-order cone constraints can be represented by a set of SOC inequalities \cite{SOCP-p-norm}, which is computationally more friendly. Box ($\infty$-norm) and ellipsoidal (2-norm) uncertainty sets are special kinds of $p$-norm ones.

In the general conic case, conic duality theory \cite{CVX-Book-Ben} is used. $K^*$ stands for the dual cone of $K$, and the polyhedral uncertainty is a special kind of this case when $K$ is the nonnegative orthant. 

In the general convex case, Fenchel duality, a basic theory in convex analysis, is needed. Notation $h^*$ stands for the convex conjugate function, i.e. $h^*(x) = \sup_y x^T y -h(y)$. The detailed proof of RC reformulations and more examples can be found in \cite{SRO-CVX-RCs}.

Above analysis focuses on the situation in which problem functions are linear in decision variables, and problem\ data are affine in some uncertain parameters, such as the form $a = \bar a + P \zeta$. For robust quadratic optimization, robust semidefinite optimization, robust conic optimization, and robust discrete optimization, in which the optimization problem is nonlinear and discontinuous, please refer to \cite{RO-Detail-1} and \cite{RO-Detail-2}; for quadratic type uncertainty, please refer to \cite{RO-Detail-1} (in Sect. 1.4) and \cite{SRO-CVX-RCs}.

\subsection{Formulation Issues}
\label{App-C-Sect01-03}

To help practitioners build a well-defined and easy-to-solve robust optimization model, some important modeling issues and deeper insights are discussed in this section. 

\vspace{12pt}
{\noindent \bf 1. Choosing the uncertainty set} 

Since a robust solution remains feasible if the uncertain data does not step outside the uncertainty set, the level of robustness mainly depends on the shape and size of the uncertainty set. The more reliable, the higher the cost. One may wish to seek a trade-off between reliability and economy. This inspires the development of smaller uncertainty sets with a certain probability guarantee that the constraint violation is unlikely to happen. Such guarantees are usually described via a chance constraint 
\begin{equation}
\label{eq:App-03-SRO-Chance-Cons}
\Pr\nolimits _\zeta [a(\zeta)^T x \le b] \ge 1 - \varepsilon
\end{equation}
For $\varepsilon = 0$, chance constraint (\ref{eq:App-03-SRO-Chance-Cons}) is protected in the traditional sense of RO. When $\varepsilon >0$, it becomes challenging to derive tractable reformulation for (\ref{eq:App-03-SRO-Chance-Cons}), especially when the probability distribution of uncertain data is unclear or inaccurate. In fact, this issue is closely related to the DRO that will be discussed later on. Here we provide some simple results which help the decision maker choose the parameter of the uncertainty set. 

It is revealed that if $\mathbb E [\zeta] = 0$, the components of $\zeta$ are independent, and the uncertainty set takes the form
\begin{equation}
\label{eq:App-03-SRO-US-BOX-ELP}
Z = \{\zeta ~|~ \|\zeta\|_2 \le {\rm \Omega},~ \|\zeta\|_\infty \le 1 \}
\end{equation}
then chance constraint (\ref{eq:App-03-SRO-Chance-Cons}) holds with a probability of at least $1-\exp(-{\rm \Omega^2}/2)$ (see  \cite{RO-Detail-1}, Proposition 2.3.3).

Moreover, if the uncertainty set takes the form 
\begin{equation}
\label{eq:App-03-SRO-US-BOX-Budget}
Z = \{\zeta ~|~ \|\zeta\|_1 \le {\rm \Gamma},~ \|\zeta\|_\infty \le 1 \}
\end{equation}
then chance constraint (\ref{eq:App-03-SRO-Chance-Cons}) holds with a probability of at least $1-\exp(-{\rm \Gamma^2}/2L)$, where $L$ is the dimension of $\zeta$ (see  \cite{RO-Detail-1}, Proposition 2.3.4, and \cite{RO-Price-Robust}).

It is proposed to construct uncertainty sets based on the
central limit theorem. If each component of $\zeta$ is independent and identically distributed with mean $\mu$ and variance $\sigma^2$, the uncertainty set can be built as \cite{SRO-US-CLT}
\begin{equation}
\label{eq:App-03-SRO-US-CLT}
Z = \left\{ \zeta ~\middle|~ \left|\sum_{i=1}^L \zeta_i - L \mu \right| \le \rho \sqrt{L} \sigma,~ \|\zeta\|_\infty \le 1 \right\}
\end{equation}
where parameter $\rho$ is used to control the probability guarantee. Variations of this formulation can take other distributional information into account, such as data correlation and long tail-effect. It is a special kind of polyhedral uncertainty, however, it is unbounded for $L > 1$, since the components can be arbitrarily large as long as their summation is relatively small. Unboundedness may prevents establishing tractable RCs. 

Additional references are introduced in further reading.

\vspace{12pt}
{\noindent \bf 2. How to solve a problem without a clear tractable reformulation?}

The existence of a tractable reformulation for a static RO problem largely depends on the type of the uncertainty set. If the robust counterpart cannot be
written as a tractable convex program, a smart remedy is to use an adaptive scenario generation procedure: first solve the problem with a smaller uncertainty set $Z_S$ which is a subset of the original one $Z$, and the problem with $Z_S$ has a known tractable reformulation. If the optimal solution $x^*$ is robust against all scenarios in $Z$, it is also an optimal solution of the original problem. Otherwise, we have to identify a scenario $\zeta^* \in Z$ which leads to the most severe violation, which can be implemented by solving
\begin{equation}
\label{eq:App-03-SRO-Scen-Gen-Sub}
\max~ \left\{ (P^T x^*)^T \zeta ~|~{\zeta \in Z} \right\}
\end{equation}
where $Z$ is a closed and convex set as validated in Assumption \ref{ap:App-03-SRO-3}, and then append a cutting plane
\begin{equation}
\label{eq:App-03-SRO-Scenario-Cut}
a(\zeta^*)^T x \le b
\end{equation}
to the  reformulation problem. (\ref{eq:App-03-SRO-Scenario-Cut}) removes $x$ that will cause infeasibility in scenario $\zeta^*$, so is called a feasibility cut. It is linear and does  not alter tractability. Then the updated problem is solved again. According to Proposition \ref{pr:App-03-Invariance-2}, the new solution $x^*$ will be robust for uncertain data in the convex hull of $Z_S \cup \zeta^*$. Above procedure continues until robustness is certified over the original uncertainty set $Z$.

This simple approach often converges quickly in a few number of iterations. Its advantage is that tractability is preserved. When we choose $Z_S = \zeta^0$, where $\zeta^0$ is the nominal scenario or forecast, it could be more efficient than using convex reformulations, because only LPs (whose sizes are almost equal to the problem without uncertainty, and grows slowly) and simple convex programs (\ref{eq:App-03-SRO-Scen-Gen-Sub}) are solved, see \cite{SRO-Cut-Generation} for a comparison. This paradigm is an essential strategy for solving the adjustable RO problems in the next section.

\vspace{12pt}
{\noindent \bf 3. How to deal with equality constraints?} 

Although the theory of static RO is relatively mature, it encounters difficulties in dealing with equality constraints. For example, consider $x + a = 1$ where $a \in [0,0.1]$ is uncertain. However, one can seldom find a solution that makes the equality hold true for multiple values of $a$. The problem remains if you write a equality into a pair of opposite inequalities. In fact, this issue is inevitable in the static setting. In addition, this limitation will lead to completely different robust counterpart formulations for originally equivalent deterministic problems. 

Consider the inequality $a x \le 1$, which is equivalent to $a x + s = 1,~ s \ge 0$. Suppose $a$ is uncertain and belongs to interval $[1,2]$, their respective robust counterparts are given by
\begin{equation}
\label{eq:App-03-SRO-Example-1}
a x \le 1,~ \forall a \in [1,2]
\end{equation}
and
\begin{equation}
\label{eq:App-03-SRO-Example-2}
a x + s = 1,~ \forall a \in [1,2],~  s \ge 0 
\end{equation}
The feasible set for (\ref{eq:App-03-SRO-Example-1}) is $x \le 1/2$, and is $x = 0$ for (\ref{eq:App-03-SRO-Example-2}). By observing this difference,  it is suggested that a static RO model should avoid using slack variables in constraints with uncertain parameters.

Sometimes, the optimization problem may contain state variables which can respond to parameter changes by adjusting their values. In such circumstance, equality constraint can be used to eliminate state variables. Nevertheless, such an action may lead to a problem that contains nonlinear uncertainties, which are challenging to solve. An example is taken from \cite{RO-Guide} to illustrate this issue. The constraints are    
\begin{equation}
\label{eq:App-03-SRO-Example-3}
\begin{lgathered}
\zeta_1 x_1 + x_2 + x_3 = 1  \\
x_1 + x_2 + \zeta_2 x_3 \le 5
\end{lgathered}
\end{equation}
where $\zeta_1$ and $\zeta_2$ are uncertain. 

If $x_1$ is a state variable and $\zeta_1 \ne 0$, substituting $x_1 = (1 - x_2 - x_3)/\zeta_1$ in the second inequality results in
\begin{equation}
\left( 1-\frac{1}{\zeta_1} \right) x_2 + \left( \zeta_2 - \frac{1}{\zeta_1} \right) x_3 \le 5 - \frac{1}{\zeta_1}  \notag
\end{equation}
in which the uncertainty becomes nonlinear in the coefficients.

If $x_2$ is a state variable, substituting $x_2 = 1 - \zeta_1 x_1 - x_3$ in the inequality yields
\begin{equation}
(1 - \zeta_1) x_1 + (\zeta_2 -1) x_3 \le 4 \notag
\end{equation}
in which the uncertainty sustains linear in the coefficients.

If $x_3$ is a state variable, substituting $x_3 = 1 - \zeta_1 x_1 - x_2$ in the inequality gives
\begin{equation}
(1 - \zeta_1 \zeta_2) x_1 + (1 - \zeta_2) x_2 \le 5 - \zeta_2 \notag
\end{equation}
in which the uncertainty is  nonlinear in the coefficients.

In conclusion, in the case that $x_2$ is a state variable, the problem is easier from a computational perspective. It is important to note that the physical interpretation of variable elimination is to determine the adjustable variable with exact information on the uncertain data. If no adjustment is allowed in (\ref{eq:App-03-SRO-Example-3}), the only robust feasible solution is $x_1 = x_3 =0$, $x_2 = 1$, which is rather restrictive. The adjustable RO will be elaborated in detail in the next section.

\vspace{12pt}
{\noindent \bf 4. Pareto efficiency of the robust solution}

The concept of Pareto efficiency in RO problems is proposed in \cite{SRO-Pareto-1}. If the optimal solution under the worst-case data realization is not unique, it is rational to compare their performances in non-worst-case scenarios: an alternative solution may give an improvement in the objective value for at least one data scenario without deteriorating the objective performances in all other scenarios. To present related concept tersely, we restrict the discussion on the following robust LP with objective uncertainty
\begin{equation}
\label{eq:App-03-RLP-Obj}
\max\{p^T x~|~ \mbox{s.t. } x \in X,~ \forall p \in W \} =
\max_{x \in X} \left\{ \min_{p \in W} p^T x \right\}
\end{equation}
where $W = \{p ~|~ D p \ge d \}$ is a polyhedral uncertainty set for the price vector $p$; $X = \{ x~|~A x \le b \}$ is the feasible region which is independent of the uncertainty. More general cases are elaborated in \cite{SRO-Pareto-1}. We consider this form because it is easy to discuss related issues, although objective uncertainty can be moved into constraints.

For a given strategy $x$, the worst-case uncertainty is
\begin{equation} 
\label{eq:App-03-RLP-Dual}
\min \{ p^T x ~|~ \mbox{s.t.}~ p \in W \} = 
\max \{ d^T y ~|~ \mbox{s.t.}~ y \in Y \} 
\end{equation}
where $Y =\{ y ~|~ D^T y = x,~ y \ge 0\}$ is the feasible set for dual variable $y$. Substituting (\ref{eq:App-03-RLP-Dual}) in (\ref{eq:App-03-RLP-Obj}) gives 
\begin{equation}
\label{eq:App-03-RLP-Obj-RC}
\max\{d^T y~|~ \mbox{s.t. } D^T y = x,~ y \ge 0,~x \in X\}
\end{equation}
which is an LP. Its solution $x$ is the robust optimal one to (\ref{eq:App-03-RLP-Obj}), and the worst-case price $p$ can be found by solving the left-hand side LP in (\ref{eq:App-03-RLP-Dual}). Let $z^{RO}$ be the optimal value of (\ref{eq:App-03-RLP-Obj-RC}), and then the set of robust optimal solutions for (\ref{eq:App-03-RLP-Obj}) can be expressed via
\begin{equation}
\label{eq:App-03-RLP-Obj-XRO}
X^{RO} = \{x~|~ x \in X: \exists y \in Y 
\mbox{ such that } y^T d \ge z^{RO} \} 
\end{equation}

If (\ref{eq:App-03-RLP-Obj-RC}) has a unique optimal solution, $X^{RO}$ is a singleton; otherwise, a Pareto optimal robust solution can be formally defined. 

\begin{definition}
\label{df:App-03-Pareto-Efficiency}
\cite{SRO-Pareto-1} $x \in X^{RO}$ is a Pareto optimal solution for problem (\ref{eq:App-03-RLP-Obj}) if there is no other $\bar x \in X$ such that $p^T \bar x \ge p^T x,~ \forall p \in W$ and $\bar p^T \bar x > \bar p^T x$ for some $\bar p \in W$.
\end{definition}

The terminology ``Pareto optimal'' is borrowed from multi-objective optimization theory: RO problem (\ref{eq:App-03-RLP-Obj}) is viewed as a multi-objective LP with infinitely many objectives, each of which corresponds to a particular $p \in W$. Some interesting problems are elaborated. 

\vspace{12pt}
{\noindent \bf a. Pareto efficiency test}

In general, it is not clear whether $X^{RO}$ contains multiple solutions, at least before a solution $x \in X^{RO}$ is found. To test whether a given solution $x$ is a robust optimal one or not, it is proposed to solve a new LP
\begin{equation}
\label{eq:App-03-PRO-Test-LP}
\begin{aligned}
\max_{y}  ~~   &  \bar p^T y   \\
\mbox{s.t.} ~~ &  y \in W^*        \\
               &  x + y \in X 
\end{aligned}
\end{equation}
where $\bar p$ is a relative interior of the polyhedral uncertainty set $W$, which is usually set to the nominal scenario, and $W^* = \{ y~|~\exists \lambda:d^T \lambda \ge 0,~ D^T \lambda = y,~ \lambda \ge 0\}$ is the dual cone of $W$. Please refer to Sect. \ref{App-A-Sect02-01} and equation (\ref{eq:App-01-Dual-Polytope-2}) for the dual cone of a polyhedral set. Since $y = 0$, $\lambda = 0$ is always feasible in (\ref{eq:App-03-PRO-Test-LP}), the optimal value is either zero or strictly positive. In the former case, $x$ is also a Pareto optimal solution; in the latter case, $\bar x = x + y^*$ dominates $x$ and itself is Pareto optimal for any $y^*$ that solves LP (\ref{eq:App-03-PRO-Test-LP}) \cite{SRO-Pareto-1}. The interpretation of (\ref{eq:App-03-PRO-Test-LP}) is clear: since $y \in W^*$, $y^T p$ must be non-negative for all $p \in W$. If we can find $y$ that leads to a strict objective improvement for $\bar p$, then $x+y$ would be Pareto optimal. 

In view of the above interpretation, it is a direct conclusion that for an arbitrary relative interior point $\bar p \in W$, the optimal solutions to the problem
\begin{equation}
\label{eq:App-03-PRO-Find-LP}
\max~ \left\{ \bar p^T x ~|~ x \in X^{RO} \right\}
\end{equation}
are Pareto optimal.

\vspace{12pt}
{\noindent \bf b. Characterizing the set of Pareto optimal solutions}

It is interesting to characterize the Pareto optimal solution set $X^{PRO}$.

After we get $z^{RO}$ and $X^{RO}$, solve the following LP
\begin{equation}
\label{eq:App-03-XPRO}
\begin{aligned}
\max_{x,y,\lambda} ~~ &  \bar p^T y   \\
\mbox{s.t.} ~~ &  d^T \lambda \ge 0,~ D^T \lambda = y,~ \lambda \ge 0 \\
& x \in X^{RO},~ x + y \in X 
\end{aligned}
\end{equation}
and we can conclude $X^{PRO}=X^{RO}$ if and only if the optimal value of (\ref{eq:App-03-XPRO}) is equal to 0 \cite{SRO-Pareto-1}. If this is true, the decision maker would not have to worry about Pareto efficiency, as any solution in $X^{RO}$ is also Pareto optimal. More broadly, the set $X^{PRO}$ is shown to be non-convex and is contained in the boundary of $X^{RO}$.

\vspace{12pt}
{\noindent \bf c. Optimization over Pareto optimal solutions}

In the case that $X^{PRO}$ is not a singleton, one may consider to optimize a linear secondary objective over $X^{PRO}$, i.e.:
\begin{equation}
\label{eq:App-03-OPTI-XPRO}
\max \{ r^T x ~|~ \mbox{s.t. } x \in X^{PRO} \}  
\end{equation}

It is demonstrated in \cite{SRO-Pareto-1} that if $r$ lies in the relative interior of $W$, the decision maker can simply replace $X^{PRO}$ with $X^{RO}$ in (\ref{eq:App-03-OPTI-XPRO}) without altering the problem solution, due to the property revealed in (\ref{eq:App-03-PRO-Find-LP}). In more general cases, problem (\ref{eq:App-03-OPTI-XPRO}) can be formulated as an MILP \cite{SRO-Pareto-1}
\begin{equation}
\label{eq:App-03-OPTI-XPRO-MILP}
\begin{aligned}
\max_{x,\mu,\eta,z}~~ & r^T x   \\
\mbox{s.t.}~~ & x \in X^{RO}    \\
& \mu \le M(1-z)  \\
& b - A x \le Mz  \\
& DA^T \mu - d \eta \ge D \bar p \\
& \mu, \eta \ge 0, z \in \{0,1\}^{m}
\end{aligned}
\end{equation}
where $M$ is a sufficiently large number, $m$ is the dimension of vector $z$. To show their equivalence, it is revealed that the feasible set of (\ref{eq:App-03-OPTI-XPRO-MILP}) depicts an optimal solution of (\ref{eq:App-03-PRO-Test-LP}) with a zero objective value \cite{SRO-Pareto-1}. In other words, the constraints of (\ref{eq:App-03-OPTI-XPRO-MILP}) contain the KKT optimality condition of (\ref{eq:App-03-PRO-Test-LP}). To see this, the binary vector $z$ imposes the complementarity and slackness condition $\mu^T (b-Ax)=0$, which ensures $\lambda, \mu, \eta$ are the optimal solution of the following primal-dual
LP pair   
\begin{equation}
\mbox{Primal : }
\begin{aligned}
\max_{\lambda}~~ & \bar p^T D^T \lambda  \\
\mbox{s.t.}~~ & \lambda \ge 0 \\
              & d^T \lambda \ge 0 \\
              & A D^T \lambda \le b - A x 
\end{aligned}
\quad 
\mbox{Dual : }
\begin{aligned}
\min_{\mu, \eta}~~ & \mu^T (b-Ax)  \\
\mbox{s.t.}~~ & \mu \ge 0  \\
              & \eta \ge 0 \\
              & D A^T \mu - d \eta \ge D \bar p
\end{aligned}
\notag
\end{equation}
The original variable $y$ in (\ref{eq:App-03-PRO-Test-LP}) is eliminated via equality $ D^T \lambda = y$ in the dual cone $W^*$. According to strong duality, the optimal value of the primal LP (\ref{eq:App-03-PRO-Test-LP}) is $\bar p^T D^T \lambda = \bar p^T y = 0$, and Pareto optimality is guaranteed. 

In practice, Pareto inefficiency is not a contrived phenomenon, see various examples in \cite{SRO-Pareto-1} and power market examples in \cite{SRO-Pareto-2,SRO-Pareto-3}.

\vspace{12pt}
{\noindent \bf 5. On max-min and min-max formulations}

In many literatures, the robust counterpart problem of (\ref{eq:App-03-SRO-LP-UA}) is written as a min-max form
\begin{equation}
\label{eq:App-03-SRO-RC-min-max}
\begin{aligned}
\mbox{Opt-1} = \min_x \max_{A \in W} \left\{c^T x ~|~ \mbox{s.t. } A x \le b \right\}
\end{aligned}
\end{equation}
which means $x$ is determined before $A$ takes a value in $W$, and the decision maker can foresee the worst consequence of deploying $x$ brought by the perturbation of $A$. To make a prudent decision that is insensitive to data perturbation, the decision maker resorts to minimizing the maximal objective.    

The max-min formulation
\begin{equation}
\label{eq:App-03-SRO-RC-max-min}
\begin{aligned}
\mbox{Opt-2} = \max_{A \in W} \min_x \left\{c^T x ~|~ \mbox{s.t. } A x \le b \right\}
\end{aligned}
\end{equation}
has a different interpretation: the decision maker can first observe the realization of uncertainty, and then recovers the constraints by deploying a corrective action $x$ as a response to the observed $A$. Certainly, this specific $x$ may not be feasible for other $A \in W$. On the other hand, the uncertainty, like a rational player, can foresee the optimal action taken by the human decision maker, and select a strategy that will yield a maximal objective value even an optimal corrective action is deployed.   

From above analysis, the feasible region of $x$ in (\ref{eq:App-03-SRO-RC-min-max}) is a subset of that in (\ref{eq:App-03-SRO-RC-max-min}), because (\ref{eq:App-03-SRO-RC-max-min})
only accounts for a special scenario in $W$. As a result, their optimal values satisfy Opt-1 $\ge$ Opt-2.

Consider the following problem in which the uncertainty is not constraint-wise
\begin{equation}
\label{eq:App-03-SRO-RC-Toy}
\begin{aligned}
\min_x ~~ & x_1 + x_2  \\
\mbox{s.t.} ~~ & x_1 \ge a_1,~ x_2 \ge a_2,~ \forall a \in W 
\end{aligned}
\end{equation}
where $W = \{a ~|~ a \ge 0,~ \|a\|_2 \le 1 \}$. 

For the min-max formulation, since $x$  should be feasible for all possible values of $a$, it is necessary to require $x_1 \ge 1$ and $x_2 \ge 1$, and Opt-1 $=2$ for problem (\ref{eq:App-03-SRO-RC-Toy}). 

As for the max-min formulation, as $x$ is determined in response to the value of $a$, it is clear that the optimal choice is $x_1 = a_1$ and $x_2 = a_2$, so the problem becomes      
\begin{equation}
\begin{aligned}
\max_a ~~ & a_1 + a_2  \\
\mbox{s.t.} ~~ & a^2_1 + a^2_2 \le 1 
\end{aligned}   \notag
\end{equation}
whose optimal value is Opt-2 $=\sqrt{2} <$ Opt-1.

As a short conclusion, static RO models discussed in this section are used to immunize against constraint violation or objective volatility caused by data perturbations, without jeopardizing computational tractability. General approaches involve reformulating the original uncertainty dependent constraints into deterministic convex ones without uncertain data, such that feasible solutions of the robust counterpart program remain feasible for all data realizations in the pre-specified uncertainty set, which interprets the meaning of robustness.

\section{Adjustable Robust Optimization}
\label{App-C-Sect02}

Several reasons call for developing new decision-making mechanisms to overcome limitations of the static RO approach: 1) Equality constraints often give rise to infeasible robust counterpart problems in the static setting; 2) real-world decision-making process may involve multiple stages, in which some decisions indeed can be made after the uncertain data has been known or can be predicted accurately. Take power system operation for an example, the on-off status of generation units must be made several hours before real-time dispatch when the renewable power is unclear; however, the output of some units (called AGC units) can change in response to the real values of system demands and renewable generations. This section will be devoted to the adjustable robust optimization (ARO) with two stages, which leverages the adaptability in the second stage. We still focus our attention on the linear case.  

\subsection{Basic Assumptions and Formulations}
\label{App-C-Sect02-01}

The essential difference between static RO and ARO approaches stems from the manner of decision making.

\begin{assumption}
\label{ap:App-03-ARO-1}
In an ARO problem, some variables are ``here-and-now'' decisions, whereas the rest are ``wait-and-see'' decisions: they can be made at a later moment according to the observed data.
\end{assumption}

In analogy to the static case, the decision-making mechanism  can be explained.

\begin{assumption}
\label{ap:App-03-ARO-2}
Once the here-and-now decisions are made, there must be at least one valid wait-and-see decision which is able to recover constraints in response to the observed data realization, if the actual data is within the uncertainty set. 
\end{assumption}

In this regard, we can say here-and-now decisions are robust against the uncertainty, and wait-and-see decisions are adaptive to the uncertainty. These terminologies are borrowed from two-stage SO models. In fact, there is a close relation between two-stage SO and two-stage RO \cite{ARO-TSSO-Relation-1,ARO-TSSO-Relation-2}. 

Now we are ready to post the compact form of a linear ARO problem with an uncertain constraint right-hand side:
\begin{equation}
\label{eq:App-03-ARO}
\min_{x \in X} \left\{ c^T x + \max_{w \in W} \min_{y(w) \in Y(x,w)} d^T y(w)  \right\} 
\end{equation}
where $x$ is the here-and-now decision variable (or the first-stage decision variable), and $X$ is the feasible region of $x$; $w$ is the uncertain parameter, and $W$ is the uncertainty set, which has been discussed in the previous section; $y(w)$ is the wait-and-see decision variable (or second-stage decision variable), which can be adjusted according to the actual data of $w$, so it is represented as a function of $w$; $Y$ is the feasible region of $y$ given the values of $x$ and $w$, because the here-and-now decision is not allowed to change in this stage, and the exact value of $w$ is known. It has a polyhedral form
\begin{equation}
\label{eq:App-03-ARO-Y(x,w)}
Y(x,w) = \{y ~|~ Ax + By + Cw \le b \}
\end{equation}
where $A$, $B$, $C$, and $b$ are constant matrices and vector with compatible dimensions. It is clear that both of the here-and-now decision $x$ and the data uncertainty $w$ can influence the feasible region $Y$ in the second stage. We define $w = 0$ the nominal scenario and assume $0$ is a relative interior of $W$. Otherwise, we can decompose the uncertainty as $w = w^0 + \Delta w$ and merge the constant term $Cw^0$ into the right-hand side as $b \to b - C w^0$, where $w^0$ is the predicted or expected value of $w$, and ${\rm \Delta} w$ is the forecast error, which is the real uncertain parameter. 

It should be pointed out that $x$, $w$, and $y$ may contain discrete decision variables. Later we will see, integer  variables in $x$ and $w$ do not significantly alter the solution algorithm of ARO. However, because integrality in $y$ prevents  the use of LP duality theory, the computation will be greatly challenged. Although we assume coefficient matrices are constants in (\ref{eq:App-03-ARO-Y(x,w)}), most results in this section can be generalized if matrix $A$ is a linear function in $w$; the situation would be complicated in matrix $B$ is uncertainty-dependent. The purpose for the specific form in (\ref{eq:App-03-ARO-Y(x,w)}) is that it is more dedicated to the problems considered in this book: uncertainties originates from renewable/load volatility can be modeled by term $Cw$ in (\ref{eq:App-03-ARO-Y(x,w)}), and the coefficients representing component and network parameters are constants.

Assumption \ref{ap:App-03-ARO-2} inspires the definition for a feasible solution of ARO (\ref{eq:App-03-ARO}).

\begin{definition}
\label{df:App-03-ARO-XR}
A first-stage decision $x$ is called robust feasible in (\ref{eq:App-03-ARO}) if the feasible region $Y(x,w)$ is non-empty for all $w \in W$, and the set of  robust feasible solutions are given by: 
\begin{equation}
\label{eq:App-03-ARO-XR}
X_R = \{x ~|~ x \in X : \forall w \in W,~ Y(x,w) \ne \emptyset \}  
\end{equation}
\end{definition}

Please be aware of the sequence in (\ref{eq:App-03-ARO-XR}): $x$ takes its value first, and then parameter $w$ chooses a value in $W$ before some $y \in Y$ does. The non-emptiness of $Y$ is guaranteed by the selection of $x$ for an arbitrary $w$. If we swap the latter two terms and write $Y(x,w) \ne \emptyset$, $\forall w \in W$, like the form in a static RO, it sometimes cause confusion that both $x$ and $y$ are here-and-now type decisions, the adaptiveness vanishes, and thus $X_R$ may become empty if uncertainty appears in an equality constraint, as analyzed in the previous section.    

The definition of an optimal solution depends on the decision maker's attitude towards the cost in the second stage. In (\ref{eq:App-03-ARO}), we adopt the following definition.
\begin{definition}
\label{df:App-03-SRO-Optimality-MMC}
(Min-max cost criterion) An optimal solution of (\ref{eq:App-03-ARO}) is a pair of here-and-now decision $x \in X_R$ and wait-and-see decision $y(w^*)$ corresponding to the worst-case scenario $w^* \in W$, such that the total cost in scenario $w^*$ is minimal, where the worst-case scenario $w^*$ means that for the fixed $x$, the optimal second-stage cost is maximized over $W$.
\end{definition}

Other criteria may give different robust formulations. For example, the minimum nominal cost formulation and min-max regret formulation.

\begin{definition}
\label{df:App-03-SRO-Optimality-MNC}
(Minimum nominal cost criterion) An optimal solution under the minimum nominal cost criterion is a pair of here-and-now decision $x \in X_R$ and wait-and-see decision $y^0$ corresponding to the nominal scenario $w^0 = 0$, such that the total cost in scenario $w^0$ is minimal.
\end{definition}

The minimum nominal cost criterion leads to the following robust formulation   
\begin{equation}
\label{eq:App-03-ARO-V1}
\begin{aligned}
\min~~ & c^T x +  d^T y  \\
\mbox{s.t.} ~~ & x \in X_R \\
       & y \in Y(x,w^0)
\end{aligned}
\end{equation}
where robustness is guranteed by $X_R$. 

To explain the concept of regret, the minimum perfect-information
total cost is 
\begin{equation*}
C_P(w) = \min \left\{ c^T x + d^T y ~|~  x \in X,~  y \in Y(x,w) \right\}
\end{equation*}
where $w$ is known to the decision maker. For a fixed first-stage decision $x$, the maximum regret is defined as 
\begin{equation*}
\mbox{Reg}(x) = \max_{w \in W} \left\{ \min_{y \in Y(x,w)} \{c^T x + d^T y\} - C_P(w) \right\}
\end{equation*}

\begin{definition}
\label{df:App-03-SRO-Optimality-MMC}
(Min-max regret criterion) An optimal solution under the min-max regret criterion is a pair of here-and-now decision $x$ and wait-and-see decision $y(w)$, such that the worst-case regret under all possible scenarios $w \in W$ is minimized.
\end{definition}

The min-max regret cost criterion leads to the following robust formulation
\begin{equation}
\label{eq:App-03-ARO-MMR}
\min_{x \in X} \left\{ c^T x + \max_{w \in W} \left\{ \min_{y \in Y(x,w)} d^T y - \min_{x^\prime \in X,y^\prime \in Y(x,w)} \left\{ c^T x^\prime + d^T y^\prime \right\} \right\} \right\} \\
\end{equation}

\vspace{12pt}

In an ARO problem, we can naturally assume that the uncertainty set is a polyhedron. To see this,  if $x$ is a robust solution under an uncertainty set consists of discrete scenarios, i.e., $W=\{w^1,w^2,\cdots w^S\}$, according to Definition \ref{df:App-03-ARO-XR}, there exist corresponding $\{y^1,y^2,\cdots y^S\}$ such that 
\begin{equation}
\begin{gathered}
B y^1 \le b - A x -C w^1  \\
B y^2 \le b - A x -C w^2  \\
\vdots  \\
B y^S \le b - A x -C w^S
\end{gathered}
\notag
\end{equation}
For non-negative weighting parameters $\lambda_1, \lambda_2,\cdots, \lambda_S \ge 0$, $\sum_{s=1}^S \lambda_s =1$, we have 
\begin{equation}
\sum_{s=1}^S \lambda_s ( B y^s ) \le \sum_{s=1}^S \lambda_s (b-Ax-Cw^s) 
\notag
\end{equation}
or equivalently
\begin{equation}
B \sum_{s=1}^S \lambda_s y^s \le b - Ax - C \sum_{s=1}^S \lambda_s w^s 
\notag
\end{equation}
indicating that for any $w = \sum_{s=1}^S \lambda_s w^s \in \mbox{conv}(\{w^1,w^2,\cdots w^S\})$, the wait-and-see decision $y = \sum_{s=1}^S \lambda_s y^s$ can recover all constraints, and thus $Y(x,w) \ne \emptyset$. This property inspires the following proposition that is in analogy to Proposition \ref{pr:App-03-Invariance-2}

\begin{proposition}
\label{pr:App-03-Invariance-3}
Suppose $x$ is a robust feasible solution for a discrete uncertainty set $\{w^1,w^2,\cdots w^S\}$, then it remains robust feasible if we replace the uncertainty set with its convex hull. 
\end{proposition}

Proposition \ref{pr:App-03-Invariance-3} also implies that in order to ensure the robustness of $x$, it is sufficient to consider the extreme points of a bounded polytope. Suppose the vertices of the polyhedral uncertainty set are $w^s$, $s=1,2,\cdots,S$. Consider the following set   
\begin{equation}
\label{eq:App-03-ARO-XR-Lift}
{\rm \Xi} = \{x,y^1,y^2,\cdots,y^S ~|~ 
A x + B y^s \le b -C w^s,~ s = 1,2,\cdots,S \} 
\end{equation}
Robust feasible region $X_R$ is the projection of polyhedron $\rm \Xi$ on $x$-space, which is also a polyhedron (Theorem B.2.5 in \cite{CVX-Book-Ben}).

\begin{proposition}
\label{pr:App-03-XR}
If the uncertainty set has a finite number of extreme points, set $X_R$ is a polytope. 
\end{proposition}

Despite the nice theoretical properties, it is still difficult to solve an ARO problem in its general form (\ref{eq:App-03-ARO}). There have been considerable efforts spent on developing different approximations and approaches to tackle the computational challenges. We leave the solution methods of ARO problems to the next subsection. Here we demonstrate the benefit from postponing some decisions to the second stage via a simple example taken from \cite{RO-Detail-1}.

Consider an uncertain LP
\begin{equation}
\begin{aligned}
\min_x ~~ x_1  \\ 
\mbox{s.t.}~~ & x_2 \ge 0.5 \xi x_1 + 1 ~~ (a_\xi)  \\
              & x_1 \ge (2 - \xi) x_2 ~~~~ (b_\xi)  \\
              & x_1 \ge 0,~ x_2 \ge 0 ~~~~ (c_\xi)
\end{aligned}   \notag
\end{equation}
where $\xi \in [0, \rho]$ is an uncertain parameter and $\rho$ is a constant (level of uncertainty) which may take a value in the open interval $(0,1)$.

In a static setting, both $x_1$ and $x_2$ must be independent of $\xi$. When $\xi = \rho$, constraint $(a_{\xi})$ suggests $x_2 \ge 0.5 \rho x_1 +1$; when $\xi = 0$, constraint $(b_\xi)$ indicates $x_1 \ge 2 x_2$; as a result, we arrive at the conclusion $x_1 \ge \rho x_1 + 2$, so the optimal value in the static case satisfies 
\begin{equation}
\mbox{Opt} \ge x_1 \ge \dfrac{2}{1-\rho} \notag
\end{equation}
Thus the optimal value tends to infinity when $\rho$ approaches 1.

Now consider the adjustable case, in which $x_2$ is a wait-and-see decision. Let $x_2 = 0.5 \xi x_1 + 1$, ($a_\xi$) is always satisfied; substituting $x_2$ in constraint ($b_\xi$) yields:
\begin{equation}
x_1 \ge (2 - \xi) (\dfrac{1}{2} \xi x_1 + 1),~ \forall \xi \in [0, \rho]
\notag 
\end{equation}
 Substituting $x_1=4$ into above inequality we have 
\begin{equation}
4 \ge 2(2 - \xi)\xi + 2-\xi, \forall \xi \in [0, \rho] \notag
\end{equation}
This inequality can be certified by the fact that $\xi \ge 0$ and $\xi(2-\xi) \le 1$, $\forall \xi \in \mathbb R$, indicating that $x_1=4$ is a robust feasible solution. Therefore, the optimal value should be no greater than 4 in the adjustable  case for any $\rho$. The difference of optimal values in two cases can go arbitrarily large, depending on the value of $\rho$.

\subsection{Affine Policy Based Approximation Model}
\label{App-C-Sect02-02}

ARO problem (\ref{eq:App-03-ARO}) is difficult to solve because the functional dependence of the wait-and-see decision on $w$ is arbitrary, and there lacks a closed-form formula to characterize the optimal solution function $y(w)$ or certify whether $Y(x,w)$ is empty or not. At this point, we consider to approximate the recurse function $y(w)$ using a simpler one, naturally, an affine function
\begin{equation}
\label{eq:App-03-ARO-Affine-Policy}
y(w) = y^0 + G w
\end{equation}
where $y^0$ is the action in the second stage for the nominal scenario $w=0$, and $G$ is the gain matrix to be designed. (\ref{eq:App-03-ARO-Affine-Policy}) is called a linear decision rule or affine policy. It explicitly characterizes the wait-and-see decisions as an affine function in the revealed uncertain data. The rationality for employing an affine policy instead of other parametric ones is that  it yields computationally tractable robust counterpart reformulations.
This finding is firstly reported in \cite{ARO-Affine-Policy}.

To validate (\ref{eq:App-03-ARO-XR}) under the linear decision rule, substituting (\ref{eq:App-03-ARO-Affine-Policy}) in (\ref{eq:App-03-ARO-Y(x,w)}) 
\begin{equation}
\label{eq:App-03-AARO-RC-1}
Ax + By^0 + (BG + C) w \le b,~ \forall w \in W
\end{equation}
In (\ref{eq:App-03-AARO-RC-1}), decision variables are $x$, $y^0$, and $G$, which should be made before $w$ is known, and thus are here-and-now decisions. The wait-and-see decision (or the incremental part) is naturally determined from (\ref{eq:App-03-ARO-Affine-Policy}) without further optimization, and cost reduction is considered in the determination of gain matrix $G$. (\ref{eq:App-03-AARO-RC-1}) is in form of (\ref{eq:App-03-SRO-RC-Single}), and hence its robust counterpart can be derived via the methods in Appendix \ref{App-C-Sect01-02}. Here we just provide the results of polyhedral uncertainty as an example. 

Suppose the uncertainty set is described by
\begin{equation}
W = \{w ~|~ S w \le h \} \notag
\end{equation}
If we assume that $y^0$ is the optimal second stage decision when $w=0$, then we have
\begin{equation}
Ax + By^0 \le b \notag
\end{equation}
Furthermore, (\ref{eq:App-03-AARO-RC-1}) must hold if 
\begin{equation}
\label{eq:App-03-AARO-RC-2}
\max_{w \in W} (BG + C)_i w \le 0,~ \forall i
\end{equation}
where $(\cdot)_i$ stands for the $i$-th row of the input matrix.
According to LP duality theory,
\begin{equation}
\label{eq:App-03-AARO-RC-3}
\max_{w \in W}~ (BG + C)_i w = \min_{{\rm \Lambda}_i \in {\rm \Pi}_i}~ 
{\rm \Lambda_i} h, ~ \forall i
\end{equation}
where $\rm \Lambda$ is a matrix consists of the dual variables, ${\rm \Lambda_i}$ is the $i$-th row of $\rm \Lambda$ and also the dual variable of the $i$-th LP in (\ref{eq:App-03-AARO-RC-3}), and the set
\begin{equation}
{\rm \Pi}_i = \{ {\rm \Lambda}_i ~|~ {\rm \Lambda}_i \ge 0,~
{\rm \Lambda}_i S = (BG + C)_i \}  \notag
\end{equation} 
is the feasible region of the $i$-th dual LP.  

The minimization operator in the right-hand side of (\ref{eq:App-03-AARO-RC-3}) can be omitted if the objective is to seek a minimum. Moreover, if we adopt the minimum nominal cost criterion, the ARO problem with a linear decision rule in the second stage can be formulated as an LP  
\begin{equation}
\label{eq:App-03-AARO-RC-LP}
\begin{aligned}
\min~~ & c^T x +  d^T y^0  \\
\mbox{s.t.} ~~ & Ax + By^0 \le b, {\rm \Lambda} h \le 0 \\
& {\rm \Lambda} \ge 0,~ {\rm \Lambda} S = BG + C  
\end{aligned}
\end{equation}
In (\ref{eq:App-03-AARO-RC-LP}), decision variables are vectors $x$ and $y^0$, gain matrix $G$ and dual matrix $\rm \Lambda$. The constraints actually constitute a lifted formulation for $X_R$ in (\ref{eq:App-03-ARO-XR}). If the min-max cost criterion is employed, the objective can be transformed into a linear inequality constraint with uncertainty via an epigraph form, whose robust form can be derived using similar procedures shown above.  

Affine policy based method is attractive because it reduces the conservatism in the static RO approach by incorporating corrective actions, and sustains computational tractability. In theory, the affine assumption more or less restricts the adaptability in the recourse stage. Nevertheless, research work in \cite{AARO-Opt-1,AARO-Opt-2,AARO-Opt-3} shows that linear decision rules are indeed optimal or near optimal for many practical problems.  

For more information on other decision rules and their reformulations, please see \cite{RO-Detail-1} (Chapter 14.3) for the quadratic decision rule, \cite{ARO-Extend-Affine-Policy} for the extended linear decision rule, \cite{ARO-Finite-Adapt-1,ARO-Finite-Adapt-2} for the piecewise constant decision rule (finite adaptability), \cite{ARO-PWL-DR-1,ARO-PWL-DR-2} for the piecewise linear decision rule, and \cite{ARO-General-DR} for generalized decision rules. The methods in \cite{ARO-Finite-Adapt-1,ARO-PWL-DR-2} can be used to cope with integer wait-and-see decision variables. See also \cite{RO-Guide}.

\subsection{Algorithms for Fully Adjustable Models}
\label{App-C-Sect02-03}

Fully adjustable models are generally NP-hard \cite{ARO-Benders-Decomposition}. To find the solution in Definition \ref{df:App-03-SRO-Optimality}, the model is decomposed into a master problem and a subproblem, which are solved iteratively, and a sequence of lower bound and upper bound of the optimal values are generated, until they get close enough to each other. To explain the algorithm for ARO problems, we discuss two instances.

\vspace{12pt}
{\noindent \bf 1. Second-stage problem is an LP}

Now we consider problem (\ref{eq:App-03-ARO}) without specific functional assumptions on the wait-and-see variables. We start from the second-stage LP with fixed $x$ and $w$:
\begin{equation}
\label{eq:App-03-ARO-Inner-LP}
\begin{aligned}
\min_{y}~~ &  d^T y  \\
\mbox{s.t.}~~ & By \le b - Ax - Cw : u  
\end{aligned} 
\end{equation}
where $u$ is the dual variable, and the dual LP of (\ref{eq:App-03-ARO-Inner-LP}) is
\begin{equation}
\label{eq:App-03-ARO-Inner-Dual}
\begin{aligned}
\max_{u}~~ &  u^T (b - Ax - Cw)  \\
\mbox{s.t.}~~ & B^T u = d,~ u \le 0   
\end{aligned} 
\end{equation}

If the primal LP (\ref{eq:App-03-ARO-Inner-LP}) has a finite optimum, the dual LP (\ref{eq:App-03-ARO-Inner-Dual}) is also feasible and has the same optimum; otherwise, if (\ref{eq:App-03-ARO-Inner-LP}) is infeasible, then (\ref{eq:App-03-ARO-Inner-Dual}) will be unbounded. Sometimes, an improper choice of $x$ indeed leads to an infeasible second-stage problem.  To detect infeasibility, consider the following LP with slack variables 
\begin{equation}
\label{eq:App-03-ARO-Inner-LP-Slack}
\begin{aligned}
\min_{y,s}~~ &  1^T s \\
\mbox{s.t.}~~ & s \ge 0     \\
              & B y - I s \le b - Ax - Cw : u 
\end{aligned} 
\end{equation}
Its dual LP is 
\begin{equation}
\label{eq:App-03-ARO-Inner-Dual-Slack}
\begin{aligned}
\max_{u}~~ &  u^T (b - Ax - Cw)  \\
\mbox{s.t.}~~ & B^T u = 0,~ -1 \le u \le 0   
\end{aligned} 
\end{equation}
(\ref{eq:App-03-ARO-Inner-LP-Slack}) and (\ref{eq:App-03-ARO-Inner-Dual-Slack}) are always feasible and have the same finite optimums. If the optimal value is equal to 0, then LP (\ref{eq:App-03-ARO-Inner-LP}) is feasible; otherwise, if the optimal value is strictly positive, then LP (\ref{eq:App-03-ARO-Inner-LP}) is infeasible. 

For notation brevity, define feasible sets for the dual variable
\begin{equation}
\begin{aligned}
U_O & = \{ u ~|~ B^T u = d,~ u \le 0 \} \\
U_F &= \{ u ~|~ B^T u = 0,~ -1 \le u \le 0 \}
\end{aligned}
\notag
\end{equation}
The former one is associated with the dual form (\ref{eq:App-03-ARO-Inner-Dual}) of the second-stage optimization problem (\ref{eq:App-03-ARO-Inner-LP}); the latter one corresponds to the dual form (\ref{eq:App-03-ARO-Inner-Dual-Slack}) of the second-stage feasibility test problem (\ref{eq:App-03-ARO-Inner-LP-Slack}).

Next, we proceed to the middle level with fixed $x$:
\begin{equation}
\label{eq:App-03-ARO-Middle-Linear-max-min}
R(x)=\max_{w \in W} \min_{y \in Y(x,w)} d^T y
\end{equation}
which is a linear max-min problem that identifies  the worst-case uncertainty. If LP (\ref{eq:App-03-ARO-Inner-LP}) is feasible for an arbitrarily given value of $w \in W$, then we conclude $x \in X_R$ defined in (\ref{eq:App-03-ARO-XR}); otherwise, if LP (\ref{eq:App-03-ARO-Inner-LP}) is infeasible for some $w \in W$, then $x \notin X_R$ and $R(x) = +\infty$.

To check whether $x \in X_R$ or not, we investigate the following problem
\begin{equation}
\label{eq:App-03-ARO-Middle-Linear-max-min-Fea}
\begin{aligned}
\max_{w}~~ & \min_{y,s}  1^T s \\
\mbox{s.t.}~~ & w \in W,~ s \ge 0     \\
              & B y - I s \le b - Ax - Cw : u 
\end{aligned} 
\end{equation}

It maximizes the minimum of (\ref{eq:App-03-ARO-Inner-LP-Slack}) over all possible values of $w \in W$. Since the minimums of (\ref{eq:App-03-ARO-Inner-LP-Slack}) and (\ref{eq:App-03-ARO-Inner-Dual-Slack}) are equal, problem (\ref{eq:App-03-ARO-Middle-Linear-max-min-Fea}) is equivalent to maximizing the optimal value of (\ref{eq:App-03-ARO-Inner-Dual-Slack}) over the uncertainty set $W$, leading to a bilinear program 
\begin{equation}
\label{eq:App-03-ARO-Middle-BLP-Fea}
\begin{aligned}
r(x) =\max_{u,w}~~ &  u^T (b - Ax - Cw)  \\
\mbox{s.t.}~~ & w \in W,~ u \in U_F   
\end{aligned} 
\end{equation}
Because both $W$ and $U_F$ are bounded, (\ref{eq:App-03-ARO-Middle-BLP-Fea})
must have a finite optimum. Clearly, $0 \in U_F$, so $r(x)$ must be non-negative. In fact, if $r(x) = 0$, then $x \in X_R$; if  $r(x)> 0$, then $x \notin X_R$. With the duality transformation, the opposite optimization operators in (\ref{eq:App-03-ARO-Middle-Linear-max-min-Fea}) come down to a traditional NLP. 

For similar reasons, by replacing the second-stage LP (\ref{eq:App-03-ARO-Inner-LP}) with its dual LP (\ref{eq:App-03-ARO-Inner-Dual}), problem (\ref{eq:App-03-ARO-Middle-Linear-max-min}) is equivalent to the following bilinear program 
\begin{equation}
\label{eq:App-03-ARO-Middle-BLP-Opt}
\begin{aligned}
r(x) =\max_{u,w}~~ &  u^T (b - Ax - Cw)  \\
\mbox{s.t.}~~ & w \in W,~ u \in U_O   
\end{aligned} 
\end{equation}

The fact that a linear max-min problem can be transformed as a bilinear program using LP duality is reported in \cite{Linear-max-min-BLP}. Bilinear programs can be locally solved by general purpose NLP solvers, but the non-convexity prevents a global optimal solution from being found easily. In what follows, we introduce some methods that exploit specific features of the uncertainty set and are widely used by the research community. In view that (\ref{eq:App-03-ARO-Middle-BLP-Fea}) and (\ref{eq:App-03-ARO-Middle-BLP-Opt}) only differ in the dual feasibility set, we will use set $U$ to refer either $U_F$ or $U_O$ in the unified solution method.  

\vspace{12pt}
{\noindent \bf a. General polytope}

Suppose that the uncertainty set is described by
\begin{equation}
W = \{w ~|~ S w \le h \} \notag
\end{equation}

An important feature in (\ref{eq:App-03-ARO-Middle-BLP-Fea}) and (\ref{eq:App-03-ARO-Middle-BLP-Opt}) is that the constraint set $W$ and $U$ are separated and there is no constraint that involves $w$ and $u$ simultaneously, so the bilinear program can be considered in the following format
\begin{equation}
\label{eq:App-03-ARO-BLP-Poly-1}
\begin{aligned}
\max_{u \in U}~~  u^T (b - A x) + \max_{w} ~~ & (- u^T C w ) \\
\mbox{s.t.}~~ & S w \le h : \xi   
\end{aligned} 
\end{equation}
The bilinear term $u^T C w$ is non-convex. If we treat the second part $\max _{w \in W} (- u^T C w )$ as an LP in $w$ where $u$ is a parameter, whose KKT optimality condition is given by
\begin{equation}
\label{eq:App-03-ARO-BLP-Poly-2}
\begin{gathered}
0 \le \xi \bot h - S w \ge 0 \\
S^T \xi + C^T u = 0   
\end{gathered} 
\end{equation}
The stationary point of LCP (\ref{eq:App-03-ARO-BLP-Poly-2}) gives the optimal primal and dual solutions simultaneously. As the uncertainty set is a bounded polyhedron, the optimal solution must be bounded, and strong duality holds, so we can replace $- u^T C w$ in the objective with a linear term $h^T \xi$ and additional constraints in (\ref{eq:App-03-ARO-BLP-Poly-2}). Moreover, the complementarity and slackness condition in (\ref{eq:App-03-ARO-BLP-Poly-2}) can be linearized via the method in Appendix \ref{App-B-Sect03-05}. In summary, problem (\ref{eq:App-03-ARO-BLP-Poly-1}) can be solved via an equivalent MILP
\begin{equation}
\label{eq:App-03-ARO-BLP-Poly-MILP}
\begin{aligned}
\max_{u,w,\xi}~~ & u^T (b - A x) + h^T \xi \\
\mbox{s.t.}~~ & u \in U,~ \theta \in \{0,1\}^m \\
              & S^T \xi + C^T u = 0  \\
              & 0 \le \xi \le M(1 - \theta) \\
              & 0 \le h - S w \le M \theta 
\end{aligned} 
\end{equation}
where $m$ is the dimension of $\theta$, and $M$ is a large enough constant. Compared with  (\ref{eq:App-03-ARO-BLP-Poly-1}), non-convexity migrates from the objective function to the constraints with binary variables. The number of binary variables in (\ref{eq:App-03-ARO-BLP-Poly-MILP}) only depends on the number of constraints in set $W$, and is independent of the dimension of $x$.  

Another heuristic method for bilinear programs in the form of (\ref{eq:App-03-ARO-Middle-BLP-Fea}) and (\ref{eq:App-03-ARO-Middle-BLP-Opt}) is the mountain climbing method in \cite{BLP-Mountain-Climbing}, which is summarized in Algorithm \ref{Ag:App-03-BLP-Mountain-Climbing}
\begin{algorithm}[!htp]
\normalsize
\caption{\bf : Mountain climbing}
\begin{algorithmic}[1]
\STATE Choose a convergence tolerance $\varepsilon>0$, and an initial $w^* \in W$;
\STATE Solve the following LP with current $w^*$
\begin{equation}
\label{eq:App-03-BLP-Mountain-Climbing-LP1}
 R_1 = \max_{u \in U} ~ u^T ( b - A x - C w^*) 
\end{equation}
The optimal solution is $u^*$ and the optimal value is $R_1$;
\STATE Solve the following LP with current $u^*$
\begin{equation}
\label{eq:App-03-BLP-Mountain-Climbing-LP2}
 R_2 = \max_{w \in W} ~ (b - A x - C w)^T u^* 
\end{equation}
The optimal solution is $w^*$ and the optimal value is $R_2$;
\STATE If $R_2 - R_1  \le \varepsilon$, report the optimal value $R_2$ as well as the optimal solution $w^* ,u^*$, and terminate; otherwise, go to step 2.
\end{algorithmic}
\label{Ag:App-03-BLP-Mountain-Climbing}
\end{algorithm} 

The optimal solutions of LPs must be found at one of the vertices of its feasible region, hence $w^* \in {\rm{vert}}(W)$ and $u^*  \in {\rm{vert}}(U)$ hold. As its name implies, the sequence of objective values generated by  Algorithm \ref{Ag:App-03-BLP-Mountain-Climbing} is monotonically increasing, until a local maximum is found \cite{BLP-Mountain-Climbing}. The convergence is guaranteed by the finiteness of $\mbox{vert}(U)$ and $\mbox{vert}(W)$. If we try multiple initial points that are chosen elaborately and pick up the best one among the returned results, the solution quality is often satisfactory.  The key point is, these initial points should span along most directions in the $w$-subspace. For example, one may search the $2m$ points on the boundary of $W$ in directions $\pm e^m_i$, $i=1,2,\cdots,m$, where $m$ is the dimension of $w$, and $e^m_i$ is the $i$-th column of an $m \times m$ identity matrix. As LPs can be solved very efficiently, Algorithm \ref{Ag:App-03-BLP-Mountain-Climbing} is especially suitable for the instances with very complicated $U$ and $W$, and usually outperforms general NLP solvers for bilinear programs with disjoint constraints.

Algorithm \ref{Ag:App-03-BLP-Mountain-Climbing} is also valid if $W$ is other convex set, say, an ellipsoid, and converges to a local optimum in a finite number of iterations for a given precision \cite{BLP-Mountain-Climbing-BCVX}.  

\vspace{12pt}
{\noindent \bf b. Cardinality constrained uncertainty set}

A continuous cardinality constrained uncertainty set in the form of (\ref{eq:App-03-SRO-US-Card}) is a special class of the polyhedral case, see the transformation in (\ref{eq:App-03-SRO-US-Card-Lift}). Therefore, the previous method can be applied, and the number of inequalities in the polyhedral form is $3m+1$, which is equal to the number of binary variables in MILP (\ref{eq:App-03-ARO-BLP-Poly-MILP}). As revealed in Proposition \ref{pr:App-03-Invariance-3}, for a polyhedral uncertainty set, we can merely consider the extreme points. 

Consider a discrete cardinality constrained uncertainty set
\begin{subequations}
\label{eq:App-03-ARO-US-Card}
\begin{equation}
\label{eq:App-03-ARO-W-Card}
W = \left\{ w \middle| \begin{gathered}
w_j = w^0_j + w^+_j z^+_j - w^-_j z^-_j, \forall j \\
\exists ~ z^+,z^- \in Z   
\end{gathered} \right\}   
\end{equation}
\begin{equation}
\label{eq:App-03-ARO-Z-Card}
Z = \left\{ z^+,z^- \middle| \begin{gathered}
z^+,~ z^- \in \{0,1\}^m  \\ 
z^+_j + z^-_j \le 1,~ \forall j \\
 1^T ( z^+ + z^-)  \le {\rm \Gamma}  
\end{gathered} \right\}   
\end{equation}
\end{subequations}
where the budget of uncertainty ${\rm \Gamma} \le m$ is an integer. In (\ref{eq:App-03-ARO-W-Card}), each element $w_j$ takes one of three possible values: $w^0_j$, $w^0_j + w^+_j$, and $w^0_j - w^-_j$, and at most $\rm \Gamma$ of the $m$ elements $w_j$ can take a value that is not equal to $w^0_j$. If the forecast error is symmetric, i.e., $w^+_j = w^-_j$, then (\ref{eq:App-03-ARO-US-Card}) is called symmetric as the nominal scenario locates at the center of $W$. We discuss this case separately because this representation allows to linearize the non-convexity in (\ref{eq:App-03-ARO-Middle-BLP-Fea}) and (\ref{eq:App-03-ARO-Middle-BLP-Opt}) with fewer binary variables.   

Expanding the bilinear term $u^T C w$ in an element-wise form
\begin{equation}
u^T C w = u^T C w^0 + \sum_i \sum_j 
( c_{ij} w^+_j u_i z^+_j -  c_{ij} w^-_j u_i z^-_j )  \notag
\end{equation}
where $c_{ij}$ is the element of matrix $C$. Let
\begin{equation}
v^+_{ij} = u_i z^+_j,~ v^-_{ij} = u_i z^-_j,~ \forall i,~ \forall j \notag
\end{equation}
the bilinear term can be expressed via a linear function. The product involving a binary variable and a continuous variable can be linearized via the method illuminated in Appendix \ref{App-B-Sect02-02}. 

In conclusion, bilinear subproblems (\ref{eq:App-03-ARO-Middle-BLP-Fea}) and (\ref{eq:App-03-ARO-Middle-BLP-Opt}) can be solved via MILP 
\begin{equation}
\label{eq:App-03-ARO-W-Card-MILP}
\begin{aligned}
\max ~~ & u^T ( b - A x ) - u^T C w^0 - \sum_i \sum_j 
( c_{ij} w^+_j v^+_{ij} -  c_{ij} w^-_j v^-_{ij} )  \\
\mbox{s.t.} ~~ & u \in U,~  \{z^+,z^- \} \in Z  \\
& 0 \le v^+_{ij} - u_j \le M (1 - z^+_j), -M z^+_j \le v^+_{ij} \le 0,~
\forall i, \forall j \\ 
& 0 \le v^-_{ij} - u_j \le M (1 - z^-_j), -M z^-_j \le v^-_{ij} \le 0,~
\forall i, \forall j  
\end{aligned}
\end{equation}
where $M=1$ for problem (\ref{eq:App-03-ARO-Middle-BLP-Fea}) since $-1 \le u \le 0$, and $M$ is a sufficiently large number for problem (\ref{eq:App-03-ARO-Middle-BLP-Opt}),
because there is no clear bounds for the dual variable $u$.  The number of binary variables in MILP (\ref{eq:App-03-ARO-W-Card-MILP}) is $2m$, which is less than that in (\ref{eq:App-03-ARO-BLP-Poly-MILP}) if the uncertainty set is replaced by its convex hull. The number of additional continuous variables $v^+_{ij}$ and $v^-_{ij}$ is also moderate since the matrix $C$ is sparse.

Finally, we are ready to give the decomposition algorithm which is proposed in \cite{ARO-CCG}. In light of Proposition \ref{pr:App-03-Invariance-3},
it is sufficient to consider the extreme points $w^1$, $w^2$, $\cdots$, $w^S$ in the uncertainty set, inspiring the following epigraph formulation which is equivalent to (\ref{eq:App-03-ARO})
\begin{equation}
\label{eq:App-03-ARO-Epigraph}
\begin{aligned}
\min_{x,y^s,\eta } ~~ & c^T x + \eta  \\
\mbox{s.t.} ~~ &   x \in X  \\
& \eta \ge d^T y^s,~ \forall s  \\
& A x + B y^s \le b - C w^s, \forall s 
\end{aligned} 
\end{equation}
Recall (\ref{eq:App-03-ARO-XR-Lift}), the last constraint is in fact a lifted formulation for $X_R$. For polytope and cardinality constrained uncertainty sets, the number of extreme points are finite, but may grow exponentially in the dimension of uncertainty. Actually, it is difficult and also unnecessary to enumerate every extreme point, because most of them actually provide redundant constraints. A smart method is to identify active scenarios which contribute binding constraints in $X_R$. This motivation has been widely used in complex optimization problems and formalized in Sect. \ref{App-C-Sect01-03}. The procedure of the adaptive scenario generation algorithm for ARO is summarized in Algorithm \ref{Ag:App-03-ARO-Adaptive-Scenario-Generation}.

\begin{algorithm}[!htp]
\normalsize
\caption{\bf : Adaptive scenario generation}
\begin{algorithmic}[1]
\STATE Choose a tolerance $\varepsilon > 0$, set $LB = -\infty$, $UB = +\infty$, iteration index $k = 0$, and the critical scenario set $O = w^0$;
\STATE Solve the following master problem
\begin{equation}
\label{eq:App-03-ARO-ASG-Master}
\begin{aligned}
\min_{x,y^s,\eta } ~~ & c^T x + \eta  \\
\mbox{s.t.} ~~ &   x \in X,~ \eta \ge d^T y^s,~ s = 0,\cdots, k \\
& A x + B y^s \le b - C w^s, \forall w^s \in O 
\end{aligned} 
\end{equation}
The optimal solution is $x^{k+1}$, $\eta^{k+1}$, and update $LB = c^T x^{k+1} + \eta^{k+1} $;
\STATE Solve bilinear feasibility testing problem (\ref{eq:App-03-ARO-Middle-BLP-Fea}) with $x^{k+1}$, the optimal solution is $w^{k+1}$, $u^{k+1}$; if the optimal value $r^{k+1} > 0$, update $O = O \cup w^{k+1}$, and add a scenario cut 
\begin{equation}
\label{eq:App-03-ARO-ASG-Cut}
\eta \ge d^T y^{k+1},~ A x + B y^{k+1} \le b - C w^{k+1} 
\end{equation}
with a new variable $y^{k+1}$ to the master problem (\ref{eq:App-03-ARO-ASG-Master}), update $k \leftarrow k+1$, and go to Step 2; 
\STATE Solve bilinear optimality testing problem (\ref{eq:App-03-ARO-Middle-BLP-Opt}) with $x^{k+1}$, the optimal solution is $w^{k+1}$, $u^{k+1}$, and the optimal value is $R^{k+1}$; update $O = O \cup w^{k+1}$ and $UB = c^T x^{k+1} + R^{k+1}$, create scenario cut (\ref{eq:App-03-ARO-ASG-Cut}) with a new variable $y^{k+1}$.
\STATE If $UB - LB \le \varepsilon$, report the optimal solution, terminate; otherwise, add the scenario cut in step 4 to the master problem (\ref{eq:App-03-ARO-ASG-Master}), update $k \leftarrow k+1$, and go to step 2;
\end{algorithmic}
\label{Ag:App-03-ARO-Adaptive-Scenario-Generation}
\end{algorithm} 

Algorithm \ref{Ag:App-03-ARO-Adaptive-Scenario-Generation} converges in a finite number of iterations, which is bounded by the  number of extreme points of the uncertainty set. In practice, this algorithm often converges in a few iterations, because problems (\ref{eq:App-03-ARO-Middle-BLP-Fea}) and (\ref{eq:App-03-ARO-Middle-BLP-Opt}) always identify the most critical scenario that should be considered. This is why we name the algorithm ``adaptive scenario generation''. It is called ``constraint-and-column generation algorithm'' in \cite{ARO-CCG}, because the numbers of decision variables (columns) and constraints increase simultaneously. Please note that the scenario cut streamlines the feasibility cut and optimality cut used in the existing literature.  

Bilinear subproblems (\ref{eq:App-03-ARO-Middle-BLP-Fea}) and (\ref{eq:App-03-ARO-Middle-BLP-Opt}) can be solved by the methods discussed previously, according to the form of the uncertainty set. In Algorithm \ref{Ag:App-03-ARO-Adaptive-Scenario-Generation}, we utilize $w$ to create scenario cuts, which are also called primal cuts. In fact, the optimal dual variable $u$ of (\ref{eq:App-03-ARO-Middle-BLP-Fea}) and (\ref{eq:App-03-ARO-Middle-BLP-Opt}) provides sensitivity information, and can be used to construct dual cuts, which is a single inequality in the first-stage variable $x$. See Benders decomposition algorithm in \cite{ARO-Benders-Decomposition}. Since scenario cuts are much tighter than Benders cuts, Algorithm  \ref{Ag:App-03-ARO-Adaptive-Scenario-Generation} is the most prevalent method for solving ARO problems.     

If matrix $A$ is uncertainty-dependent, the scenario constraints in the master problem (\ref{eq:App-03-ARO-ASG-Master}) becomes $A(w^s) x + B y^s \le b - C w^s$, $\forall w^s$, where $A(w^s)$ is constant but varies in different scenarios; the objective function of bilinear subproblems changes to $u^T [b - A(w)x - Cw]$, where $x$ is given in the subproblem. If $A$ can be expressed as a linear function in $w$, the problem structure remains the same, and previous methods are still valid. Even if the second-stage problem is an SOCP, the adaptive scenario generation framework remains applicable, and the key procedure is to solve a max-min SOCP. Such a model originates from the robust operation of a power distribution network with uncertain generation and demand. By dualizing the inner-most SOCP, the max-min SOCP is cast as a bi-convex program, which can be globally or locally solved via an MISOCP or the mountain climbing method.

Recently, the duality theory of fully-adjustable robust optimization problem has been proposed in \cite{Duality-ARO}. It has been shown that this kind of problem is self-dual, i.e., the dual problem remains an ARO. However, solving the dual problem may enjoy better efficiency. An extended CCG algorithm which always produces a feasible fist-stage decision (if one exists) is proposed in \cite{Ext-CCG-ARO}.

\vspace{12pt}
{\noindent \bf 2. Second-stage problem is an MILP}

Now we consider the case in which some of the wait-and-see decisions are discrete. As what can be observed from the previous case, the most important tasks in solving an ARO problem is to validate feasibility and optimality, which can boil down to solving a linear max-min problem. When the  wait-and-see decisions are continuous and the second-stage problem is linear, LP duality theory is applied such that the linear max-min problem is cast as a  traditional bilinear program. However, discrete variables appearing in the second stage make the recourse problem a mixed-integer linear max-min problem with a non-convex inner level, preventing the use of LP duality theory. As a result, validating feasibility and optimality becomes more challenging.

The compact form of an ARO problem with integer wait-and-see decisions can be written as  

\begin{equation}
\label{eq:App-03-ARO-MIP-Recourse}
\min_{x \in X} \left\{ c^T x + \max_{w \in W} \min_{y,z \in Y(x,w)} d^T y + g^T z \right\} 
\end{equation}
where $z$ is binary and depends on the exact value of $w$; the feasible region
\begin{equation}
Y(x,w) = \left\{ y,z ~\middle|~ \begin{gathered}
By + Gz \le b - Ax - Cw   \\
y \in \mathbb{R}^{m_1},~ z \in {\rm \Phi}
\end{gathered}  \right\} \notag
\end{equation}
where feasible set ${\rm \Phi} =  \{z | z \in \mathbb{B}^{m_2}, T z \le v\}$; $m_1$ and $m_2$ are dimensions of $y$ and $z$; $T$ and $v$ are constant coefficients; all coefficient matrices have compatible dimensions. We assume that the uncertainty set $W$ can be represented by a finite number of extreme points. This kind of problem is studied in \cite{ARO-MIP-Nested-CCG}. A nested constraint-and-column generation algorithm is proposed. 

Different from the mainstream idea that directly solves a linear max-min program as a bilinear program, the mixed-integer max-min program in (\ref{eq:App-03-ARO-MIP-Recourse}) is expanded to a tri-level problem 
\begin{equation}
\label{eq:App-03-MILP-max-min}
\begin{aligned}
\max_{w \in W} \min_{z \in {\rm \Phi}} g^T z + \min_{y}~~ & d^T y  \\
\mbox{s.t.}~~ &  By \le b - Ax - Cw -G z
\end{aligned}
\end{equation}
For the ease of discussion, we assume all feasible sets are bounded, because decision variables of practical problems have physical bounds. By replacing the innermost LP in variable $y$ with its dual LP, problem (\ref{eq:App-03-MILP-max-min}) becomes
\begin{equation}
\label{eq:App-03-MILP-Recoure-Trilevel}
\max_{w \in W} \left\{ \min_{z \in {\rm \Phi}} \left\{ g^T z + 
\max_{u \in U} u^T (b - Ax - Cw -G z) \right\} \right\}
\end{equation}
where $u$ is the dual variable, and set $U = \{ u ~|~ u \le 0,~ B^T u = d \}$. Because both $w$ and $z$ are expressed via binary variables, bilinear terms $u^T C w$ and $u^T G z$ have linear representations by using the method in Appendix \ref{App-B-Sect02-02}. Since $\rm \Phi$ has a countable number of elements, problem (\ref{eq:App-03-MILP-Recoure-Trilevel}) (in its linearized version) has the same form as ARO problem (\ref{eq:App-03-ARO}), and can be solved by Algorithm \ref{Ag:App-03-ARO-Adaptive-Scenario-Generation}. More exactly, write (\ref{eq:App-03-MILP-Recoure-Trilevel}) into an epigraph form by enumerating all possible elements $z \in {\rm \Phi}$, then perform Algorithm \ref{Ag:App-03-ARO-Adaptive-Scenario-Generation} and identify binding elements. In this way, the minimization operator in the middle level is eliminated.

The nested adaptive scenario generation algorithm for ARO problem (\ref{eq:App-03-ARO-MIP-Recourse}) with mixed-integer recourses is summarized in Algorithm \ref{Ag:App-03-ARO-Nested-ASG}. Because both $W$ and $\rm \Phi$ are finite sets with countable elements, Algorithm \ref{Ag:App-03-ARO-Nested-ASG} converges in a finite number of iterations. Notice that we do not distinguish feasibility and optimality subproblems in above algorithm due to their similarities. One can also introduce slack here-and-now variables in the second stage and penalty terms in the objective function, such that the recourse problem is always feasible. It should be pointed out that Algorithm \ref{Ag:App-03-ARO-Nested-ASG} incorporates double loops, and an MILP should be solved in each iteration in the inner loop, so we'd better not expect too much on its efficiency. Nonetheless, it is the first systematic method to solve an ARO problem with integer variables in the second stage. Another concept which should be clarified is that although the second-stage discrete variable $z$ is treated as scenario and enumerated on the fly when solving problem (\ref{eq:App-03-MILP-Recoure-Trilevel}) in step 3 (the inner loop), it is a decision variable of the master problem (\ref{eq:App-03-ARO-Nested-ASG-Master}) in the outer loop.

\begin{algorithm}[!htp]
\normalsize
\caption{\bf : Nested adaptive scenario generation}
\begin{algorithmic}[1]
\STATE Choose a tolerance $\varepsilon > 0$, set $LB = -\infty$, $UB = +\infty$, iteration index $k = 0$, and the critical scenario set $O = w^0$;
\STATE Solve the following master problem
\begin{equation}
\label{eq:App-03-ARO-Nested-ASG-Master}
\begin{aligned}
\min_{x,y,z,\eta } ~~ & c^T x + \eta  \\
\mbox{s.t.} ~~ & x \in X  \\
& \eta \ge d^T y^s + g^T z^s,~ z^s \in {\rm \Phi},~ s = 0,\cdots,k \\
& A x + B y^s + G z^s  \le b - C w^s,~ \forall w^s \in O 
\end{aligned} 
\end{equation}
The optimal solution is $x^{k+1}$, $\eta^{k+1}$, and update $LB = c^T x^{k+1} + \eta^{k+1} $;
\STATE Solve problem (\ref{eq:App-03-MILP-Recoure-Trilevel}) with $x^{k+1}$, the optimal solution is $(z^{k+1},w^{k+1},u^{k+1})$, and optimal value is $R^{k+1}$; update $O = O \cup w^{k+1}$, $UB = \min \{UB, c^T x^{k+1} + R^{k+1}\}$, create new variables $(y^{k+1},z^{k+1})$ and scenario cuts 
\begin{equation}
\label{eq:App-03-ARO-Nested-ASG-Cut}
\begin{gathered}
 \eta \ge d^T y^{k+1} + g^T z^{k+1},~ z^{k+1} \in {\rm \Phi} \\
 A x + B y^{k+1} + G z^{k+1} \le b - C w^{k+1} 
\end{gathered} 
\end{equation}
\STATE If $UB - LB \le \varepsilon$, terminate and report the optimal solution and optimal value; otherwise, add scenario cuts (\ref{eq:App-03-ARO-Nested-ASG-Cut}) to the master problem (\ref{eq:App-03-ARO-Nested-ASG-Master}), update $k \leftarrow k+1$, and go to step 2;
\end{algorithmic}
\label{Ag:App-03-ARO-Nested-ASG}
\end{algorithm}

As a short conclusion, to overcome the limitation of  traditional static RO approaches which require all decisions should be made without exact information on the underlying uncertainty, ARO employs a two-stage decision-making framework and allows a subset of decision variables to be made after the uncertain data are revealed. Under some special decision rules, computational tractability can be preserved. In fully adjustable cases, the ARO problem can be solved by a decomposition algorithm. The subproblem comes down to a (mixed-integer) linear max-min problem, which is generally challenging to solve. We introduce MILP reformulations for special classes of uncertainty sets, which are compatible with commercial solvers, and help solve an engineering optimization problem in a systematic way.

\section{Distributionally Robust Optimization}
\label{App-C-Sect03}

Static and adjustable RO models presented in Sect. \ref{App-C-Sect01} and Sect. \ref{App-C-Sect02} do not rely on specifying probability distributions of the uncertain data, which are used in SO approaches for generating scenarios, evaluating probability of constraint violation, or deriving analytic solutions for some specific problems. Instead, RO design principle aims to cope with the worst-case scenario in a pre-defined uncertainty set in the space of uncertain variables, which is a salient distinction between these two approaches. If the exact probability distribution is precisely known, optimal solutions to SO models would be less conservative than the robust ones from the statistical perspective. However, the optimal solution to SO models could have poor statistical performances if the actual distribution is not identical to the designated one \cite{Bertsimas-2006}. As for the RO approach, as it hedges against the worst-case scenario, which rarely happens in reality, the robust strategy could be conservative thus suboptimal in most cases. 

A method which aims to build a bridge connecting SO and RO approaches is the DRO, whose optimal solutions are designed for the worst-case probability distribution within a family of candidate distributions, which are described by statistic information, such as moments, and structure properties, including symmetry, unimodality, and so on.  This approach is generally less conservative than the traditional RO because dispersion effect of uncertainty is taken into account, i.e., the probability of an extreme event is low. Meanwhile,  the statistic performances of the solution is less sensitive to the perturbation in probability distributions than that of an SO model,  as it hedges against the worst distribution. Publications on this method have been proliferating rapidly in the past few years. This section only sheds light on some most representative methods which have been used in energy system studies. 

\subsection{Static Distributionally Robust Optimization}
\label{App-C-Sect03-01}

In analogy with the terminology used in Sect \ref{App-C-Sect01}, ``static'' means that all decision variables are here-and-now type. Theoretical outcomes in this part mainly come from \cite{Static-DRO}. A static DRO problem can be formulated as  
\begin{equation}
\label{eq:App-03-DRO-Model-1}
\begin{aligned}
\min_x ~~ & c^T x   \\
\mbox{s.t.} ~~ & x \in X  \\
& \Pr \left( a_i (\xi)^T x \le b_i(\xi),~ i=1,\cdots,m \right) \ge 1 - \varepsilon,~ \forall f(\xi) \in {\mathcal P}
\end{aligned}
\end{equation}
where $x$ is the decision variable, $X$ is a closed and convex set that is independent of the uncertain parameter, $c$ is a deterministic vector, and $\xi$ is the uncertain data, whose probability density function $f(\xi)$ is not known exactly, and belongs to $\mathcal P$, a set comprised of candidate distributions. Robust chance constraint in (\ref{eq:App-03-DRO-Model-1}) requires a finite number of linear inequalities depending on $\xi$ to be met with a probability of at least $1-\varepsilon$, regardless of the true probability density function of $\xi$. We assume uncertain coefficients $a_i$ and $b_i$ are linear functions in $\xi$, i.e.
\begin{equation}
\begin{gathered}
a_i(\xi) = a^0_i + \sum_{j=1}^k a^j_i \xi_j  \\
b_i(\xi) = b^0_i + \sum_{j=1}^k b^j_i \xi_j
\end{gathered}  \notag
\end{equation}
where $a^0_i$, $a^j_i$ are constant vectors and $b^0_i$, $b^j_i$ are constant scalars. Define
\begin{equation}
y^j_i (x) = (a^j_i)^T x - b^j_i,~\forall i,~ \forall j  \notag
\end{equation}
the chance constraint in (\ref{eq:App-03-DRO-Model-1}) can be expressed via
\begin{equation}
\label{eq:App-03-DRO-RCC}
\Pr \left( y^0_i (x) + y_i(x)^T \xi \le 0,~ i=1,\cdots,m \right) \ge 1 - \varepsilon,~ \forall f(\xi) \in {\mathcal P}
\end{equation}
where vector $y_i(x)=[y^1_i(x),\cdots,y^k_i(x)]^T$ is affine in $x$. Since the objective is certain and constraint violation is bounded by a small probability, problem (\ref{eq:App-03-DRO-Model-1}) is also called a robust chance-constrained program.

Chance constraints can be transformed into tractable ones that are convex in variable $x$ only for a few special cases. For example, if $\xi$ follows a Gaussian distribution, $\varepsilon \le 0.5$, and $m=1$, then the individual chance constraint without distribution uncertainty is equivalent to a single SOC constraint \cite{CCO-Gauss}. For $m >1$, joint chance constraints form  convex feasible region when the right-hand side terms $b_i(\xi)$ are uncertain and follow a log-concave distribution \cite{Static-DRO,CCP-RHS-Log-Concave}, while coefficients $a_i$, $i=1,\cdots,m$ are deterministic.

Constraint (\ref{eq:App-03-DRO-RCC}) is even more challenging at first sight: not only the random vector $\xi$, but also the probability distribution function $f(\xi)$ itself is uncertain. Because in many practical situations, probability distribution must be estimated from enough historical data, which may not be available at hand. Typically, one may only have access to some statistical indicators about $f(\xi)$, e.g. its mean value, covariance, and support set. Using a specific $f(\xi) \in \mathcal P$ may lead to over-optimistic solutions which fail to satisfy the probability guarantee under the true distribution.

Similar to the paradigm in static RO, a prudent way to immunize a chance constraint against uncertain probability distribution is to investigate the situation in the worst case, inspiring the following distributionally robust chance constraint, which is equivalent to (\ref{eq:App-03-DRO-RCC})
\begin{equation}
\label{eq:App-03-DRO-DRCC}
\inf_{f(\xi) \in \mathcal P} \Pr \left( y^0_i (x) + y_i(x)^T \xi \le 0,~ i=1,\cdots,m \right) \ge 1 - \varepsilon
\end{equation}
Clearly, if $x$ satisfies (\ref{eq:App-03-DRO-DRCC}), the probability of constraint violation is upper bounded by $\varepsilon$ for the true probability distribution of $\xi$. 

This section introduces convex optimization models for approximating
robust chance constraints under uncertain probability distributions, whose first- and second-order moments as well as the support set (or equivalently the feasible region) of random variable are known. More precisely, we let $\mathbb E_P (\xi) = \mu \in \mathbb R^k$ be the mean value and $\mathbb E_P ((\xi-\mu)(\xi-\mu)^T) = {\rm \Sigma} \in \mathbb S^k_{++}$ be the covariance matrix of random variable $\xi$ under the true distribution $P$. We define the moment matrix
\begin{equation}
{\rm \Omega} = \begin{bmatrix}
{\rm \Sigma} + \mu \mu^T & \mu \\ \mu^T & 1
\end{bmatrix}
\notag
\end{equation}
for ease of notation. 

To help readers understand the fundamental ideas in DRO, we briefly introduce the worst-case expectation problem, which will be used throughout this section. Recall that $\mathcal P$ represents the set of all probability distributions on $\mathbb R^k$ with mean vector $\mu$ and covariance matrix ${\rm \Sigma} \succ 0$,. the problem is formulated by 
\begin{equation}
\theta^m_P = \sup_{f(\xi) \in \mathcal P} \mathbb E 
\left[ (g(\xi))^+ \right]   \notag
\end{equation}
where $g:\mathbb R^k \to \mathbb R$ is a function of $\xi$; $(g(\xi))^+$ means the maximum between 0 and $g(\xi)$. Write the problem into an integral format
\begin{equation}
\label{eq:App-03-DRO-Worst-Expectation-Primal}
\begin{aligned}
\theta^m_P =  \sup_{f(\xi) \in \mathcal P} ~~ & 
\int_{\xi \in \mathbb R^k} \max \{0, g(\xi) \} f(\xi) \mbox{d} \xi \\
\mbox{s.t.} ~~ &  f(\xi) \ge 0,~ \forall \xi \in \mathbb R^k     \\
& \int_{\xi \in \mathbb R^k} f(\xi) \mbox{d} \xi = 1 :~ \lambda_0    \\
& \int_{\xi \in \mathbb R^k} \xi f(\xi) \mbox{d} \xi = \mu :~ \lambda \\
& \int_{\xi \in \mathbb R^k} \xi \xi^T f(\xi) \mbox{d} \xi 
= {\rm \Sigma} +  \mu \mu^T :~ {\rm \Lambda}
\end{aligned}    
\end{equation}

In problem (\ref{eq:App-03-DRO-Worst-Expectation-Primal}), the decision variables are the values of $f(\xi)$ over all possible $\xi \in \mathbb R^k$, so there are infinitely many decision variables, and problem (\ref{eq:App-03-DRO-Worst-Expectation-Primal}) is an infinite-dimensional LP. The former two constraints enforce $f(\xi)$ to be a valid distribution function; the latter two ensure consistent first- and second-order moments. The optimal solution gives the worst-case distribution. However, it is difficult to solve (\ref{eq:App-03-DRO-Worst-Expectation-Primal}) in its primal form. We now associate dual variables $\lambda_0 \in \mathbb R$, $\lambda \in \mathbb R^k$, and ${\rm \Lambda} \in \mathbb S^k$ with each integral constraint, and the dual problem of (\ref{eq:App-03-DRO-Worst-Expectation-Primal}) can be constructed following the duality theory of conic LP, which is given by
\begin{equation}
\label{eq:App-03-DRO-Worst-Expectation-Dual}
\begin{aligned}
\theta^m_D =  \inf_{\lambda_0,\lambda,{\rm \Lambda} } ~~ & \lambda_0 + 
\mu^T \lambda + \mbox{tr} [{\rm \Lambda}^T({\rm \Sigma} + \mu \mu^T)]  \\
\mbox{s.t.} ~~ & \lambda_0 + \xi^T \lambda + \mbox{tr} [{\rm \Lambda}^T 
(\xi \xi^T)] \\ 
& \quad \ge \max \{ 0,g(\xi) \}, \forall \xi \in \mathbb R^k
\end{aligned} 
\end{equation}
To understand this dual form in (\ref{eq:App-03-DRO-Worst-Expectation-Dual}), we can image a discrete version of (\ref{eq:App-03-DRO-Worst-Expectation-Primal}), in which $\xi_1$, $\cdots$, $\xi_n$ are sampled scenarios of the uncertain parameter, and their associated probabilities $f(\xi_1)$, $\cdots$, $f(\xi_n)$ are decision variables of (\ref{eq:App-03-DRO-Worst-Expectation-Primal}). Moreover, if we replace the integral arithmetic in the constraints with the summation arithmetic, (\ref{eq:App-03-DRO-Worst-Expectation-Primal}) comes down to a traditional LP, and its dual is also an LP, where the constraint becomes  
\begin{equation*}
\lambda_0 + \xi^T_i \lambda + \mbox{tr} [{\rm \Lambda}^T 
(\xi_i \xi^T_i)] \ge \max \{ 0,g(\xi_i) \},~ i = 1, \cdots, n 
\end{equation*}
Let $n \to +\infty$ and $\xi$ spread over $\mathbb R^k$, we can get the dual problem (\ref{eq:App-03-DRO-Worst-Expectation-Dual}).

Unlike the primal problem (\ref{eq:App-03-DRO-Worst-Expectation-Primal}) that has infinite decision variables, the dual problem (\ref{eq:App-03-DRO-Worst-Expectation-Dual}) has finite variables and an infinite number of constraints. In fact, we are optimizing over the coefficients of a polynomial in $\xi$. Because $\rm \Sigma \succ 0$, Slater condition is met, and thus strong duality holds (this conclusion can be found in many other literatures, such as \cite{Zero-Gap-GPI}), i.e., $\theta^m_P = \theta^m_D$. In the following, we will eliminate $\xi$ and reduce the constraint into convex ones in dual variables $\lambda_0$, $\lambda$, and $\rm \Lambda$. Recall the definition of matrix $\rm \Omega$, the compact form of problem (\ref{eq:App-03-DRO-Worst-Expectation-Dual}) can be expressed as
\begin{equation}
\label{eq:App-03-DRO-Worst-Expectation-Dual-Comp}
\begin{aligned}
\inf_{M \in \mathbb S^{k+1}} ~~ & \mbox{tr} [ {\rm \Omega}^T M ]  \\
\mbox{s.t.} ~~ & 
\begin{bmatrix} \xi^T & 1 \end{bmatrix} M 
\begin{bmatrix} \xi^T & 1 \end{bmatrix}^T \ge 0,~ 
\forall \xi \in \mathbb R^k \\ 
& \begin{bmatrix} \xi^T & 1 \end{bmatrix} M 
\begin{bmatrix} \xi^T & 1 \end{bmatrix}^T \ge g(\xi),~ 
\forall \xi \in \mathbb R^k \\ 
\end{aligned} 
\end{equation}
where the matrix decision variable is
\begin{equation}
M = \begin{bmatrix}
{\rm \Lambda}  &  \dfrac{\lambda}{2}  \\  
\dfrac{\lambda^T}{2}  &  \lambda_0
\end{bmatrix} \notag
\end{equation}
and the first constraint is equivalent to an LMI $M \succeq 0$.

A special case of  the worst-case expectation problem is 
\begin{equation}
\label{eq:App-03-DRO-GPI-Primal}
\theta^m_P = \sup_{f(\xi) \in \mathcal P} \Pr [ \xi \in S ] 
\end{equation}
which quantifies the maximum probability of the event $\xi \in S$, where $S$ is a Borel measurable set. This problem has a close relationship with generalized probability inequalities discussed in \cite{Zero-Gap-GPI} and the generalized moments problem studied in \cite{Moment-Book}. 
By defining an indicator function as
\begin{equation}
\mathbb I_{S} (\xi) = \left\{ 
\begin{gathered} 1 \\ 0 \end{gathered} \quad 
\begin{lgathered}
\mbox{if } \xi \in S \\
\mbox{otherwise} 
\end{lgathered}  \right. \notag
\end{equation}
The dual problem of (\ref{eq:App-03-DRO-GPI-Primal}) can be written as
\begin{equation}
\label{eq:App-03-DRO-GPI-Dual}
\begin{aligned}
\inf_{M \in \mathbb S^{k+1}} ~~ & \mbox{tr} [ {\rm \Omega}^T M ]  \\
\mbox{s.t.} ~~ & M \succeq 0,~  
\begin{bmatrix} \xi^T & 1 \end{bmatrix} M 
\begin{bmatrix} \xi^T & 1 \end{bmatrix} \ge 1,~ 
\forall \xi \in S \\ 
\end{aligned} 
\end{equation}
which is a special case of (\ref{eq:App-03-DRO-Worst-Expectation-Dual-Comp}) when $g(\xi) = \mathbb I_S(\xi)$.

Next we present how to formulate a robust chance constraint
(\ref{eq:App-03-DRO-DRCC}) as convex constraints that can be recognized by convex optimization solvers.

\vspace{12pt}
{\noindent \bf 1. Individual chance constraints}

Consider a single robust chance constraint
\begin{equation}
\label{eq:App-03-DRO-DRCC-Single}
\inf_{f(\xi) \in \mathcal P} \Pr \left( y^0 (x) + y(x)^T \xi \le 0 \right) \ge 1 - \varepsilon
\end{equation}
The feasible set in $x$ is denoted by $X^S_R$. 

To eliminate the optimization over function $f(\xi)$, we leverage the concept of conditional value-at-risk (CVaR) introduced by \cite{CVaR}. For a given loss function $L(\xi)$ and tolerance $\varepsilon \in (0,1)$, the CVaR at level $\varepsilon$ is defined as
\begin{equation}
\label{eq:App-03-CVaR}
\mbox{CVaR}(L(\xi), \varepsilon) = \inf_{\beta \in \mathbb R} \beta + 
\frac{1}{\varepsilon} \mathbb E_{f(\xi)} \left( \left[ L(\xi) - \beta \right]^+ \right)
\end{equation}
where the expectation is taken over a given probability distribution $f(\xi)$. CVaR is the conditional expectation of loss greater than the $(1-\varepsilon)$-quantile of the loss distribution. Indeed, condition
\begin{equation}
\Pr \left[ L(\xi) \le \mbox{CVaR}(L(\xi), \varepsilon) \right] \ge 1- \varepsilon
\notag
\end{equation}
holds regardless of the probability distribution and loss function $L(\xi)$ \cite{Static-DRO}. Therefore, to certify $\Pr(L(\xi) \le 0) \ge 1-\varepsilon$, a sufficient condition without probability evaluation is $\mbox{CVaR}(L(\xi), \varepsilon) \le 0$, or more precisely: 
\begin{equation}
\label{eq:App-03-CVaR-CC}
\begin{aligned}
& \sup_{f(\xi) \in \mathcal P} \mbox{CVaR} \left( y^0(x) + y(x)^T \xi, \varepsilon \right) \le 0  \\  \Longrightarrow
& \inf_{f(\xi) \in \mathcal P} \Pr \left( y^0(x) + y(x)^T \xi \le 0 \right) \ge 1 - \varepsilon
\end{aligned}
\end{equation}

According to (\ref{eq:App-03-CVaR}), above worst-case CVaR can be expressed by
\begin{equation}
\begin{lgathered}
\sup_{f(\xi) \in \mathcal P} \mbox{CVaR} \left( y^0(x) + y(x)^T \xi, \varepsilon \right)  \\
= \sup_{f(\xi) \in \mathcal P} \inf_{\beta \in \mathbb R} \left\{ \beta +
\frac{1}{\varepsilon} \mathbb E_{f(\xi)} \left( \left[ y^0(x) + y(x)^T \xi - \beta \right]^+ \right) \right\}  \\
= \inf_{\beta \in \mathbb R} \left\{ \beta + \frac{1}{\varepsilon} \sup_{f(\xi) \in \mathcal P} \mathbb E_{f(\xi)} \left( \left[ y^0(x) + y(x)^T \xi - \beta \right]^+ \right) \right\}
\end{lgathered}
\end{equation}
The maximization and minimization operators are interchangeable because of the saddle point theorem in \cite{Saddle-Point}. Recall previous analysis; the worst-case expectation can be computed from problem
\begin{equation}
\begin{aligned}
\inf_{\beta,M \in \mathbb S^{k+1}} ~~ & \mbox{tr} [ {\rm \Omega}^T M ]  \\
\mbox{s.t.} ~~ & M \succeq 0,~  \\
& 
\begin{bmatrix} \xi^T & 1 \end{bmatrix} M 
\begin{bmatrix} \xi^T & 1 \end{bmatrix} \ge 
y^0(x) + y(x)^T \xi - \beta,~ 
\forall \xi \in \mathbb R^k \\ 
\end{aligned}    \notag
\end{equation}
The semi-infinite constraint has a matrix quadratic form
\begin{equation}
\begin{bmatrix} \xi \\ 1 \end{bmatrix}^T 
\left( M - \begin{bmatrix}
0 & \dfrac{y(x)}{2} \\
\dfrac{y(x)^T}{2}  & y^0(x) - \beta
\end{bmatrix}  \right) 
\begin{bmatrix} \xi \\ 1 \end{bmatrix} \ge 
0,~ \forall \xi \in \mathbb R^k  \notag 
\end{equation}
which is equivalent to
\begin{equation}
M - \begin{bmatrix}
 0  &  \dfrac{y(x)}{2}  \\
\dfrac{y(x)^T}{2}  & y^0(x) - \beta
\end{bmatrix} \succeq 0  \notag 
\end{equation}
As a result, the worst-case CVaR can be calculated from an SDP
\begin{equation}
\label{eq:App-03-Worst-CVaR-SDP}
\begin{aligned}
\sup_{f(\xi) \in \mathcal P} ~~ & \mbox{CVaR} \left( y^0(x) + y(x)^T \xi, \varepsilon \right)  \\
= \inf_{\beta, M} ~~ & \beta + \frac{1}{\varepsilon} 
\mbox{tr}({\rm \Omega}^T M)   \\
\mbox{s.t.} ~~ &  M \succeq 0 \\
& M \succeq \begin{bmatrix}
0 & \dfrac{y(x)}{2}  \\
\dfrac{y(x)^T}{2} & y^0(x) - \beta
\end{bmatrix} 
\end{aligned}
\end{equation}

It is shown that the indicator $\Rightarrow$ in (\ref{eq:App-03-CVaR-CC}) is in fact an equivalence $\Leftrightarrow$ \cite{Static-DRO} in static DRO. In conclusion, robust chance constraint (\ref{eq:App-03-DRO-DRCC-Single}) can be written as a convex set in variable $x$, $\beta$, and $M$ as follows
\begin{equation}
\label{eq:App-03-DRO-XSR}
X^S_R = \left\{ x ~\middle|~ \begin{lgathered}
\exists  \beta \in \mathbb R,~ M \succeq 0 \mbox{ such that} \\
\beta + \frac{1}{\varepsilon} 
\mbox{tr}({\rm \Omega}^T M) \le 0  \\
M \succeq \begin{bmatrix}
0 & \dfrac{y(x)}{2}  \\
\dfrac{y(x)^T}{2}  & y^0(x) - \beta
\end{bmatrix}
\end{lgathered} \right\}
\end{equation}

\vspace{12pt}
{\noindent \bf 2. Joint chance constraints}

Now consider the joint robust chance constraints
\begin{equation}
\label{eq:App-03-DRO-DRCC-Joint}
\inf_{f(\xi) \in \mathcal P} \Pr \left( y^0_i (x) + y_i(x)^T \xi \le 0, 
i = 1,\cdots,m \right) \ge 1 - \varepsilon
\end{equation}
The feasible set in $x$ is denoted by $X^J_R$. 

Let $\alpha$ be the vector of strictly positive scaling parameters, and $\mathcal A = \{ \alpha ~|~ \alpha > 0 \}$. It is clear that constraint
\begin{equation}
\label{eq:App-03-DRO-DRCC-Joint-Para}
\inf_{f(\xi) \in \mathcal P} \Pr \left[ \max_{i=1,\cdots, m} \left\{ \alpha_i \left( y^0_i (x) + y_i(x)^T \xi \right) \right\} \le 0 \right] \ge 1 - \varepsilon
\end{equation}
imposes the same feasible region in variable $x$ as (\ref{eq:App-03-DRO-DRCC-Joint}). Nonetheless, it turns out that parameter $\alpha_i$ can be co-optimized to improve the quality of the convex approximation for $X^J_R$. (\ref{eq:App-03-DRO-DRCC-Joint-Para}) is a single robust chance constraint, and can be conservatively approximated by a worst-case CVaR constraint  
\begin{equation}
\label{eq:App-03-DRO-DRCC-Joint-Para-Worst-CVaR}
\sup_{f(\xi) \in \mathcal P} \mbox{CVaR} \left[ \max_{i=1,\cdots, m} \left\{ \alpha_i \left( y^0_i (x) + y_i(x)^T \xi \right) \right\}, \varepsilon \right] \le 0
\end{equation}
It defines a feasible region in variable $x$ with auxiliary parameter $\alpha \in \mathcal A$, which is denoted by $X^J_R(\alpha)$. Clearly, $X^J_R(\alpha) \subseteq X^J_R$, $\forall \alpha \in \mathcal A$. Unlike (\ref{eq:App-03-Worst-CVaR-SDP}), condition (\ref{eq:App-03-DRO-DRCC-Joint-Para-Worst-CVaR}) is $\alpha$-dependent.

By observing the fact that
\begin{equation}
\begin{aligned}
& \begin{bmatrix} \xi^T & 1 \end{bmatrix} M 
  \begin{bmatrix} \xi^T & 1 \end{bmatrix} \ge 
\max_{i=1,\cdots,m} \left\{ \alpha_i \left( y^0_i (x) + y_i(x)^T \xi \right) \right\} - \beta,~ \forall \xi \in \mathbb R^k  \\ \Longleftrightarrow
& \begin{bmatrix} \xi^T & 1 \end{bmatrix} M 
  \begin{bmatrix} \xi^T & 1 \end{bmatrix} \ge 
  \alpha_i \left( y^0_i (x) + y(x)_i^T \xi \right) - \beta,~ 
\forall \xi \in \mathbb R^k,~ i = 1,\cdots,m \\  \Longleftrightarrow
& M - \begin{bmatrix}
0  &  \dfrac{\alpha_i y_i(x)}{2}  \\
\dfrac{\alpha_i y_i(x)^T}{2} & \alpha_i y^0_i(x) - \beta
\end{bmatrix} \succeq 0, ~  i = 1,\cdots, m
\end{aligned} \notag
\end{equation}
and employing the optimization formulation of the worst-case expectation problem, the worst-case CVaR in (\ref{eq:App-03-DRO-DRCC-Joint-Para-Worst-CVaR}) can be calculated by
\begin{equation}
\label{eq:App-03-DRO-DRCC-Joint-Para-Worst-CVaR-SDP}
\begin{aligned}
& J(x,\alpha) = \sup_{f(\xi) \in \mathcal P} \mbox{CVaR} \left[ \max_{i=1,\cdots, m} \left\{ \alpha_i \left( y^0_i (x) + y_i(x)^T \xi \right) \right\}, \varepsilon \right]    \\  =  &
\inf_{\beta \in \mathbb R} \left\{ \beta + \frac{1}{\varepsilon} \sup_{f(\xi) \in \mathcal P}  \mathbb E_{f(\xi)}  \left( \left[ \max_{i=1,\cdots,m} \left\{ \alpha_i \left( y^0_i(x) + y_i(x)^T \xi \right) \right\} - \beta \right]^+ \right) \right\}   \\  =  &
\inf_{\beta, M} \left\{ \beta + \frac{1}{\varepsilon} \mbox{tr}({\rm \Omega}^T M) ~\middle|~ \mbox{s.t.} ~   M \succeq 0,~ M \succeq 
\begin{bmatrix}
0 & \dfrac{\alpha_i y_i(x)}{2} \\
\dfrac{\alpha_i y_i(x)^T}{2} & y^0_i(x) - \beta
\end{bmatrix}, \forall i  \right\}
\end{aligned}
\end{equation}

In conclusion, for any fixed $\alpha \in \mathcal A$, the worst-case CVaR constraint (\ref{eq:App-03-DRO-DRCC-Joint-Para-Worst-CVaR}) can be written as a convex set in variables $x$, $\beta$, and $M$ as follows
\begin{equation}
\label{eq:App-03-DRO-XJR-alpha}
X^J_R(\alpha) = \left\{ x ~\middle|~ \begin{lgathered}
\exists  \beta \in \mathbb R,~ M \succeq 0 \mbox{ such that} \\
\beta + \frac{1}{\varepsilon} 
\mbox{tr}({\rm \Omega}^T M) \le 0  \\
M \succeq 
\begin{bmatrix}
0 & \dfrac{\alpha_i y_i(x)}{2} \\
\dfrac{\alpha_i y_i(x)^T}{2} & y^0_i(x) - \beta
\end{bmatrix}, \forall i 
\end{lgathered} \right\}
\end{equation}

Moreover, it is revealed in \cite{Static-DRO} that the union $\bigcup_{\alpha \in \mathcal A} X^J_R(\alpha)$ gives an exact description of $X^J_R$, which indicates that the original robust chance constrained program 
\begin{equation}
\label{eq:App-03-DRO-Model-2}
\min_x \left\{ c^T x  ~\middle|~ \mbox{s.t. } x \in X \cap X^J_R \right \}
\end{equation}
and the worst-case CVaR formulation
\begin{equation}
\min_{x,\alpha} \left\{ c^T x  ~\middle|~ \mbox{s.t. } 
x \in X \cap X^J_R(\alpha),~ \alpha \in \mathcal A \right \} \notag
\end{equation}
or equivalently 
\begin{equation}
\label{eq:App-03-DRO-Model-3}
\min_{x,\alpha} \left\{ c^T x  ~\middle|~ \mbox{s.t. } 
x \in X,~ \alpha \in \mathcal A,~ J(x,\alpha) \le 0 \right \}
\end{equation}
have the same optimal value. The constraints of (\ref{eq:App-03-DRO-Model-3}) contain bilinear matrix inequalities, which means that if either $x$ or $\alpha$ is fixed, $J(x,\alpha) \le 0$ in (\ref{eq:App-03-DRO-Model-3}) can come down to LMIs, however, when both $x$ and $\alpha$ are variables, the constraint is non-convex, making problem (\ref{eq:App-03-DRO-Model-3}) difficult to solve. In view of the biconvex feature \cite{{BLP-Mountain-Climbing-BCVX}}, a sequential convex optimization procedure is presented to find an approximated solution.
\begin{algorithm}[!htp]
\normalsize
\caption{\bf }
\begin{algorithmic}[1]
\STATE Choose a convergence tolerance $\varepsilon > 0$; Let the iteration counter $k=1$, $x^0 \in X \cap X^J_R(\alpha)$ be a feasible solution for some $\alpha$ and $f^0 = c^T x^0$;
\STATE Solve the following subproblem with input $x^{k-1}$
\begin{equation}
\label{eq:App-03-BMI-Sub}
\min_\alpha ~ \left\{ J(x,\alpha) ~|~ 
\mbox{s.t. } \alpha \ge \delta \bf 1 \right\} 
\end{equation}
where $\bf 1$ denotes the all-one vector with a compatible dimension, and $\delta > 0$ is a small constant; the worst-case CVaR functional is defined in (\ref{eq:App-03-DRO-DRCC-Joint-Para-Worst-CVaR-SDP}). The optimal solution is $\alpha^k$;
\STATE Solve the following master problem with input $\alpha^k$
\begin{equation}
\label{eq:App-03-BMI-Master}
\min_x ~ \left\{ c^T x ~|~ \mbox{s.t. } x \in X,~ J(x,\alpha^k) \le 0\right\}
\end{equation}
The optimal solution is $x^k$ and the optimal value is $f^k$;
\STATE If $|f^k - f^{k-1}| / |f^{k-1}| \le \varepsilon$, terminate and report the optimal solution $x^k$; otherwise, update $k \leftarrow k+1$, and go to step 2.
\end{algorithmic}
\label{Ag:App-03-BMI-Mountain-Climbing}
\end{algorithm}  

The main idea of this algorithm is to identify the best feasible region $X^J_R(\alpha)$ through successively solving the subproblem (\ref{eq:App-03-BMI-Sub}), and therefore improving the objective value. The performance of Algorithm \ref{Ag:App-03-BMI-Mountain-Climbing} is intuitively explained below. 

Because parameter $\alpha$ is optimized in the subproblem (\ref{eq:App-03-BMI-Sub}) given the value $x^k$, there must be $J(x^k,\alpha^{k+1}) \le J(x^k,\alpha^k) \le 0$, $\forall k$, demonstrating that $x^k$ is a feasible solution of the master problem (\ref{eq:App-03-BMI-Master}) in iteration $k+1$; therefore, the optimal values of  (\ref{eq:App-03-BMI-Master}) in two consecutive iterations satisfy $c^T x^{k+1} \le c^T x^k$, as the objective evaluated at the optimal solution $x^{k+1}$ in iteration $k+1$ deserves a value no greater than that is incurred at any feasible solution. In this regard, the optimal value sequence $f^k$, $k=1,2,\cdots$ is monotonically decreasing. If $X$ is bounded, the optimal solution sequence $x^k$ is also bounded, and the optimal value converges. Algorithm \ref{Ag:App-03-BMI-Mountain-Climbing} does not necessarily find the global optimum of problem (\ref{eq:App-03-DRO-Model-3}). Nevertheless, it is desired by practical problems due to its robustness since it involves only convex optimization.

In many practical applications, the uncertain data $\xi$ is known to be within a strict subset of $\mathbb R^k$, which is called the support set. We briefly outline how to incorporate the support set in the distributionally robust chance constraints. We assume the support set $\rm \Xi$ is the intersection of a finite number of ellipsoids, i.e. 
\begin{equation}
\label{eq:App-03-Support-Set}
{\rm \Xi} = \left\{ \xi \in \mathbb R^k ~\middle|~ 
\xi^T  W_i \xi  \le 1,~ i = 1,\cdots,l  \right\}   
\end{equation}
where $W_i \in \mathbb S^k_+$, $i=1,\cdots,l$, and we have $\Pr(\xi \in {\rm \Xi}) = 1$. Let $\mathcal P_{\rm \Xi}$ be the set of all candidate probability
distributions supported on $\rm \Xi$ which have identical first- and second-order moments.   

Consider the worst-case expectation problem (\ref{eq:App-03-DRO-Worst-Expectation-Primal}). If we replace $\mathcal P$ with $\mathcal P_{\rm \Xi}$, the constraints of the dual problem (\ref{eq:App-03-DRO-Worst-Expectation-Dual}) become
\begin{alignat}{2}
\begin{bmatrix} \xi^T & 1 \end{bmatrix} M 
\begin{bmatrix} \xi^T & 1 \end{bmatrix}^T & \ge 0,~ & \quad
& \forall \xi  \in {\rm \Xi} \label{eq:App-03-S-Lemma-1} \\ 
 \begin{bmatrix} \xi^T & 1 \end{bmatrix} M 
\begin{bmatrix} \xi^T & 1 \end{bmatrix}^T & \ge g(\xi),~ &
& \forall \xi \in {\rm \Xi}  \label{eq:App-03-S-Lemma-2} 
\end{alignat}
According to (\ref{eq:App-03-Support-Set}), $1-\xi^T  W_i \xi$ must be non-negative if and only if $\xi \in {\rm \Xi}$, and hence a sufficient condition for (\ref{eq:App-03-S-Lemma-1}) is the existence of constants $\tau_i \ge 0$, $i = 1,\cdots,l$, such that
\begin{equation}
\label{eq:App-03-S-Lemma-3} 
\begin{bmatrix} \xi^T & 1 \end{bmatrix} M 
\begin{bmatrix} \xi^T & 1 \end{bmatrix}^T - 
\sum_{i=1}^l \tau_i \left( 1-\xi^T  W_i \xi \right) \ge 0
\end{equation}

Under this condition, as long as $\xi \in {\rm \Xi}$, we have 
\begin{equation*}
\begin{bmatrix} \xi^T & 1 \end{bmatrix} M 
\begin{bmatrix} \xi^T & 1 \end{bmatrix}^T \ge  
\sum_{i=1}^l \tau_i \left( 1-\xi^T  W_i \xi \right) \ge 0 
\end{equation*}

Arrange (\ref{eq:App-03-S-Lemma-3}) as a matrix quadratic form
\begin{equation*} 
\begin{bmatrix} \xi^T & 1 \end{bmatrix} 
\left(M - \sum_{i=1}^l \tau_i \begin{bmatrix}
-W_i & {\bf 0}  \\
{\bf 0}^T & 1
\end{bmatrix} \right)
\begin{bmatrix} \xi \\ 1 \end{bmatrix}  \ge 0,~
\forall \xi \in \mathbb R^k 
\end{equation*}

As a result, (\ref{eq:App-03-S-Lemma-1}) can be reduced to an LMI in variables $M$ and $\tau$
\begin{equation}
\label{eq:App-03-S-Lemma-4}
M - \sum_{i=1}^l \tau_i 
\begin{bmatrix}
-W_i & {\bf 0}  \\
{\bf 0}^T & 1
\end{bmatrix}   \succeq 0
\end{equation}
For similar reasons, by letting $g(\xi) = y^0(x) + y(x)^T \xi - \beta$, (\ref{eq:App-03-S-Lemma-2}) can be conservatively approximated by the following LMI
\begin{equation}
\label{eq:App-03-S-Lemma-5}
M - \sum_{i=1}^l \tau_i 
\begin{bmatrix}
-W_i & {\bf 0}  \\
{\bf 0}^T & 1
\end{bmatrix}  \succeq 
\begin{bmatrix}
0 & \frac{1}{2} y(x) \\
\frac{1}{2} y(x)^T & y^0(x) - \beta
\end{bmatrix} 
\end{equation}
In fact, (\ref{eq:App-03-S-Lemma-4}) and (\ref{eq:App-03-S-Lemma-5}) are special cases of S-Lemma. Based upon these outcomes, most formulations in this section can be extended to consider the bounded support set $\rm \Xi$ in the form of (\ref{eq:App-03-Support-Set}). For polyhedral and some special classes of convex support sets, one may utilize the nonlinear Farkas lemma (Lemma 2.2 in \cite{Static-DRO}) to derive tractable reformulations.

\subsection{Adjustable Distributionally Robust Optimization}
\label{App-C-Sect03-02}

As explained in Appendix \ref{App-C-Sect01}, the traditional static RO encounters difficulties in dealing with equality constraints. This plight remains in the DRO approach following a static setting. Consider $x + \xi = 1$ where $\xi \in [0,0.1]$ is uncertain, while its mean and variance are known. For any given $x^*$, the worst-case probability $\inf_{f(\xi) \in \mathcal P} \Pr[x^* + \xi = 1] =0$, because one can always find a feasible probability distribution function $f(\xi)$ that satisfies the first- and second-order moment constraints, whereas $f(1-x^*) = 0$.   

To vanquish this difficulty, it is necessary to incorporate wait-and-see decisions. A simple remedy is to impose an affine recourse policy without involving optimization in the second stage, giving rise to an affine-adjustable RO with distributional uncertainty and linear decision rule, which can be solved by the method in Appendix \ref{App-C-Sect03-01}. 

This section aims to investigate the following adjustable DRO with completely flexible wait-and-see decisions 
\begin{equation}
\label{eq:App-03-ADRO-Model-1}
\min_{x \in X} \left\{ c^T x + \sup_{f(w) \in \mathcal P} 
\mathbb E_{f(w)} Q(x,w) \right\}
\end{equation}
where $x$ is the first-stage (here-and-now) decision, and $X$ is its feasible set; the uncertain parameter is denoted by $w$; the probability distribution $f(w)$ belongs to the Chebyshev ambiguity set (whose first- and second-order moments are known)
\begin{equation}
\label{eq:App-03-ADRO-Ambiguity-Set}
\mathcal P =\left\{ f(w) \middle| \begin{gathered}
 f(w) \ge 0,~\forall w \in W   \\
 \int_{w \in W} f(w) \mbox{d} w = 1    \\
 \int_{w \in W} w f(w)\mbox{d} w = \mu  \\
 \int_{w \in W} w w^T f(w)\mbox{d}w = {\rm \Theta}
\end{gathered}  \right\}   
\end{equation}
supported on $W = \{w~|~ (w - \mu)^T Q (w-\mu) \le {\rm \Gamma} \}$, where matrix ${\rm \Theta} = {\rm \Sigma} + \mu \mu^T$ represents the second-order moment; $\mu$ is the mean value and $\rm \Sigma$ is the covariance matrix. The expectation in (\ref{eq:App-03-ADRO-Model-1})  is taken over the worst-case $f(w)$ in $\mathcal P$, and the second-stage problem under fixed $x$ and $w$ is an LP
\begin{equation}
\label{eq:App-03-ADRO-Model-2}
 Q(x,w) = \min_{y \in Y(x,w)} d^T y 
\end{equation}
$Q(x,w)$ is its optimal value function under fixed $x$ and $w$. The feasible set of the second-stage problem is 
\begin{equation}
Y(x,w) = \{ y ~|~ B y \le b - A x -C w \} \notag
\end{equation}
Matrices $A$, $B$, $C$ and vectors $b$, $c$, $d$ are constant coefficients in the model. We assume that the second-stage problem is always feasible, i.e.,  $\forall x \in X$, $\forall w \in W:$ $Y(x,w) \ne \emptyset$ and is bounded, and thus $Q(x,w)$ has a finite optimal value. This can be implemented by introducing wait-and-see type slack variables and adding penalties in the objective of (\ref{eq:App-03-ADRO-Model-2}).

The difference between problems (\ref{eq:App-03-ARO}) and (\ref{eq:App-03-ADRO-Model-1}) stems from the descriptions of uncertainty and the criteria in the objective function: more information of the dispersion effect, such as the covariance matrix, is taken into account in the latter one, and the objective function in (\ref{eq:App-03-ADRO-Model-1}) is an expectation reflecting the statistical behavior of the second-stage cost, rather than the one in (\ref{eq:App-03-ARO}) which is associated with only a single worst-case scenario, and leaves the performances in all other scenarios un-optimized. Because the probability distribution is uncertain, it is prudent to investigate the worst-case outcome in which the expected cost of the second stage is maximized. This formulation is advantageous in several ways: first, the requirement on the exact probability distribution is not necessary, and the optimal solution is insensitive to the family of distributions with common mean and covariance; second, the dispersion of the uncertainty is also taken into account, which helps reduce model conservatism: since the variance is fixed, a scenario that leaves far away from the forecast would have a low probability; finally, it is often important to tackle the tail effect, which indicates that the occurrence of a rare event may induce heavy losses in spite of its low probability. Such phenomenon is naturally taken into account in (\ref{eq:App-03-ADRO-Model-1}). In what follows, we outline the method proposed in \cite{App03-Sect3-ADRO-1} to solve the adjustable DRO problem (\ref{eq:App-03-ADRO-Model-1}). A slight modification is that an ellipsoid support set is considered.

\vspace{12pt}
{\noindent \bf 1. The worst-case expectation problem}

We consider the following worst-case expectation problem with a fixed $x$
\begin{equation}
\label{eq:App-03-ADRO-Sub-Worst-Expectation}
\sup_{f(w) \in \mathcal P} 
\mathbb E_{f(w)} Q(x,w)
\end{equation}

According to the discussions for problem (\ref{eq:App-03-DRO-Worst-Expectation-Primal}), the dual problem of (\ref{eq:App-03-ADRO-Sub-Worst-Expectation}) is
\begin{equation}
\label{eq:App-03-ADRO-Sub-Worst-Expectation-Dual}
\begin{aligned}
\min_{H,h,h_0}~~ & \mbox{tr}(H^T {\rm \Theta}) + \mu^T h + h_0 \\ 
\mbox{s.t.}~~ & w^T H w + h^T w + h_0 \ge Q(x,w),~\forall w \in W 
\end{aligned} 
\end{equation}
where $H$, $h$, $h_0$ are dual variables. Nevertheless, the optimal value function $Q(x,w)$ is not given in a closed form. From the LP duality theory
\begin{equation}
Q(x,w) = \max_{u \in U} ~ u^T (b - A x - C w)  \notag
\end{equation}
where $u$ is the dual variable of LP (\ref{eq:App-03-ADRO-Model-2}), and its feasible set is given by
\begin{equation}
U = \{ u ~|~ B^T u = d,~ u \le 0 \} \notag
\end{equation}
Because we have assumed that $Q(x,w)$ is bounded, the optimal solution of the dual problem can be found at one of the extreme points of
$U$, i.e.,
\begin{equation}
\label{eq:App-03-ADRO-Sub-Critical-Vertex}
\exists u^* \in \mbox{vert}(U): \ Q(x,w) = (b-Ax-Cw)^T u^* 
\end{equation}
where $\mbox{vert}(U) = \{u^1,u^2,\cdots,u^{N_E }\}$ stands for the vertices of  polyhedron $U$, and $N_E = |\mbox{vert}(U)|$ is the cardinality of vert($U$). In view of this, the constraint of (\ref{eq:App-03-ADRO-Sub-Worst-Expectation-Dual}) can be expressed as 
\begin{equation}
w^T H w + h^T w + h_0 \ge (b-Ax-Cw)^T u^i,~
\forall w \in W,~ i=1,\cdots,N_E  \notag
\end{equation}
Recall the definition of $W$; a certification for above condition is
\begin{equation}
\begin{aligned}
w^T H w & + h^T w + h_0 - (b-Ax-Cw)^T u^i \\ 
  & \ge \lambda [{\rm \Gamma} - (w-\mu)^T Q (w-\mu)] \ge 0,~ 
  \forall w \in \mathbb R^k,~ i=1,\cdots,N_E  
\end{aligned}    \notag
\end{equation}
which has the following compact matrix form
\begin{equation}
\label{eq:App-03-ADRO-Sub-Cons-LMI}
\begin{bmatrix} w \\  1  \end{bmatrix}^T  M^i 
\begin{bmatrix} w \\  1  \end{bmatrix} \ge 0, 
\forall w \in \mathbb R^k,~ i=1,\cdots,N_E  
\end{equation}
where
\begin{equation}
\label{eq:App-03-ADRO-Sub-Cons-LMI-MI}
M^i = \begin{bmatrix}
 H + \lambda Q  & \dfrac{h-C^T u^i}{2} - \lambda Q \mu  \\
\dfrac{h^T- (u^i)^T C}{2}-\lambda \mu^T Q & h_0 -(b-Ax)^T u^i - \lambda ({\rm \Gamma} - \mu^T Q \mu) 
\end{bmatrix}  
\end{equation}
and (\ref{eq:App-03-ADRO-Sub-Cons-LMI}) simply reduces to $M^i \succeq 0$, $i=1,\cdots,N_E$.

Finally, problem (\ref{eq:App-03-ADRO-Sub-Worst-Expectation-Dual}) comes down to the following SDP
\begin{equation}
\label{eq:App-03-ADRO-Sub-Worst-Expectation-Dual-SDP}
\begin{aligned}
\min_{H,h,h_0,\lambda}~~ & \mbox{tr}(H^T {\rm \Theta}) + \mu^T h + h_0 \\
\mbox{s.t.}~~ & M^{i} (H,h,h_0,\lambda) \succeq 0,~ i=1,\cdots,N_E   \\
& \lambda \in \mathbb R^+ 
\end{aligned} 
\end{equation}
where $M^{i} (H,h,h_0,\lambda)$ is defined in (\ref{eq:App-03-ADRO-Sub-Cons-LMI-MI}). Above results can be readily extended if the support set is the intersection of ellipsoids.

\vspace{12pt}
{\noindent \bf 2. Adaptive constraint generation algorithm}

Due to the positive semi-definiteness of the covariance matrix $\rm \Sigma$, the duality gap between problems (\ref{eq:App-03-ADRO-Sub-Worst-Expectation}) and (\ref{eq:App-03-ADRO-Sub-Worst-Expectation-Dual}) is zero \cite{App03-Sect3-ADRO-1}, and hence we can replace the worst-case expectation in (\ref{eq:App-03-ADRO-Model-1}) with its dual form, yielding  
\begin{equation}
\label{eq:App-03-ADRO-Model-3}
\begin{aligned}
\min ~~ & c^T x +  \mbox{tr}(H^T {\rm \Theta}) + \mu^T h + h_0 \\
\mbox{s.t.}~~ & M^{i} (H,h,h_0,\lambda) \succeq 0,~ i=1,\cdots,N_E   \\
& x \in X,~ \lambda \in \mathbb R^+ 
\end{aligned}
\end{equation}

Problem (\ref{eq:App-03-ADRO-Model-3}) is an SDP. However, the number of vertices in set $U$ ($|\mbox{vert}(U)|$) may increase exponentially in the dimension of $U$. It is non-trivial to enumerate all of them. However, because of weak duality, only the one which is optimal in the dual problem provides an active constraint, as shown in (\ref{eq:App-03-ADRO-Sub-Critical-Vertex}), and the rest are redundant inequalities. To identify the critical vertex in (\ref{eq:App-03-ADRO-Sub-Critical-Vertex}), we solve problem (\ref{eq:App-03-ADRO-Model-3}) in iterations: in the master problem, a subset of $\mbox{vert}(U)$ is used to formulate a relaxation, then check whether the following constraint 
\begin{equation}
\label{eq:App-03-ADRO-Sub-Double-Enumeration}
w^T H w + h^T w + h_0 \ge (b-Ax-Cw)^T u,~
\forall w \in W,~ \forall u \in U
\end{equation}
is fulfilled. If yes, the relaxation is exact and the optimal solution is found; otherwise, find a new vertex of $U$ at which constraint (\ref{eq:App-03-ADRO-Sub-Double-Enumeration}) is violated, and then add a cut to the master problem so as to tighten the relaxation, till constraint (\ref{eq:App-03-ADRO-Sub-Double-Enumeration}) is satisfied. The flowchart is summarized in Algorithm \ref{Ag:App-03-ADRO-Delayed-Vertex-Generation}.
\begin{algorithm}[!htp]
\normalsize
\caption{\bf }
\begin{algorithmic}[1]

\STATE Choose a convergence tolerance $\epsilon > 0$ and an initial vertex set $V_E \subseteq \mbox{vert}(U)$.

\STATE Solve the following master problem 
\begin{equation}
\label{eq:App-03-ADRO-Model-AVG-Master}
\begin{aligned}
\min ~~ & c^T x +  \mbox{tr}(H^T {\rm \Theta}) + \mu^T h + h_0 \\
\mbox{s.t.}~~ & M^{i} (H,h,h_0,\lambda) \succeq 0,~ \forall u^i \in V_E \\
& x \in X,~ \lambda \in \mathbb R^+ 
\end{aligned}
\end{equation}
 The optimal value is $R^*$, and the optimal solution is $(x^*, H, h, h_0)$.

\STATE Solve the following sub-problem with obtained $(x^*, H, h, h_0)$
\begin{equation}
\label{eq:App-03-ADRO-Model-AVG-Sub}
\begin{aligned}
\min_{w,u}~~ & w^T H w + h^T w+h_0 - (b-Ax^* -Cw)^T u  \\
\mbox{s.t.}~~ & w \in W,~ u \in U
\end{aligned}
\end{equation}
The optimal value is $r^*$, and the optimal solution is $u^* $and $w^*$.

\STATE If $r^* \ge - \varepsilon$, terminate and report the optimal solution $x^*$and the optimal value $R^*$; otherwise, $V_E = V_E \cup u^*$, add an LMI cut $M (H,h,h_0,\lambda) \succeq 0$ associated with the current $u^*$ to the master problem (\ref{eq:App-03-ADRO-Model-AVG-Master}), and go to step 2.
\end{algorithmic}
\label{Ag:App-03-ADRO-Delayed-Vertex-Generation}
\end{algorithm}  

Algorithm \ref{Ag:App-03-ADRO-Delayed-Vertex-Generation} terminates in a finite number of iterations which is bounded by $|\mbox{vert}(U)|$. Actually, it will converge within a few iterations, because the sub-problem (\ref{eq:App-03-ADRO-Model-AVG-Sub}) in step 3 always identifies the most critical vertex in $\mbox{vert}(U)$. It is worth mentioning that the subproblem (\ref{eq:App-03-ADRO-Model-AVG-Sub}) is a non-convex program. Despite that it can be solved by general NLP solvers, we suggest three approaches with different computational complexity and optimality guarantees.

1. If the support set $W = \mathbb R^k$, it can be verified that matrix $M_i$ becomes
\begin{equation}
\begin{bmatrix}
 H & \dfrac{h + C^T u^i}{2} \\
\dfrac{(h + C^T u^i)^T}{2} & h_0 - (b - Ax)^T u^i
\end{bmatrix} \succeq 0 \notag
\end{equation}
Then there must be $H \succeq 0$, and non-convexity appears in the bilinear term $u^T C w$. In such circumstance, problem (\ref{eq:App-03-ADRO-Model-AVG-Sub}) can be solved via a mountain climbing method similar to Algorithm \ref{Ag:App-03-BLP-Mountain-Climbing} (but here the mountain is actually a pit because the objective is to be minimized).  

2. In the case that $W$ is an ellipsoid, above iterative approach is still applicable; however, the $w$-subproblem in which $w$ is to be optimized may become non-convex because $H$ may be indefinite. Since $u$ is fixed in the $w$-subproblem, non-convex term $w^T H w$ can be decomposed as the difference of two convex functions as $w^T (H + \alpha I) w - \alpha w^T w$, where $\alpha$ is a constant such that $H + \alpha I$ is positive-definite, and the $w$-subproblem can be solved by the convex-concave procedure elaborated in \cite{CCP-Boyd}, or any existing NLP solver.   

3. As a non-convex QP, problem (\ref{eq:App-03-ADRO-Model-AVG-Sub}) can be globally solved by the MILP method presented in Appendix \ref{App-A-Sect04}. This method could be time consuming with the growth in problem sizes.

\section{Data-driven Robust Stochastic Program}
\label{App-C-Sect04}

Most classical SO methods assume that the probability distribution of uncertain factors is exactly known, which is an input of the problem. However, such information heavily relies on historical data, and may not be available at hand or accurate enough. Using an inaccurate distribution in a classical SO model could lead to biased results. To cope with ambiguous probability distributions, a natural way is to consider a set of possible candidates derived from available data, instead of a single distribution, just as the moment-inspired ambiguity set used in DRO. In this section, we investigate some useful SO models with distributional uncertainty described by divergence ambiguity sets, which is referred to as robust SO. When the distribution is discrete, the distributional uncertainty is interpreted by the perturbation of probability value associated with each scenario; when the distribution is continuous, the distance of two density functions should be specified first. In this section, we consider $\rm \Phi$-divergence and Wasserstein metric based ambiguity sets. 

\subsection{Robust Chance Constrained Stochastic Program}
\label{App-C-Sect04-01}

We introduce robust chance-constrained stochastic programs with distributional robustness. The ambiguous PDF is modeled based on $\phi$-divergence, and the optimal solution provides constraint feasibility guarantee with desired probability even in the worst-case distribution.  In short, the underlying problem possesses the following features:  

1) The PDF is continuous and the constraint violation probability is a functional.

2) Uncertain parameters do not explicitly appear in the objective function.

Main results of this section come from \cite{App03-Sect4-RCCP}.

\vspace{12pt}
{\noindent \bf 1. Problem formulation}

In a traditional chance-constrained stochastic
linear program, the decision maker seeks a cost-minimum solution at which some certain constraints can be met with a given probability, yielding:
\begin{equation}
\label{eq:App-03-CCP-Model}
\begin{aligned}
\min ~~ & c^T x  \\
\mbox{s.t.}~~ & \Pr [C(x,\xi)] \ge 1-\alpha  \\
& x \in X 
\end{aligned}
\end{equation}
where $x$ is the vector of decision variables; $\xi$ is the vector of uncertain parameters, and the exact (joint) probability distribution is apparent to the decision maker; vector $c$ represents the cost coefficients; $X$ is a polyhedron that is independent of $\xi$; $\alpha$ is the risk level or the maximum allowed probability of constraint violation; $C(x,\xi)$ collects all uncertainty dependent constraints, whose general form is given by
\begin{equation}
\label{eq:App-03-CCP-Cons-Uncertain}
C(x,\xi) = \{\xi~|~ \exists y : A(\xi) x + B(\xi) y \le b(\xi) \} 
\end{equation}
where $A$, $B$, $b$ are constant coefficient matrices that may contain uncertain parameters; $y$ is a recourse action that can be made after $\xi$ is known. In the presence of $y$, we call (\ref{eq:App-03-CCP-Model}) a two-stage problem; otherwise, it is a single-stage problem if $y$ is null. We don't consider the cost of recourse actions in the objective function in its current form. In case of need, we can add the second-stage cost $d^T y(\xi)$ in the objective function, and $\xi$ is a specific scenario  which $y(\xi)$ corresponds to; for instance, robust optimization may consider a max-min cost scenario or a max-min regret scenario; traditional SO often tackles the expected second-stage cost $\mathbb E[d^T y(\xi)]$. We leave it to the end of this section to discuss how to deal with the second-stage cost in the form of worst-case expectation like (\ref{eq:App-03-ADRO-Sub-Worst-Expectation}), and show that the problem can be convexified under some technical assumptions.       

In the chance constraint, for a given $x$, the probability of constraint satisfaction can be evaluated for a particular probability distribution of $\xi$. Traditional studies on chance-constrained programs often assume that the distribution of $\xi$ is perfectly known. However, this assumption can be very strong because it requires a lot of historical data. Moreover, the optimal solution may be sensitive to the true distribution and thus highly suboptimal in practice. To overcome these difficulties, a prudent method is to consider a set of probability distributions  belonging to a pre-specified ambiguity set $D$, and require that the chance constraint should be satisfied under all possible distributions in $D$, resulting in the following robust chance-constrained programming problem:
\begin{equation}
\label{eq:App-03-RCCP-Model}
\begin{aligned}
\min ~~ & c^T x  \\
\mbox{s.t.}~~ & \inf_{f(\xi) \in D}  \Pr[C(x,\xi)] 
\ge 1-\alpha  \\
& x \in X 
\end{aligned}
\end{equation}
where $f(\xi)$ is the probability density function of random variable $\xi$.

The ambiguity set $D$ in (\ref{eq:App-03-RCCP-Model}) which includes distributional information can be constructed in a data-driven fashion, such as the moment based ones used in Appendix \ref{App-C-Sect03}. Please see \cite{Am-Set-Overview} for more information on establishing $D$ based on moment data and other structural properties, such as symmetry and unimodality. The tractability of (\ref{eq:App-03-RCCP-Model}) largely depends on the form of $D$. For example: if $D$ is built on the mean value and covariance matrix (which is called a Chebyshev ambiguity set), a single robust chance constraint can be reformulated as an LMI and a set of joint robust chance constraints can be approximated by BMIs \cite{Static-DRO}; probability of constraint violation under more general moment based ambiguity sets can be evacuated by solving conic optimization problems \cite{Am-Set-Overview}. 

A shortcoming of moment description is that it does not provide a direct measure on the distance between the candidate PDFs in $D$ and a reference distribution. Two PDFs with the same moments may differ a lot in other aspects. Furthermore, the worst-case distribution corresponding to a Chebyshev ambiguity set always puts more weights away from the mean value, subject to the variance. As such, the long-tail effect is a source of conservatism. In this section, we consider the confidence set built around a reference distribution. The motivation is:  the decision maker may have some knowledge on what distribution the uncertainty follows, although such a distribution could be inexact, and the true density function would not deviate far away from it. 

To describe distributional ambiguity in term of a PDF, the first problem is how to characterize the distance between two functions.  One common measure on the distance between density functions is the $\phi$-divergence, which is defined as \cite{Am-Set-Phi-Divergence} 
\begin{equation}
\label{eq:App-03-RCCP-Phi-Div}
D_\phi(f \| f_0) = \int_{\rm \Omega}  
\phi \left( \dfrac{f(\xi)}{f_0(\xi)} \right) f_0(\xi) d\xi
\end{equation} 
where $f$ and $f_0$ stand for the particular density function and the estimated one (or the reference distribution), respectively; function $\phi$ satisfies:
\begin{equation*}
\begin{lgathered}
\mbox{(C1)}~~ \phi(1) = 0 \\
\mbox{(C2)}~~ 0\phi(x/0) = \begin{cases}
x \lim_{p \to +\infty} \phi(p)/p & \mbox{if } x > 0  \\
0 &  \mbox{if } x = 0 
\end{cases} \\
\mbox{(C3)}~~ \phi(x)= +\infty~ \mbox{for } x<0 \\
\mbox{(C4)}~~ \phi(x)~ \mbox{is a convex function on}~ \mathbb R^+ 
\end{lgathered}
\end{equation*}
It is proposed in \cite{Am-Set-Phi-Divergence} that the ambiguity set can be built as:
\begin{equation}
\label{eq:App-03-Conf-Set-Phi-Div}
D = \{P: D_\phi(f \| f_0) \le d, f={\rm d} P/{\rm d} \xi\} 
\end{equation} 
where the tolerance $d$ can be adjusted by the decision maker according to their attitudes towards risks. The ambiguity set in (\ref{eq:App-03-Conf-Set-Phi-Div}) can be denoted as $D_{\phi}$, without causing confusion with the definition of $\phi$-divergence $D_\phi(f \| f_0) $.  Compared to the moment-based ambiguity sets, especially the Chebyshev ambiguity set, where only the first- and second-order moments are involved, the density based description captures the overall profile of the ambiguous distribution, so may hopefully provide less conservative solutions. However, it hardly guarantees consistent moments. Which one is better depends on data availability: if we are more confident on the reference distribution, (\ref{eq:App-03-Conf-Set-Phi-Div}) may be better; otherwise, if we only have limited statistic information such as mean and variance, then the moment-based ones are  more straightforward.

\begin{table}[!htp]
\footnotesize
\renewcommand{\arraystretch}{1.3}
\renewcommand{\tabcolsep}{1em}
\caption{Instances of $\phi$-divergences}
\centering
\begin{tabular}{cc}
\toprule 
     Divergence         &  function $\phi(x)$    \\
\midrule
   KL-divergence        &   $x \log x - x + 1$   \\
reverse KL-divergence   &   $- \log x$     \\
Hellinger distance      &   $({\sqrt x} -1)^2$   \\
Variation distance      &   $|x-1|$     \\
    J-divergence        &   $(x-1)\log x$  \\
 $\chi^2$ divergence    &   $(x-1)^2$      \\
 $\alpha$-divergence    &   
 $\begin{cases}
  \dfrac{4}{1-\alpha^2} \left( 1-x^{(1+\alpha)/2} \right) 
  & \mbox{If}~  \alpha \ne \pm 1 \\
  x \ln x  &  \mbox{If}~  \alpha = 1  \\
  - \ln x  &  \mbox{If}~  \alpha = -1 
 \end{cases}$ \\
\bottomrule
\end{tabular}  
\label{tab:App03-RCCP-01}
\end{table}

Many commonly seen divergence measures are special cases of $\phi$-divergence, coinciding with a particular choice of function $\phi$. Some examples are given in Table \ref{tab:App03-RCCP-01} \cite{App03-Sect4-Example-Phi-Div}. In what follows, we will use the KL-divergence. According to its corresponding function $\phi$, the KL-divergence is given by
\begin{equation}
\label{eq:App-03-RCCP-Phi-Div}
D_\phi(f \| f_0) = \int_{\rm \Omega}  
\log \left( \dfrac{f(\xi)}{f_0(\xi)} \right) f(\xi) d\xi
\end{equation} 

Before presenting the main results in \cite{App03-Sect4-RCCP}, the definition of conjugate duality is given. For a univariate function $g: \mathbb R \to \mathbb R \cup \{+\infty\}$, its conjugate function $g^*: \mathbb R \to \mathbb R \cup \{+\infty\}$ is defined as
\begin{equation*}
g^*(t) = \sup_{x \in \mathbb R}  \{tx - g(x)\} 
\end{equation*}

For a valid function $\phi$ for $\phi$-divergence satisfying (C1)-(C4), its conjugate function $\phi^*$ is convex, nondecreasing, and the following condition holds \cite{App03-Sect4-RCCP} 
\begin{equation}
\label{eq:App-03-DRSO-Conjugate-Ineq}
\phi^*(x) \ge x
\end{equation}
Besides, if $\phi^*$ is a finite constant on a closed interval $[a, b]$, then it is a finite constant on the interval $(-\infty,b]$.

\vspace{12pt}
{\noindent \bf 2. Equivalent formulation}

It is revealed in \cite{App03-Sect4-RCCP} that when the confidence set $D$ is constructed based on $\phi$-divergence, robust chance constrained program (\ref{eq:App-03-RCCP-Model}) can be easily transformed into a traditional chance-constrained program (\ref{eq:App-03-CCP-Model}) at the reference distribution by calibrating the confidence tolerance $\alpha$.

\begin{theorem}
\label{thm:App03-RCCP-Phi-Div}
{\rm \cite{App03-Sect4-RCCP}}  Let $\mathbb P_0$ be the cumulative distribution function generated by density function $f_0$, then the robust chance constraint
\begin{equation}
\label{eq:App-03-Thm-1}
\inf_{\mathbb P(\xi) \in \{D_{\phi}(f\|f_0) \le d\} }  
\Pr[C(x,\xi)] \ge 1-\alpha
\end{equation}
constructed based on $\phi$-divergence is equivalent to a traditional chance constraint
\begin{equation}
\label{eq:App-03-Thm-2}
\Pr\nolimits_0[C(x,\xi)] \ge 1-\alpha^\prime_+
\end{equation}
where $\Pr_0$ means that the probability is evaluated at the reference distribution $\mathbb P_0$, $\alpha^\prime_+ = \max\{\alpha^\prime,0\}$, and $\alpha^\prime$ can be computed by 
\begin{equation*}
\alpha^\prime = 1 - \inf_{z \in Z} \left\{
\dfrac{\phi^*(z_0+z)-z_0-\alpha z +d}{\phi^*(z_0+z)-\phi^*(z_0)} \right\}
\end{equation*}
where 
\begin{equation*}
Z = \left\{ z \middle| \begin{lgathered}
z>0,~ z_0 + \pi z \le l_\phi \\
\underline m(\phi^*) \le z+z_0 \le \overline m(\phi^*)
\end{lgathered} \right\}
\end{equation*}
In above formula, constants $l_\phi = \lim_{x \to +\infty} \phi(x)/x$, $\overline m(\phi^*) = \inf \{ m: \phi^*(m)=+\infty\}$,  $\underline m(\phi^*) = \sup \{ m: \phi^*$ is a finite constant on $(-\infty,m]\}$, Table \ref{tab:App03-RCCP-02} summarizes the values of these parameters for typical $\phi$-divergence measures, and  
\begin{equation*}
\pi = \begin{cases}
-\infty &  \mbox{if Leb}\{[f_0=0]\}=0  \\
0       &  \mbox{if Leb $\{[f_0=0]\}>0$ and Leb$\{[f_0=0] \backslash C(x,\xi) \}=0$} \\
1       &  \mbox{otherwise}
\end{cases}
\end{equation*}
where Leb$\{\cdot\}$ is the Lebesgue measure on $\mathbb R^{Dim(\xi)}$.
\end{theorem}

\begin{table}[!htp]
\footnotesize
\renewcommand{\arraystretch}{1.3}
\renewcommand{\tabcolsep}{1em}
\caption{Values of $l_\phi$, $\underline m(\phi^*)$, and $\overline m(\phi^*)$ for $\phi$-divergences}
\centering
\begin{tabular}{cccc}
\toprule 
$\phi$-Divergence   & $l_\phi$  & $\underline m(\phi^*)$ & $\overline m(\phi^*)$\\
\midrule
   KL-divergence        & $+\infty$ & $-\infty$ & $+\infty$ \\
Hellinger distance      &    $1$    & $-\infty$ &    $1$    \\
Variation distance      &    $1$    &   $-1$    &    $1$    \\
    J-divergence        & $+\infty$ & $-\infty$ & $+\infty$ \\
 $\chi^2$ divergence    & $+\infty$ &   $-2$    & $+\infty$ \\

\bottomrule
\end{tabular}  
\label{tab:App03-RCCP-02}
\end{table}

The values of $\alpha^\prime$ for the Variation distance and the $\chi^2$ divergence have analytical expressions; for the KL divergence, $\alpha^\prime$ can be computed from one-dimensional line search. Results are shown in Table \ref{tab:App03-RCCP-03}.

\begin{table}[!htp]
\footnotesize
\renewcommand{\arraystretch}{1.3}
\renewcommand{\tabcolsep}{1em}
\caption{Values of $\alpha^\prime$ for some $\phi$-divergences}
\centering
\begin{tabular}{cc}
\toprule 
   $\phi$-Divergence    &  $\alpha^\prime$    \\
\midrule
 $\chi^2$ divergence    &  $\alpha^\prime = \alpha- \dfrac{\sqrt{d^2 + 4d (\alpha-\alpha^2)}-(1-2\alpha)d}{2d+2}$      \\
                        &                     \\
 Variation distance     &  $\alpha^\prime = \alpha - \dfrac{1}{2}d$      \\
                        &                     \\
   KL-divergence        &  $\alpha^\prime = 1- \inf_{x \in (0,1)} \left\{ \dfrac{{\rm e}^{-d}x^{1-\alpha}-1}{x-1} \right\}$     \\
\bottomrule
\end{tabular}  
\label{tab:App03-RCCP-03}
\end{table}

For the KL divergence, calculating $\alpha^\prime$ entails
solving $\inf_{x\in (0,1)} h(x)$ where
\begin{equation*}
h(x) = \dfrac{{\rm e}^{-d}x^{1-\alpha}-1}{x-1} 
\end{equation*}
Its first-order derivative is given by 
\begin{equation*}
h^\prime (x) = \dfrac{1-\alpha{\rm e}^{-d}x^{1-\alpha}-(1-\alpha){\rm e}^{-d} x^{-\alpha}}{(x-1)^2},~~ \forall x \in (0,1) 
\end{equation*}

To claim the convexity of $h(x)$, we need to show that $h^\prime (x)$ is an increasing function in $x \in (0,1)$. To this end, first notice that the denominator $(x-1)^2$ is a decreasing function in $x$ on the open interval $(0,1)$; then we can show the numerator is an increasing function in $x$, because its first-order derivative gives
\begin{equation*}
(1-\alpha{\rm e}^{-d}x^{1-\alpha}-(1-\alpha){\rm e}^{-d} x^{-\alpha})^\prime _x = \alpha (1-\alpha){\rm e}^{-d} (x^{-\alpha-1}-x^{-\alpha}) >0,~ \forall x \in (0,1)
\end{equation*}
Hence $h^\prime(x)$ is monotonically increasing, and $h(x)$ is a convex function in $x$. Moreover, because $h^\prime(x)$ is continuous in $(0,1)$, and $\lim_{x \to 0^+} h^\prime(x) = -\infty$, $\lim_{x \to 1^-} h^\prime(x) = +\infty$, there must be some $x^* \in [\delta, 1-\delta]$ such that $h^\prime(x^*)=0$, i.e., the infimum of $h(x)$ is attainable. The minimum of $h(x)$ can be calculated by solving a nonlinear equation $h^\prime(x)=0$ via Newton's method, or a derivative-free line search, such as the golden section search algorithm. Either scheme is computationally inexpensive.  

Finally, we discuss the connection between the modified tolerance $\alpha^\prime$ and its original value $\alpha$. Because a set of distributions are considered in (\ref{eq:App-03-Thm-1}), the threshold  in (\ref{eq:App-03-Thm-2}) should be greater than the original one, i.e., $1-\alpha^\prime \ge 1-\alpha$ must hold. To see this, recall inequality (\ref{eq:App-03-DRSO-Conjugate-Ineq}) of conjugate function, we have
\begin{equation*}
\alpha \phi^* (z_0 + z) + (1 - \alpha) \phi^∗ (z_0) \ge 
\alpha (z_0 + z) + (1 - \alpha) z_0
\end{equation*}
The right-hand side gives $\alpha z + z_0$; in the ambiguity set (\ref{eq:App-03-RCCP-Phi-Div}), $d$ is strictly positive, therefore
\begin{equation*}
\alpha \phi^* (z_0 + z) + (1 - \alpha) \phi^∗ (z_0) \ge 
\alpha z + z_0 -d
\end{equation*}
which gives
\begin{equation*}
\phi^* (z_0 + z) - z_0 - az + d  \ge
(1 - \alpha) (\phi^∗ (z_0 + z) - \phi^∗ (z_0) )
\end{equation*}
Recall the expression of $\alpha^\prime$ in Theorem \ref{thm:App03-RCCP-Phi-Div}, we arrive at
\begin{equation*}
1-\alpha^\prime = \dfrac{\phi^* (z_0+z)-z_0-az+d}{\phi^∗ (z_0+z)-\phi^∗(z_0)}
\ge 1 - \alpha 
\end{equation*}
which is the desired conclusion.

Theorem \ref{thm:App03-RCCP-Phi-Div} concludes that the complexity of handling
a robust chance constraint is almost the same as that of tackling a traditional chance constraint associated with the reference distribution $\mathbb P_0$, except for the efforts on computing $\alpha^\prime$. If $\mathbb P_0$ belongs to the family of log-concave distributions, then the chance constraint is convex. As a special case, if $\mathbb P_0$ is the Gaussian distribution or a uniform distribution on ellipsoidal support, a single chance constraint can boil down to   a second-order cone \cite{App03-Sect4-Example-Q-Distribution}. For more general cases, the chance constraint is non-convex in $x$. In such circumstance, we will use  risk based reformulation and the sampling average approximation (SAA) approach. 

\vspace{12pt}
{\noindent \bf 3. Risk and SAA based reformulation}

Owing to the different descriptions on dispersion ambiguity and presence of the wait-and-see decision $y$, unlike DRO problem (\ref{eq:App-03-DRO-Model-1}) with static robust chance constraint (\ref{eq:App-03-DRO-DRCC}) which can be transformed into an SDP, constraint (\ref{eq:App-03-Thm-1}) is treated in a different way, as demonstrated in Theorem \ref{thm:App03-RCCP-Phi-Div}: it comes down to a traditional chance constraint (\ref{eq:App-03-Thm-2}) while the dispersion ambiguity is taken into account by a modification in the confidence level. The remaining task is to express (\ref{eq:App-03-Thm-2}) as a solver-compatible form. 

\vspace{12pt}
{\noindent \bf 1) Loss function}

For given $x$ and $\xi$, constraints in $C(x,\xi)$ cannot be met if no $y$ satisfying $A(\xi) x + B(\xi) y \le b(\xi)$ exists. To quantify the constraint violation under scenario $\boldsymbol{\xi}$ and first-stage decision $x$, define the following loss function $L(x,\xi)$ 
\begin{equation}
\label{eq:App-03-RCCP-Loss-Function}
\begin{aligned}
L(x,\xi) = \min_{y,\sigma} ~~ & \sigma \\
\text{s.t.}~~ & A(\xi) x + B(\xi) y \le b(\xi) + \sigma {\bf 1}
\end{aligned}
\end{equation}
where {\bf 1} is an all-one vector with compatible dimension. If $L(x,\xi) \ge 0$, the minimum of slackness $\sigma$ under the joint efforts of the recourse action $y$ is defined as the loss; otherwise, demands are satisfiable after the uncertain parameter is known. As we assume $C(x,\xi)$ is a bounded polytope, problem (\ref{eq:App-03-RCCP-Loss-Function}) is always feasible and bounded below. Therefore, the loss function $L(\boldsymbol{x}, \boldsymbol{\xi})$ is well-defined, and the chance constraint (\ref{eq:App-03-Thm-2}) can be written as  
\begin{equation}
\label{eq:App-03-RCC-Loss-Fun}
\Pr\nolimits_0 [ L(x,\xi) \leq 0 ]  \geq 1 - \alpha^\prime_+
\end{equation}

In this way, the joint chance constraints are consolidated into a single one, just like what has been done in (\ref{eq:App-03-DRO-DRCC-Joint}) and (\ref{eq:App-03-DRO-DRCC-Joint-Para}).

\vspace{12pt}
{\noindent \bf 2) VaR based reformulation: An MILP}

For a given probability tolerance $\beta$ and a first-stage decision $x$, the $\beta$-VaR for loss function $L(x,\xi)$ under the reference distribution PDF $\mathbb P_0$ is defined as
\begin{equation}
\label{eq:App-03-RCC-VaR-Def}
\beta \mbox{-VaR}(x) = \min \left\{ a \in \mathbb {R} \middle| \int_{L(x,\xi) \leq a} f_0(\xi) \mbox{d} \xi \ge \beta \right\}
\end{equation}
which interprets the threshold $a$ such that the loss is no greater than $a$ will hold with a probability no less than $\beta$. According to (\ref{eq:App-03-RCC-VaR-Def}), an equivalent expression of chance constraint (\ref{eq:App-03-RCC-Loss-Fun}) is
\begin{equation}
\label{eq:App-03-RCC-CC-VaR}
(1 - \alpha^\prime_+) \mbox{-VaR} (x) \le 0
\end{equation}
So that probability evaluation is obviated. Furthermore, if SAA is used, (\ref{eq:App-03-RCC-Loss-Fun}) and (\ref{eq:App-03-RCC-CC-VaR}) indicate that the scenarios which will lead to $L(x,\xi) > 0$ account for a fraction of $\alpha_{1+}$ among all sampled data.

Let $\xi_1,\xi_2, \cdots, \xi_q$ be $q$ scenarios sampled from random variable $\boldsymbol{\xi}$. We use $q$ binary variables $z_1, z_2, \cdots, z_q$ to identify possible infeasibility: $z_k = 1$ implies that constraints cannot be satisfied in scenario $\xi_k$. To this end, let $M$ be a large enough constant, consider inequality
\begin{equation}
\label{eq:App-03-RCC-Loss-SSA} 
A(\xi_k) x + B(\xi_k) y_k \le b(\xi_k) + M z_k
\end{equation}
In (\ref{eq:App-03-RCC-Loss-SSA}), if $z_k = 0$, recourse action $y_k$ will recover all constraints in scenario $\xi_k$, and thus $C(x,\xi_k)$ is non-empty; otherwise, if no such a recourse action $y_k$ exists, then constraint violation will take place. To reconcile infeasibility, $z_k = 1$ so that (\ref{eq:App-03-RCC-Loss-SSA}) becomes redundant, and there is actually no constraint for scenario $\boldsymbol{\xi_k}$. The fraction of sampled scenarios which will incur inevitable constraint violations is counted by $\sum_{k=1}^q z_k/q$. So we can write out the following MILP reformulation for robust chance-constrained program (\ref{eq:App-03-RCCP-Model}) based on VaR and SAA
\begin{equation}
\label{eq:App-03-RCCP-VaR-MILP}
\begin{aligned}
\min ~~ & c^T x  \\
\mbox{s.t.}~~ & x \in X  \\
&  A(\xi_k) x + B(\xi_k) y_k \le b(\xi_k) + M z_k,~ k=1,\cdots,q \\
& \sum_{k = 1}^{q} z_k \leq q\alpha^\prime_+,~ z_k \in \{ 0, 1 \},~ k=1,\cdots,q
\end{aligned}
\end{equation}

In MILP (\ref{eq:App-03-RCCP-VaR-MILP}), constraint violation can happen in at most $q \alpha_{1+}$ out of $q$ scenarios in the reference distribution, according to Theorem \ref{thm:App03-RCCP-Phi-Div}, and the reliability requirement (\ref{eq:App-03-Thm-1}) under all possible distributions in ambiguity set $D_\phi$ can be guaranteed by the selection of $\alpha^\prime_+$. Improved MILP formulations of chance constraints which do not rely on the specific big-M parameter are comprehensively studied in \cite{App03-Sect4-CC-SAA-MILP}, and some structure properties of the feasible region are revealed. 

\vspace{12pt}
{\noindent \bf 3) CVaR based reformulation: An LP}

The number of binary variables in MILP (\ref{eq:App-03-RCCP-VaR-MILP}) is equal to the number of sampled scenarios. To guarantee the accuracy of SAA, a large number of scenarios are required, preventing MILP (\ref{eq:App-03-RCCP-VaR-MILP}) from being solved efficiently. To ameliorate this plight, we provide a conservative LP approximation for problem (\ref{eq:App-03-RCCP-Model}) based on the properties of CVaR revealed in \cite{CVaR}.

The $\beta$-CVaR for the loss function $L(x,\xi)$ is defined as  
\begin{equation}
\label{eq:App-03-VaR-Def}
\beta \mbox{-CVaR} (x) = \frac{1} {1 - \beta} \int_{L(x,\xi) \ge \beta \text{-VaR} (x)} L(x,\xi) f(\boldsymbol{\xi}) d \xi
\end{equation}
which interprets  the conditional expectation of loss that is no less than $\beta$-VaR; therefore, relation 
\begin{equation}
\label{eq:App-03-VaR-CVaR}
\beta \mbox{-VaR} \leq \beta \mbox{-CVaR}
\end{equation} 
always holds, and a conservative approximation of constraint (\ref{eq:App-03-RCC-CC-VaR}) is
\begin{equation}
(1 - \alpha^\prime_+) \mbox{-CVaR} (x) \le 0
\label{eq:App-03-RCC-CC-CVaR}
\end{equation}
Inequality  (\ref{eq:App-03-RCC-CC-CVaR}) is a sufficient condition for (\ref{eq:App-03-RCC-CC-VaR}) and (\ref{eq:App-03-RCC-Loss-Fun}). This conservative replacement is apposite to the spirit of robust optimization. In what follows, we will reformulate (\ref{eq:App-03-RCC-CC-CVaR}) in a solver-compatible form.

According to \cite{CVaR}, the left-hand side of (\ref{eq:App-03-RCC-CC-CVaR}) is equal to the optimum of the following minimization problem 
\begin{equation}
\min_{\gamma} \left\{ \gamma + \frac{1}{\alpha^\prime_+} \int_{\xi \in \mathbb {R}^K} \max \{ L(x,\xi) - \gamma, 0 \} f(\xi) \mbox{d} \xi  \right\}
\label{eq:App-03-CVaR-Opt-Int}
\end{equation}

By performing SAA, the integral in (\ref{eq:App-03-CVaR-Opt-Int}) renders a summation over discrete sampled scenarios $\xi_1,\xi_2, \cdots, \xi_q$, resulting in
\begin{equation}
\label{eq:App-03-CVaR-Opt-Discrete}
\min_{\gamma} \left\{ \gamma + \frac{1}{q \alpha^\prime_+} \sum_{k=1}^q \max \left\{ L(x,\xi_k) - \gamma, 0 \right\} \right\}
\end{equation}

By introducing auxiliary variable $s_k$, the feasible region defined by (\ref{eq:App-03-RCC-CC-CVaR}) can be expressed via 
\begin{gather*}
\exists \gamma \in \mathbb R, s_k \in \mathbb R^+,~\sigma_k \in \mathbb R,~ k = 1,\cdots,q \\
\sigma_k - \gamma  \le s_k,~ k = 1,\cdots,q  \\  
A(\xi_k) x + B(\xi_k) y_k \le b(\xi_k) + \sigma_k {\bf 1},
~ k = 1,\cdots,q  \\
\gamma + \frac{1}{q \alpha^\prime_+} \sum_{k=1}^q s_k \leq 0 
\end{gather*}

Now we can write out the the conservative LP reformulation for robust chance constrained program (\ref{eq:App-03-RCCP-Model}) based on CVaR and SAA
\begin{equation}
\label{eq:App-03-RCCP-CVaR-LP}
\begin{aligned}
\min_{x,y,s,\gamma} ~~ & c^T x  \\
\mbox{s.t.}~~ & x \in X,~ \gamma + \frac{1}{q \alpha^\prime_+} \sum_{k=1}^q s_k \leq 0,~s_k \geq 0,~ k = 1,\cdots,q  \\
&  A(\xi_k) x + B(\xi_k) y_k - b(\xi_k) \le (\gamma + s_k) {\bf 1},
~ k = 1,\cdots,q 
\end{aligned}
\end{equation}
where $\sigma_k$ is eliminated.

According to (\ref{eq:App-03-VaR-CVaR}), condition (\ref{eq:App-03-RCC-CC-CVaR}) guarantees (\ref{eq:App-03-RCC-CC-VaR}) as well as (\ref{eq:App-03-RCC-Loss-Fun}), so chance constraint in (\ref{eq:App-03-Thm-1}) holds with a probability no less (usually higher) than $1-\alpha$, regardless of the true distributions in confidence set $D_\phi$. Since (\ref{eq:App-03-VaR-CVaR}) is usually a strict inequality, this fact will introduce some extent of conservatism in the CVaR based LP model (\ref{eq:App-03-RCCP-CVaR-LP}).

\begin{figure}[!htp]
  \centering
  \includegraphics[scale = 0.40]{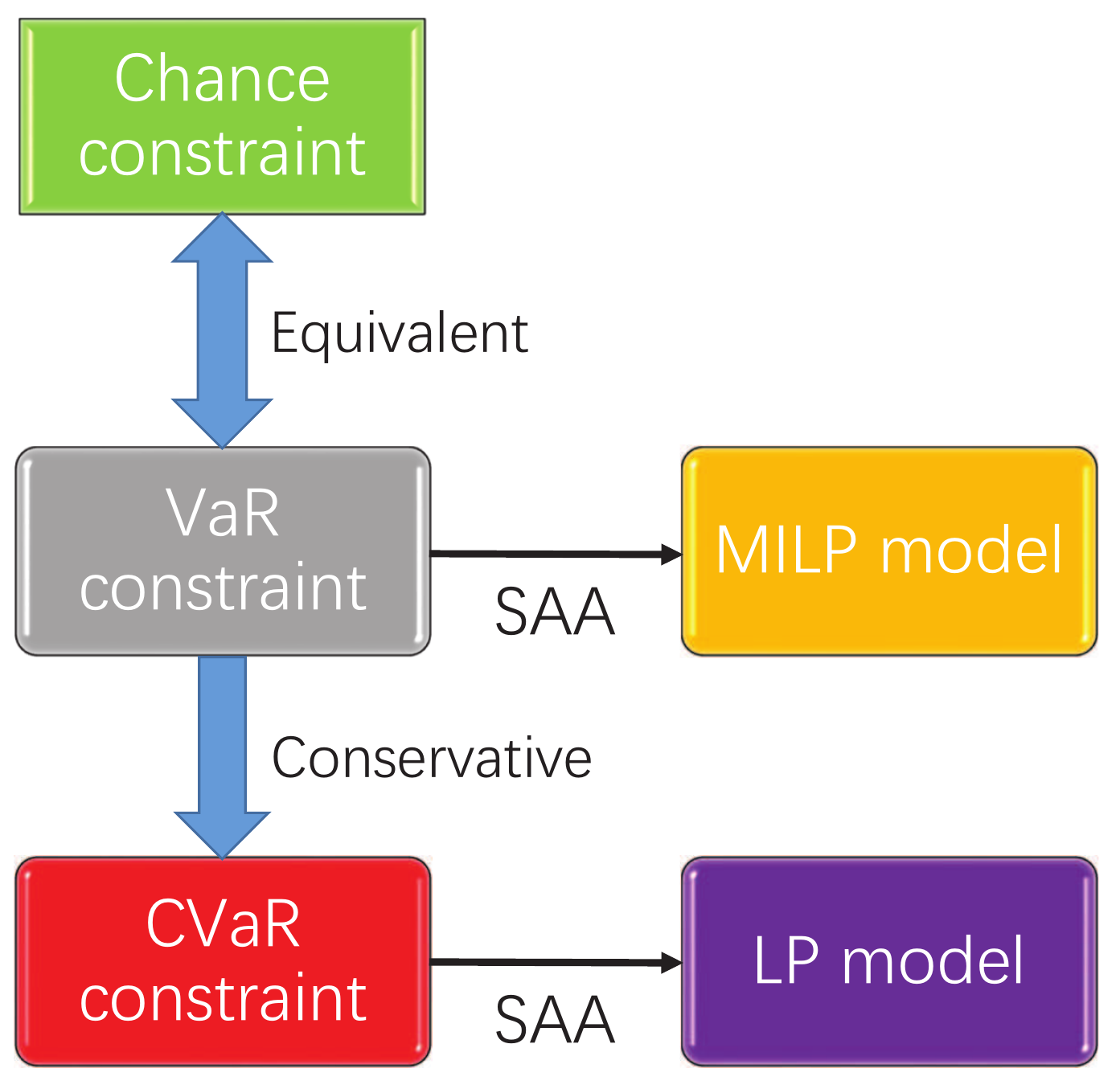}
  \caption{Relations of the models discussed in this section.}
  \label{fig:Fig-App03-01}
\end{figure}

Relations among different mathematical models discussed in this section are summarized in Fig. \ref{fig:Fig-App03-01}.

\vspace{12pt}
{\noindent \bf 4. Considering second-stage cost}

Finally, we elaborate how to solve problem (\ref{eq:App-03-RCCP-Model}) with a second-stage cost in the sense of worst-case expectation, i.e. 
\begin{equation}
\label{eq:App-03-RCCP-MEdy-1}
\begin{aligned}
\min_x~ & \left\{c^T x + \max_{P(\xi) \in D_{KL}} \mathbb E_P [Q(x,\xi)] \right\}  \\
\mbox{s.t.}~~ & x \in X  \\
& \sup_{P(\xi) \in D^\prime}  \Pr[C(x,\xi)] 
\ge 1-\alpha
\end{aligned}
\end{equation}
where $Q(x,\xi)$ is the optimal value function of the second-stage problem
\begin{equation*}
\begin{aligned}
Q(x,\xi) = \min ~~ & q^T y   \\
\mbox{s.t.}~~ & B(\xi) y \le b(\xi) - A(\xi) x   
\end{aligned}
\end{equation*}
which is an LP for a fixed first-stage decision $x$ and a given parameter $\xi$;
\begin{equation*}
D_{KL}=\{P(\xi)~|~ D^{KL}_\phi(f \| f_0) \le d_{KL}(\alpha^*),~ f = {\rm d} P/ {\rm d} \xi\}
\end{equation*} 
is the KL-divergence based ambiguity set, and $d_{KL}$ is an $\alpha$-dependent threshold which determines the size of the ambiguity set, and $\alpha^*$ reflects the confidence level: the real distribution is contained in $D_{KL}$ with a probability no less than $\alpha^*$. For discrete distributions, the KL-divergence measure has the form of 
\begin{equation*}
D^{KL}_\phi(f \parallel f_0) = \sum_s \rho_s \log \dfrac{\rho_s}{\rho^0_s}
\end{equation*}
In either case, there are infinitely many PDFs satisfying the inequality in the ambiguity set $D_{KL}$ when $d_{KL} > 0$. Otherwise, when $d_{KL}=0$, the ambiguity set $D_{KL}$ becomes a singleton, and the model (\ref{eq:App-03-RCCP-MEdy-1}) degenerates to a traditional SO problem. In practice, the user can specify the value of $d_{KL}$ according to the attitude towards risks. Nevertheless, the proper value of $d_{KL}$ can be obtained from probability theory. Intuitively, the more historical data we possess, the closer the reference PDF $f_0$ leaves from the true one, and the smaller $d_{KL}$ should be set. 

 Suppose we have totally $M$ samples with equal probabilities to fit in $N$ bins, and there are $M_1$, $M_2$, $\cdots$, $M_N$ samples fall into each bin, then  the discrete reference PDF for the histogram is $\{\pi_1, \cdots,\pi_N\}$, where $\pi_i = M_i/M$, $i = 1, \cdots, N$. Let $\pi^r_1$, $\cdots$, $\pi^r_N$ be the real probability of each bin, according to the discussions in \cite{Am-Set-Phi-Divergence},  random variable $2M \sum_{i=1}^N \pi^r_i \log (\pi^r_i/\pi_i)$ follows $\chi^2$ distribution with $N-1$ degrees of freedom. Therefore, the confidence threshold can be calculated from
\begin{equation*}
d_{KL}(\alpha^*) = \dfrac{1}{2M}\chi_{N-1,\alpha^*}^2
\end{equation*} 
where $\chi_{N-1,\alpha^*}^2$ stands for the $\alpha^*$ upper quantile of $\chi^2$ distribution with $N-1$ degrees of freedom. For other divergence based ambiguity sets, please see more discussions in \cite{Am-Set-Phi-Divergence}. Robust chance constraints in (\ref{eq:App-03-RCCP-MEdy-1}) are tackled using the method presented previously, and the objective function will be treated independently. The ambiguity sets in the objective function and chance constraints could be the same one or different ones, and thus are distinguished by $D_{KL}$ and $D^\prime$.   

Sometimes, it is imperative to coordinately optimize the costs in both stages. For example, in the facility planning problem, the first stage represents the investment decision and the second stage describes the operation management. If we only optimize the first-stage cost, then the facilities with lower investment costs will be preferred, but they may suffer from higher operating costs, and not be the optimal choice from the long-term aspect.   

To solve (\ref{eq:App-03-RCCP-MEdy-1}), we need a tractable reformulation for the worst-case expectation problem under KL-divergence ambiguity set  
\begin{equation}
\label{eq:App-03-MEdy-KL-Div}
\max_{P(\xi) \in D_{KL}} \mathbb E_P [Q(x,\xi)] 
\end{equation}
under fixed $x$. It is proved in \cite{App03-Sect4-RCCP-mEdy,Am-Set-Phi-Divergence} that problem (\ref{eq:App-03-MEdy-KL-Div}) is equivalent to 
\begin{equation}
\label{eq:App-03-MEdy-KL-Div-Dual}
\min_{\alpha \ge 0} ~ \alpha  \log \mathbb E_{P_0} [{\rm e}^{Q(x,\xi)/\alpha}] + \alpha d_{KL}
\end{equation}
where $\alpha$ is the dual variable. Formulation (\ref{eq:App-03-MEdy-KL-Div-Dual}) has two advantages: first, the expectation is evaluated associated with the reference distribution $P_0$, which is much easier than optimizing over the ambiguity set $D_{KL}$; second, the maximum operator switches to a minimum  operator, which is consistent with the objective function of the decision making problem. We will use SAA to express the expectation, giving rise to a discrete version of problem (\ref{eq:App-03-MEdy-KL-Div-Dual}). In fact, in discrete cases, (\ref{eq:App-03-MEdy-KL-Div-Dual}) can be derived from (\ref{eq:App-03-MEdy-KL-Div}) using Lagrange duality. The following interpretation is given in \cite{App03-Sect4-KL-Div-UC}.

Denote by $\xi_1,\cdots,\xi_s$  the representative scenarios in the discrete distribution; their corresponding probabilities in the reference PDF and the actual PDF are given by $P_0 = \{p^0_1,\cdots,p^0_s\}$ and $P = \{p_1,\cdots,p_s\}$, respectively. Then problem (\ref{eq:App-03-MEdy-KL-Div}) can be written in a discrete form as
\begin{equation}
\label{eq:App-03-MEdy-KL-Div-Discrete}
\begin{aligned}
\max_p~~ & \sum_{i=1}^s   p_i Q(x,\xi_i)  \\
\mbox{s.t.} ~~ &  \sum_{i=1}^s p_i  \log \left( \dfrac{p_i}{p^0_i} \right) \le d_{KL}  \\
& p \ge 0,~ 1^T p =1
\end{aligned}
\end{equation}
where vector $p=[p_1,\cdots,p_s]^T$ is the decision variable. According to Lagrange duality theory, the objective function of the dual problem is
\begin{equation}
\label{eq:App-03-MEdy-KL-Div-Discrete-Dual-Obj-1}
g(\alpha,\mu) = \alpha d_{KL} + \mu + \sum_{i=1}^s \max_{p_i \ge 0}  p_i \left( Q(x,\xi_i) - \mu - \alpha \log \left( \dfrac{p_i}{p^0_i} \right)  \right)
\end{equation}
where $\mu$ is the dual variable associated with equality constraint $1^T p =1$, and $\alpha$ with the KL-divergence inequality. Substituting $t_i = p_i / p^0_i$ into (\ref{eq:App-03-MEdy-KL-Div-Discrete-Dual-Obj-1}) and eliminating $p_i$, we get
\begin{equation*}
g(\alpha,\mu) = \alpha d_{KL} + \mu + \sum_{i=1}^s \max_{t_i \ge 0}~ p^0_i t_i \left( Q(x,\xi_i) - \mu - \alpha \log t_i  \right)
\end{equation*}

Calulating the first-order derivative of $t_i ( Q(x,\xi_i) - \mu - \alpha \log t_i )$ with respect to $t_i$, the optimal solution is 
\begin{equation*}
t_i = {\rm e} ^{\frac{ Q(x,\xi_i) - \mu - \alpha}{\alpha}} > 0
\end{equation*}
and the maximum is
\begin{equation*}
\alpha {\rm e} ^{\frac{ Q(x,\xi_i) - \mu - \alpha}{\alpha}} 
\end{equation*}
As a result, the dual objective reduces to 
\begin{equation}
\label{eq:App-03-MEdy-KL-Div-Discrete-Dual-Obj-2}
g(\alpha,\mu) = \alpha d_{KL} + \mu + \alpha \sum_{i=1}^s 
p^0_i {\rm e} ^{\frac{ Q(x,\xi_i) - \mu - \alpha}{\alpha}}
\end{equation}
and the dual problem of (\ref{eq:App-03-MEdy-KL-Div-Discrete}) can be rewritten as
\begin{equation}
\label{eq:App-03-MEdy-KL-Div-Discrete-Dual-1}
\min_{\alpha \ge 0,\mu}~ g(\alpha,\mu)
\end{equation}

The optimal solution $\mu^*$ must satisfy $\partial g / \partial \mu = 0$, yielding
\begin{equation*}
\sum_{i=1}^s p^0_i {\rm e} ^{\frac{ Q(x,\xi_i) - \mu^* - \alpha}{\alpha}}=1
\end{equation*}
or
\begin{equation*}
\mu^* = \alpha \log \sum_{i=1}^s p^0_i ~{\rm e}^{Q(x,\xi_i)/\alpha} - \alpha
\end{equation*}
Substituting above relations into $g(\alpha,\mu)$ results in the following dual problem
\begin{equation}
\label{eq:App-03-MEdy-KL-Div-Discrete-Dual-2}
\min_{\alpha \ge 0}~ \left\{ \alpha d_{KL} + \alpha \log \sum_{i=1}^s p^0_i ~{\rm e}^{Q(x,\xi_i)/\alpha} \right\}  
\end{equation}
which is a discrete form of (\ref{eq:App-03-MEdy-KL-Div-Dual}). 

In (\ref{eq:App-03-RCCP-MEdy-1}), replacing the inner problem (\ref{eq:App-03-MEdy-KL-Div}) with its Lagrangian dual form (\ref{eq:App-03-MEdy-KL-Div-Discrete-Dual-2}), we can obtain an equivalent mathematical program
\begin{equation}
\label{eq:App-03-RCCP-MEdy-2}
\begin{aligned}
\min~ & \left\{c^T x + \alpha d_{KL} + \alpha \log \sum_{i=1}^s p^0_i ~{\rm e}^{\theta_i/\alpha} \right\}  \\
\mbox{s.t.}~~ & x \in X,~ \alpha \ge 0,~ \theta_i = q^T y_i,~ \forall i  \\
& A(\xi_i) x + B(\xi_i) y_i \le b(\xi_i),~ \forall i \\
& \mbox{Cons-RCC}
\end{aligned}
\end{equation}
where Cons-RCC stands for the LP based formulation of robust chance constraints, so the constraints in problem (\ref{eq:App-03-RCCP-MEdy-2}) are all linear, and the only nonlinearity rests in the last term of the objective function. In what follows, we will show it is actually a convex function in $\theta_i$ and $\alpha$.

In the first step, we claim that the following function is convex (\cite{CVX-Book-Boyd}, page 87, in Example 3.14)
\begin{equation*}
h_1(\theta) = \log \left( \sum_{i=1}^s {\rm e}^{\theta_i} \right)
\end{equation*}
Since the composition with an affine  mapping preserves convexity (\cite{CVX-Book-Boyd}, Sect. 3.2.2), a new function
\begin{equation*}
h_2(\theta) = h_1 (A \theta + b)
\end{equation*}
remains convex under linear mapping $\theta \to A \theta + b$. Let $A$ be an identity matrix, and  
\begin{equation*}
b = \begin{bmatrix}
\log p^0_1 \\ \vdots  \\  \log p^0_s
\end{bmatrix}
\end{equation*}
then we have 
\begin{equation*}
h_2(\theta) = \log \left( \sum_{i=1}^s p^0_i {\rm e}^{\theta_i} \right)
\end{equation*}
is a convex function; at last, function
\begin{equation*}
h_3(\alpha,\theta) = \alpha h_2(\theta/\alpha) 
\end{equation*}
is the perspective of $h_2(\theta)$, so is also convex (\cite{CVX-Book-Boyd}, page 89, Sect. 3.2.6).

In view of this convex structure, (\ref{eq:App-03-RCCP-MEdy-2}) essentially gives rise to a convex program, and the local minimum is also the global one. However, according to our experiments, general purpose NLP solvers still have difficulty to solve (\ref{eq:App-03-RCCP-MEdy-2}). Therefore, we employ the outer approximation method \cite{App03-Sect4-OA-1,App03-Sect4-OA-2}. The motivation is to solve the epigraph form of (\ref{eq:App-03-RCCP-MEdy-2}), in which nonlinearity is moved into the constraints; then linearize the feasible region with an increasing number of cutting planes generated in an iteration algorithm, until certain convergence criterion is met. In this way, the hard problem (\ref{eq:App-03-RCCP-MEdy-2}) can be solved via a sequence of LPs. The outer approximation algorithm is outlined in Algorithm \ref{Ag:App-03-DRO-Outer-Approximation}. Because (\ref{eq:App-03-RCCP-MEdy-2}) is a convex program, the cutting planes will not remove any feasible point, and Algorithm \ref{Ag:App-03-DRO-Outer-Approximation} finds the global optimal solution in finite steps, regardless of the initial point. But for sure, the number of iterations is  affected by the quality of initial guess. A proper initiation could be obtained by solving a traditional SO problem without considering distribution uncertainty.  

\begin{algorithm}[!htp]
\normalsize
\caption{\bf Outer Approximation}
\begin{algorithmic}[1]

\STATE Choose an initial point $(\theta^1,\alpha^1)$ and convergence tolerance $\epsilon > 0$, the initial objective value is $R^1=0$, and  iteration index $k=1$. 

\STATE Solve the following master problem which is an LP 
\begin{equation}
\label{eq:App-03-RCCP-MEdy-OA-Lin}
\begin{aligned}
\min_{\alpha,\theta,\gamma,x}~~ & c^T x + \alpha d_{KL} + \gamma   \\
\mbox{s.t.}~~ & h_3(\alpha^j,\theta^j) + \nabla h_3(\alpha^j,\theta^j) 
\begin{bmatrix}
\alpha - \alpha^j \\
\theta - \theta^j  
\end{bmatrix} 
\le \gamma,~ j=1,\cdots,k \\ 
& x \in X,~ \alpha \ge 0,~ \theta_i = q^T y_i,~ \forall i  \\
& A(\xi_i) x + B(\xi_i) y_i \le b(\xi_i),~ \forall i \\
& \mbox{Cons-RCC}
\end{aligned}
\end{equation}
 The optimal value is $R^{k+1}$, and the optimal solution is $(x^{k+1},\theta^{k+1},\alpha^{k+1})$.

\STATE If $R^{k+1} - R^k \le \varepsilon$, terminate and report the optimal solution $(x^{k+1},\theta^{k+1},\alpha^{k+1})$; otherwise, update $k \leftarrow k+1$, calculate the gradient $\nabla h_3$ at the obtained solution $(\alpha^k,\theta^k)$, add the following cut to problem (\ref{eq:App-03-RCCP-MEdy-OA-Lin}), and go to step 2.
\begin{equation}
\label{eq:App-03-RCCP-MEdy-OA-Cut}
h_3(\alpha^k,\theta^k) + \nabla h_3(\alpha^k,\theta^k) 
\begin{bmatrix}
\alpha - \alpha^k \\
\theta - \theta^k  
\end{bmatrix} 
\le \gamma 
\end{equation}
\end{algorithmic}
\label{Ag:App-03-DRO-Outer-Approximation}
\end{algorithm}  

\begin{figure}[!htp]
  \centering
  \includegraphics[scale = 0.60]{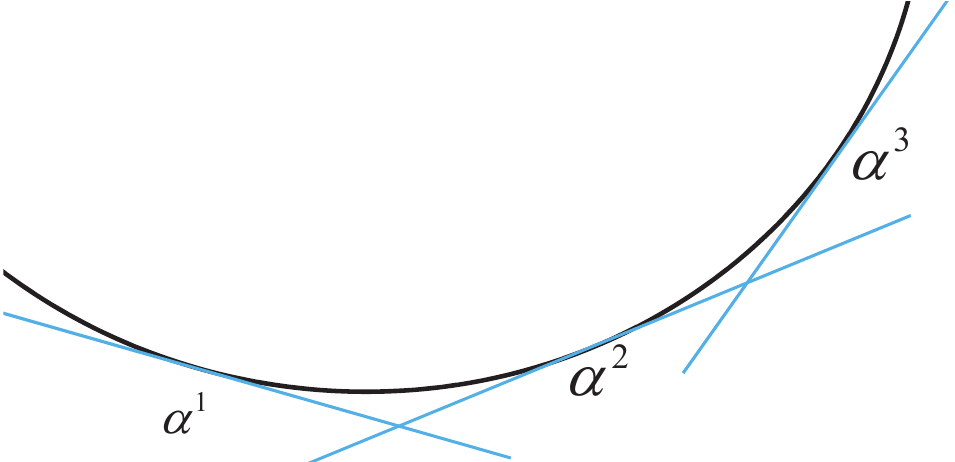}
  \caption{Illustration of the outer approximation algorithm.}
  \label{fig:Fig-App03-02}
\end{figure}

The motivation of Algorithm \ref{Ag:App-03-DRO-Outer-Approximation} is illustrated in \ref{fig:Fig-App03-02}. The original objective function is nonlinear but convex. In the epigraph form (\ref{eq:App-03-RCCP-MEdy-OA-Lin}), we generate a set of linear cuts (\ref{eq:App-03-RCCP-MEdy-OA-Cut}) dynamically according to the optimal solution found in step 2, then the convex region can be approximated with arbitrarily high accuracy around the optimal solution. The convergence of the very basic version of outer approximation method has been analyzed in \cite{App03-Sect4-OA-3,App03-Sect4-OA-4}. In fact, Algorithm \ref{Ag:App-03-DRO-Outer-Approximation} is very efficient to solve problem (\ref{eq:App-03-RCCP-MEdy-2}), because problem (\ref{eq:App-03-RCCP-MEdy-OA-Lin}) is an LP, the objective function is smooth, and the algorithm often converges in a few number of iterations.

\subsection{Stochastic Program with Discrete Distributions}
\label{App-C-Sect04-02}

In ARO discussed in Appendix \ref{App-C-Sect02}, the uncertain parameter is assumed to reside in the so-called uncertainty set. Every element in this set is treated equally, so the scenario in the worst case must be one of the extreme points of the uncertainty set, which is the main source of conservatism in the traditional RO paradigm. In contrast, in the classic two-stage SO, uncertain parameter $\xi$ is modeled through a certain probability distribution $P$, and the expected cost is minimized, giving rise to 
\begin{equation}
\label{eq:App-03-TTSP}
\begin{aligned}
\min ~~ & c^T x + \mathbb E_P [Q(x,\xi)]  \\
\mbox{s.t.}~~ & x \in X   
\end{aligned}
\end{equation}
where the bounded polyhedron $X$ is the feasible region of first-stage decision $x$, $\xi$ is the uncertain parameter, and $Q(x,\xi)$ is the optimal value function of the second-stage problem, which is an LP  for fixed $x$ and $\xi$
\begin{equation}
\label{eq:App-03-TTSP-Stage2}
\begin{aligned}
Q(x,\xi) = \min ~~ & q^T y   \\
\mbox{s.t.}~~ & B(\xi) y \le b(\xi) - A(\xi) x   
\end{aligned}
\end{equation}
where $q$ is the cost coefficients, $A(\xi)$, $B(\xi)$, and $b(\xi)$ are constant matrices affected by uncertain data, $y(\xi)$ is the second-stage decision, which is the reaction to the realization of uncertainty.

Since the true PDF of $\xi$ is difficult to obtain in some circumstances, in this section, we do not require perfect knowledge on the probability distribution $\mathbb P$ of random variable $\xi$, and let it be ambiguous around a reference distribution and reside in an ambiguity set $D$, which can be constructed from limited historical data. We take all possible distributions in the  ambiguity set into consideration, so as to minimize the expected cost in the worst-case distribution, resulting in the following model
\begin{equation}
\label{eq:App-03-RTSP}
\begin{aligned}
\min ~~ & c^T x + \max_{P(\xi) \in D} \mathbb E_P [Q(x,\xi)]  \\
\mbox{s.t.}~~ & x \in X   
\end{aligned}
\end{equation}

Compared with (\ref{eq:App-03-RCCP-Model}), constraint violation is not allowed in problem (\ref{eq:App-03-RTSP}), and the second-stage expected cost in the worst-case distribution is considered. It is a particular case of (\ref{eq:App-03-RCCP-MEdy-1}) without chance constraints. Specifically, we will utilize discrete distributions in this section. This formulation enjoys several benefits. One is the easy exposition of the density function. In previous sections, the candidate in the moment or divergence based ambiguity sets is not given in an analytical form, and vanishes during the dual transformation. As a result, we don't have clear knowledge on the worst-case distribution. For discrete distributions, the density function is a vector of real entries associated with the probability of each representative scenario. We can easily construct the ambiguity set and optimize an expectation over discrete distributions. The other originates from the computational perspective, which can be seen later. Main results in this section come from \cite{App03-Sect4-DRSO-Zhao,App03-Sect4-DRSO-Ding}.

\vspace{12pt}
{\noindent \bf 1. Modeling the confidence set}

For a given set of historical data with $M$ elements, which can be regarded as $M$ samples of the random variable, we can draw a histogram with $K$ bins as an estimation of the reference distribution. Suppose that the numbers of samples fall in each bin is $M_1,M_2,\cdots,M_K$, where $\sum^K_{i=1} M_i=M$,  then the reference (empirical) distribution of the uncertain data is given by $\mathbb P_0 = [p^0_1,\cdots,p^0_K]$, where $p^0_i = M_i/M$, $i=1,\cdots,K$. Since the data may not be enough to fit a PDF with high accuracy, the actual distribution should be close to but might be different from its reference. It is proposed in \cite{App03-Sect4-DRSO-Zhao} to  construct the ambiguity set using statistical inference corresponding to a given tolerance. Two types of ambiguity sets are suggested based on $L_1$ norm and $L_{\infty}$ norm
\begin{equation}
\label{eq:App-03-RTSP-D1}
D_1 = \left\{ \mathbb P \in \mathbb R^K_+ \middle| \| \mathbb P - \mathbb P_0 \|_1 \le \theta \right\} =  \left\{ p \in {\rm \Delta}_K \middle| \sum_{i=1}^K \left| p_i - p^0_i \right| \le \theta \right\}
\end{equation}
\begin{equation}
\label{eq:App-03-RTSP-D2}
D_\infty = \left\{ \mathbb P \in \mathbb R^K_+ \middle| \| \mathbb P - \mathbb P_0 \|_\infty \le \theta \right\} =  \left\{ p \in {\rm \Delta}_K \middle| \max_{1\le i \le K} \left| p_i - p^0_i \right| \le \theta \right\}
\end{equation}
where ${\rm \Delta}_K = \{p \in [0,1]^K:{\bf 1}^T p = 1\}$. These two ambiguity sets can be easily expressed by polyhedral sets as follows  
\begin{equation}
\label{eq:App-03-RTSP-D1-Poly}
D_1 = \left\{ p \in {\rm \Delta}_K  \middle| 
\begin{lgathered}
\exists t \in \mathbb R^K_+ : 
\sum\nolimits_{k=1}^K t_k \le \theta \\
t_k \ge p_k - p^0_k,~ k = 1,\cdots,K  \\
t_k \ge p^0_k - p_k,~ k = 1,\cdots,K  
\end{lgathered}
\right\}
\end{equation}
\begin{equation}
\label{eq:App-03-RTSP-D2-Poly}
D_\infty = \left\{ p \in {\rm \Delta}_K \middle| 
\begin{lgathered}
\theta \ge p_k - p^0_k,~ k=1,\cdots, K \\ 
\theta \ge p^0_k - p_k,~ k=1,\cdots, K 
\end{lgathered} \right\}
\end{equation}
where $p=[p_1,\cdots,p_K]^T$ is the variable in the ambiguity set; $t=[t_1,\cdots,t_K]^T$ is the lifting (auxiliary) variable in $D_1$; parameter $\theta$ reflects decision maker's confidence level on the distance between the reference distribution and the true one. Apparently, the more historical data we utilize, the smaller their distance will be. Provided with $M$ observations and $K$ bins, the quantitative relation between the value of $\theta$ and the number of samples are given by \cite{App03-Sect4-DRSO-Zhao}
\begin{equation}
\label{eq:App-03-RTSP-D1-Conf}
\Pr \{ \|\mathbb P - \mathbb P_0\|_1 \le \theta \} \ge 
1-2K {\rm e}^{-2M \theta/K}
\end{equation}
\begin{equation}
\label{eq:App-03-RTSP-D2-Conf}
\Pr \{ \|\mathbb P - \mathbb P_0\|_\infty \le \theta \} \ge 
1-2K {\rm e}^{-2M \theta}
\end{equation}

According to (\ref{eq:App-03-RTSP-D1-Conf}) and (\ref{eq:App-03-RTSP-D2-Conf}), if we want to maintain (\ref{eq:App-03-RTSP-D1}) and (\ref{eq:App-03-RTSP-D2}) with a confidence level of $\beta$, parameter $\theta$ should be selected as 
\begin{equation}
\label{eq:App-03-RTSP-D1-Theta}
\mbox{For } D_1:~ \theta_1 = \dfrac{K}{2M} \ln \dfrac{2K}{1-\beta}
\end{equation}
\begin{equation}
\label{eq:App-03-RTSP-D2-Theta}
\mbox{For } D_\infty:~ \theta_\infty = \dfrac{1}{2M} \ln \dfrac{2K}{1-\beta}
\end{equation}

As the size of sampled data approaches infinity, $\theta_1$ and $\theta_\infty$ decrease to 0, and the reference distribution converges to the true one. Accordingly, problem (\ref{eq:App-03-RTSP}) becomes a traditional two-stage SO.

\vspace{12pt}
{\noindent \bf 2. CCG based decomposition algorithm}

Let $\xi^k$ denote the representative scenario of the $k$-th bin, $p_k$ be the corresponding probability, and $P = [p_1,\cdots,p_K]$ belongs to the ambiguity set in form of (\ref{eq:App-03-RTSP-D1}) or (\ref{eq:App-03-RTSP-D2}), then problem (\ref{eq:App-03-RTSP}) can be written as 
\begin{equation}
\label{eq:App-03-RTSP-Discrete}
\begin{aligned}
\min ~~ & c^T x + \max_{\mathbb P} \sum_{k=1}^K p_k \min q^T y^k \\
\mbox{s.t.}~~ & x \in X,~ \mathbb P \in D  \\
& A(\xi^k) x + B(\xi^k) y^k \le b(\xi^k), \forall k 
\end{aligned}
\end{equation}

Problem (\ref{eq:App-03-RTSP-Discrete}) has a min-max-min structure and can be solved by the Benders decomposition method \cite{App03-Sect4-DRSO-Zhao} or the CCG method \cite{App03-Sect4-DRSO-Ding}. The latter one will be introduced in the rest of this section. It decomposes problem (\ref{eq:App-03-RTSP-Discrete}) into a lower bounding master problem and an upper bounding subproblem, which are solved iteratively until the gap between the upper bound and lower bound gets smaller than a convergence tolerance. The basic idea has been explained in Appendix \ref{App-C-Sect02-03}. As we can see in \cite{App03-Sect4-DRSO-Ding}, the second-stage problem can be a broader class of convex programs, such as an SOCP.

\vspace{12pt}
{\noindent \bf 1) Subproblem}

For a given first-stage decision $x$, the subproblem aims to find the worst-case distribution, which comes down to a max-min program shown below 
\begin{equation}
\label{eq:App-03-RTSP-SP-1}
\max_{\mathbb P \in D} \sum_{k=1}^K p_k \min_{y^k \in Y_k(x)} q^T y^k  
\end{equation}
where 
\begin{equation}
\label{eq:App-03-RTSP-SP-Yk}
Y_k = \{ y^k ~|~ B(\xi^k) y^k \le b(\xi^k) - A(\xi^k) x\},~ \forall k  
\end{equation}

Problem (\ref{eq:App-03-RTSP-SP-1}) has some unique features
 that  facilitate the computation: 

(1) Feasible sets $Y_k$ are decoupled. 

(2) The probability variables $p_k$ do not affect feasible sets $Y_k$. 

(3) The ambiguity set $D$ and feasible sets $Y_k$ are decoupled.

Although  (\ref{eq:App-03-RTSP-SP-1}) seems nonlinear due to the production of scalar variable $p_k$ and vector variable $y^k$ in the objective function, as we can see in the following discussion, it is equivalent to an LP or can be decomposed into several LPs, and thus can be solved efficiently.  

\vspace{6pt}
{\emph{An equivalent LP} }

Because $p_k \ge 0$, we can exchange the summation operator and the minimization operator, and problem  (\ref{eq:App-03-RTSP-SP-1}) can be written as
\begin{equation}
\label{eq:App-03-RTSP-SP-2}
\max_{\mathbb P \in D} \min_{y^k \in Y_k(x)}  \sum_{k=1}^K p_k q^T y^k  
\end{equation}

For the inner minimization problem, $p_k$ is constant, so it is an LP, whose dual problem is       
\begin{equation*}
\begin{aligned}
\max_{\mu^k}~~ & \sum_{k=1}^K \left(b(\xi^k)-A(\xi^k) x \right)^T \mu^k \\
\mbox{s.t.}~~ & \mu_k \le 0,~ B^T (\xi^k) \mu^k = p_k q,~ \forall k 
\end{aligned}
\end{equation*}
where $\mu^k$ are dual variables. Substituting it into (\ref{eq:App-03-RTSP-SP-2}), and combining two maximization operators, we obtain 
\begin{equation}
\label{eq:App-03-RTSP-SP-3}
\begin{aligned}
\max_{p_k,\mu^k}~~ & \sum_{k=1}^K \left(b(\xi^k)-A(\xi^k) x \right)^T \mu^k\\
\mbox{s.t.}~~ & \mu_k \le 0,~ B^T (\xi^k) \mu^k = p_k q,~ \forall k  \\
& (p_1,\cdots,p_k) \in D
\end{aligned}
\end{equation}
Since $D$ is polyhedral, problem (\ref{eq:App-03-RTSP-SP-3}) is in fact an LP. The optimal solution offers the worst-case distribution $[p^*_1,\cdots,p^*_K]$, which will be used to generate cuts in the master problem. The recourse actions $y^k$ in each scenario will be provided by the optimal solution of the master problem. 

Despite of the fact that LP is acknowledged as the most tractable mathematical programming problem, however, when $K$ is extremely large, it is still challenging to solve (\ref{eq:App-03-RTSP-SP-3}) or even store it in a computer. Nevertheless, the separability of feasible regions allows solving (\ref{eq:App-03-RTSP-SP-1}) in a decomposition manner. 
 
\vspace{6pt}
{\emph{A decomposition method} }

As mentioned above, $p_k$ has no impact on $Y_k$, which are decoupled; moreover, because $p_k$ is a scalar in the objective function of each inner minimization problem, it does not affect the optimal solution $y^k$. In view of this convenience, problem (\ref{eq:App-03-RTSP-SP-1}) can be decomposed into $K+1$ smaller LPs, and can be solved in parallel. To this end, for each $\xi^k$, solve the following LP:
\begin{equation*}
h^*_k = \min_{y^k \in Y_k(x)} q^T y^k,~ k=1,\cdots,K
\end{equation*}
The optimal value is $h^*_k$; after obtaining optimal values $(h^*_1,\cdots,h^*_K)$ of the $K$ LPs, we can retrieve the worst-case distribution through solving an additional LP 
\begin{equation*}
\max_{\mathbb P \in D}~~ \sum_{k=1}^K  p_k h^*_k
\end{equation*}

In fact, if the second-stage problem is a conic program (in \cite{App03-Sect4-DRSO-Ding}, it is  an SOCP), above discussions are still valid, as long as the strong duality holds. 

It is interesting to notice that in the ARO problem in Sect. \ref{App-C-Sect02-03}, the subproblem comes down to a non-convex bilinear program after dualizing the inner minimization problem, and is generally NP-hard; in this section, the subproblem actually gives rise to LPs, whose complexity is polynomial in problem sizes. The reason accounting for this difference is that the uncertain parameter in (\ref{eq:App-03-RTSP-Discrete}) is expressed by sampled scenarios and thus is constant; the distributional uncertainty appearing in the objective function does not influence the constraints of the second stage problem, and thus the linear max-min problem (\ref{eq:App-03-RTSP-SP-2}) reduces to an LP after a dual transformation.

\vspace{12pt}
{\noindent \bf 2) The CCG algorithm}

The motivation of CCG algorithm has been thoroughly discussed in Appendix \ref{App-C-Sect02-03}. In this section, for a fixed $x$, the optimal value of subproblem  (\ref{eq:App-03-RTSP-SP-1}) is denoted by $Q(x)$, and $c^T x + Q(x)$ gives an upper bound of the optimal solution of (\ref{eq:App-03-RTSP-Discrete}), because the first-stage variable is un-optimized. Then a set of new variables and optimality cuts are generated and added into master problem. If the subproblem is infeasible in some scenario, then a set of feasibility cuts are assigned to the master problem. The master problem starts from a subset of $D$, which is updated by including the worst-case distribution identified by the subproblem. Forasmuch, the master problem is a relax version of the original problem (\ref{eq:App-03-RTSP-Discrete}), and provides a lower bound on the optimal value. The flowchart of the CCG procedure for problem (\ref{eq:App-03-RTSP-Discrete}) is given in Algorithm \ref{Ag:App-03-RTSO-CCG}. This algorithm will terminate in a finite number of iterations, as the confidence set $D$ has finite extreme points.

\begin{algorithm}[!htp]
\normalsize
\caption{\bf }
\begin{algorithmic}[1]
\STATE Choose a convergence tolerance $\varepsilon>0$, and an initial probability vector $p^0 \in D$; Set LB $=-\infty$, UB $=+\infty$, and iteration index $s=0$. 
\STATE Solve the master problem
\begin{equation}
\label{eq:App-03-RTSO-CCG-MP}
\begin{aligned}
\min_{x,\eta,y^{k,m}} ~~ & c^T x + \eta   \\
\mbox{s.t.} ~~ & x \in X,~~ \eta \ge \sum_{k=1}^K p^m_k q^T y^{k,m},~ 
m \in \mbox{Opt}\{ 0,1,\cdots,s\} \\
& A(\xi^k) x + B(\xi^k) y^{k,m} \le b(\xi^k),~
m \in \mbox{Opt}\{ 0,1,\cdots,s\},~ \forall k  \\
& A(\xi^k) x + B(\xi^k) y^{k,m} \le b(\xi^k),~
m \in \mbox{Fea}\{ 0,1,\cdots,s\},~k \in I(s)
\end{aligned}
\end{equation}
where Opt$\{*\}/$Fea$\{*\}$ selects the iterations in which an optimality (feasibility) cut is generated; $I(s)$ depicts the index of scenarios in which the second-stage problem is infeasible in iteration $s$. The optimal solution is $(x^*,\eta^*)$; update LB $=c^T x^* + \eta^*$;
\STATE Solve subproblem (\ref{eq:App-03-RTSP-SP-1}) with current $x^*$. If there exists some $\xi^k$ such that $Y_k(x^*)=\emptyset$, then generate new variable $y^{k,s}$, update $I(s)$, and add the following feasibility cut to the master problem
\begin{equation}
\label{eq:App-03-RTSO-CCG-Fea-Cut}
 A(\xi^k) x + B(\xi^k) y^{k,s} \le b(\xi^k),~ k \in I(s)
\end{equation}
Otherwise, if $Y_k(x^*) \ne \emptyset, \forall k$, subproblem (\ref{eq:App-03-RTSP-SP-1}) can be solved. The optimal solution is $p^{s+1}$, and the optimal value is $Q(x^*)$; update UB $=$ min$\{\mbox{UB}, c^T x^* + Q(x^*)\}$, create new variables $(y^{1,s+1},\cdots,y^{k,s+1})$, and add the following optimality cut to the master problem
\begin{equation}
\label{eq:App-03-RTSO-CCG-Opt-Cut}
\begin{gathered}
\eta \ge \sum_{k=1}^K p^{s+1}_k q^T y^{k,s+1}  \\
A(\xi^k) x + B(\xi^k) y^{k,s+1} \le b(\xi^k),~ \forall k
\end{gathered}
\end{equation}
\STATE If UB$-$LB$< \varepsilon$, terminate and report the optimal first-stage solution $x^*$ as well as the worst-case distribution $p^{s+1}$; otherwise, update $s \leftarrow s+1$, and go to step 2.
\end{algorithmic}
\label{Ag:App-03-RTSO-CCG}
\end{algorithm}

\subsection{Formulations based on Wasserstein Metric}

Up to now, the KL-divergence based ambiguity set based formulations have received plenty of research, because it enjoys some convenience when deriving the robust counterpart. For example, it has already known in Sect. \ref{App-C-Sect04-01} that robust chance constraints under KL-divergence ambiguity set can reduce to a traditional chance constraints under the empirical distribution with a rescaled confidence level, and the worst-case expectation problem under KL-divergence ambiguity set is equivalent to a convex program. However, according to its definition, KL-divergence ambiguity set may encounter theoretical difficulty to represent confidence sets for continuous distribution \cite{Am-Set-Wasserstein-1}, because the empirical distribution calibrated from finite data must be discrete, and any distribution in the KL-divergence ambiguity set must assign positive probability mass to each sampled scenario. As a continuous distribution has a density function, it must reside outside the KL-divergence ambiguity set regardless of the sampled scenarios. In contrast, Wasserstein metric based ambiguity sets contain both discrete and continuous distributions. It offers an explicit confidence level for the unknown distribution belonging to the set, and enables the decision maker more informative guidance to control the model conservativeness. This section introduces state-of-the-art results in robust SO with Wasserstein metric based ambiguity sets. The most critical problem is the robust counterparts of the worst-case expectation problem and robust chance constraints, which will be discussed respectively. They can be embedded in single- and two-stage robust SO problems without substantial barriers. The materials in this section mainly come from \cite{Am-Set-Wasserstein-1}.

\vspace{12pt}
{\noindent \bf 1. Wasserstein metric based ambiguity set}

Let $\rm \Xi$ be the support set of multi-dimensional random variable $\xi \in \mathbb R^m$. $M({\rm \Xi})$ represent all probability distributions $\mathbb Q$ supported on $\rm \Xi$, and $\mathbb E_{\mathbb Q}[\|\xi\|]=\int_{\rm \Xi} \|\xi\| \mathbb Q({\rm d}\xi) < \infty$, where $\|\cdot\|$ stands for an arbitrary norm on $\mathbb R^m$. 
\begin{definition}
\label{def:App-03-Wass-Metric}
Wasserstein metric $d_W: M({\rm \Xi}) \times M({\rm \Xi}) \to \mathbb R_+$ is defined as
\begin{equation*}
d_W(\mathbb Q,\mathbb Q_0) = \inf \left( \int_{\rm \Xi^2} \left\| \xi - \xi^0 \right\| {\rm \Pi} ({\rm d} \xi, {\rm d} \xi^0) \middle| 
\begin{gathered}
\mbox{ $\rm \Pi$ is a joint distribution of $\xi$ and} \\
\mbox{ $\xi^0$ with marginals $\mathbb Q$ and $\mathbb Q_0$} 
\end{gathered} \right)
\end{equation*}
for two probability distributions $\mathbb Q, \mathbb Q_0 \in M({\rm \Xi})$.
\end{definition}

As a special case, for two discrete distributions, Wasserstein metric is given by
\begin{equation}
d_W(\mathbb Q,\mathbb Q_0) = \inf_{\pi \ge 0} \left( \sum_i \sum_j \pi_{ij} \left\| \xi_j - \xi^0_i \right\| ~\middle|~ 
\begin{gathered}
\sum\nolimits_j \pi_{ij} = p^0_i,~ \forall i \\
\sum\nolimits_i \pi_{ij} = p_j,~   \forall j 
\end{gathered}~~ \right)
\label{eq:App-03-Wass-Def-Dis}
\end{equation}
where $p^0_i$ and $p_j$ denote the probability of representative scenario $\xi^0_i$ and $\xi_j$.

In either case, the decision variable $\rm \Pi$ (or $\pi_{ij}$) represents the probability mass transported from $\xi^0_i$ to $\xi_j$, therefore, the Wasserstein metric can be viewed as the minimal cost of a transportation plan, where the distance $\| \xi_j - \xi^0_i \|$ encodes the transportation cost of unit mass.

Sometimes, the Wasserstein metric can be represented in the dual form 
\begin{equation}
d_W(\mathbb Q,\mathbb Q_0) = \sup_{f \in L} \left( 
\int_{\rm \Xi} f(\xi) \mathbb Q({\rm d} \xi) - 
\int_{\rm \Xi} f(\xi) \mathbb Q_0({\rm d} \xi) \right)
\label{eq:App-03-Wass-Def-Dual}
\end{equation}
where $L=\{f:|f(\xi)-f(\xi^0)|\le \| \xi - \xi^0 \|,\forall \xi,\xi^0 \in {\rm \Xi} \}$ (Theorem 3.2, \cite{Am-Set-Wasserstein-1}, which was firstly discovered by Kantorovich and Rubinstein \cite{Kantorovich-Rubinstein} for distributions
with a bounded support). 

With above definition, the Wasserstein ambiguity set is the ball of radius $\epsilon$ centered at the empirical distribution $\mathbb Q_0$
\begin{equation}
\label{eq:App-03-Wass-Ambiguity-Set}
D_W = \left\{ \mathbb Q \in M({\rm \Xi}):d_W(\mathbb Q,\mathbb Q_0) \le \epsilon \right\}
\end{equation}
where $\mathbb Q_0$ is constructed with $N$ independent data samples
\begin{equation*}
\mathbb Q_0 = \frac{1}{N} \sum_{i=1}^N \delta_{\xi^0_i}
\end{equation*}
where $\delta_{\xi^0_i}$ stands for Dirac distribution concentrating unit mass at $\xi^0_i$.

Particularly, we require the unknown distribution $\mathbb Q$ follow a light tail assumption, i.e., there exists $a > 1$ such that 
\begin{equation*}
\int_{\rm \Xi} {\rm e}^{\|\xi\|^a} \mathbb Q({\rm d} \xi) < \infty
\end{equation*}

This assumption indicates that the tail of distribution $\mathbb Q$ decays at an exponential rate. If $\rm \Xi$ is bounded and compact, this assumption trivially holds. Under this assumption, modern measure concentration theory provides the following finite sample guarantee for the unknown distribution belonging to Wasserstein ambiguity set
\begin{equation}
\Pr\left[ d_W (\mathbb Q, \mathbb Q_0) \ge \epsilon \right] \le 
\begin{dcases}
c_1 {\rm e}^{-c_2 N \epsilon^{\max\{m,2\}}}  & \mbox{ if }\epsilon \le 1\\
c_1 {\rm e}^{-c_2 N \epsilon^a}  &  \mbox{ if } \epsilon > 1
\end{dcases}
\label{eq:App-03-Wass-Conf-Level}
\end{equation}
where $c_1, c_2$ are positive constants depending on $a$, $A$, and $m$ and $m \ne 2$.

Equation (\ref{eq:App-03-Wass-Conf-Level}) provides a priori estimate of the confidence level for $\mathbb Q \notin D_W$. On the other hand, we can utilize (\ref{eq:App-03-Wass-Conf-Level}) to select parameter $\epsilon$ of the Wasserstein ambiguity set such that $D_W$ contains the uncertain distribution $\mathbb Q$ with probability $1-\beta$ for some prescribed $\beta$. This requires solving $\epsilon$ from the right-hand side of (\ref{eq:App-03-Wass-Conf-Level}) with a given left-hand side $\beta$, resulting in
\begin{equation}
\epsilon = \begin{dcases}
\left( \dfrac{\ln(c_1 \beta^{-1})}{c_2 N} \right)^{1/\max\{m,2\}}  & \mbox{ if } N \ge \dfrac{\ln(c_1 \beta^{-1})}{c_2} \\
\left( \dfrac{\ln(c_1 \beta^{-1})}{c_2 N} \right)^{1/a}  & \mbox{ if } N < \dfrac{\ln(c_1 \beta^{-1})}{c_2}
\end{dcases}
\label{eq:App-03-Wass-Radius}
\end{equation}
Wasserstein ambiguity set with above radius can be regarded as a confidence
set for the unknown distribution $\mathbb Q$ as in statistical testing.

\vspace{12pt}
{\noindent \bf 2. Worst-case expectation problem}

A robust SO problem under Wasserstein metric naturally requests to minimize the worst-case expected cost:
\begin{equation}
\label{eq:App-03-Wass-RSO}
\inf_{x \in X} \sup_{\mathbb Q \in D_W} \mathbb E_{\mathbb Q} [h(x,\xi)]
\end{equation}
We demonstrate how to solve the core problem: the worst-case expectation
\begin{equation}
\label{eq:App-03-Wass-SupE}
\sup_{\mathbb Q \in D_W} \mathbb E_{\mathbb Q} [l(\xi)]
\end{equation}
where $l(\xi) = \max_{1 \le k \le K} l_k(\xi)$ is the payoff function, consisting of the point-wise maximum of $K$ elementary functions. For notation brevity, the dependence on $x$ is suppressed and will be recovered later on when necessary. We further assume that the support set $\rm \Xi$ is closed and convex, and specific $l(\xi)$ will be discussed.

Problem (\ref{eq:App-03-Wass-SupE}) renders an infinite-dimensional optimization
problem for continuous distribution. Nonetheless, the inspiring work in \cite{Am-Set-Wasserstein-1} show that (\ref{eq:App-03-Wass-SupE}) can be reformulated as a finite-dimensional convex program for various payoff functions. To see this, expand the worst-case expectation as 
\begin{equation*}
\sup_{\mathbb Q \in D_W} \mathbb E_{\mathbb Q} [l(\xi)] = \left\{
\begin{aligned}
\sup_{\rm \Pi}~ & \int_{\rm \Xi} l(\xi) \mathbb Q ({\rm d} \xi)   \\
\mbox{s.t.}~ &   \int_{\rm \Xi^2} \left\| \xi - \xi^0 \right\| {\rm \Pi} ({\rm d} \xi, {\rm d} \xi^0) \le \epsilon  \\
& \mbox{ $\rm \Pi$ is a joint distribution of $\xi$} \\
& \mbox{ $\xi^0$ with marginals $\mathbb Q$ and $\mathbb Q_0$} 
\end{aligned} \right.
\end{equation*}

According to the law of total probability, $\rm \Pi$ can be decomposed as the marginal distribution $\mathbb Q_0$ of $\xi^0$ and the conditional distributions $\mathbb Q_i$ of $\xi$ given $\xi^0 = \xi^0_i$:  
\begin{equation*}
{\rm \Pi} = \frac{1}{N} \sum_{i=1}^N \delta_{\xi^0_i} \otimes \mathbb Q_i
\end{equation*}
and the worst-case expectation evolves into a generalized moment problem in conditional distributions $\mathbb Q_i$, $i \le N$
\begin{equation*}
\sup_{\mathbb Q \in D_W} \mathbb E_{\mathbb Q} [l(\xi)] = \left\{
\begin{aligned}
\sup_{\mathbb Q_i \in M({\rm \Xi})}  & \frac{1}{N} \sum_{i=1}^N 
\int_{\rm \Xi} l(\xi) \mathbb Q_i ({\rm d} \xi)   \\
\mbox{s.t.} ~~~ &  \frac{1}{N} \sum_{i=1}^N \int_{\rm \Xi} 
\left\| \xi - \xi^0_i \right\|  \mathbb Q_i ({\rm d} \xi) \le \epsilon
\end{aligned}   \right.
\end{equation*}
Using standard Lagrangian duality, we obtain
\begin{equation*}
\begin{aligned}
\sup_{\mathbb Q \in D_W} \mathbb E_{\mathbb Q} [l(\xi)] =
\sup_{\mathbb Q_i \in M({\rm \Xi})} \inf_{\lambda \ge 0} &
\frac{1}{N} \sum_{i=1}^N  \int_{\rm \Xi} l(\xi) \mathbb Q_i ({\rm d} \xi)\\
& + \lambda \left(\epsilon -\frac{1}{N} \sum_{i=1}^N \int_{\rm \Xi} 
\left\| \xi - \xi^0_i \right\|  \mathbb Q_i ({\rm d} \xi) \right)
\end{aligned}
\end{equation*}
\begin{equation*}
\begin{aligned}
& \le \inf_{\lambda \ge 0} \sup_{\mathbb Q_i \in M({\rm \Xi})}  \lambda \epsilon + \frac{1}{N} \sum_{i=1}^N \int_{\rm \Xi} \left( l(\xi) - \lambda \left\| \xi - \xi^0_i \right\| \right) \mathbb Q_i ({\rm d} \xi)  \\
& = \inf_{\lambda \ge 0} \lambda \epsilon + \frac{1}{N} \sum_{i=1}^N \sup_{\xi \in {\rm \Xi}} \left( l(\xi) - \lambda \left\| \xi - \xi^0_i \right\| \right)
\end{aligned}
\end{equation*}
Decision variables $\lambda$ and $\xi$ have finite dimensions. The last problem can be reformulated as  
\begin{equation}
\begin{aligned}
\inf_{\lambda, s_i}~& \lambda \epsilon +  \frac{1}{N} \sum_{i=1}^N s_i  \\
\mbox{s.t.} ~ & \sup_{\xi \in {\rm \Xi}} \left( l_k(\xi) - \lambda \left\| \xi - \xi^0_i \right\| \right)  \le s_i \\
& i=1,\cdots,N,~k=1,\cdots,K  \\
& \lambda \ge 0
\end{aligned}
\label{eq:App-03-Wass-SupE-Reduce-1}
\end{equation}
From the definition of dual norm, we know $\lambda \left\| \xi - \xi^0_i \right\| = \max_{\|z_{ik}\|_* \le \lambda} \langle z_{ik}, \xi-\xi^0_i \rangle$, so the constraints give rise to
\begin{equation*}
\begin{aligned}
\sup_{\xi \in {\rm \Xi}} \left( l_k(\xi) - \lambda \left\| \xi - \xi^0_i \right\|\right)
= & \sup_{\xi \in {\rm \Xi}} \left( l_k(\xi) - \max_{\|z_{ik} \|_* \le \lambda} \langle z_{ik}, \xi-\xi^0_i \rangle \right)  \\ 
= &  \sup_{\xi \in {\rm \Xi}} \min_{\|z_{ik}\|_* \le \lambda} l_k(\xi) - \langle z_{ik}, \xi-\xi^0_i \rangle  \\
\le & \min_{\|z_{ik}\|_* \le \lambda} \sup_{\xi \in {\rm \Xi}}~ l_k(\xi) - \langle z_{ik}, \xi-\xi^0_i \rangle 
\end{aligned}
\end{equation*}
Substituting it into problem (\ref{eq:App-03-Wass-SupE-Reduce-1}) leads to a more restricted feasible set and a larger objective value, yielding 
\begin{equation}
\begin{aligned}
\inf_{\lambda, s_i}~& \lambda \epsilon +  \frac{1}{N} \sum_{i=1}^N s_i  \\
\mbox{s.t.} ~ & \min_{\|z_{ik}\|_* \le \lambda} \sup_{\xi \in {\rm \Xi}}~ l_k(\xi) - \langle z_{ik}, \xi-\xi^0_i \rangle  \le s_i \\
& i=1,\cdots,N,~ k = 1,\cdots, K  \\
& \lambda \ge 0
\end{aligned}
\label{eq:App-03-Wass-SupE-Reduce-2}
\end{equation}
The constraints of (\ref{eq:App-03-Wass-SupE-Reduce-2}) trivially suggests the feasible set of $\lambda$ is $\lambda \ge \|z_{ik}\|_*$, and the min operator in constraints can be omitted because it is in compliance with the objective function. Therefore, we arrive at 
\begin{equation}
\begin{aligned}
\inf_{\lambda, s_i}~& \lambda \epsilon +  \frac{1}{N} \sum_{i=1}^N s_i  \\
\mbox{s.t.} ~ & \sup_{\xi \in {\rm \Xi}}~ \left( l_k(\xi) - \langle z_{ik}, \xi \rangle \right) + \langle z_{ik}, \xi^0_i \rangle \le s_i,~ \lambda \ge \|z_{ik}\|_* \\
& i=1,\cdots,N, ~ k = 1,\cdots K    
\end{aligned}
\label{eq:App-03-Wass-SupE-Reduce-3}
\end{equation}
It is proved in \cite{Am-Set-Wasserstein-1} that problems (\ref{eq:App-03-Wass-SupE}) and (\ref{eq:App-03-Wass-SupE-Reduce-3}) are actually equivalent. Next, we will derive the concrete forms of (\ref{eq:App-03-Wass-SupE-Reduce-3}) under specific payoff function $l(\xi)$ and uncertainty set $\rm \Xi$. Unlike \cite{Am-Set-Wasserstein-1} which relies on conjugate functions in convex analysis, we mainly exploit LP duality theory, which is more friendly to readers with engineering background.

{\bf Case 1:} Convex PWL payoff function  $l(\xi) = \max_{1 \le k \le K} \{a^T_k \xi + b_k\}$ and bounded polyhedral uncertainty set ${\rm \Xi} = \{\xi \in \mathbb R^m:C \xi \le d\}$. The key point is the supremum regarding $\xi$ in the following constraint
\begin{equation*}
\sup_{\xi \in {\rm \Xi}} \left( a^T_k \xi - \langle z_{ik}, \xi \rangle \right) + b_k + \langle z_{ik}, \xi^0_i \rangle \le s_i
\end{equation*}
For each $k$, the supremum is an LP 
\begin{equation*}
\begin{aligned}
\max~~ & (a_k - z_{ik})^T \xi  \\
\mbox{s.t.} ~~& C \xi \le d
\end{aligned}
\end{equation*}
Its dual LP reads
\begin{equation*}
\begin{aligned}
\min~~ & d^T \gamma_{ik}  \\
\mbox{s.t.} ~~& C^T \gamma_{ik} = a_k - z_{ik} \\
& \gamma_{ik} \ge 0
\end{aligned}
\end{equation*}
Therefore, $z_{ik} = a_k - C^T \gamma_{ik}$. Because of strong duality, we can replace the  supremum by the objective of the dual LP, which gives rise to: 
\begin{equation*}
d^T \gamma_{ik} + b_k + \langle a_k - C^T \gamma_{ik}, \xi^0_i \rangle \le s_i,~ i=1,\cdots,N,~ k = 1,\cdots,K
\end{equation*}
Arrange all constraints together, we obtain a convex program which is equivalent to problem (\ref{eq:App-03-Wass-SupE-Reduce-3}) in Case 1:
\begin{equation}
\begin{aligned}
\inf_{\lambda, s_i}~& \lambda \epsilon +  \frac{1}{N} \sum_{i=1}^N s_i  \\
\mbox{s.t.} ~ & b_k + a^T_k \xi^0_i + \gamma^T_{ik} (d - C^T \xi^0_i) \le s_i,~ i=1,\cdots,N,~ k = 1,\cdots,K \\
& \lambda \ge \|a_k - C^T \gamma_{ik}\|_*,~ \gamma_{ik} \ge 0,~ i=1,\cdots,N,~ k = 1,\cdots,K
\end{aligned}
\label{eq:App-03-Wass-SupE-Conic-Case1}
\end{equation}

In the absence of distributional uncertainty, or $\epsilon = 0$ which implies that Wasserstein ambiguity set $D_W$ is a singleton, $\lambda$ can take any non-negative value without changing the objective function. Because all sampled scenarios must belong to the support set, i.e. $d - C^T \xi^0_i \ge 0$, $\forall i$ holds, so there must be $\gamma_{ik}=0$ at the optimal solution, leading to an optimal value of $\sum_{i=1}^N s_i/N$, where $s_i=\max_{1 \le k \le K} \{a^T_k \xi^0_i + b_k\}$, which represents the sample
average of the payoff function under the empirical distribution.

{\bf Case 2:} Concave PWL payoff function  $l(\xi) = \min_{1 \le k \le K} \{a^T_k \xi + b_k\}$ and bounded polyhedral uncertainty set ${\rm \Xi} = \{\xi \in \mathbb R^m:C \xi \le d\}$. In such circumstance, the supremum regarding $\xi$ in the constraint becomes
\begin{equation*}
\max_{\xi \in {\rm \Xi}}~ \left\{ -z^T_i \xi + \min_{1 \le k \le L} \left\{ a^T_k \xi + b_k \right\}  \right\}
\end{equation*}
which is equivalent to an LP 
\begin{equation*}
\begin{aligned}
\max~~ & -z^T_i \xi + \tau_i \\
\mbox{s.t.} ~~ & A \xi + b \ge \tau_i {\bf 1} \\
& C \xi \le d
\end{aligned}
\end{equation*}
where the $k$-th row of $A$ is $a^T_k$; the $k$-th entry of $b$ is $b_k$; $\bf 1$ is all-one vector with a compatible dimension. Its dual LP reads
\begin{equation*}
\begin{aligned}
\min~~ & b^T \theta_i + d^T \gamma_i \\
\mbox{s.t.} ~~& - A^T \theta_i + C^T \gamma_i= -z_i \\
& {\bf 1}^T \theta_i = 1,~ \theta_i \ge 0,~ \gamma_i \ge 0 
\end{aligned}
\end{equation*}
Therefore, $z_i = A^T \theta_i - C^T \gamma_i$. Because of strong duality, we can replace the  supremum by the objective of the dual LP, which gives rise to: 
\begin{equation*}
b^T \theta_i + d^T \gamma_i + \langle A^T \theta_i - C^T \gamma_i, \xi^0_i \rangle \le s_i,~ i=1,\cdots,N
\end{equation*}
Arrange all constraints together, we obtain a convex program which is equivalent to problem (\ref{eq:App-03-Wass-SupE-Reduce-3}) in Case 2:
\begin{equation}
\begin{aligned}
\inf_{\lambda, s_i}~~& \lambda \epsilon +  \frac{1}{N} \sum_{i=1}^N s_i \\
\mbox{s.t.} ~~ & \theta^T_i (b+A\xi^0_i) + \gamma^T_i ( d - C \xi^0_i) \le s_i,~ i=1,\cdots,N \\
& \lambda \ge \|A^T \theta_i - C^T \gamma_i \|_*,~ i=1,\cdots,N  \\
&\gamma_i \ge 0,~ \theta_i \ge 0,~ {\bf 1}^T \theta_i = 1,~ i=1,\cdots,N
\end{aligned}
\label{eq:App-03-Wass-SupE-Conic-Case2}
\end{equation}
There will be no $k$ index for the constraints, because it is packaged in $A$ and $b$.

An analogous analysis shows that if $\epsilon=0$, there must be $\gamma_i = 0$ and 
\begin{equation*}
s_i = \min \{\theta^T_i (b+A\xi^0_i): \theta_i \ge 0,~ {\bf 1}^T \theta_i = 1 \} = \min_{1 \le k \le K} \{a^T_k \xi + b_k\}
\end{equation*}
implying $\sum_{i=1}^N s_i/N$ is the sample average of the payoff function
under the empirical distribution.

Now we focus our attention on the min-max problem (\ref{eq:App-03-Wass-RSO}) which frequently arises in two-stage robust SO, which entails evaluation of the expected recourse cost from an LP parameterized in $\xi$. We investigate two cases depending on where $\xi$ appears.

{\bf Case 3:} Uncertain cost coefficients: $l(\xi) = \min_y \{y^T Q \xi: Wy \ge h - Ax\}$ where $x$ is the first-stage decision variable, $y$ represents the recourse action, and the feasible region is always non-empty. In this case, the supremum regarding $\xi$ in the constraint becomes
\begin{equation*}
\begin{aligned}
& \max_{\xi \in {\rm \Xi}}~ \left\{ -z^T_i \xi + \min_{y} \left\{ y^T Q \xi: Wy \ge h - Ax \right\}  \right\}  \\
=& \min_{y} \left\{\max_{\xi \in {\rm \Xi}} \left\{ \left(Q^T y-z_i\right)^T \xi \right\} : Wy \ge h - Ax\right\}
\end{aligned}
\end{equation*}
Replace the inner LP with its dual, we get an equivalent LP  
\begin{equation*}
\begin{aligned}
\min_{\gamma_i,y_i}~~ & d^T \gamma_i \\
\mbox{s.t.} ~~ & C^T \gamma_i = Q^T y_i - z_i,~ \gamma_i \ge 0 \\
& Wy_i \ge h - Ax
\end{aligned}
\end{equation*}
Here we associated variable $y$ with a subscript $i$ to highlight its dependence on the value of $\xi$. Therefore, $z_i = Q^T y_i - C^T \gamma_i$, and we can replace the  supremum by the objective of the dual LP, which gives rise to: 
\begin{equation*}
d^T \gamma_i + \langle Q^T y_i - C^T \gamma_i, \xi^0_i \rangle \le s_i,~ i=1,\cdots,N
\end{equation*}
Arrange all constraints together, we obtain a convex program which is equivalent to problem (\ref{eq:App-03-Wass-SupE-Reduce-3}) in Case 3:
\begin{equation}
\begin{aligned}
\inf_{\lambda, s_i}~~& \lambda \epsilon +  \frac{1}{N} \sum_{i=1}^N s_i \\
\mbox{s.t.} ~~ & y^T_i Q \xi^0_i + \gamma^T_i(d - C^T \xi^0_i) \le s_i,~ i=1,\cdots,N\\
& \lambda \ge \|Q^T y_i - C^T \gamma_i \|_*, ~\gamma_i \ge 0,~ i=1,\cdots,N\\
& W y_i \ge h-Ax,~  i=1,\cdots,N 
\end{aligned}
\label{eq:App-03-Wass-SupE-Conic-Case3}
\end{equation}

Without distributional uncertainty, $\epsilon = 0$, $\lambda$ can be arbitrary nonnegative value; for similar reason, we have $\gamma_i=0$ and $s_i =y^T_i Q \xi^0_i$ at optimum. So problem (\ref{eq:App-03-Wass-SupE-Conic-Case3}) is equivalent to the SAA problem under the empirical distribution
\begin{equation*}
\min_{y_i}~ \left\{ \frac{1}{N} \sum_{i=1}^N y^T_i Q \xi^0_i: W y_i \ge h-Ax \right\}
\end{equation*}

{\bf Case 4:} Uncertain constraint right-hand side: 
\begin{equation*}
\begin{aligned}
l(\xi) & = \min_y ~\{q^T y: Wy \ge H \xi + h - Ax\}  \\
& = \max_{\theta} \left\{ \theta^T(H\xi + h - Ax):W^T \theta = q,
\theta \ge 0 \right\}  \\
& = \max_k ~ v^T_k (H\xi + h - Ax) = \max_k \left\{ v^T_k H \xi + v^T_k(h-Ax) \right\}
\end{aligned}
\end{equation*}
where $v_k$ is the vertices of polyhedron $\{\theta:W^T \theta = q, \theta \ge 0\}$. In this way, $l(\xi)$ is expressed as a convex PWL function. Applying the result in Case 1, we obtain a convex program which is equivalent to problem (\ref{eq:App-03-Wass-SupE-Reduce-3}) in Case 4:
\begin{equation}
\begin{aligned}
\inf_{\lambda, s_i}~& \lambda \epsilon +  \frac{1}{N} \sum_{i=1}^N s_i  \\
\mbox{s.t.} ~ & v^T_k(h-Ax) + v^T_k H \xi^0_i + \gamma^T_{ik} (d - C^T \xi^0_i) \le s_i,~ i=1,\cdots,N,~ \forall k  \\
& \lambda \ge \|H^T v_k - C^T \gamma_{ik}\|_*,~ \gamma_{ik} \ge 0,~ i=1,\cdots,N,~ \forall k
\end{aligned}
\label{eq:App-03-Wass-SupE-Conic-Case4}
\end{equation}

For similar reason, without distributional uncertainty, we have $\gamma_{ik}=0$ and $s_i = v^T_k(h-Ax) + v^T_k H \xi^0_i = q^Ty_i$ at optimum, where the last equality is because of strong duality. So problem (\ref{eq:App-03-Wass-SupE-Conic-Case3}) is equivalent to the SAA problem under the empirical distribution
\begin{equation*}
\min_{y_i}~ \left\{ \frac{1}{N} \sum_{i=1}^N q^T y_i: W y_i \ge H \xi + h - Ax \right\}
\end{equation*}

The following discussions are devoted to the computational tractability.
\begin{itemize}
\item If the 1-norm or $\infty$-norm is used to define Wasserstein metric, their dual norms are $\infty$-norm and 1-norm respectively, then problems (\ref{eq:App-03-Wass-SupE-Conic-Case1})-(\ref{eq:App-03-Wass-SupE-Conic-Case4}) reduce to LPs whose sizes grow with the number $N$ of sampled data. If the Euclidean norm is used, the resulting problems will be SOCP.

\item For Case 1, Case 2 and Case 3, the remaining equivalent LPs scale polynomially and can be therefore readily solved. As for Case 4, the number of vertices may grow exponential in the problem size. However, one can adopt a decomposition algorithm similar to CCG which iteratively identifies critical vertices without enumerating all of them.

\item The computational complexity of all equivalent convex programs is independent of the size of the Wasserstein ambiguity set.

\item  It is shown in \cite{Am-Set-Wasserstein-1} that the worst-case expectation can also be computed from the following problem
\begin{equation}
\label{eq:App-03-Wass-SupE-Extreme-Q}
\begin{aligned}
\sup_{\alpha_{ik},q_{ik}} ~~& \frac{1}{N} \sum_{i=1}^N \sum_{k=1}^K 
\alpha_{ij} l_k \left( \xi^0_i - \frac{q_{ik}}{\alpha_{ik}} \right) \\
\mbox{s.t.}~~ &  \frac{1}{N} \sum_{i=1}^N \sum_{k=1}^K \|q_{ik}\| \le \epsilon\\
& \alpha_{ik} \ge 0, \forall i,\forall k,~ 
\sum_{k=1}^K \alpha_{ik} =1, \forall i  \\
& \xi^0_i - \frac{q_{ik}}{\alpha_{ik}} \in {\rm \Xi},~ 
\forall i, \forall k
\end{aligned}
\end{equation}
Non-convex term arise from the fraction $q_{ik}/\alpha_{ik}$. In fact, problem (\ref{eq:App-03-Wass-SupE-Extreme-Q}) is convex following the definition of extended perspective function \cite{Am-Set-Wasserstein-1}. Moreover, if $[\alpha_{ik}(r),q_{ik}(r)]_{r \in \mathbb N}$ is a sequence of feasible solutions and the corresponding objective values converge to the supremum of (\ref{eq:App-03-Wass-SupE-Extreme-Q}), then the discrete distribution
\begin{equation*}
\mathbb Q_r = \frac{1}{N} \sum_{i=1}^N \sum_{k=1}^K \alpha_{ik}(r)\delta_{\xi_{ik}(r)},~ \xi_{ik}(r) = \xi^0_i - \frac{q_{ik}(r)}{\alpha_{ik}(r)}
\end{equation*}
approaches the worst-case distribution in $D_W$ \cite{Am-Set-Wasserstein-1}.
\end{itemize}

\vspace{12pt}
{\noindent \bf 3. Static robust chance constraints}

Another important issue in SO is chance constraint. Here we discuss robust joint chance constraints in the following form
\begin{equation}
\label{eq:App-03-Wass-RCC-Def}
\inf_{\mathbb Q \in D_W} \Pr [a(x)^T \xi_i \le b_i(x), i=1,\cdots,I] \ge 1-\beta
\end{equation}
where $x$ is the decision variable; the chance constraint involves $I$ inequalities with uncertain parameter $\xi_i$ supported on set ${\rm \Xi}_i \subseteq \mathbb R^n$ for each $i$. The joint probability distribution $\mathbb Q$ belongs to the Wasserstein ambiguity set. $a(x) \in \mathbb R^n$ and $b(x) \in \mathbb R$ are affine mappings of $x$, where $a(x)=\eta x + (1-\eta) {\bf 1}$, $\eta \in \{0,1\}$, and $b_i(x) = B^T_ix + b^0_i$. When $\eta=1$ ($\eta=0$), (\ref{eq:App-03-Wass-RCC-Def}) involves left-hand (right-hand) uncertainty. ${\rm \Xi}=\prod_i {\rm \Xi_i}$ is the support set of $\xi = [\xi^T_1,\cdots,\xi^T_I]^T$. The robust chance constraint (\ref{eq:App-03-Wass-RCC-Def}) requires that all inequalities be met for all possible distributions in Wasserstein ambiguity set $D_W$ with a probability of at least $1-\beta$, where $\beta \in (0,1)$ denotes a prescribed risk tolerance. The feasible region stipulated by (\ref{eq:App-03-Wass-RCC-Def}) is $X$. We will introduce main results from \cite{Am-Set-Wasserstein-2} while avoiding rigorous mathematical proofs.

\begin{assumption}
\label{ap:App-03-Wass-RCC}
The support set $\rm \Xi$ is an $n \times I$-dimensional vector space, and the distance metric in Wasserstein ambiguity set is $d(\xi,\zeta)=\|\xi - \zeta\|$.
\end{assumption}

\begin{theorem}
\cite{Am-Set-Wasserstein-2} Under Assumption (\ref{ap:App-03-Wass-RCC}), $X=Z_1 \cup Z_2$, where
\begin{equation}
\label{eq:App-03-Wass-RCC-Z1}
Z_1 = \left\{x \in \mathbb R^n ~\middle|~ \begin{lgathered}
\epsilon v - \beta \gamma \le \frac{1}{N} \sum_{j=1}^N z_j  \\
z_j + \gamma \le \max \left\{ b_i(x) - a(x)^T \zeta^j_i,0 \right\} \\
i = 1, \cdots I,~ j = 1,\cdots, N  \\
z_j \le 0,~ j = 1,\cdots,N  \\
\|a(x)\|_* \le v,~ \gamma \ge 0
\end{lgathered}  \right\}
\end{equation}
where $\epsilon$ is the radius of the Wasserstein ambiguity set, $N$ is the number of sampled scenarios in the empirical distribution, and
\begin{equation}
\label{eq:App-03-Wass-RCC-Z2}
Z_2 = \{x \in \mathbb R^n ~|~ a(x) = 0,~ b_i(x) \ge 0,~i=1,\cdots I \}
\end{equation}
\label{th:App-03-Wass-RCC-Exact}
\end{theorem}

In Theorem \ref{th:App-03-Wass-RCC-Exact}, $Z_2$ is trivial: If $ \eta = 1$, then $Z_2=\{x \in \mathbb R^n~|~ x=0,~b_i \ge 0,~ \forall i\}$; If $ \eta = 0$, then $Z_2=\emptyset$. $Z_1$ can be reformulated as an MILP compatible form if it is bounded. By linearizing the second constraint, we have
\begin{equation}
Z_1 = \left\{x \in \mathbb R^n ~\middle|~ \begin{lgathered}
\epsilon v - \beta \gamma \le \frac{1}{N} \sum_{j=1}^N z_j  \\
z_j + \gamma \le s_{ij},~ \forall i, \forall j  \\
b_i(x) - a(x)^T \zeta^j_i \le s_{ij} \le M_{ij} y_{ij},~\forall i, \forall j \\
s_{ij} \le b_i(x) - a(x)^T \zeta^j_i+M_{ij}(1-y_{ij}),~\forall i, \forall j \\
\|a(x)\|_* \le v,~ \gamma \ge 0,~z_j \le 0,~ \forall j  \\
s_{ij} \ge 0,~ y_{ij} \in \{0,1\},~ \forall i,~ \forall j
\end{lgathered}  \right\}
\label{eq:App-03-Wass-RCC-Z1-MILP}
\end{equation}
where $\forall i$ and $\forall j$ are short for $i=1,\cdots,I$ and $j=1,\cdots,N$, respectively;
\begin{equation*}
M_{ij} \ge \max_{x \in Z_1} \left| b_i(x) - a(x)^T \zeta^j_i \right|
\end{equation*}

It is easy to see that if $b_i(x) - a(x)^T \zeta^j_i <0$, then $y_{ij}=0$ (otherwise $s_{ij} \le b_i(x) - a(x)^T \zeta^j_i <0$), hence $s_{ij} = 0=\max \{ b_i(x) - a(x)^T \zeta^j_i,0 \}$. If  $b_i(x) - a(x)^T \zeta^j_i >0$, then $y_{ij}=1$ (otherwise $b_i(x) - a(x)^T \zeta^j_i \le M_{ij}y_{ij} = 0$), hence $s_{ij} = b_i(x) - a(x)^T \zeta^j_i =\max\{ b_i(x) - a(x)^T \zeta^j_i,0\}$. If $b_i(x) - a(x)^T \zeta^j_i =0$, then we have $s_{ij} = 0$ regardless of the value of $y_{ij}$. In conclusion, (\ref{eq:App-03-Wass-RCC-Z1}) and (\ref{eq:App-03-Wass-RCC-Z1-MILP}) are equivalent.

 For right-hand uncertainty in which $\eta=0$, $a(x)={\bf 1}$, $X = Z_1$ because $Z_2 = \emptyset$. Moreover, variable $v$ in (\ref{eq:App-03-Wass-RCC-Z1-MILP}) is equal to 1 if  1-norm is used in Wasserstein ambiguity set $D_W$, indicating $v \ge \|{\bf 1}\|_\infty =1$ in $Z_1$.

In (\ref{eq:App-03-Wass-RCC-Z1-MILP}), a total number of $I \times N$ binary variables are introduced to linearize the $\max\{a,b\}$ function, making the problem challenging to solve. An inner approximation of $Z$ is to simply replace $\max \{ b_i(x) - a(x)^T \zeta^j_i,0 \}$ with its first input, yielding a parameter-free approximation
\begin{equation}
\label{eq:App-03-Wass-RCC-CVaR}
Z = \left\{x \in \mathbb R^n ~\middle|~ \begin{lgathered}
\epsilon v - \beta \gamma \le \frac{1}{N} \sum_{j=1}^N z_j  \\
z_j + \gamma \le b_i(x) - a(x)^T \zeta^j_i,~ \forall i,\forall j \\
z_j \le 0,~ \forall j,~ \|a(x)\|_* \le v,~ \gamma \ge 0
\end{lgathered}  \right\}
\end{equation}
This formulation can be derived from CVaR model, and enjoys better computational tractability.

\vspace{12pt}
{\noindent \bf 4. Adaptive robust chance constraints}

Robust chance constraint program with Wasserstein metric is studied in \cite{App03-Sect4-DRSO-Was-4} in a different but more general form. The problem is as follows
\begin{equation}
\label{eq:App-03-DRCC-Wass-1}
\begin{aligned}
\min_{x \in X}~~ & c^T x \\
\mbox{s.t.} ~~ & \inf_{\mathbb Q \in D_W} \Pr [F(x,\xi) \le 0] \ge 1-\beta
\end{aligned}
\end{equation}
where $X$ is a bounded polyhedron, $F: \mathbb R^n \times {\rm \Xi} \to \mathbb R$ is a scalar function that is convex in $x$ for every $\xi$. This formulation is general enough to capture joint chance constraints. To see this, suppose $F$ contains $K$ individual constraints, then $F$ can be defined as the component-wise maximum as in (\ref{eq:App-03-DRO-DRCC-Joint-Para}). 

Here we develop a technique to solve two-stage problems where $F(x,\xi)$ is the optimal value of another LP parameterized in $x$ and $\xi$. More precisely, we consider 
\begin{subequations}
\label{eq:App-03-TSDRCC-Wass}
\begin{equation}
\label{eq:App-03-TSDRCC-Wass-1}
\begin{aligned}
\min_{x \in X}~~ & c^T_1 x \\
\mbox{s.t.} ~~ & \sup_{\mathbb Q \in D_W} \Pr [f(x,\xi) \ge c^T_2 \xi] \le \beta
\end{aligned}
\end{equation}
\begin{equation}
\label{eq:App-03-TSDRCC-Wass-2}
\begin{aligned}
f(x,\xi) = \min ~~ & c^T_3 y \\
\mbox{s.t.} ~~ & Ax + By + C\xi \le d 
\end{aligned}
\end{equation}
\end{subequations}
where in (\ref{eq:App-03-TSDRCC-Wass-1}), the robust chance constraint can be regarded a risk limiting requirement, and the threshold value depends on uncertain parameter $\xi$. We assume LP (\ref{eq:App-03-TSDRCC-Wass-2}) is always feasible (relatively complete recourse) and has finite optimum. Second-stage cost can be considered in the objective function of (\ref{eq:App-03-TSDRCC-Wass-1}) in form of worst-case expectation which has been discussed in previous sections and is omitted here for the sake of brevity. Here we focus on coping with second-stage LP in robust chance constraint.  

Define loss function  
\begin{equation}
\label{eq:App-03-TSDRCC-Wass-Loss-Fun}
g(x,\xi) = f(x,\xi) - c^T_2 \xi
\end{equation}

Recall the relation between chance constraint and CVaR discussed in Sect. \ref{App-C-Sect03-01}, a sufficient condition of robust chance constraint in (\ref{eq:App-03-TSDRCC-Wass-1}) is CVaR$(g(x,\xi),\beta) \le 0$, $\forall \mathbb Q \in D_W$, or equivalently
\begin{equation}
\label{eq:App-03-TSDRCC-CVaR-1}
\sup_{\mathbb Q \in D_W} \inf_{\gamma \in \mathbb R} \beta \gamma +  \mathbb E_{\mathbb Q} (\max \{g(x,\xi)-\gamma,0\}) \le 0
\end{equation}
According to \cite{App03-Sect4-DRSO-Was-4}, constraint (\ref{eq:App-03-TSDRCC-CVaR-1}) can be conservatively approximated by
\begin{equation}
\label{eq:App-03-TSDRCC-CVaR-2}
\epsilon L + \inf_{\gamma \in \mathbb R} \left\{ \beta \gamma + \frac{1}{N} \sum_{i=1}^N \max \{g(x,\xi^i)-\gamma,0\}  \right\} \le 0
\end{equation}
where $\epsilon$ is the parameter in Wasserstein ambiguity set $D_W$, $L$ is a constant satisfying $g(x,\xi) \le L\| \xi \|_1$, and $\xi^i$, $i=1,\cdots,N$ are samples of uncertain data. Substituting (\ref{eq:App-03-TSDRCC-Wass-2}) and (\ref{eq:App-03-TSDRCC-Wass-Loss-Fun}) into (\ref{eq:App-03-TSDRCC-CVaR-2}), we obtain an LP that is equivalent to problem (\ref{eq:App-03-TSDRCC-Wass})
\begin{equation}
\label{eq:App-03-TSDRCC-CVaR-Eqv-LP}
\begin{aligned}
\min~~ & c^T_1 x \\
\mbox{s.t.} ~~ & x \in X,~ \epsilon L +  \beta \gamma + \frac{1}{N} \sum_{i=1}^N s_{i}   \le 0 \\
& s_i \ge 0,~ s_i \ge c_3^T y^i - c^T_2 \xi^i - \gamma,~ i = 1,\cdots,N \\
& Ax + By^i + C\xi^i \le d,~ i = 1,\cdots,N  \\
\end{aligned}
\end{equation}
where $y^i$ is the second-stage decision associated with $\xi^i$. This formulation could be very conservative due to three reasons. First, worst-case distribution is considered; second, CVaR constraint (\ref{eq:App-03-TSDRCC-CVaR-1}) is a pessimistic approximation of chance constraints; finally, sampling constraint (\ref{eq:App-03-TSDRCC-CVaR-2}) is a pessimistic approximation of (\ref{eq:App-03-TSDRCC-CVaR-1}). 

More discussions on robust chance constraints with Wasserstein metric under various settings can be found in \cite{App03-Sect4-DRSO-Was-4}.

\vspace{12pt}
{\noindent \bf 5. Use of forecast data}

Wasserstein metric enjoys many advantages, such as finite-sample performance guarantee and existence of tractable reformulation. However, moment information is not used, especially the first-order moment reflecting the prediction, which can be updated with time rolling on, so the worst-case distribution generally has a mean value different from the forecast (if available). To incorporate forecast data, we propose the following Wasserstein ambiguity set with fixed-mean 
\begin{equation}
\label{eq:App-03-Wass-Mom}
D^M_W = \left\{ \mathbb  Q \in D_W  \middle| \mathbb E_{\mathbb Q} [\xi] = \hat \xi  \right\}
\end{equation}
and the worst-case expectation problem can be expressed as
\begin{subequations}
\label{eq:App-03-MaxE-Wass-Mom}
\begin{align}  
\sup_{\mathbb Q \in D^M_W} ~~ & \mathbb E_{\mathbb Q} [l(\xi)] \\
 = \sup_{f^n(\xi)} ~~ & \frac{1}{N}  \sum_{n=1}^N  \int_{\rm \Xi} l(\xi) f^n(\xi) {\rm d} \xi  \\
\mbox{s.t.} ~~ &  \frac{1}{N}  \sum_{n=1}^N  \int_{\rm \Xi} \| \xi - \xi^n \|_p f^n(\xi) {\rm d} \xi  \le \epsilon : \lambda  \\
& \int_{\rm \Xi} f^n(\xi) {\rm d} \xi = 1 : \theta_n,~ n = 1,\cdots, N  \\
& \frac{1}{N}  \sum_{n=1}^N  \int_{\rm \Xi} \xi f^n(\xi) {\rm d} \xi = \hat \xi : \rho 
\end{align}
\end{subequations}
where $l(\xi)$ is a loss function similar to that in (\ref{eq:App-03-Wass-SupE}), $f^n(\xi)$ is the conditional density function under historical data sample $\xi^n$, dual variables $\lambda$, $\theta_n$, and $\rho$ are listed following a colon. Similar to the discussions for problem (\ref{eq:App-03-DRO-Worst-Expectation-Primal}), the dual problem of (\ref{eq:App-03-MaxE-Wass-Mom}) is
\begin{subequations}
\label{eq:App-03-Dual-MaxE-Wass-Mom}
\begin{align}  
\min_{\lambda \ge 0,\theta_n,\rho} ~~ &  (\lambda \epsilon + \rho^T \hat \xi) N + \sum_{n=1}^N \theta_n  \label{eq:App-03-Dual-MaxE-Wass-Mom-Obj}\\
\mbox{s.t.} ~~ &  \theta_n + \lambda \| \xi - \xi^n \|_p + \rho^T \xi \ge l(\xi), \forall \xi \in {\rm \Xi},~ \forall n  \label{eq:App-03-Dual-MaxE-Wass-Mom-Cons}
\end{align}
\end{subequations}

For $p=2$, polyhedral $\rm \Xi$ and PWL $l(\xi)$, constraint (\ref{eq:App-03-Dual-MaxE-Wass-Mom-Cons}) can be transformed into the intersection of PSD cones, and problem (\ref{eq:App-03-Dual-MaxE-Wass-Mom}) gives rise to an SDP; some examples can be found in Sect. \ref{App-C-Sect03-01}.
If $\rm \Xi$ is described by a single quadratic constraint,  constraint (\ref{eq:App-03-Dual-MaxE-Wass-Mom-Cons}) can be reformulated by using the well-known S-Lemma, which has been discussed in Sect. \ref{App-A-Sect02-04}, and problem (\ref{eq:App-03-Dual-MaxE-Wass-Mom}) still comes down to an SDP. For $p=1$ or $p=+\infty$, polyhedral $\rm \Xi$ and PWL $l(\xi)$, constraint (\ref{eq:App-03-Dual-MaxE-Wass-Mom-Cons}) can be transformed into a polyhedron using duality theory, and problem (\ref{eq:App-03-Dual-MaxE-Wass-Mom}) gives rise to an LP. Because the ambiguity set is more restrictive, problem (\ref{eq:App-03-Dual-MaxE-Wass-Mom}) would be less conservative than problem (\ref{eq:App-03-Wass-SupE-Reduce-3}) in which the mean value of uncertain data is free. 

A Wasserstein-moment metric with variance is exploited in \cite{App03-Sect4-DRSO-Was-5} and applied to wind power dispatch. Nevertheless, the ambiguity set neglects first-order moment and considers second-order moment. This formulation is useful when little historical data is available at hand.

\vspace{12pt}

As a short conclusion, distributionally robust optimization and data-driven robust stochastic optimization leverage statistical information on the uncertain data and overcome the conservatism of traditional robust optimization approaches which are built upon the worst-case scenario. The core issue is the equivalent convex reformulation of the worst-case expectation problem or the robust chance constraint over the uncertain probability distribution restricted in the ambiguity set. Optimization over a moment based  ambiguity set can be formulated as a semi-infinite LP, whose dual problem gives rise to SDPs, and hence can be readily solved. When additional structure property is taken into account, such as unimodality, more sophisticated treatment is need. As for the robust stochastic programming, tractable reformulation of the worst-case expectation and robust chance constraints is the central issue. Robust chance constraint under a $\phi$-divergence based ambiguity set are equivalent to traditional chance constraint under the empirical distribution but with a modified confidence level, and it can be transformed into an MILP or approximated by LP based on risk theory under the help of sampling average approximation technique, so does a robust chance constraint  under a Wasserstein metric based ambiguity set, following somewhat different expressions. The worst-case expectation under $\phi$-divergence based ambiguity set boils down to a convex program with linear constraints and a nonlinear objective function, which can be efficiently solved via outer approximation algorithm. The worst-case expectation under Wasserstein ambiguity set comes down to a conic program which is convex and readily solvable. Unlike the max-min problem in traditional robust optimization method identifying the worst-case scenario which the decision maker wishes to avoid, the worst-case expectation problem in distributionally robust optimization and robust stochastic programming is solved in its dual form, whose solution is less intuitive to the decision maker; moreover, it may not be easy to recover the primal optimal solution, i.e., the worst-case probability. The worst-case distribution in the robust chance constrained stochastic programming is discussed in \cite{App03-Sect4-RCCP,Am-Set-Wasserstein-1}; the worst-case discrete  distribution in a two-stage stochastic program with min-max expectation can be computed via a polynomial complexity algorithm. Nonetheless, from a practical perspective, what the human decision makers actually need to deploy is merely the here-and-now decision, and the worst probability distribution is usually not very important, since corrective actions can be postponed to a later stage when the uncertain data have been observed or can be predicted with high accuracy.

\section{Further Reading}
\label{App-C-Sect05}

Uncertainty is ubiquitous in real-life decision-making problems, and the decision maker usually has limited information and statistic data on the uncertain factors,  which makes robust optimization very attractive in practice, as it is tailored to the available information at hand, and often gives rise to computationally tractable reformulations. Although the original idea can date back to \cite{RO-Soyster} in 1970s, it is during the past two decades that the fundamental theory of robust optimization has been systematically developed. This research field is even more active during the past five years. This chapter aims to help beginners get an overview on this method and understand how to apply robust optimization in practice. We provide basic models and tractable reformulations, called the robust counterparts, for various robust optimization models under different assumptions on the uncertainty and decision-making manner. Basic theory of robust optimization is provided in \cite{RO-Detail-1,RO-Detail-2}. Comprehensive surveys can be found in \cite{RO-Guide,RO-Survey}. Here we shed more light on several important topics in robust optimization.  

Uncertainty sets play a decisive role on the performance of a robust solution. A larger set could protect the system against a higher level of uncertainty, and increase the cost as well. However, the probability that uncertain data take their wort-case values is usually small. The decision-maker needs to make a trade-off between reliability and economy. Ambiguous chance constraints and their approximations are discussed in Chapter 2 of \cite{RO-Detail-1}, based on which the parameter in the uncertainty set can be selected. It is proposed in \cite{Un-Set-Data-Driven} to construct uncertainty sets from historical data and statistical tests. The connection of uncertainty sets and coherent risk measures are revealed in \cite{Un-Set-Risk-Measure}. It is shown that the distortion risk measure leads to a polyhedral uncertainty set. Specifically, the connection of CVaR and uncertainty sets is discussed in \cite{Un-Set-CVaR}. A reverse correspondence is reported in \cite{Risk-Measure-Un-Set}, demonstrating that robust optimization could generalize the concepts of risk measures. A data-driven approach is proposed in \cite{Un-Set-Data-Driven} to construct uncertainty sets for robust optimization based on statistical hypothesis tests. The counterpart problems are shown to be tractable, and optimal solutions satisfy constraints with finite-sample probabilistic guarantee.

Distributionally robust optimization integrates statistic information, worst-case expectation, and robust probability guarantee in a holistic optimization framework, in which the uncertainty is modeled via an ambiguous probability distribution. The choice of ambiguity sets for candidate distributions affects not only the model conservatism, but also the existence of tractable reformulations. Various ambiguity sets have been proposed in the literature, which can be roughly classified into two categories:

1) Moment ambiguity sets. All PDFs share the same moment data, usually the first- and second-order moments, and  structured properties, such as symmetry and unimodality. For example, Markov ambiguity set contains all distributions with the same mean and support, and the worst-case expectation is shown to be equivalent to LPs \cite{Am-Set-Markov}. Chebyshev ambiguity set is composed of all distributions with known expectation and covariance matrix, and usually leads to SDP counterparts \cite{Static-DRO,Am-Set-Chebyshev-1,Am-Set-Chebyshev-2}; the Gauss ambiguity set contains all unimodal distributions in the Chebyshev ambiguity set, and also gives rise to SDP reformulations \cite{Am-Set-Gauss-1}.    

2) Divergence ambiguity sets. All PDFs are close to a reference distribution in term of a specified measure. For example, the Wasserstein ambiguity quantifies the divergence via Wasserstein metric \cite{Am-Set-Wasserstein-1,Am-Set-Wasserstein-2,Am-Set-Wasserstein-3}; the $\phi$-divergence ambiguity \cite{Am-Set-Phi-Divergence,Am-Set-Phi-Div-1} characterizes the divergence of two probability density functions through the distance of special non-negative weights (for discrete distributions) or integrals (for continuous distributions). 

More information on the types of ambiguity sets and reformulations of their distributionally robust counterparts can be found in \cite{Am-Set-Overview}. According to the latest research progress, the moment based   distributionally robust optimization is relatively mature and has been widely adopted in engineering, because the semi-infinite LP formulation and its dual for the worst-case expectation problem offer a systematic approach to analyze the impact of uncertain distributions. However, when more complicated ambiguity sets are involved, such as the Gauss ambiguity set, deriving a tractable reformulation needs more sophisticated approaches. The study on the latter category, which directly imposes uncertainty  on the distributions is attracting growing attentions in the past two or three years, because it makes full use of historical data, which can better capture the unique feature of uncertain factors under investigation.  

Data-driven robust stochastic programming, conceptually the same as distributionally robust optimization but preferred by some researchers, has been studied using $\phi$-divergence in \cite{App03-Sect4-DRSO-Phi-1,App03-Sect4-DRSO-Phi-2}, and Wasserstein metric in \cite{Am-Set-Wasserstein-1,Am-Set-Wasserstein-2,Am-Set-Wasserstein-3,App03-Sect4-DRSO-Was-1,App03-Sect4-DRSO-Was-2,App03-Sect4-DRSO-Was-3,App03-Sect4-DRSO-Was-4}, because a tractable counterpart problem can be derived under such ambiguity sets. 

Many decision-making problems in engineering and finance often require that a certain risk measure associated with random variables should be limited below a threshold. However, the probability distribution of random variables is not exactly known; therefore, the risk limiting constraint must be able to withstand perturbations of distribution in a reasonable range. This entails a tractable reformulation of a risk measure under distributional uncertainty. This problem has been comprehensively discussed in \cite{Am-Set-Wasserstein-3}. In more recent publications, CVaR under moment ambiguity set with unimodality is studied in \cite{App03-Sect4-DR-Risk-1}; VaR and CVaR under moment ambiguity set are discussed in  \cite{App03-Sect4-DR-Risk-2}; distortion risk measure under Wasserstein ambiguity set is considered in \cite{App03-Sect4-DR-Risk-3}.

In multi-stage decision making, causality is a pivotal issue for practical implementation, which means that the wait-and-see decisions in the current stage cannot depend on the information of uncertainty in future stages. For example, in a unit commitment problem with 24 periods, the wind power output is observed period-by-period. It is shown in \cite{App03-Sect5-Causal-1} that the two-stage robust model in \cite{ARO-Benders-Decomposition} offers non-causal dispatch strategies, which are in fact not robust. A multi-stage causal unit commitment model is suggested in \cite{App03-Sect5-Causal-1,App03-Sect5-Causal-2} based on affine policy. Causality is put to effect by imposing block diagonal constraints on the gain matrix of affine policy. Causality is also called non-anticipativity in some literature, such as \cite{App03-Sect5-Causal-3}, which is attracting attention from practitioners \cite{App03-Sect5-Causal-4,App03-Sect5-Causal-5}.  

For some other interesting topics on robust optimization, such as the connection with stochastic optimization, connection with risk theory, and applications in engineering problems other than those in power systems, readers can refer to \cite{RO-Survey}. Nonlinear issues have been addressed in \cite{SRO-CVX-RCs,App03-Sect5-RNLP-1}. Optimization models with uncertain SOC and SDP constraints are discussed in \cite{App03-Sect5-RSDP-1,App03-Sect5-RSDP-2}. The connection among robust optimization, data utilization, and machine learning has been reviewed in \cite{App03-Sect5-Opt-Data-ML}.

%
%
%

\motto{Life is not a game. Still, in this life, we choose the games we live to play.}

\chapter{Equilibrium Problems}
\label{App-D} 

The concept of an equilibrium describes a state that the system has no incentive to change. These incentives can be profit-driven in the case of competitive markets or a reflection of physical laws such as energy flow equations. In this sense, equilibrium is encompasses broader concepts than the solution of a game.  Equilibrium is a fundamental notation appearing in various disciplines in economics and engineering. Identifying the equilibria allows eligible authorities to predict the system state at a future time or design reasonable policies for regulating a system or a market. This is not saying that an equilibrium state must appear sooner or later, partly because decision makers in reality have only limited rationality and information. Nevertheless, the awareness of such an equilibrium could be helpful for system design and operation. In this chapter, we restrict our attention in the field of game theory, which entails simultaneously solving multiple interactive optimization problems. We review the notions of some quintessential equilibrium problems and show how they can be solved via traditional optimization methods. These problems can be roughly categorized into two classes: the first one contains only one level: all players must make a decision simultaneously, which is referred to as a Nash-type game; the second one has two levels: decisions are made sequentially by two groups of players, called the leaders and the followers. This category is widely known as Stackelberg-type games, or multi-leader-follower games, or equilibrium programs with equilibrium constraints (EPEC). Unlike a traditional mathematical programming problem where the decision maker is unique, in an equilibrium problem or a game, multiple decision makers seek optimums of individual optimization problems parameterized in the optimal solutions of others.

General notations used throughout this chapter are defined as follows. Specific symbols are explained in the individual sections. In the game theoretic language, a decision maker is called a player.  Vector $x=(x_1,\cdots,x_n)$ refers to the joint decisions of all upper-level players or the so-called leaders in a bilevel setting, where $x_i$ stands for the decisions of leader $i$; $x_{-i}=(x_1,\cdots,x_{i-1},x_{i+1},\cdots,x_n)$ refers to the rivals' actions for leader $i$. Similarly, $y=(y_1,\cdots,y_m)$ refers to the joint decisions of all lower-level players or the so-called followers, where $y_j$ stands for the decisions of follower $j$; $y_{-j}=(y_1,\cdots,y_{j-1},y_{j+1},\cdots,y_m)$ refers to the rivals' actions for follower $j$. $\lambda$ and $\mu$ are Lagrangian dual multipliers associated with inequality and equality constraints.

\section{Standard Nash Equilibrium Problem}
\label{App-D-Sect01}

After J. F. Nash published his work on the equilibrium of $n$-person non-cooperative games in early 1950s \cite{App-04-Nash-1,App-04-Nash-2}, game theory quickly became a new branch of operational research. Nash equilibrium problem (NEP) captures the interactive behaviors of strategic players, in which each player's utility depends on the actions of other players. During decades of wonderful research, a variety of new concepts and algorithms of Nash equilibriums have been proposed and applied to almost every area of knowledge. This section just reviews some basic concepts and the most prevalent best-response algorithms.

\subsection{Formulation and Optimality Condition}
\label{App-D-Sect01-01}

In a standard $n$-person non-cooperative game, each player minimizes his payoff function $f_i(x_i,x_{-i})$ which depends on all players' actions.
The strategy set $X_i=\{x_i \in \mathbb R^{k_i} ~|~ g_i(x_i) \le 0\}$ of player $i$ is independent of $x_{-i}$. The joint strategy set of the game is the Cartesian product  of $X_i$, i.e., $X = \prod_{i=1}^n X_i$, and $X_{-i} = \prod_{j \ne i} X_j$. Roughly speaking, the non-cooperative game is a collection of coupled optimization problems, where player $i$ chooses $x_i \in X_i$ that minimizes his payoff $f_i(x_i,x_{-i})$ given his rivals' strategies $x_{-i}$, or mathematically
\begin{equation}
\label{eq:App-04-NE-Problem}
\left. 
\begin{aligned}
\min_{x_i} ~~ &  f_i(x_i,x_{-i}) \\
\mbox{s.t.}~~ &  g_i(x_i) \le 0 : \lambda_i  
\end{aligned}
\right\},~ i = 1,\cdots,n
\end{equation}
In the problem of player $i$, the decision variable is $x_i$, and $x_{-i}$ is regarded as parameters; $\lambda_i$ is the dual variable. 

The Nash equilibrium consists of a strategy profile such that every player's strategy constitutes the best response to all other players' strategies, or in other words, no player can further reduce his payoff by changing his action unilaterally. Therefore, the Nash equilibrium is a stable state which can sustain spontaneously. The mathematical definition is formally given below.

\begin{definition}
\label{df:App-04-NE}
A strategy vector $x^* \in X$ is a Nash equilibrium if the condition 
\begin{equation}
\label{eq:App-04-NE-Definition}
f_i(x^*_i,x^*_{-i}) \le f_i(x_i,x^*_{-i}),~ \forall x_i \in X_i
\end{equation}
holds for all players.
\end{definition}

Condition (\ref{eq:App-04-NE-Definition}) naturally interprets the fact that at a Nash equilibrium, if any player choose an alternative strategy, his payoff may grow, which is undesired. To depict a Nash equilibrium, a usual approach is the fixed-point of best-response mapping. Let $B_i(x_{-i})$ be the set of optimal strategies of player $i$ given the strategies $x_{-i}$ of others, then set $B(x) = \prod_{i=1}^n B_i(x_{-i})$ is the best-response mapping of the game. It is clear that $x^*$ is a Nash equilibrium if and only if $x^* \in B(x^*)$, i.e., $x^*$ is a fixed point of $B(x)$. This fact establishes the foundation for analyzing Nash equilibria using the well-developed fixed-point theory. However, conducting the fixed-point analysis usually requires the best-response mapping $B(x)$ in a closed form. Moreover, to declare the existence and uniqueness of a Nash equilibrium, the mapping should be contractive \cite{App-04-Fixed-Point-1}. These strong assumptions inevitably limit the applicability of fixed-point method. For example, in many instances, the best-response mapping $B(x)$ is neither contractive nor continuous, but Nash equilibria may still exist.   

Another way to characterize the Nash equilibrium is the KKT system approach. Generally speaking, in a standard Nash game, each player is facing an NLP parameterized in the rivals' strategies. If we consolidate the KKT optimality conditions of all these NLPs in (\ref{eq:App-04-NE-Problem}), we get the following KKT system 
\begin{equation}
\label{eq:App-04-NE-KKT}
\left.  \begin{gathered}
\nabla_{x_i} f_i (x_i ,x_{-i}) + \lambda_i^T \nabla_{x_i} g_i (x_i) = 0 \\
\lambda_i \ge 0,~g(x_i) \le 0,~ \lambda_i^T g_i (x_i) = 0 
\end{gathered}  \right\}~i = 1,\cdots,n  
\end{equation}
If $x^*$ is a Nash equilibrium that satisfies (\ref{eq:App-04-NE-Definition}), and any standard constraint qualification holds for every player's problem in (\ref{eq:App-04-NE-Problem}), then $x^*$ must be a stationary point of the concentrated KKT system (\ref{eq:App-04-NE-KKT}) \cite{App-04-GNEP-KKT-1}; and vice versa: if all problems in (\ref{eq:App-04-NE-Problem}) meet a standard constraint qualification, and a point $x^*$ together with a proper vector of dual multipliers $\lambda = (\lambda_1,\cdots,\lambda_n)$ solves KKT system (\ref{eq:App-04-NE-KKT}), then $x^*$ is also a Nash equilibrium that satisfies (\ref{eq:App-04-NE-Definition}).

Problem (\ref{eq:App-04-NE-KKT}) is an NCP and is the optimality condition of Nash equilibrium. It is a natural attempt to retrieve an equilibrium by solving NCP (\ref{eq:App-04-NE-KKT}) without deploying an iterative algorithm, which may suffer from divergence. To obviate the computational challenges brought by the complementarity and slackness constraints in KKT system (\ref{eq:App-04-NE-KKT}), a merit function approach and an interior-point method are comprehensively discussed in \cite{App-04-GNEP-KKT-1}.

\subsection{Variational Inequality Formulation}
\label{App-D-Sect01-02}

An alternative perspective to study the NEP is to formulate it as a variational inequality (VI) problem. This approach is pursued in \cite{App-04-GNEP-VI-1}. The advantage of variational inequality approach is that it permits an easy access to existence and uniqueness results without the best-response mapping. From a computational point of view, it naturally leads to easily implementable algorithms along with provable convergence performances.

Given a closed and convex set $X \in \mathbb R^n$ and a mapping $F:X \to \mathbb R^n$, a variational inequality problem, denoted by VI($X,F$), is to determine a point $x^* \in X$ satisfying \cite{App-04-VI-1}
\begin{equation}
\label{eq:App-04-VI}
(x-x^*)^T F(x^*)  \ge 0,~  \forall x \in X
\end{equation}

To see the connection between a VI problem and a traditional convex optimization problem that seeks a minimum of a convex function $f(x)$ over a convex set $X$, let us assume that the optimal solution is $x^*$, then the feasible region must not lie in the half space where $f(x)$ decreases; geometrically, the line segment connecting any $x \in X$ with $x^*$  must form an acute angle with the gradient of $f$ at $x^*$, which can be mathematically described as $(x-x^*)^T \nabla f(x^*) \ge 0,~  \forall x \in X$. This condition can be concisely expressed by VI($X,\nabla f$) \cite{App-04-VI-2}.

However,  when the Jacobian matrix of $F$ is not symmetric, $F$ cannot be written as the gradient of another scalar function, and hence the variational inequality problem encompasses broader classes of problems than traditional mathematical programs. For example, when $X = \mathbb R^n$, problem (\ref{eq:App-04-VI}) degenerates into a system of equations $F(x^*)=0$; when $X = \mathbb R^n_+$, problem (\ref{eq:App-04-VI}) comes down to an NCP $0 \le x^* \bot F(x^*) \ge0$. 

To see the later case, $x^* \ge 0$ because it belongs to $X$; if any element of $F(x^*)$ is negative, say the first element $[F(x^*)]_1 < 0$, we let $x_1 = x^*_1 +1$, and $x_i=x^*_i$, $i=2,\cdots$, then $(x-x^*)^T F(x^*) = [F(x^*)]_1 < 0$, which is contradictive to (\ref{eq:App-04-VI}). Hence $ F(x^*) \ge 0$ must hold. Let $x=0$ in (\ref{eq:App-04-VI}), we have $(x^*)^T F(x^*) \le 0$. Because $x^* \ge 0$ and $F(x^*) \ge 0$, there must be $(x^*)^T F(x^*) = 0$, resulting in the target NCP.

The monotonicity of $F$ plays a central role in the theoretical analysis of VI problems, just like the role of convexity in mathematical programming. It has a close relationship with the Jacobian matrix $\nabla F$ \cite{App-04-GNEP-VI-1,App-04-VI-3}:
\begin{svgraybox}
\normalsize
\renewcommand{\arraystretch}{1.5}
\renewcommand{\tabcolsep}{1em}
\centering
\begin{tabular}{lll}
$F(x)$ is monotone on $X$     & $\Leftrightarrow$  &  $\nabla F(x) \succeq 0$, $\forall x \in X$  \\ 
$F(x)$ is strictly monotone on $X$ & $\Leftarrow$  &  $\nabla F(x) \succ 0$, $\forall x \in X$  \\ 
$F(x)$ is strongly monotone on $X$ & $\Leftrightarrow$  &  $\nabla F(x) -c_{m}I \succeq 0$, $\forall x \in X$  \\ 
\end{tabular}
\end{svgraybox}
\noindent where $c_m$ is a strictly positive constant. As a correspondence to convexity, a differentiable function $f$ is convex (strictly convex, strongly convex) on $X$ if and only if $\nabla f$ is monotone (strictly monotone, strongly monotone) on $X$.

Conceptually, monotonicity (convexity) is the weakest, since the matrix $\nabla F(x)$ can have zero eigenvalues; strict monotonicity (strict convexity) is stronger, as all eigenvalues of matrix $\nabla F(x)$ are strictly positive; strong monotonicity (strong convexity) is the strongest, because the smallest eigenvalue of matrix $\nabla F(x)$ should be greater than a given positive number. Intuitively, a strong convex function must be more convex than a given convex quadratic function; for example, $f(x)=x^2$ is strongly convex on $\mathbb R$; $f(x)=1/x$ is convex on $\mathbb R_{++}$ and strongly convex on $(0,1]$.

To formulate an NEP as a VI problem and establish the existence and uniqueness result of Nash equilibria, we list some assumptions on the convexity and smoothness of each player's problem.
\begin{assumption}
\label{ap:App-04-NE-Convex-Smooth} \cite{App-04-GNEP-VI-1}

1) The strategy set $X_i$ is non-empty, closed and convex;

2) Function $f_i(x_i,x_{-i})$ is convex in $x_i \in X_i$ for fixed $x_{-i} \in X_{-i}$;

3) Function $f_i(x_i,x_{-i})$ is continuously differentiable in $x_i \in X_i$ for fixed $x_{-i} \in X_{-i}$;

4) Function $f_i(x_i,x_{-i})$ is twice continuously differentiable in $x \in X$ with bounded second derivatives.
\end{assumption}

\begin{proposition}
\label{pr:App-04-NE-VI} \cite{App-04-GNEP-VI-1}
In a standard NEP NE($X,f$), where $f = (f_1,\cdots,f_n)$, if conditions 1)-3) in Assumption \ref{ap:App-04-NE-Convex-Smooth} are met, then the game is equivalent to a variational inequality problem VI($X,F$) with
\begin{equation*}
X = X_1 \times \cdots \times X_n
\end{equation*}
and 
\begin{equation*}
F(x) = (\nabla_{x_1} f_1 (x),\cdots,\nabla_{x_n} f_n (x))
\end{equation*}
\end{proposition}

In the VI problem corresponding to a traditional mathematical program, the   Jacobian matrix $\nabla F$ is symmetric, because it is the Hessian matrix of a scalar function. However, in Proposition \ref{pr:App-04-NE-VI}, the Jacobian matrix $\nabla F$ for an NEP is generally non-symmetric.  Building upon the VI reformulation, the standard results on solution properties of VI problems \cite{App-04-VI-1} can be extended to standard NEPs. 

\begin{proposition}
\label{pr:App-04-NE-Property}
Given an NEP NE($X,f$), all conditions in Assumption \ref{ap:App-04-NE-Convex-Smooth} are met, then we have the following statements:

1) If $F(x)$ is strictly monotone, then the game has at most one Nash equilibrium.

2) If $F(x)$ is strongly monotone, then the game has a unique Nash equilibrium.

\end{proposition}

Some sufficient guarantees for $F(x)$ to be (strictly, strongly)
monotone are given in \cite{App-04-GNEP-VI-1}. It should be pointed out that the equilibrium concept in the sense of Proposition \ref{pr:App-04-NE-Property} is termed the pure-strategy Nash equilibrium, so as to distinguish it from the mixed-strategy Nash equilibrium which will appear later on. 

\subsection{Best Response Algorithms}
\label{App-D-Sect01-03}

A major benefit of the VI reformulation is that it leads to easily implementable solution algorithms. Here we list two of them. Readers who are interested in the proofs on their performances can consult \cite{App-04-GNEP-VI-1}.

\vspace{12pt}
{\noindent \bf 1. Algorithms for strongly convex cases}

The first algorithm is a totally asynchronous-iterative one, in which players may update their strategies with different frequencies. Let $T=\{0,1,2,\cdots\}$ be the indices of iteration steps, and $T_i \subseteq T$ be the set of steps in which player $i$ updates his own strategy $x_i$. The notation $x^k_i$ implies that at step $k \notin T_i$, $x^k_i$ remains unchanged. Let $t^i_j(k)$ be the latest step at which the strategy of player $j$ is received by player $i$ at step $k$. Therefore, if player $i$ updates his strategy at step $k$, he uses the following strategy profile offered by other players:
\begin{equation}
\label{eq:App-04-NE-Strategy-Update}
x^{t^i(k)}_{-i} = \left(x^{t^i_1(k)}_1,\cdots,x^{t^i_{i-1}(k)}_{i-1},
x^{t^i_{i+1}(k)}_{i+1},\cdots,x^{t^i_n(k)}_{n}  \right)
\end{equation}

Using above definitions, the totally asynchronous-iterative algorithm is summarized in Algorithm \ref{Ag:App-04-NE-Asy-BR}. Some technique conditions for which the schedules $T_i$ and $t^i_j(k)$ should satisfy in order to be implementable in practice are discussed in \cite{App-04-Update-Sequence-1,App-04-Update-Sequence-2}, which are assumed to be satisfied without particular mention. 

\begin{algorithm}[!htp]
\normalsize
\caption{\bf : Asynchronous best-response algorithm}
\begin{algorithmic}[1]
\STATE Choose a convergence tolerance $\varepsilon>0$ and a feasible initial point $x^0 \in X$; the iteration index is $k=0$;
\STATE For player $i=1,\cdots,n$, update the strategy $x^{k+1}_i$ as
\begin{equation}
\label{eq:App-04-NE-Asy-BR}
x^{k+1}_i = \begin{cases}
x^*_i \in \arg \min_{x_i} \left\{ f_i\left( x_i,x^{t^i(k)}_{-i} \right) ~\middle|~ x_i \in X_i \right\}  &  \mbox{if } k \in T_i  \\
x^n_i  &  \mbox{otherwise}
\end{cases}
\end{equation}

\STATE If $\| x^{k+1} - x^k \|_2 \le \varepsilon$, terminate and report $x^{k+1}$ as the Nash equilibrium; otherwise, update $k \leftarrow k+1$, and go to step 2.
\end{algorithmic}
\label{Ag:App-04-NE-Asy-BR}
\end{algorithm} 

A sufficient condition which guarantees the convergence of Algorithm \ref{Ag:App-04-NE-Asy-BR}  is provided in \cite{App-04-GNEP-VI-1}. Roughly speaking, Algorithm \ref{Ag:App-04-NE-Asy-BR} would converge if $f_i(x)$ is strongly convex in $x_i$. However, this is a strong assumption, which cannot be satisfied even if there is only one point where the partial Hessian matrix $\nabla^2_{x_i} f_i(x)$ of player $i$ is singular. 

Algorithm \ref{Ag:App-04-NE-Asy-BR} reduces to some classic algorithms by enforcing a special updating procedure, i.e., a particular selection of $T_i$ and $t^i_j(k)$.
For example,  if players update their strategies simultaneously (sequentially), Algorithm \ref{Ag:App-04-NE-Asy-BR} becomes the Jacobi (Gauss-Seidel) type iterative scheme. Interestingly, the asynchronous best-response algorithm is robust against data missing or delay, and is guaranteed to find the unique Nash equilibrium. This feature greatly relaxes the requirement on data synchronization and simplifies the design of communication systems, and makes this class of algorithm very appealing in distributed system operations.

\vspace{12pt}
{\noindent \bf 2. Algorithms for convex cases}

To relax the strong monotonicity assumption on $F(x)$, the second algorithm has been proposed in \cite{App-04-GNEP-VI-1}, which only uses the monotonicity property and is summarized below. Algorithm \ref{Ag:App-04-NE-Proximal} converges to a Nash equilibrium, if each player's optimization problem is convex (or $F(x)$ is monotone), which significantly improves its applicability. 
\begin{algorithm}[!htp]
\normalsize
\caption{\bf : Proximal Decomposition Algorithm}
\begin{algorithmic}[1]
\STATE Given $\{ \rho_n\}_{n=0}^\infty$, $\varepsilon >0$, and $\tau > 0$, choose a feasible initial point $x^0 \in X$; 
\STATE Find an equilibrium $z^0$ of the following NEP using Algorithm \ref{Ag:App-04-NE-Asy-BR}
\begin{equation}
\label{eq:App-04-NE-Proximal}
\left. \begin{aligned}
\min_{x_i} ~~ & f_i( x_i,x_{-i}) + \tau \| x_i - x^0_i  \|^2_2  \\
\mbox{s.t.}~~ & x_i \in X_i 
\end{aligned}  \right\}, i=1,\cdots,n
\end{equation}

\STATE If $\| z^0 - x^0 \|_2 \le \varepsilon$, terminate and report $x^0$ as the Nash equilibrium; otherwise, update $x^0 \leftarrow (1-\rho_n) x^0 + \rho_n z^0$, and go to step 2.
\end{algorithmic}
\label{Ag:App-04-NE-Proximal}
\end{algorithm} 

Algorithm \ref{Ag:App-04-NE-Proximal} is a double-loop method: the inner loop identifies a Nash equilibrium of the regularized game (\ref{eq:App-04-NE-Proximal}) with $x^0$ being a parameter, which is updated in each iteration, and the outer loop updates $x^0$ by selecting a new point along the line connecting $x^0$ and $z^0$. Notice that in step 2, as long as $\tau$ is large enough, the Hessian matrix $\nabla^2_{x_i} f_i (x) + 2 \tau I$ must be positive definite, and thus the best-response algorithm applied to (\ref{eq:App-04-NE-Proximal}) is guaranteed to converge to the unique Nash equilibrium. See \cite{App-04-GNEP-VI-1} for more details about parameter selection.

The penalty term $\tau \| x_i - x^0_i \|^2$ limits the change of optimal strategies in two consecutive iterations, and can be interpreted as a damping factor that attenuates possible oscillations during the computation. It is worth mentioning that the penalty parameter $\tau$ significantly impacts the convergence rate of Algorithm \ref{Ag:App-04-NE-Proximal} and should be carefully selected. If it is too small, the damping effect of the penalty term is limited, and the oscillation may still take place; if it is too large, the increment of $x$ in each step is very small, and Algorithm \ref{Ag:App-04-NE-Proximal} may suffer from a slow convergence rate. The optimal value of $\tau$ is problem-dependent. There is not a universal way to determine its best value.    

Recently, single-loop distributed algorithms for monotone Nash games are proposed in \cite{App-04-GNEP-ITR}, which authors believe to be promising in practical applications. In these two schemes, the regularization parameter is updated at once after each iteration is completed, rather than when the regularized problem is approximately solved, and players can select their parameter independently.

\subsection{Nash Equilibrium of Matrix Games}
\label{App-D-Sect01-04}

As explained before, not all Nash games have an equilibrium, especially when the strategy set and the payoff function are non-convex or discrete. To widen the equilibrium notion and reveal deeper insights on the behaviors of players in such instances, it is instructive to revisit some simple games, called the matrix game, which is the primary research object of game theorists. 

The bimatrix game refers to a matrix game involving two players P1 and P2. The numbers of possible strategies of P1 and P2 are $m$ and $n$, respectively. $A = \{a_{ij}\} \in \mathbb M^{m \times n}$ is the payoff matrix of P1:  when P1 chooses strategy $i$ and P2 selects strategy $j$, the payoff of P1 is $a_{ij}$. The payoff matrix $B \in \mathbb M^{m \times n}$ of P2 can be defined in the same way. In a matrix game, each player is interested to determine a probability distribution of his actions, such that his expected payoff is minimized. Let $x_i$ ($y_j$) be the probability that P1 (P2) will use strategy $i$ ($j$), vectors $x=[x_1,\cdots,x_m]^T$ and $y = [y_1,\cdots,y_n]^T$ are called mixed strategies, clearly,  
\begin{equation}
\label{eq:App-04-MG-Mixed-Strategy}
\begin{gathered}
x \ge 0,~ \sum_{i=1}^m x_i = 1 \quad \mbox{or} \quad x \in {\rm \Delta}_m\\
y \ge 0,~ \sum_{j=1}^n y_j = 1 \quad \mbox{or} \quad y \in {\rm \Delta}_n\\
\end{gathered}
\end{equation}
where ${\rm \Delta}_m$ and ${\rm \Delta}_n$ are simplex slices in $\mathbb R^m$ and $\mathbb R^n$.

\vspace{12pt}
{\noindent \bf 1. Two-person zero-sum games}

The zero-sum game represents a totally competitive situation: P1's gain is P2's loss, so the sum of their payoff matrices is $A+B = 0$, as its name suggests. Such type of game has been well studied in vast literature since von Neumann found the famous Minimax Theorem in 1928. The game is revisited from a mathematical programming perspective in \cite{App-04-Minimax-LP}. The proposed linear programming method is especially powerful for instances with a high-dimensional payoff matrix. Next, we briefly introduce this method. 

Let us begin with a payoff matrix $A = \{a_{ij}\}$, $a_{ij} > 0$, $\forall i,j$ with strictly positive entries (otherwise, we can add a constant to every entry, such that the smallest entry becomes positive, and the equilibrium strategy remains the same). The expected payoff of P1 is given by  
\begin{equation}
\label{eq:App-04-ZSG-Payoff}
V_A = \sum_{i=1}^m \sum_{j=1}^n  x_i a_{ij} y_j  = x^T A y
\end{equation}
which must be positive because of the element-wise positivity assumption on $A$. 

Since $B=-A$ and minimizing $-x^T A y$ is equivalent to maximizing $x^T A y$, the two-person zero-sum game has a min-max form as
\begin{equation}
\label{eq:App-04-ZSG-1}
\min_{x \in {\rm \Delta}_m} \max_{y \in {\rm \Delta}_n} ~~ x^T A y 
\end{equation}
or 
\begin{equation}
\label{eq:App-04-ZSG-2}
\max_{y \in {\rm \Delta}_n} \min_{x \in {\rm \Delta}_m} ~~ x^T A y 
\end{equation}

The solution to the two-person zero-sum matrix game (\ref{eq:App-04-ZSG-1}) or (\ref{eq:App-04-ZSG-2}) is called a mixed-strategy Nash equilibrium, or the saddle point of a min-max problem. It satisfies 
\begin{gather}
(x^*)^T A y^* \le x^T A y^* ,~ \forall x \in {\rm \Delta}_m \notag \\
(x^*)^T A y^* \ge y^T A x^* ,~ \forall y \in {\rm \Delta}_n \notag
\end{gather}

To solve this game, consider (\ref{eq:App-04-ZSG-1}) in the following format
\begin{equation}
\label{eq:App-04-ZSG-3}
\min_x \left\{ v_1(x) ~\middle|~ x \in {\rm \Delta}_m \right\}  
\end{equation}
where $v_1(x)$ is the optimal value function of the problem faced by P2 with the fixed strategy $x$ of P1
\begin{equation}
v_1(x) = \max_y \left\{ x^T A y ~\middle|~ y \in {\rm \Delta}_n \right\}\notag
\end{equation}

In view of the feasible region defined in (\ref{eq:App-04-MG-Mixed-Strategy}), $v_1(x)$ is equal to the maximal element of vector $x^T A$, which is strictly positive, and the inequality
\begin{equation}
A^T x \le {\bf 1}^n v_1(x)  \notag
\end{equation}
holds. Furthermore, introducing a normalized vector $\bar x = x/v_1(x)$, we have  
\begin{equation}
\begin{gathered}
\bar x  \ge 0, ~~ A^T \bar x \le {\bf 1}^n \\
v_1(x) = ({\bar x}^T {\bf 1}^m)^{-1}   
\end{gathered}  \notag
\end{equation}
Taking these relations into account, problem (\ref{eq:App-04-ZSG-3})
becomes
\begin{equation}
\label{eq:App-04-ZSG-4}
\begin{aligned}
\min_{\bar x} ~~ & ({\bar x}^T {\bf 1}^m)^{-1}  \\
\mbox{s.t.}~~    &  A^T \bar x \le {\bf 1}^n       \\
                 &  \bar x  \ge 0
\end{aligned} 
\end{equation}

Because the objective is strictly positive and monotonic, the optimal solution of (\ref{eq:App-04-ZSG-4}) keeps unchanged if we choose to maximize ${\bar x}^T {\bf 1}^m$ under the same constraints, giving rise to the following LP 
\begin{equation}
\label{eq:App-04-ZSG-5}
\begin{aligned}
\max_{\bar x} ~~ & {\bar x}^T {\bf 1}^m  \\
\mbox{s.t.}~~    &  A^T \bar x \le {\bf 1}^n       \\
                 &  \bar x  \ge 0
\end{aligned} 
\end{equation}
Let $\bar x^*$ and $\bar v^*_1$ be the optimal solution and optimal value of LP (\ref{eq:App-04-ZSG-5}). According to the analysis of variable transformation, the optimal expected payoff $v^*_1$ and the optimal mixed strategy $x^*$ of P1 in this game are given by 
\begin{equation}
\label{eq:App-04-ZSG-6}
v^*_1 = 1 / \bar v^*_1,~~ x^* = \bar x^* / \bar v^*_1
\end{equation}

Consider (\ref{eq:App-04-ZSG-2}) in the same way, we obtain the following LP for P2:
\begin{equation}
\label{eq:App-04-ZSG-7}
\begin{aligned}
\min_{\bar y} ~~ & {\bar y}^T {\bf 1}^n  \\
\mbox{s.t.}~~    &  A \bar y \ge {\bf 1}^m    \\
                 &  \bar y  \ge 0
\end{aligned} 
\end{equation}
Denote by $\bar y^*$ and $\bar v^*_2$ the optimal solution and optimal value of LP (\ref{eq:App-04-ZSG-7}), and then the optimal expected payoff $v^*_2$ and the optimal mixed strategy $y^*$ of P2 in this game can be posed as
\begin{equation}
\label{eq:App-04-ZSG-8}
v^*_2 = 1 / \bar v^*_2,~~ y^* = \bar y^* / \bar v^*_2
\end{equation}

In summary, the mixed-strategy Nash equilibrium of two-person zero-sum matrix game (\ref{eq:App-04-ZSG-1}) is $(x^*,y^*)$, and the payoff of P1 is $v^*_1$. Interestingly, we notice that problems (\ref{eq:App-04-ZSG-5}) and (\ref{eq:App-04-ZSG-7}) constitute a pair of dual LPs, implying that their optimal values are equal, and the optimal solution $y^*$ in (\ref{eq:App-04-ZSG-8}) also solves the inner LP of (\ref{eq:App-04-ZSG-1}). This observation
leads to two important conclusions: 

1) The Nash equilibrium of a two-person zero-sum matrix game, or the saddle point,  can be computed by solving a pair of dual LPs. In fact, if one player's strategy, say $x^*$, have been obtained from (\ref{eq:App-04-ZSG-5}) and (\ref{eq:App-04-ZSG-6}), the rival's strategy can be retrieved by solving (\ref{eq:App-04-ZSG-1}) with $x = x^*$. 

2) The decision sequence of a two-person zero-sum game is interchangeable without influencing the saddle point. 

\vspace{12pt}
{\noindent \bf 2. General bimatrix games}

In more general two-person matrix games, the sum of payoff matrices is not equal to zero, and each player wishes to minimize its own expected payoff taking the other player's strategy as given. In the setting of mixed strategies, players are selecting the probability distribution among available strategies rather than a single action (the pure strategy), and the respective optimization problems are
as follows
\begin{equation}
\label{eq:App-04-BMG-1}
\begin{aligned}
\min_x ~~   &  x^T A y  \\
\mbox{s.t.}~~ &  x^T {\bf 1}^m =1 : \lambda   \\
            &  x \ge 0
\end{aligned}
\end{equation}
\begin{equation}
\label{eq:App-04-BMG-2}
\begin{aligned}
\min_y ~~   &  x^T B y  \\
\mbox{s.t.}~~ &  y^T {\bf 1}^n =1 : \gamma   \\
            &  y \ge 0
\end{aligned}
\end{equation}

The pair of probability distributions $(x^*,y^*)$ is called a mixed-strategy Nash equilibrium if
\begin{gather}
(x^*)^T A y^* \le x^T A y^* ,~ \forall x \in {\rm \Delta}_m \notag \\
(x^*)^T B y^* \le (x^*)^T B y ,~ \forall y \in {\rm \Delta}_n \notag
\end{gather}

Unlike the zero-sum case, there is not an equivalent LP that can extract the Nash equilibrium. Performing the KKT system approach, we write out the KKT condition for (\ref{eq:App-04-BMG-1})
\begin{gather*}
 A y - \lambda {\bf 1}^m - \mu = 0  \\
 0 \le \mu \bot x \ge 0   \\
 x^T {\bf 1}^m = 1 
\end{gather*}
where $\mu$ is the dual variable associated with the non-negative constraint, and can be eliminated from the first equality. Concentrating the KKT conditions of LPs (\ref{eq:App-04-BMG-1}) and (\ref{eq:App-04-BMG-2}) gives 
\begin{equation}
\label{eq:App-04-BMG-3}
\begin{gathered}
0 \le A y   - \lambda {\bf 1}^m \bot x \ge 0,~ x^T {\bf 1}^m = 1  \\
0 \le B^T x - \gamma  {\bf 1}^n \bot y \ge 0,~ y^T {\bf 1}^n = 1  \\
\end{gathered}
\end{equation}

Complementarity condition (\ref{eq:App-04-BMG-3}) can be solved by setting  $\lambda = \gamma = 1$ and omitting equality constraints, and recovering them at a later normalization step, i.e., we first solve 
\begin{equation}
\label{eq:App-04-BMG-4}
\begin{gathered}
0 \le A y   - {\bf 1}^m \bot x \ge 0  \\
0 \le B^T x - {\bf 1}^n \bot y \ge 0  \\
\end{gathered}
\end{equation}
Suppose that the solution is ($\bar x, \bar y$), then the Nash equilibrium is
\begin{equation}
\label{eq:App-04-BMG-5}
\begin{gathered}
x^* = \bar x / \bar x^T {\bf 1}^m  \\
y^* = \bar y / \bar y^T {\bf 1}^n  
\end{gathered}
\end{equation}
and the corresponding multipliers are derived from (\ref{eq:App-04-BMG-3}) as
\begin{equation}
\label{eq:App-04-BMG-6}
\begin{gathered}
\lambda^* = (x^*)^T A y^*  \\
\gamma^*  = (x^*)^T B y^*  \\
\end{gathered}
\end{equation}

On the other hand, if ($x^*, y^*$) is a Nash equilibrium and solves (\ref{eq:App-04-BMG-3}) with multipliers ($\lambda^*,\gamma^*$), we can observe that ($x^*/\gamma^*,y^*/\lambda^*$) solves (\ref{eq:App-04-BMG-4}), therefore
\begin{equation}
\label{eq:App-04-BMG-7}
\begin{gathered}
\bar x = \dfrac{x^*}{(x^*)^T B y^*}  \\
\bar y = \dfrac{y^*}{(x^*)^T A y^*}  
\end{gathered}
\end{equation}

Now we can see that identifying the mixed-strategy Nash equilibrium of a bimatrix game entails solving KKT system (\ref{eq:App-04-BMG-3}) or (\ref{eq:App-04-BMG-4}), which is called a linear complementarity problem (LCP). A classical algorithm for LCP is the Lemke's method \cite{App-04-LCP-Lemke-1,App-04-LCP-Lemke-2}. Another systematic way to solve an LCP is to reformulate it as an MILP using the method described in Appendix \ref{App-B-Sect03-05}. Nonetheless, there are more tailored MILP models for LCPs, which will be detailed in Sect. \ref{App-D-Sect04-02}.

Unlike the pure-strategy Nash equilibrium, whose existence relies on some assumptions on convexity, the mixed-strategy Nash equilibrium for matrix games, which is the discrete probability distribution among available actions, always exists \cite{App-04-Nash-2}. If a game with two players has no pure-strategy Nash equilibrium, and each player can choose actions from a finite strategy set, we can then calculate the payoff matrices as well as the mixed-strategy Nash equilibrium, which informs the likelihood that the player will adopt each corresponding pure strategy.

\subsection{Potential Games}
\label{App-D-Sect01-05}

Despite that a direct certification of the existence and uniqueness of a pure-strategy Nash equilibrium for a general game model is non-trivial, when the game possesses some special structures, such a certification becomes axiomatic. One of these guarantees is the existence of an exact potential function, and the associated problem is known as the potential game \cite{App-04-Potential-Game-1}. Four types of potential games are listed in \cite{App-04-Potential-Game-1}, categorized by the type of the potential function. Other extensions of the potential game have been studied as well. For a complete introduction, we recommend \cite{App-04-Potential-Game-2}.

\begin{definition}
\label{pr:App-04-PG-1} 
(Exact potential game) A game is an exact potential game if there is a potential function $U(x)$ such that:
\begin{equation}
\label{eq:App-04-PG-1}
\begin{gathered}
f_i(x_i,x_{-i}) - f_i(y_i,x_{-i}) = U(x_i,x_{-i}) - U(y_i,x_{-i})  \\
\forall x_i,y_i \in X_i,~ \forall x_{-i} \in X_{-i},~ i = 1,\cdots,n
\end{gathered}
\end{equation}
\end{definition}

In an exact potential game, the change in the utility/payoff of any single player due to the unilateral strategy deviation leads to the same amount of change in the potential function. Among various variations of potential games which are defined by relaxing the strict equality (\ref{eq:App-04-PG-1}), the exact potential game is the most fundamental one and has attracted the majority of research interests. Throughout this section, the term potential game means the exact one without particular mention.  

The condition for a game being a potential game and the method for constructing the potential function are given in the following proposition.

\begin{proposition}
\label{pr:App-04-PG-1}
\cite{App-04-Potential-Game-1} Suppose the payoff functions $f_i$, $i=1,\cdots,n$ in a game are twice continuously differentiable, then a potential function exists if and only if
\begin{equation}
\label{eq:App-04-PG-2}
\dfrac{\partial^2 f_i}{\partial x_i \partial x_j} = 
\dfrac{\partial^2 f_j}{\partial x_i \partial x_j},~~ \forall i,j = 1,\cdots,n
\end{equation}
and the potential function can be constructed as
\begin{equation}
\label{eq:App-04-PG-3}
U(v) - U(z) = \sum_{i=1}^n \int_0^1 (x^\prime_i(t))^T \dfrac{\partial f_i}{\partial x_i} (x(t))  d t
\end{equation}
where $x(t): [0,1] \to X$ is a continuously differentiable path in $X$ connecting strategy profile $v$ and a fixed strategy profile $z$, such that $x(0)=z$, $x(1)=v$.
\end{proposition}

To obtain (\ref{eq:App-04-PG-3}), first, a direct consequence of (\ref{eq:App-04-PG-1}) is 
\begin{equation}
\label{eq:App-04-PG-4}
\dfrac{\partial f_i}{\partial x_i} = \dfrac{\partial U}{\partial x_i},~ i=1,\cdots,n
\end{equation}

For any smooth curve $C(t): [0,1] \to X$ and any function $U$ with a continuous gradient $\nabla U$, the gradient theorem in calculus tells us
\begin{equation}
U(C_{end}) - U(C_{start}) = \int_C \nabla U(s) d s  \notag
\end{equation}
where vector $s$ represents points along the integral trajectory $C$ parameterized in a scalar variable. Introducing $s=x(t)$: when $t=0$, $s=C_{start}=z$; when $t=1$, $s=C_{end}=v$. By the chain rule, $ds = x^\prime(t) dt$, and hence we get
\begin{equation}
\begin{aligned} 
U(v) - U(z) & = \int_0^1 (x^\prime(t))^T \nabla U(x(t)) d t  \\
            & = \sum_{i=1}^n \int_0^1 (x_i^\prime(t))^T 
              \dfrac{\partial U}{\partial x_i} (x(t)) d t 
\end{aligned}  \notag
\end{equation}

Then, if $U$ is a potential function, substituting (\ref{eq:App-04-PG-4}) into above equation gives equation (\ref{eq:App-04-PG-3}). 

In summary, for a standard NEP with continuous payoff functions, we can check whether it is a potential game, and further construct its potential function, if the right-hand side of (\ref{eq:App-04-PG-3}) has a closed form expression. Nevertheless, in some particular cases, the potential function can be observed without calculating an integral.  

1. The payoff functions of the game can be decomposed as  
\begin{equation}
f_i(x_i,x_{-i}) = p_i(x_i) + Q(x),~ \forall i = 1,\cdots,n \notag
\end{equation}
where the first term only depends on $x_i$, and the second term that couples all players' strategies and appears in every utility function is identical. In such circumstance, the potential function is instantly posed as
\begin{equation}
U(x) = Q(x) + \sum_{i=1}^n p_i(x_i)   \notag
\end{equation}
which can be verified through its definition in (\ref{eq:App-04-PG-1}).

2. The payoff functions of the game can be decomposed as  
\begin{equation}
f_i(x_i,x_{-i}) = p_i(x_{-i}) + Q(x),~ \forall i = 1,\cdots,n \notag
\end{equation}
where the first term only depends on the joint actions of opponents $x_{-i}$, and the second term is common and identical to all players. In such circumstance, the potential function is $Q(x)$. This is easy to understand because $x_{-i}$ is constant in the decision-making problem of player $i$ and thus the first term $p_i(x_{-i})$ can be omitted from the objective function.

3. The payoff function of each player has a form of
\begin{equation}
f(x_i,x_{-i}) = \left( a+b \sum_{j=1}^n x_j \right) x_i + c_i(x_i),~ 
i=1,\cdots,n  \notag
\end{equation}
Obviously,
\begin{equation}
\dfrac{\partial^2 f_i}{\partial x_i \partial x_j} = 
\dfrac{\partial^2 f_j}{\partial x_i \partial x_j} = b  \notag
\end{equation}
therefore, a potential function exists and is given by
\begin{equation}
U(x) = a \sum_{i=1}^n x_i + \dfrac{b}{2} \sum_{i=1}^n \sum_{j \ne i} x_ix_j + b \sum_{i=1}^n  x^2_i  + \sum_{i=1}^n c_i(x_i)  \notag
\end{equation}

The potential function provides a convenient way to analyze the Nash equilibria of potential games, since the function coincides with incentives of all players. 

\begin{proposition}
\label{pr:App-04-PG-2}  \cite{App-04-Potential-Game-2}
If game $G_1$ is a potential game with potential function $U(x)$; $G_2$ is another game with the same number of players and their payoff functions are $F_1(x_1,x_{-1}) = \cdots = F_n(x_n,x_{-n}) = U(x)$. Then $G_1$ and $G_2$ have the same set of Nash equilibria.
\end{proposition}

This is easy to understand because an equilibrium of $G_1$ satisfies
\begin{equation}
f_i(x^*_i,x^*_{-i}) \le f_i(x_i,x^*_{-i}), \forall x_i \in X_i,~
i = 1,\cdots,n  \notag
\end{equation}
By the definition of potential function (\ref{eq:App-04-PG-1}), this gives
\begin{equation}
\label{eq:App-04-PG-5}
U(x^*_i,x^*_{-i}) \le U(x_i,x^*_{-i}), \forall x_i \in X_i,~ i = 1,\cdots,n
\end{equation}
So any equilibrium of $G_1$ is an equilibrium of $G_2$. Similarly, the reverse holds, too. 

In Proposition \ref{pr:App-04-PG-2}, the identical interest game $G_2$ is actually an optimization problem. The potential function builds a bridge between an NEP and a mathematical programming problem. Let $X = X_1 \times \cdots \times X_n$, (\ref{eq:App-04-PG-5}) can be written as 
\begin{equation}
\label{eq:App-04-PG-6}
U(x^*) \le U(x), \forall x \in X
\end{equation}
On this account, we have

\begin{proposition}
\label{pr:App-04-PG-3}   \cite{App-04-Potential-Game-1}
Every minimizer of the potential function $U(x)$ in $X$ is a (pure-strategy) Nash equilibrium of the potential game.
\end{proposition}

Proposition \ref{pr:App-04-PG-3} is very useful. It reveals the fact that computing a Nash equilibrium of a potential game is equivalent to solving a traditional mathematical program. Meanwhile, the existence and uniqueness results of Nash equilibrium for potential games can be understood from the solution property of NLPs.   

\begin{proposition}
\label{pr:App-04-PG-4}  \cite{App-04-Potential-Game-2}
Every potential game with a continuous potential function $U(x)$ and a compact strategy space $X$ has at least one (pure-strategy) Nash equilibrium. If $U(x)$ is strictly convex, then the Nash equilibrium is unique.
\end{proposition}

Propositions \ref{pr:App-04-PG-3}-\ref{pr:App-04-PG-4} make no reference on the convexity of individual payoff functions of players. Moreover, if the potential function $U(x)$ is non-convex and has multiple local minimums, then each local optimizer corresponds to a local Nash equilibrium where $X_i$ in (\ref{eq:App-04-PG-5}) is replaced with the intersection of $X_i$ with a neighborhood region of $x^*_i$.

\section{Generalized Nash Equilibrium Problem}
\label{App-D-Sect02}

In above developments for standard NEPs, we have assumed that the strategy sets are decoupled: the available strategies of each player do not depend on other players' choices. However, there are indeed many practical cases where the strategy sets are interactive. For example, when players consume a common resource, the total consumption should not exceed the inventory quantity. The generalized Nash equilibrium problem (GNEP), invented in \cite{App-04-GNEP-1}, relaxes the strategy independence assumption in classic NEPs and allows the feasible set of each player's actions to depend on the rivals' strategies. For a comprehensive review, we recommend \cite{App-04-GNEP-2}.

\subsection{Formulation and Optimality Condition}
\label{App-D-Sect02-01}

Denote by $X_i(x_{-i})$ the strategy set of player $i$ when others select $x_{-i}$. In a GNEP, given the value of $x_{-i}$, each player $i$ determines a strategy $x_i \in X_i(x_{-i})$ which minimizes a payoff function $f_i(x_i,x_{-i})$. In this regard, a GNEP with $n$ players is the joint solution of $n$ coupled optimization problems
\begin{equation}
\label{eq:App-04-GENP-MP}
\left. \begin{aligned}
\min_{x_i} ~~ & f_i(x_i,x_{-i})  \\
\mbox{s.t.}~~ & x_i \in X_i(x_{-i})
\end{aligned}  \right\},~~ i=1,\cdots,n
\end{equation}

In (\ref{eq:App-04-GENP-MP}), correlation
takes place not only in the objective function, but also in the constraints.

\begin{definition}
\label{df:App-04-GENP}
A generalized Nash equilibrium (GNE), or the solution of a
GNEP, is a feasible point $x^*$ such that
\begin{equation}
f(x^*_i,x^*_{-i})  \le f(x_i,x^*_{-i}),~ \forall x_i \in X_i (x_{-i})
\end{equation}
holds for all players.
\end{definition}

In its full generality, the GNEP is much more difficult than an NEP due to the variability of strategy sets. In this section, we restrict our attention to a particular class of GNEP: the so-called GNEP with shared convex constraints. In such a problem, the strategy sets can be expressed as
\begin{equation}
\label{eq:App-04-GENP-Convex-Xi}
X_i(x_{-i}) = \left\{x_i ~\middle|~ x_i \in Q_i,~ 
g(x_i,x_{-i}) \le 0 \right\},~ i = 1,\cdots,n
\end{equation}
where $Q_i$ is a closed and convex set which involves only $x_i$; $g(x_i,x_{-i}) \le 0$ represents the shared constraints. They consist of a set of convex inequalities coupling all players' strategies and are identical in $X_i(x_{-i})$, $i=1,\cdots,n$. Sometimes, $Q_i$ and $g(x_i,x_{-i}) \le 0$ are also mentioned as local and global constraints, respectively. 

\begin{figure}[!htp]
\centering
\includegraphics[scale=0.50]{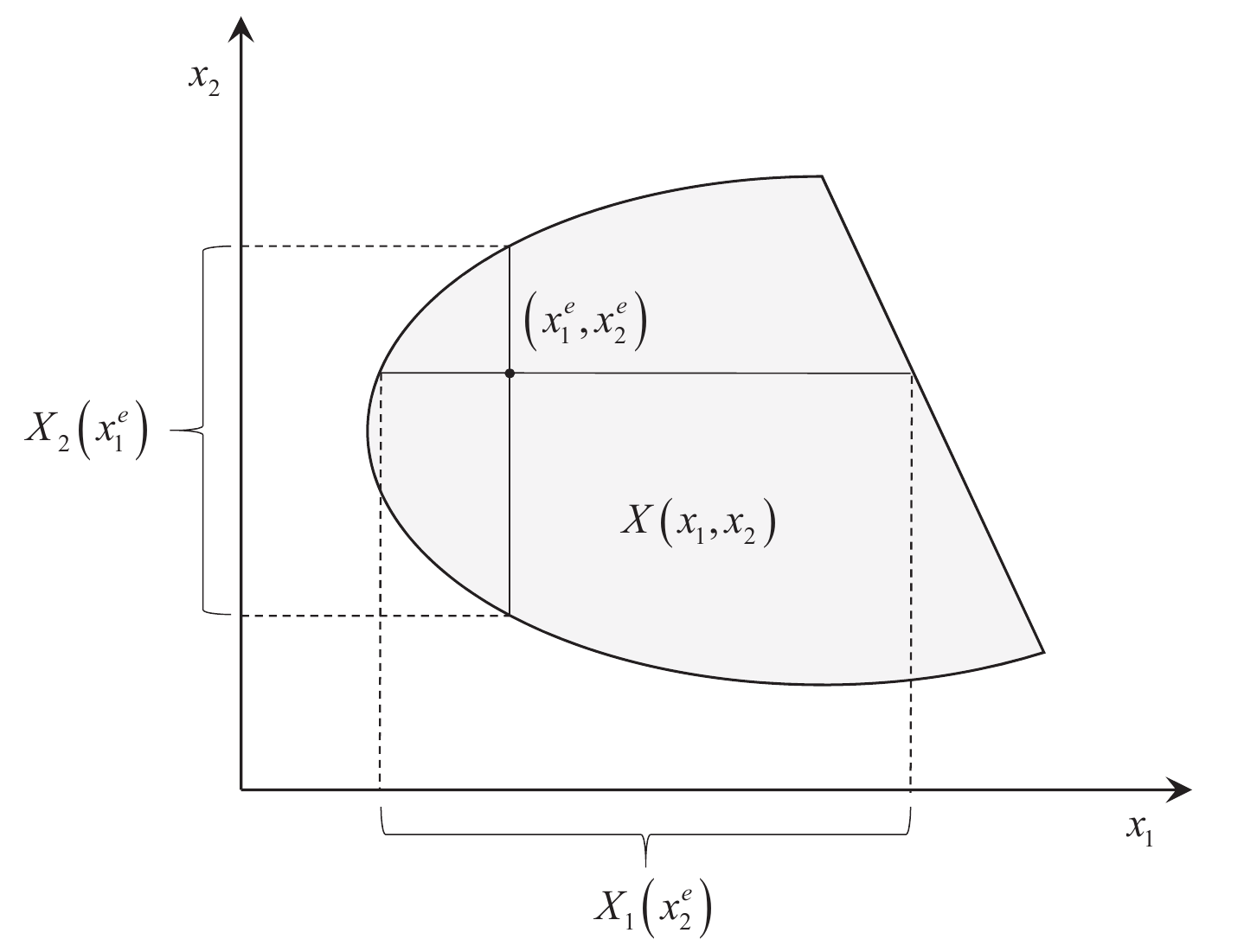}
\caption{Relations of $X$ and the individual strategy sets.}
\label{fig:App-04-01}
\end{figure}

In the absence of shared constraints, the GNEP reduces to a standard NEP. Define the feasible set of strategy profile $x=(x_1,\cdots,x_n)$ in a GNEP
\begin{equation}
\label{eq:App-04-GENP-Convex-X}
X = \left\{x ~\middle|~ x \in \prod_{i=1}^n Q_i,~ g(x) \le 0 \right\}
\end{equation}
It is easy to see that $X_i(x_{-i})$ is a slice of $X$. A geometric interpretation of (\ref{eq:App-04-GENP-Convex-Xi}) is illustrated in Fig. \ref{fig:App-04-01}. It is seen that the choice of $x_1$  influences the feasible interval $X_2(x_1)$ of Player 2. 

We make some assumptions on the smoothness and convexity for a GNEP with shared constraints.

\begin{assumption}
\label{ap:App-04-GNEP-Convex-Smooth} \

1) Strategy set $Q_i$ of each player is nonempty, closed, and convex.

2) Payoff function $f_i(x_i,x_{-i})$ of each player is twice continuously differentiable in $x$ and convex in $x_i$ for every fixed $x_{-i}$.

3) Functions $g(x) = (g_1(x),\cdots,g_m(x))$ are differentiable and convex in $x$.
\end{assumption}

In analogy with the NEP, concatenating the KKT optimality condition of each optimization problem in (\ref{eq:App-04-GENP-MP}) gives us what is called the KKT condition of the GNEP. For notation brevity, we omit local constraints ($Q_i=\mathbb R^{n_i}$) and assume that $X_i(x_{-i})$ contains only global constraints. Write out the KKT condition of GNEP (\ref{eq:App-04-GENP-MP})
\begin{equation}
\label{eq:App-04-GNEP-KKT}
\left.  \begin{gathered}
\nabla_{x_i} f_i (x_i ,x_{-i}) + \lambda_i^T \nabla_{x_i} g (x) = 0 \\
\lambda_i \ge 0,~g(x) \le 0,~ \lambda_i^T g (x) = 0   
\end{gathered}  \right\}~i = 1,\cdots,n  
\end{equation}
where $\lambda_i$ is the Lagrange multiplier vector associated with the global constraints in the $i$-th player's problem. 

\begin{proposition}
\label{pr:App-04-GNEP-KKT}  \cite{App-04-GNEP-2} \

1) Let $\bar x = (\bar x_1,\cdots,\bar x_n)$ be the equilibrium of a GNEP, then a multiplier vector $\bar \lambda = (\bar \lambda_1,\cdots,\bar \lambda_n)$ exists, such that the pair $(\bar x, \bar \lambda)$ solves KKT system (\ref{eq:App-04-GNEP-KKT}).

2) If  $(\bar x, \bar \lambda)$ solves KKT system (\ref{eq:App-04-GNEP-KKT}), and Assumption \ref{ap:App-04-GNEP-Convex-Smooth} holds, then $\bar x$ is an equilibrium of GNEP (\ref{eq:App-04-GENP-MP}) with shared convex constraints.

\end{proposition}

However, in contrast to an NEP, the solutions of an GNEP may be non-isolated and constitute a low dimensional manifold, because $g(x)$ is a common constraint shared by all, and the Jacobian of the KKT system may appear to be singular. A meticulous explanation is provided in \cite{App-04-GNEP-2}. We give a graphic interpretation for this phenomenon. 

Consider a GNEP with two players:
\begin{gather}
\mbox{Player 1:} \quad  
\left\{ \begin{aligned}
\max_{x_1} ~~  &  x_1  \\
\mbox{s.t.}~~  &  x_1 \in X_1(x_2)
\end{aligned}  \right.  \notag \\
\mbox{Player 2:} \quad  
\left\{ \begin{aligned}
\max_{x_2} ~~  &  x_2  \\
\mbox{s.t.}~~  &  x_2 \in X_2(x_1)
\end{aligned}  \right.  \notag
\end{gather}
where 
\begin{gather}
X_1(x_2) = \{x_1 ~|~ x_1 \ge 0,~ g(x_1,x_2) \le 0 \} \notag  \\
X_2(x_1) = \{x_2 ~|~ x_2 \ge 0,~ g(x_1,x_2) \le 0 \} \notag  
\end{gather}
and the global constraint set is
\begin{equation*}
\left\{ x ~\middle|~ \begin{gathered}
g_1 = 2 x_1 + x_2 \le 0  \\   
g_2 = x_1 + 2 x_2 \le 0
\end{gathered}  \right\}    
\end{equation*}

\begin{figure}[!htp]
\centering
\includegraphics[scale=0.40]{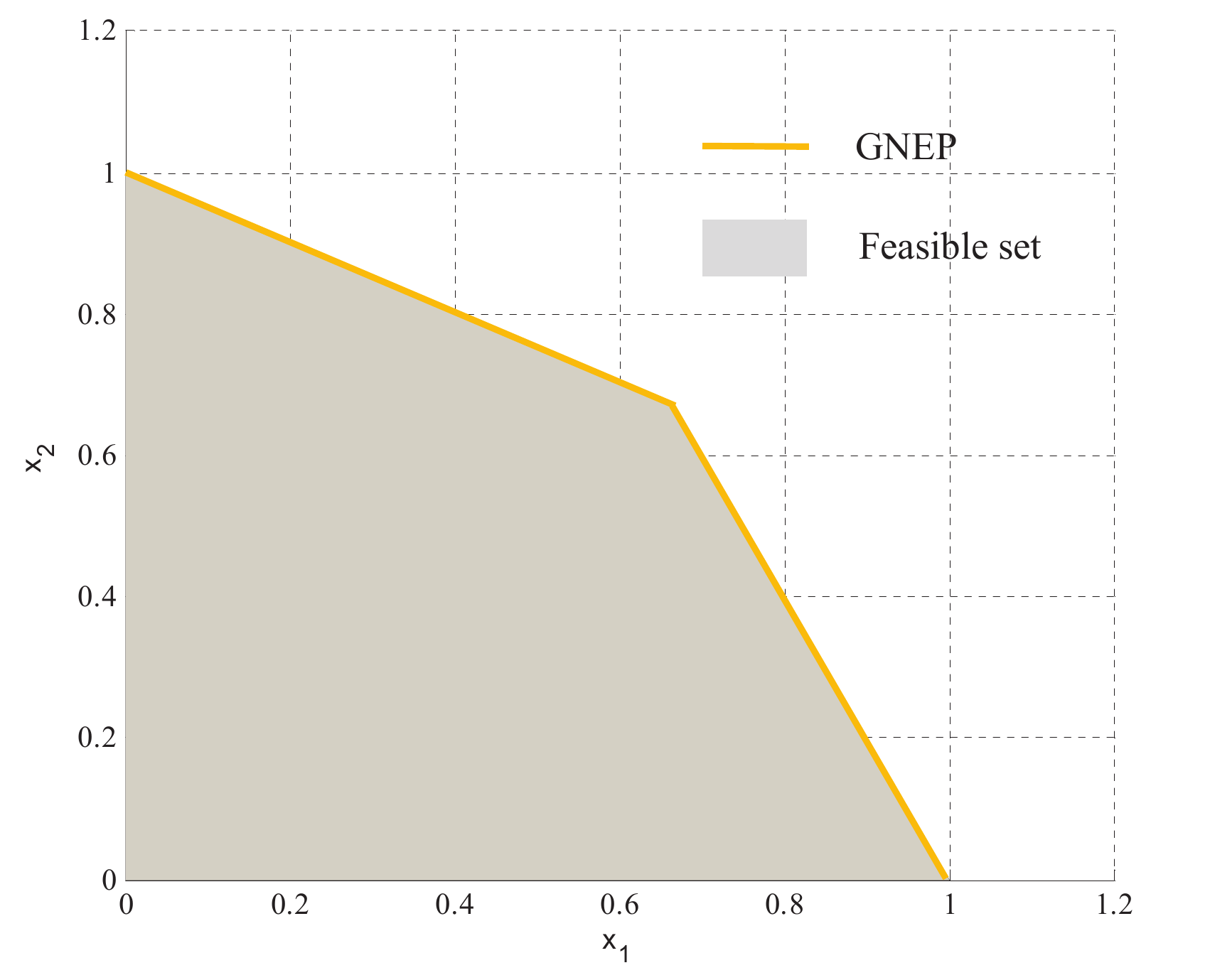}
\caption{Illustration of the equilibria of a simple GNEP.}
\label{fig:App-04-02}
\end{figure}

The feasible set $X$ of the strategy profile is plotted in Fig. \ref{fig:App-04-02}. It can be verified that any point on the line segments
\begin{equation}
L_1 = \left\{ (x_1,x_2) ~\middle|~ 0 \le x_1 \le \frac{2}{3},~
\frac{2}{3} \le x_2 \le 1,~ x_1 + 2 x_2 = 2 \right\}  \notag
\end{equation}
and
\begin{equation}
L_2 = \left\{ (x_1,x_2) ~\middle|~ \frac{2}{3} \le x_1 \le 1,~
0 \le x_2 \le \frac{2}{3},~ 2 x_1 + x_2 = 2 \right\}  \notag
\end{equation}
is an equilibrium point that satisfies Definition \ref{df:App-04-GENP}.

To refine a meaningful equilibrium from the infinitely many candidates, it is proposed to impose additional conditions on the Lagrange multipliers associated with shared constraints \cite{App-04-GNEP-3}. The outcome is called a restricted Nash equilibrium. Two special cases are discussed here. 

\vspace{12pt}
{\noindent \bf 1. Normalized Nash equilibrium}

The normalized Nash equilibrium
is firstly introduced in \cite{App-04-NNE-Rosen}. It incorporates a cone constraint on the dual multipliers
\begin{equation}
\label{eq:App-04-NNE}
\lambda_i = \beta_i \lambda_0,~ \beta_i > 0,~ i = 1,\cdots,n 
\end{equation}
where $\lambda_0 \in \mathbb R^m_+$. Solving KKT system (\ref{eq:App-04-GNEP-KKT}) with constraint (\ref{eq:App-04-NNE}) gives an  equilibrium solution. It is shown that for any given $\beta \in \mathbb R^n_{++}$, a normalized Nash equilibrium exists as long as the game is feasible. Moreover, if the mapping \begin{equation}
F(\beta): \mathbb R^n \to \mathbb R^n = \left(
\begin{gathered}
\frac{1}{\beta_1} \nabla_{x_1}  f_1(x_1,x_{-1})  \\  \vdots   \\
\frac{1}{\beta_n} \nabla_{x_n}  f_n(x_n,x_{-n})
\end{gathered}  \right)  \notag
\end{equation}
parameterized in $\beta$ is strictly monotone (by assuming convexity of payoff functions), then the normalized Nash equilibrium is unique.

The relation given in (\ref{eq:App-04-NNE}) indicates that the dual variables $\lambda_i$ associated with the shared constraints are a constant vector scaled by different scalars. From an economic perspective, this means that the shadow prices of common resources at any normalized Nash equilibrium are proportional among each player. 

\vspace{12pt}
{\noindent \bf 2. Variational equilibrium}

Recall the variational inequality formulation for the NEP in Proposition \ref{pr:App-04-NE-VI}, a GNEP with shared convex constraints can be treated in the same way: Let  $F(x) = (\nabla_{x_1} f_1 (x),\cdots,\nabla_{x_n} f_n (x))$ be a mapping, and the feasible region $X$ is defined in (\ref{eq:App-04-GENP-Convex-X}), then every solution of variational inequality problem VI($X,F$) gives an equilibrium solution of the GNEP, which is called the variational equilibrium (VE). 

However, unlike an NEP and its associated VI problem which have the same solutions, not all equilibria of the GNEP are preserved when it is passed to a corresponding VI problem; see \cite{App-04-GNEP-VI-2,App-04-GNEP-VI-3} for examples and further details. In fact, a solution $x^*$ of a GNEP is a VE if and only if it solves KKT system (\ref{eq:App-04-GNEP-KKT}) with the following constraints on the Lagrange dual multipliers \cite{App-04-GNEP-VI-1,App-04-GNEP-2,App-04-GNEP-VI-2}:
\begin{equation}
\label{eq:App-04-VE}
\lambda_1 = \cdots = \lambda_n = \lambda_0 \in \mathbb R^m_+
\end{equation}
implying that all players perceive the same shadow prices of common resource at a VE. The VI approach has two important implications. First, it allows us analyze a GNEP using well-developed VI theory, such as conditions which could guarantee the existence and uniqueness of the equilibrium point; second, condition (\ref{eq:App-04-VE}) gives an interesting economic interpretation of the VE, and inspires pricing-based distributed algorithms to compute an equilibrium solution, which will be discussed in the next section.

The concept of potential game for NEPs directly applies to GNEPs. If a GNEP with shared convex constraints possesses a potential function $U(x)$ which satisfies (\ref{eq:App-04-PG-1}), an equilibrium can be retrieved from a mathematical program which minimizes the potential function over the feasible set $X$ defined in (\ref{eq:App-04-GENP-Convex-X}). To reveal the connection of the optimal solution and the VE, we omit constraints in the local strategy sets $Q_i$, $i=1,\cdots,n$ for notation simplicity, and write out the mathematical program as follows
\begin{equation}
\label{eq:App-04-Potential-GNEP-1}
\begin{aligned}
\min_x  ~~ &  U(x)  \\
\mbox{s.t.}  ~~ &  g(x) \le 0
\end{aligned}   
\end{equation}
whose KKT optimality condition is given by
\begin{equation}
\label{eq:App-04-Potential-GNEP-2}
\begin{gathered}
\nabla_{x} U (x) + \lambda^T \nabla_{x} g (x) = 0  \\
\lambda \ge 0,~g(x) \le 0,~ \lambda^T g (x) = 0
\end{gathered}
\end{equation}

The first equality can be decomposed into $n$ sub-equations
\begin{equation}
\label{eq:App-04-Potential-GNEP-3}
\nabla_{x_i} U (x) + \lambda^T \nabla_{x_i} g (x) = 0,~ i=1,\cdots,n 
\end{equation}
Recall (\ref{eq:App-04-PG-4}), $\nabla_{x_i} U (x) = \nabla_{x_i} f_i (x_i ,x_{-i})$, substituting it into (\ref{eq:App-04-Potential-GNEP-2}) we have
\begin{gather}
\nabla_{x_i} f_i (x_i ,x_{-i}) + \lambda^T \nabla_{x_i} g (x) = 0,~ i=1,\cdots,n \notag  \\
\lambda \ge 0,~g(x) \le 0,~ \lambda^T g (x) = 0   \notag
\end{gather}
which is exactly KKT system (\ref{eq:App-04-GNEP-KKT}) with identical shadow price constraint (\ref{eq:App-04-VE}). In this regard, we can see

\begin{proposition} 
\label{pr:App-04-Potential-GNEP}
Optimizing the potential function of a GNEP with shared convex constraints gives a variational equilibrium. 
\end{proposition}

\begin{figure}[!htp]
\centering
\includegraphics[scale=0.40]{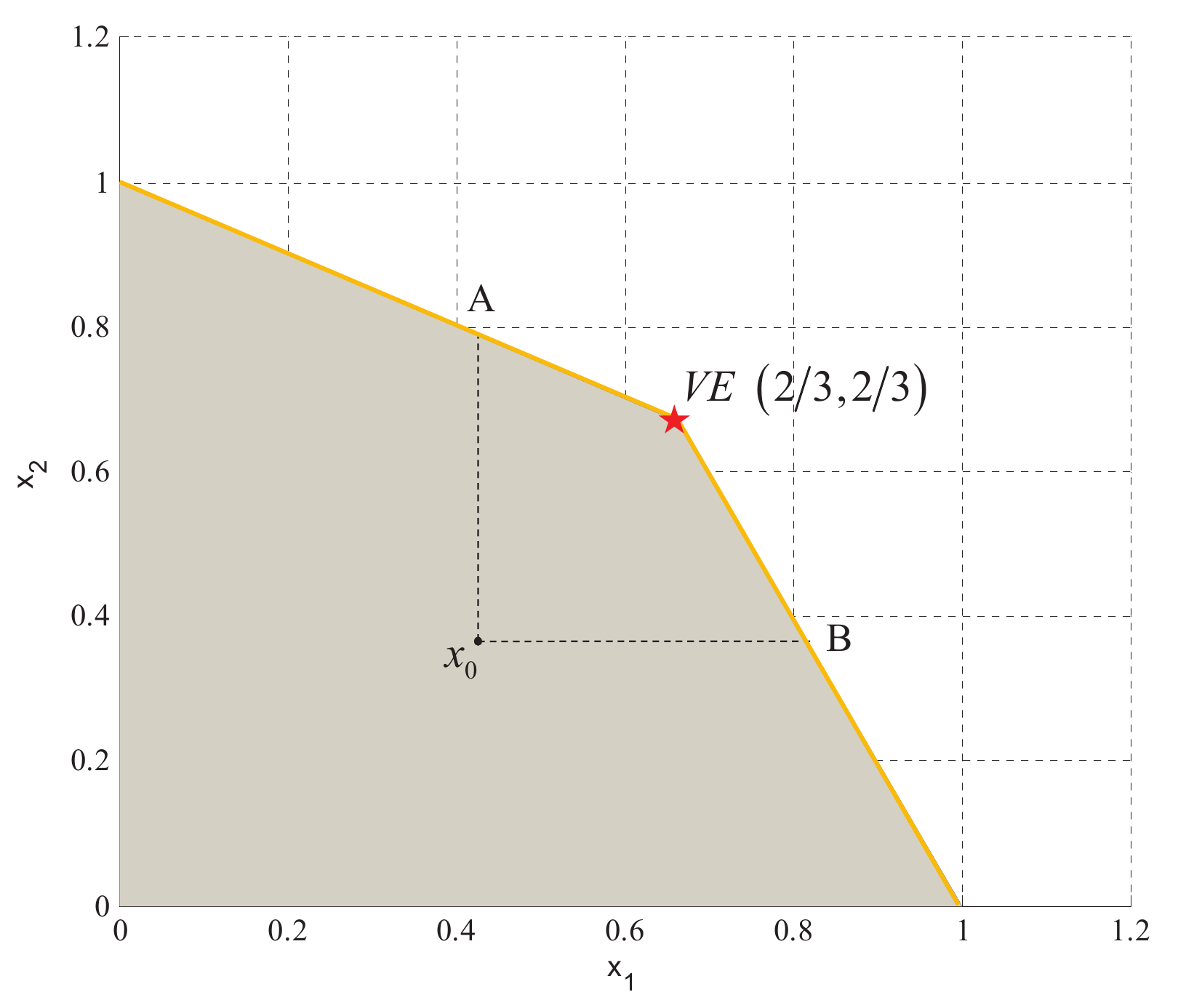}
\caption{Illustration of the variational equilibrium.}
\label{fig:App-04-03}
\end{figure}

Consider the example shown in Fig. \ref{fig:App-04-02} again, $(2/3,2/3)$ is the unique VE of the GNEP, which is plotted in Fig. \ref{fig:App-04-03}. The corresponding dual variables of global constraints $g_1 \le 0$ and $g_2 \le 0$ are $(1/3,1/3)$.  

\subsection{Best-Response Algorithm}
\label{App-D-Sect02-02}

The presence of shared constraints wrecks the Cartesian structure of $\prod_{i=1}^n Q_i$ in a standard Nash game, and prevents a direct application of the best response methods presented in Appendix \ref{App-D-Sect01-03} to solve an NEP. Moreover, even if an equilibrium can be found, it may depend on the initial point as well as the optimization sequence, because solutions of a GNEP are non-isolated. To illustrate this pitfall, take Fig. \ref{fig:App-04-03} for an example. Suppose we pick up an arbitrary point $x_0 \in X$ as the initial value. If we first maximize $x_1$ ($x_2$), the point moves to $B$ ($A$), and then in the second step, $x_2$ ($x_1$) does not change, because it is already an equilibrium solution in the sense of Definition \ref{df:App-04-GENP}. In view of this, fixed-point iteration may give any outcome on the line segments connecting $(2/3,2/3)$ and $(0,1)/(1,0)$, depending on the initiation.

This section introduces the distributed algorithms proposed in \cite{App-04-GNEP-VI-1} which identify a VE of GNEP (\ref{eq:App-04-GENP-MP}) with shared convex constraints. Motivated by the Lagrange decomposition framework, we can rewrite problem (\ref{eq:App-04-GENP-MP}) in a more convenient form. Consider finding a pair ($x,\lambda$), where $x$ is the equilibrium of the following standard NEP $\mathcal G(\lambda)$ with a given vector $\lambda$ of Lagrange multipliers 
\begin{equation}
\label{eq:App-04-GNEP-DIS-1}
\mathcal G(\lambda): \quad  \left\{ 
\begin{aligned}
\min_{x_i} ~~  &  f_i(x_i,x_{-i}) + \lambda ^T g(x)  \\
\mbox{s.t.}~~  &  x_i \in  Q_i
\end{aligned}  \right\},~ i=1,\cdots,n
\end{equation}
and furthermore, a complementarity constraint
\begin{equation}
\label{eq:App-04-GNEP-DIS-2}
0 \le \lambda \bot -g(x)  \ge 0
\end{equation}
Problem (\ref{eq:App-04-GNEP-DIS-1})-(\ref{eq:App-04-GNEP-DIS-2}) has a clear economic interpretation: suppose the shared constraints represent the availability of some common resources, vector $\lambda$ can be viewed as the prices paid by players for consuming these resources. Actually, when a resource is adequate, the inequality constraint is not binding and the Lagrange dual multiplier is zero; the dual multiplier or shadow price is positive only if a resource   becomes scarce, indicated by a binding inequality constraint. This relation has been imposed in constraint (\ref{eq:App-04-GNEP-DIS-2}). 

The KKT conditions of $\mathcal G(\lambda)$ (\ref{eq:App-04-GNEP-DIS-1}) in conjunction with condition (\ref{eq:App-04-GNEP-DIS-2}) turn out to be the VE condition of GNEP (\ref{eq:App-04-GENP-MP}). In view of this connection, a VE can be found by solving (\ref{eq:App-04-GNEP-DIS-1})-(\ref{eq:App-04-GNEP-DIS-2}) in a distributed manner based on previous algorithms developed for NEPs. Likewise, we discuss strongly convex cases and convex cases separately, due to their different convergence guarantees.

\vspace{12pt}
{\noindent \bf 1. Algorithms for strongly convex cases}

Suppose that the game $\mathcal G(\lambda)$ in (\ref{eq:App-04-GNEP-DIS-1}) is strongly convex and has a unique Nash equilibrium $x(\lambda)$ for any given $\lambda \ge 0$. This uniqueness condition allows defining the map 
\begin{equation}
\label{eq:App-04-GNEP-ALG-1}
{\rm \Phi} (\lambda): \lambda \to -g(x(\lambda))
\end{equation}
which quantifies the negative violation of the shared constraints at $x(\lambda)$. Based on (\ref{eq:App-04-GNEP-DIS-1})-(\ref{eq:App-04-GNEP-DIS-2}), the distributed algorithm is provided as follows. 

\begin{algorithm}[!htp]
\normalsize
\caption{\bf }
\begin{algorithmic}[1]
\STATE Choose an initial price vector $\lambda^0 \ge 0$. The iteration index is $k = 0$. 
\STATE Given $\lambda^k$, find the unique equilibrium $x(\lambda^k)$ of $\mathcal G(\lambda^k)$ using Algorithm \ref{Ag:App-04-NE-Asy-BR}.

\STATE If $0 \le \lambda^k \bot {\rm \Phi} (\lambda^k) \ge 0$ is satisfied, terminate and report $x(\lambda^k)$ as the VE; otherwise, choose $\tau_k > 0$, and update the price vector according to
\begin{equation}
\lambda^{k+1}=\left[\lambda^k - \tau_k {\rm \Phi} (\lambda^k) \right]^+ \notag
\end{equation}
set $k \leftarrow k+1$, and go to step 2.
\end{algorithmic}
\label{Ag:App-04-GNE-1}
\end{algorithm} 

Algorithm \ref{Ag:App-04-GNE-1} is a double-loop method. The range of parameter $\tau_n$ and convergence proof have been thoroughly discussed in \cite{App-04-GNEP-VI-1} based on the monotonicity of the mapping $F+ \nabla g(x) \lambda$, where $F=(\nabla_{x_i} f_i)_{i=1}^n$, and $\nabla g(x)$ is a matrix whose $i$-th column is equal to $\nabla g_i$.

\vspace{12pt}
{\noindent \bf 2. Algorithms for convex cases}

Now we consider the case in which the VI associated with problem (\ref{eq:App-04-GNEP-DIS-1})-(\ref{eq:App-04-GNEP-DIS-2}) is
merely monotone (at least one problem in (\ref{eq:App-04-GNEP-DIS-1}) is not strongly convex). In such circumstance, the convergence of Algorithm \ref{Ag:App-04-GNE-1} is no longer guaranteed. This is not only because Algorithm \ref{Ag:App-04-NE-Asy-BR} for the inner loop game $\mathcal G(\lambda)$ may not converge, but also because the outer loop has to be complicated. To circumvent this difficulty, we try to convexify the game using regularization terms as what has been done in Algorithm \ref{Ag:App-04-NE-Proximal}. To this end, we have to explore an optimization reformulation for the complementarity constraint (\ref{eq:App-04-GNEP-DIS-2}), which is given by 
\begin{equation}
\lambda \in \arg \min_{\bar \lambda} \left\{ - \bar \lambda^T g(x) ~\middle|~ \bar \lambda \ge 0 \right\}  \notag
\end{equation}

Then, consider the following ordinary NEP with $n+1$ players in which the last player controls the price vector $\lambda$:
\begin{equation}
\label{eq:App-04-GNEP-n+1}
\begin{aligned}
\min_{x_i} ~~  & \left\{ f_i(x_i,x_{-i}) + \lambda ^T g(x) 
~\middle|~  x_i \in  Q_i  \right\},~ i = 1,\cdots,n  \\
\min_{\lambda} ~~ &  \left\{ -\lambda^T g(x) ~\middle|~ \lambda \ge 0 \right\}
\end{aligned}  
\end{equation}
where the last player solves an LP in variable $\lambda$ parameterized in $x$. At the equilibrium, $g(x) \le 0$ is implicitly satisfied. To see this, because $Q_i$ is bounded, problems of the first $n$ players must have a finite optimum for arbitrary $\lambda$. If $g(x) \nleq 0$, the last problem has an infinite optimum, imposing a large penalty on the constraint that is violated, and thus the first $n$ players will alter their strategies accordingly. Whenever $g(x) \le 0$ is met, the last LP must have a zero minimum, which satisfies (\ref{eq:App-04-GNEP-DIS-2}). In summary, this extended game (\ref{eq:App-04-GNEP-n+1}) has the same equilibria as problem (\ref{eq:App-04-GNEP-DIS-1})-(\ref{eq:App-04-GNEP-DIS-2}). Since the strategy sets of (\ref{eq:App-04-GNEP-n+1}) have a Cartesian structure, Algorithm \ref{Ag:App-04-NE-Proximal} can be applied to find an equilibrium.

\begin{algorithm}[!htp]
\normalsize
\caption{\bf }
\begin{algorithmic}[1]
\STATE Given $\{ \rho_n\}_{n=0}^\infty$, $\varepsilon >0$, and $\tau > 0$, choose a feasible initial point $x^0 \in X$ and an initial price vector $\lambda^0$; the iteration index is $k=0$.
\STATE Given $z^k = (x^k,\lambda^k)$, find a Nash equilibrium $z^{k+1} = (x^{k+1},\lambda^{k+1})$ of the following regularized NEP using Algorithm \ref{Ag:App-04-NE-Asy-BR}
\begin{equation}
\begin{aligned}
\min_{x_i} ~~  & \left\{ f_i(x_i,x_{-i}) + \lambda ^T g(x) + \tau \left\| x_i - x^k_i  \right\|^2_2 ~\middle|~  x_i \in  Q_i  \right\},~ i = 1,\cdots,n\\
\min_{\lambda} ~~ &  \left\{ -\lambda^T g(x) + \tau \left\| \lambda - \lambda^k  \right\|^2_2 ~ ~\middle|~ \lambda \ge 0 \right\}
\end{aligned}   \notag
\end{equation}

\STATE If $\| z^{k+1} - z^k \|_2 \le \varepsilon$, terminate and report $x^{k+1}$ as the variational equilibrium; otherwise, update $k \leftarrow k+1$, $z^k \leftarrow (1-\rho_k) z^{k-1} + \rho_k z^k$, and go to step 2.
\end{algorithmic}
\label{Ag:App-04-GNE-2}
\end{algorithm} 

The convergence of Algorithm \ref{Ag:App-04-GNE-2} is guaranteed under a sufficiently large $\tau$. More quantitative discussions on parameter selection and convergence conditions can be found in \cite{App-04-GNEP-VI-1}. In practice, the value of $\tau$ should be carefully chosen to achieve satisfactory computational performances.

\vspace{12pt}

In NEPs and GNEPs, players make simultaneous decisions. In real-life decision making problems, there are many situations in which players can move sequentially. In the rest of this chapter, we consider three kinds of bilevel games, in which the upper-level (lower-level) players are called leaders (followers), and leaders make decisions prior to follower’s. The simplest one is the Stackelberg game, or the single-leader-single-follower game, or just the bilevel program; Stackelberg game can be generalized by incorporating multiple players in the upper and lower levels. Players at the same level make decisions simultaneously, whereas followers' actions are subject to leaders' movements, forming an NEP parameterized in the leaders' decisions. When there is only one leader, the problem is called a mathematical program with equilibrium constraints (MPEC); when there are multiple leaders, the problem is referred to as an equilibrium  program with equilibrium constraints (EPEC). It is essentially a bilevel GNEP among the leaders.

\section{Bilevel Programs}
\label{App-D-Sect03}

Bilevel program is a special mathematical program with another optimization problems nested in the constraints. The main problem is called the upper-level problem, and the decision maker is the leader; the one nested in constraints is called the lower-level problem, and the decision maker is the follower. In game theory, a bilevel program is usually referred to as the Stackelberg game, which arises in many economic and engineering design problems.

\subsection{Bilevel Programs with a Convex Lower Level}

{\bf 1. Mathematic model and single-level equivalence}

A bilevel program is the most basic instance of bilevel games. The leader moves first and chooses a decision $x$; then the follower selects its strategy $y$  solving the lower-level problem parameterized in $x$
\begin{equation}
\label{eq:App-04-BLP-LLP}
\begin{aligned}
\min_y ~~ & f(x,y)   \\
\mbox{s.t.}~~ & g(x,y) \le 0: \lambda  \\
          & h(x,y) = 0: \mu 
\end{aligned}
\end{equation}
where $\lambda$ and $\mu$ following the colon are dual variables associated with inequality and equality constraints, respectively. We assume that problem (\ref{eq:App-04-BLP-LLP}) is convex and the KKT condition is necessary and sufficient for a global optimum
\begin{equation}
\label{eq:App-04-BLP-LLP-KKT}
\mbox{Cons-KKT} = \left\{(x,y,\lambda,\mu) ~\middle|~
\begin{gathered}
\nabla_y f(x,y) + \lambda^T \nabla_y g(x,y) + \mu^T \nabla_y h(x,y) = 0 \\
0 \le \lambda \bot - g(x,y) \ge 0, ~ h(x,y) = 0 
\end{gathered}  \right\}
\end{equation}
The set of optimal solutions of problem (\ref{eq:App-04-BLP-LLP}) is denoted by $S(x)$. If (\ref{eq:App-04-BLP-LLP}) is strictly convex, the optimal solution is unique, and $S(x)$ reduces to a singleton. 

When the leader minimizes its payoff function $F(x,y)$, the best response $y(x) \in S(x)$ is taken into account.  The leader's problem is formally
described as
\begin{equation}
\label{eq:App-04-BLP-ULP}
\begin{aligned}
\min_{x, \bar y} ~~ & F(x,\bar y)   \\
\mbox{s.t.}~~ & x \in X  \\
          & \bar y \in S(x)  
\end{aligned}
\end{equation}
Notice that although $\bar y$ acts as a decision variable of the leader, it is actually controlled by the follower through the best response mapping $S(x)$. When the leader makes decisions, it will take the response from the follower into account. When $S(x)$ is a singleton, qualifier $\in$ reduces to $=$; otherwise, if $S(x)$ contains more than one elements, (\ref{eq:App-04-BLP-ULP}) assumes that the follower will choose the one which is preferred by the leader. Therefore,  (\ref{eq:App-04-BLP-ULP}) is called an optimistic equivalence. On the contrary, the pessimistic equivalence assumes that  the follower will choose the one which is unfavorable for the leader, which is more difficult to solve. As for the optimistic case, replacing $\bar y \in S(x)$ with KKT condition (\ref{eq:App-04-BLP-LLP-KKT}) leads to the NLP formulation of the bilevel program, or more exactly, a mathematical program with complementarity constraints (MPCC)
\begin{equation}
\label{eq:App-04-BLP-MPCC}
\begin{aligned}
\min_{x,\bar y,\lambda,\mu} ~~ & F(x,\bar y)   \\
\mbox{s.t.}~~ & x \in X,~ 
(x,\bar y,\lambda,\mu) \in \mbox{Cons-KKT}
\end{aligned}
\end{equation} 

Although the lower-level problem (\ref{eq:App-04-BLP-LLP}) is convex, the best reaction map of the follower characterized by Cons-KKT is non-convex, so a bilevel program is intrinsically non-convex and generally difficult to solve. 

\vspace{12pt}
{\noindent \bf 2. Why bilevel programs are difficult to solve?}

Two difficulties prevent an MPCC from being solved reliably and efficiently.

1) The feasible region of (\ref{eq:App-04-BLP-ULP}) is non-convex: even if objective functions and constraints of the leader and the follower are linear, the complementarity and slackness condition in (\ref{eq:App-04-BLP-LLP-KKT}) is still non-convex. An NLP solver only finds a local solution for non-convex problems, if succeeds, and global optimality can hardly be guaranteed.

2) Despite of its non-convexity, the failure to meet ordinary constraint qualifications creates another barrier for solving an MPCC. NLP algorithms generally stop when a stationary point of the KKT conditions is found; however,
due to the presence of the complementarity and slackness condition, the dual multipliers may not be well-defined because of the violation of  standard constraint qualifications. Therefore, NLP solvers may fail to find a local optimum without particular treatment on the complementarity constraints. To see how constraint qualifications are violated, consider the following simplest linear complementarity constraint 
\begin{equation}
x \ge 0,~ y \ge 0, ~ x^T y =0,~ x \in \mathbb R^5,~ y \in \mathbb R^5  \notag
\end{equation}
The Jacobian matrix of the active constraints at point ($\bar x, \bar y$) is
\begin{equation}
J = \left[ ~ 
\begin{gathered} 
e_{\bar x}  \\  0  \\  \bar y 
\end{gathered} ~ 
\begin{gathered} 
 0 \\ e_{\bar y}   \\  \bar x 
\end{gathered} ~
 \right]  \notag
\end{equation}
where $e_{\bar x}$ and $e_{\bar y}$ are zero-one matrices corresponding to the active constraints $x_i =0$, $i \in I$, $y_j = 0$, $j \in J$, where $I \bigcup J = \{1,2,3,4,5\}$, and $I \bigcap J$ is not necessarily empty. Suppose that $I = \{1,2,4\}$ and $J = \{3,5\}$, then  
\begin{equation}
J = \begin{bmatrix}
1  &  0  &  0  &  0  &  0  &  0  &  0  &  0  &  0  &  0 \\
0  &  1  &  0  &  0  &  0  &  0  &  0  &  0  &  0  &  0 \\
0  &  0  &  0  &  1  &  0  &  0  &  0  &  0  &  0  &  0 \\
0  &  0  &  0  &  0  &  0  &  0  &  0  &  1  &  0  &  0 \\
0  &  0  &  0  &  0  &  0  &  0  &  0  &  0  &  0  &  1 \\
\bar y_1  & \bar y_2  & \bar y_3  & \bar y_4 & \bar y_5 &
\bar x_1  & \bar x_2  & \bar x_3  & \bar x_4 & \bar x_5
\end{bmatrix}  \notag
\end{equation}
Since $\bar x_1 = \bar x_2 = \bar x_4 = 0$ and $\bar y_3 = \bar y_5 = 0$, it is apparent that the row vectors of $J$ are linearly dependent at point ($\bar x, \bar y$). The same applies to any ($\bar x, \bar y$) regardless of the indices $I$ and $J$ of active constraints, because whenever $y_j > 0$, complementarity will enforce $x_i=0,i=j$, creating a binding inequality in $x \ge 0$ and a row in matrix $J$ whose $i$-th  element is 1; whenever $x_i > 0$, complementarity will enforce $y_j=0,j=i$, creating a binding inequality in $y \ge 0$ and a row in matrix $J$ whose $(i+5)$-th  element is 1. Therefore, the last row of $J$ can be represented by a linear combination of the other rows. 

Above discussion and conclusion on linear complementarity constraints also apply to the nonlinear case, because the Jacobian matrix $J$ has the same structure. In this regard, above difficulty is an intrinsic phenomenon in MPCCs.     

\begin{proposition}
\label{pr:App-04-MPCC-MFCQ}
Complementarity and slackness conditions violate the linear independent constraint qualification at any feasible solution.
\end{proposition}

From a geometric perspective, the feasible region of complementarity constraints consists of slices like $x_i=0$, $y_j = 0$; there is no strictly feasible point and the Slater's condition does not hold. In conclusion, general purpose NLP solvers are not numerically reliable for solving MPCCs, although they
were once used to carry out such tasks.   

\vspace{12pt}
{\noindent \bf 3. Methods for solving MPCCs}

In view of the
limitations of standard NLP algorithms, new constraint qualifications are proposed to define stationary solutions so as to solve MPCCs through conventional NLP methods, such as the Bouligand-, Clarke-, Mordukhovich-, weakly-, and Strongly-stationary constraint qualifications. See \cite{App-04-BLP-CQ-1,App-04-BLP-CQ-2,App-04-BLP-CQ-3,App-04-BLP-CQ-4} for further information. Through some proper transformation, MPCCs can be solved via standard NLP algorithms as well. Several approaches are available for this task. 

\vspace{12pt}
{\noindent \bf a. Regularization method \cite{App-04-MPCC-Reg-1,App-04-MPCC-Reg-2,App-04-MPCC-Reg-3}}

In this approach, the non-negativity and complementarity requirements
\begin{equation}
\label{eq:App-04-MPCC-Reg-1}
x \ge 0,~ y \ge 0,~ xy =0
\end{equation} 
are approximated by
\begin{equation}
\label{eq:App-04-MPCC-Reg-2}
x \ge 0,~ y \ge 0, ~ xy \le \varepsilon
\end{equation}
Please note that $xy \ge 0$ is a natural result of non-negativity requirements on $x$ and $y$. When $\varepsilon=0$, (\ref{eq:App-04-MPCC-Reg-2}) is equivalent to (\ref{eq:App-04-MPCC-Reg-1}); when $\varepsilon > 0$, (\ref{eq:App-04-MPCC-Reg-2}) defines a larger feasible region than (\ref{eq:App-04-MPCC-Reg-1}), so this approach is sometimes called a relaxation method. The smaller $\varepsilon$ is, the closer any feasible point $(x,y)$ is to achieve complementarity.
if $x$ and $y$ are vectors with non-negative elements, $x^T y =0$ is the same as $x_i y_i =0$, $i=1,2,\cdots$. The same procedure can be applied if $x$ and $y$ are replaced by nonlinear functions.
 
Since Slater's condition holds for the feasible set defined by (\ref{eq:App-04-MPCC-Reg-2}) with $\epsilon > 0$, NLP solvers can be used to solve related optimization problem. In a regularization procedure for solving an MPCC, the relaxation (\ref{eq:App-04-MPCC-Reg-2}) is applied with gradually decreased value of $\varepsilon$ for implementation issues. If the initial value of $\varepsilon$ is too small, the solver may be numerically unstable and fail to find a feasible solution.

\vspace{12pt}
{\noindent \bf b. Penalization method \cite{App-04-MPCC-Pen-1,App-04-MPCC-Pen-2,App-04-MPCC-Pen-3}}

In this approach, the complementarity condition $xy=0$ is removed from the set of constraints; instead, an associated penalty term $xy/\varepsilon$ is added to the objective function to create an extra cost whenever complementarity is not satisfied. Since $x$ and $y$ are non-negative, as indicated by (\ref{eq:App-04-MPCC-Reg-1}), the penalty term would never take a negative value. In this way, the feasible region becomes much simpler.

In a penalization procedure for solving an MPCC, a sequence of NLPs are solved iteratively with gradually decreased value of $\varepsilon$, and the violation of complementarity condition gradually approaches to 0 as iterations proceed. If $\varepsilon$ is initiated too small, the penalty coefficient $1/\varepsilon$ is very large which may cause an ill-conditioned problem and numeric instability. One advantage of this iterative procedure is that  the optimal solution in iteration $k$ can be used as the initial guess in iteration $k+1$, since the feasible region does not change, and the solution in every iteration is feasible in the next one. A downside of this approach is that the NLP solver generally identifies a local optimum. In consequence, a smaller $\epsilon$ may not necessarily lead to a solution that gets closer to the feasible region.

\vspace{12pt}
{\noindent \bf c. Smoothing method \cite{App-04-MPCC-Smooth-1,App-04-MPCC-Smooth-2}}

This approach employs the perturbed Fischer-Burmeister function 
\begin{equation}
\label{eq:App-04-Fischer-Burmeister-1}
\phi(x,y,\varepsilon) = x + y - \sqrt{x^2 + y^2 + \varepsilon}
\end{equation}
which is firstly introduced in \cite{App-04-MPCC-SQP-1} for LCPs, and shown particularly useful in SQP methods for solving MPCCs in \cite{App-04-MPCC-SQP-2}. Clearly, when $\varepsilon= 0$, the function $\phi$ reduces to the standard Fischer-Burmeister function
\begin{equation}
\label{eq:App-04-Fischer-Burmeister-2}
\phi(x,y,0) = 0 ~\Longleftrightarrow~ x \ge 0,~ y \ge 0,~ xy=0
\end{equation}
$\phi(x,y,0)$ is not smooth at the origin $(0,0)$. When $\varepsilon > 0$, the function $\phi$ satisfies
\begin{equation}
\label{eq:App-04-Fischer-Burmeister-2}
\phi(x,y,\varepsilon) = 0 ~\Longleftrightarrow~ x \ge 0,~ y \ge 0,~ 
xy = \varepsilon/2 
\end{equation}
and is smooth in $x$ and $y$.

In view of this, complementarity and slackness condition (\ref{eq:App-04-MPCC-Reg-1}) can be replaced by $\phi(x,y,\varepsilon) = 0$ and further embedded in NLP models. When $\varepsilon$ tends to 0, (\ref{eq:App-04-MPCC-Reg-1}) is enforced approximately.

\vspace{12pt}
{\noindent \bf d. Sequential quadratic programming (SQP) \cite{App-04-MPCC-SQP-1,App-04-MPCC-SQP-2,App-04-MPCC-SQP-3}}

SQP is a general purpose NLP method. In each iteration of SQP, the quadratic functions in complementarity constraints are approximated by a linear one, and the nonlinear objective function is replaced with their second-order Taylor series, constituting a quadratic program with linear constraints (maybe in conjunction with trust region bounds). At the optimal solution, nonlinear constraints are linearized, the objective function is approximated again, and then the SQP algorithm proceeds to the next iteration. 

When applied to an MPCC, the SQP method is often capable of finding a local optimal solution, without a sequence of user-specified $\varepsilon_k$ approaching to 0, probably because the SQP solver itself is endowed with some softening ability, e.g., when a quadratic program encounters numeric issues, the SQP solver SNOPT automatically relaxes some hard constraints and penalizes violations in the objective function.

The aforementioned classical methods are discussed in \cite{App-04-BLP-NLP-1}, and numeric experiences are reported in \cite{App-04-MPCC-NLP-Test}.

\vspace{12pt}
{\noindent \bf e. MINLP methods}

Due to the wide applications in various engineering disciplines, solution methods of MPCCs continue to be an active research area. Recall that the complementarity constraints in form of $g(x) \ge 0$, $h(x) \ge 0$, $g(x)h(x)=0$ is equivalent to
\begin{equation}
0 \le g(x) \le M z,~ 0 \le h(x) \le M(1-z) \notag
\end{equation}
where $z$ is a binary variable, $M$ is a sufficiently large constant. Therefore, an MPCC can be converted to a mixed integer nonlinear program
(MINLP). MINLP removes the numeric difficulty in MPCC; however, the computation complexity remains. If all functions in (\ref{eq:App-04-BLP-MPCC}) are linear, or there are only a few complementarity constraints, the resulting MILP or MINLP model may be solved within reasonable time; otherwise, the branch-and-bound algorithm could offer upper and lower bounds on the optimal value.

\vspace{12pt}
{\noindent \bf f. Convex relaxation/approximation methods}

If all functions in MPCC (\ref{eq:App-04-BLP-MPCC}) are linear, it is a non-convex QCQP in which non-convexity originates from the complementarity constraints. When the problem scale is large, the MILP method may be time-consuming. Inspired by the success of convex relaxation methods in non-convex QCQPs, there have been increasing interests for developing convex relaxation methods for MPCCs. An SDP relaxation method is proposed in \cite{App-04-MPCC-SDP-1}, which is embedded in a branch-and-bound algorithm to solve the MPCC. For the MPCC derived from a bilevel polynomial program, it is proposed to solve a sequence of SDPs with increasing problem sizes, so as to solve the original problem globally \cite{App-04-MPCC-SDP-2,App-04-MPCC-SDP-3}. Convex relaxation methods have been applied to power market problems in \cite{App-04-MPCC-SDP-4,App-04-MPCC-SDP-5}. Numerical experiments show that the combination of MILP and SDP relaxation can greatly reduce the computation time. Nonetheless, please bear in mind that in the SDP relaxation model, the decision variable is a matrix with a dimension of $n \times n$, so solving the SDP model may still be a challenging task, although it is convex.  

Recently, a DC programming approach is proposed in \cite{App-04-LPCC-DCP} to solve LPCC in the penalized version. In this approach, the quadratic penalty term is decomposed into the difference of two convex quadratic functions, and the concave part is then linearized. Computational performances reported in \cite{App-04-LPCC-DCP} are very promising.

\subsection{Special Bilevel Programs }

Although general bilevel programs are difficult, there are special cases which can be solved relatively easily. One of such classes of programs is the linear bilevel program, in which objective  functions are linear and constraints are polyhedra. The linear max-min problem is a special case of the linear bilevel program, in which the leader and the follower have completely opposite targets. Furthermore, two special market models are studied.

\vspace{12pt}
{\noindent \bf 1. Linear bilevel program}

A linear bilevel program can be written as 
\begin{equation}
\label{eq:App-04-Linear-Bilevel-1}
\begin{aligned}
\max_x ~~ & c^T x + d^T y(x)   \\
\mbox{s.t.} ~~& C x \le d  \\
& \begin{aligned}
y(x) \in \arg \min_y ~~ & f^T y \\
\mbox{s.t.} ~~ & By \le b - Ax    
\end{aligned}
\end{aligned}
\end{equation}
In problem (\ref{eq:App-04-Linear-Bilevel-1}), the follower makes a decision $y$ after the leader deploys its action $x$, which influences the feasible region of $y$. Meanwhile, the leader can predict the follower's optimal response $y(x)$, and choose a strategy that finally optimizes $c^T x + d^T y(x)$. Other matrices and vectors are constant coefficients.

Given the upper level decision $x$, the follower is facing an LP, whose KKT optimality condition is given by
\begin{equation}
\begin{gathered}
B^T u = f  \\
0 \ge u ~\bot~ Ax + By - b \le 0  
\end{gathered}  \notag
\end{equation}
The last constraint is equivalent to the following linear constraints
\begin{equation}
\begin{gathered}
-M(1-z) \le u \le 0 \\
-Mz \le Ax+By-b \le 0  
\end{gathered}  \notag
\end{equation}
where $z$ is a vector consisting of binary variables, and $M$ is a large enough constant. 

In problem (\ref{eq:App-04-Linear-Bilevel-1}), replacing follower's LP with its KKT condition gives rise to an MILP 
\begin{equation}
\label{eq:App-04-Linear-Bilevel-MILP}
\begin{aligned}
\max_{x,y,u,z} ~~ & c^T x + d^T y   \\
\mbox{s.t.} ~~& C x \le d, ~ B^T u = f  \\
& -M(1-z) \le u \le 0 \\
&-Mz \le Ax + By - b \le 0  
\end{aligned}
\end{equation}
If the number of complementarity constraints is moderate, MILP (\ref{eq:App-04-Linear-Bilevel-MILP}) can be often solved efficiently, despite of its NP-hard complexity in the worst-case. Since MILP solvers and computation hardware keep improving nowadays, it is always worthy of bearing this technique in mind. Please also be aware that the big-M parameter notably impacts the performance of solving MILP (\ref{eq:App-04-Linear-Bilevel-MILP}). A heuristic method to determine such a parameter in linear bilevel programs is proposed in \cite{App-04-LPCC-BigM}. This method firstly solves two LPs and generates a feasible solution of the equivalent MPCC; then solves a regularized version of the MPCC model using NLP solvers and identifies a local optimal solution near the obtained feasible point; finally, the big-M parameter and the binary variables are initiated according to the local optimal solution. In this way, no manually-supplied parameter is needed, and the MILP model is properly strengthened.

Another optimality certification of follower's LP is the following primal-dual optimality condition
\begin{equation}
\begin{gathered}
B^T u = f,~ u \le 0,~ Ax + By \le b  \\
u^T (b-Ax) = f^T y  
\end{gathered}  \notag
\end{equation}
The first line summarizes feasible regions of the primal and dual variables. The last equation enforces equal values on the optimums of the primal and the dual problems, which is known as the strong duality condition. 

Replacing follower's LP with the primal-dual optimality condition gives an NLP:
\begin{equation}
\label{eq:App-04-Linear-Bilevel-NLP}
\begin{aligned}
\max_{x,y,u} ~~ & c^T x + d^T y   \\
\mbox{s.t.} ~~& C x \le d, ~ Ax + By \le b \\
              & u \le 0,~ B^T u = f \\  
              & u^T (b-Ax) = f^T y
\end{aligned}
\end{equation}
The following discussion are divided in two categories based on the type of variable $x$. 

a. $x$ is {\bf continuous}. In such a general situation, there is no effective way to solve problem (\ref{eq:App-04-Linear-Bilevel-NLP}), due to the last bilinear equality. Notice the fact that $f^T y \ge u^T (b-Ax) $ always holds on the feasible region because of the weak duality, the last constraint can be relaxed and penalized in the objective function, resulting in a bilinear program over a polyhedron \cite{App-04-LBLP-Pen-1,App-04-LBLP-Pen-2,App-04-LBLP-Pen-3} 
\begin{equation}
\label{eq:App-04-Linear-Bilevel-NLP-Pen}
\begin{aligned}
\max_{x,y,u} ~~ & c^T x + d^T y - \sigma [f^T y - u^T (b-Ax)]  \\
\mbox{s.t.} ~~& C x \le d, ~ Ax + By \le b \\
              & u \le 0,~ B^T u = f   
\end{aligned}
\end{equation}
where $\sigma > 0$ is a penalty parameter. In problem (\ref{eq:App-04-Linear-Bilevel-NLP-Pen}), the constraints on $u$ and $(x,y)$ are decoupled, so this problem can be solved by Algorithm \ref{Ag:App-03-BLP-Mountain-Climbing} (mountain climbing) in Appendix \ref{App-C-Sect02-03}, if global optimality is not mandatory.

In some problems, the upper-level decision influences the lower-level cost function, and has no impact on the feasible region in the lower level. For example, the tax rate design or a retail market pricing belongs to such category. The same procedure can be performed to solve this kind of bilevel problem. We recommend the MILP model, because in the penalized model, both $f^T y$ and $u^T Ax$ are non-convex. A tailored retail market model will be introduced later.

b. $x$ is {\bf binary}. In such circumstance, the bilinear term $ u^T A x =\sum_{ij} A_{ij} u_i x_j$ can be linearized by replacing $u_i x_j$ with a new continuous variable $v_{ij}$ together with auxiliary linear inequalities enforcing $v_{ij}=u_i x_j$. In this way, the last inequality translates into 
\begin{equation*}
\begin{gathered}
u^T b - \sum\nolimits_{ij} A_{ij} v_{ij} = f^T y  \\
u^l_i x_j \le v_{ij} \le 0,~ 
u^l_i (1-x_j) \le u_i-v_{ij} \le 0,~ \forall i,j 
\end{gathered}
\end{equation*}
where $u^l_i$ is a proper bound that does not discard the original optimal solution. As we can see, a bilevel linear program with binary upper-level variables is not necessarily harder than all continuous instances. This formulation is very useful to model interdiction problems in which $x$ mimics attack strategy. 

c. $x$ can be {\bf discretized}. Even if $x$ is continuous, we can approximate it via binary expansion 
\begin{equation*}
x_i = x^l_i + {\rm \Delta}_i \sum_{k=0}^{K} 2^k z_{ik},~ z_{ik} \in \{0,1\}
\end{equation*} 
where $x^l_i$ ($x^m_i$) is the lower (upper) bound of $x_i$, and ${\rm \Delta}_i = (x^m_i-x^l_i)/2^{K+1}$ is the step size. With this transformation, the bilinear term $u^T A x$ becomes
\begin{equation*}
\sum_{ij} A_{ij} u_i x^l_i + \sum_{ij} A_{ij} u_i  {\rm \Delta}_j \sum_{k=0}^{K} 2^k z_{jk}
\end{equation*}
The first term is linear, and $u_i z_{jk}$ in the second term can be linearized in a similar way. However, this entails introducing continuous variable with respect to indices $i$, $k$ and $k$. A low-complexity linearization method is suggested in \cite{App-04-BLLP-MILP-Sim}. It re-orders the summations in the second term as 
\begin{equation*}
\sum_j \sum_{k=0}^{K} {\rm \Delta}_j 2^k z_{jk} \sum_i A_{ij} u_i   
\end{equation*}
which can be linearized through defining an auxiliary continuous variable $v_{jk} = z_{jk} \sum_i A_{ij} u_i$ and stipulating 
\begin{equation*}
- M z_{jk} \le v_{jk} \le M z_{jk},~
- M (1-z_{jk}) \le \sum_i A_{ij} u_i - v_{jk} \le M (1-z_{jk})
\end{equation*}
where $M$ is a large enough constant.  

The core idea behind this trick is to treat $u^T A$ as a whole vector which has the same dimension as $x$, because for bilinear form $x^T v = \sum_i x_i v_i$, the dimension of summation is one, while for $x^T Q v = \sum_{ij} Q_{ij} x_i v_j$, the dimension of summation is two. This observation inspires us to conform vector dimensions while deploying such linearization.

\vspace{12pt}
{\noindent \bf 2. Linear max-min problem}

A linear max-min problem is a special case of the linear bilevel program, which  can be written as 
\begin{equation}
\label{eq:App-04-Linear-Max-Min-1}
\begin{aligned}
\max_x ~~ & c^T x + d^T y(x)   \\
\mbox{s.t.} ~~& x \in X  \\
& \begin{aligned}
y(x) \in \arg \min_y ~~ & c^T x + d^T y \\
\mbox{s.t.} ~~ & y \in Y,~ By \le b -Ax    
\end{aligned}
\end{aligned}
\end{equation}
In problem (\ref{eq:App-04-Linear-Max-Min-1}), the follower seeks an objective that is completely opposite to that of the leader. This kind of problem frequently arises in robust optimization and has been discussed in Appendix \ref{App-C-Sect02-03} from the computational perspective. Here we revisit it from a game theoretical point of view. 

Problem (\ref{eq:App-04-Linear-Max-Min-1}) can be expressed as a two-person zero-sum game       
\begin{equation}
\label{eq:App-04-Linear-Max-Min-2}
\max_{x \in X} \min_{y \in Y} \left\{  c^T x + d^T y ~\middle|~ 
Ax + By \le b \right\}   \\
\end{equation}

However, the coupled constraints make it different from a saddle point problem in the sense of a Nash game or a matrix game. Indeed, it is a Stackelberg game. Let us investigate the interchangeability of the max and min operators (decision sequence). We have already shown in Appendix \ref{App-D-Sect01-04} that swapping the order of max and min operators in a two-person zero-sum matrix game does not influence the equilibrium. However, this is not the case of (\ref{eq:App-04-Linear-Max-Min-2}) \cite{App-04-Linear-max-min}, because

\begin{equation}
\begin{aligned}
 & \max_{x \in X} \min_{y \in Y} \{ c^T x + d^T y ~|~ Ax + By \le b \} \\
=& \max_{x \in X} \left\{ c^T x + \min_{y \in Y} \{d^T y ~|~ 
By \le b - Ax\} \right\}  \\ 
\ge & \max_{x \in X} \left\{ c^T x + \min_{y \in Y}~ d^T y \right\}  \\ 
=& \max_{x \in X} ~ c^T x + \min_{y \in Y} ~ d^T y \\
=& \min_{y \in Y} \left\{ d^T y + \max_{x \in X} ~ c^T x \right\}  \\ 
\ge & \min_{y \in Y} \left\{ d^T y + \max_{x \in X} \{ c^T x ~|~ 
Ax \le b - By \} \right\}  \\
=& \min_{y \in Y} \max_{x \in X} \{ c^T x + d^T y ~|~ Ax + By \le b \}
\end{aligned}   \notag
\end{equation}
In fact, strict inequality usually holds in the third and sixth line. This result implies that owing to the presence of strategy coupling, the leader rests in a superior status, which is different from the Nash game in which players possess the same positions.

To solve linear max-min problem (\ref{eq:App-04-Linear-Max-Min-2}), there is no doubt that the aforementioned MILP transformation for general linear bilevel programs gives a possible mean for this task. Nevertheless, the special structure of (\ref{eq:App-04-Linear-Max-Min-2}) allows several alternatives which are more dedicated and effective. To this end, we will transform it into an equivalent optimization problem using LP duality theory. For the ease of notation, we merge polytope $Y$ into the coupled constraint, and the dual of lower-level LP in (\ref{eq:App-04-Linear-Max-Min-1}) (or the inner LP in (\ref{eq:App-04-Linear-Max-Min-2})) reads
\begin{equation}
\max_u ~ \{ u^T (b - Ax) ~|~ u \in U \}  \notag
\end{equation}
where $U = \{u ~|~ B^T u = d,~ u \le 0\}$ is the feasible region of dual variable $u$. As strong duality always holds for LPs, we have $d^T y = u^T (b-Ax)$. Substituting it into (\ref{eq:App-04-Linear-Max-Min-1}) we obtain
\begin{equation}
\label{eq:App-04-Linear-Max-Min-Bilinear}
\begin{aligned}
\max ~~ & c^T x + u^T b - u^T A x \\ 
\mbox{s.t.} ~~ & x \in X,~ u \in U  
\end{aligned}
\end{equation}
Problem (\ref{eq:App-04-Linear-Max-Min-Bilinear}) is a bilinear program due to the product term $u^T A x$ in variables $u$ and $x$. Several methods for solving such a problem locally or globally have been set forth in Appendix \ref{App-C-Sect02-03}, as a fundamental methodology in robust optimization. Although variable $y$ of the follower does not appear in (\ref{eq:App-04-Linear-Max-Min-Bilinear}), it can be easily recovered from the lower level of (\ref{eq:App-04-Linear-Max-Min-1}) with the obtained leader's strategy $x$.

\vspace{12pt}
{\noindent \bf 3. A retail market problem}

In a retail market, a retailer releases the prices of some goods; according to the retail prices, the customer decides on the optimal purchasing strategy subject to the demands on each goods as well as production constraints;   finally, the retailer produces or trades with a higher level market to manage the inventory, and delivers the goods to customers. This retail market can be modeled through a bilevel program. In the upper level
\begin{subequations}
\label{eq:App-04-RMarket-Retailer} 
\begin{align}
\max_{x,z} ~~  &  x^T  D_C y(x)   -  p^T  D_M  z \label{eq:App-04-RMarket-UP-1}\\
\mbox{s.t.}~~  &  A x \le a   \label{eq:App-04-RMarket-UP-2} \\
               &  B_1 y(x) + B_2 z \le b  \label{eq:App-04-RMarket-UP-3}
\end{align}
\end{subequations}
(\ref{eq:App-04-RMarket-UP-1})-(\ref{eq:App-04-RMarket-UP-3})  form retailer's problem, where vector $x$ denotes the prices of goods released by the retailor; vector $y(x)$ stands for the amounts of goods purchased by the customer, which is determined from an optimal production planning problem; $p$ is the production cost or the price in the higher level market; $z$   represents the production/purchase strategy of the retailer. Other matrices and vectors are constant coefficients. The first term in objective function (\ref{eq:App-04-RMarket-UP-1}) is the income paid by the customer, and the second term is the payoff of the retailer. The objective function is the total profit to be maximized. Because there is no competition and the retailer has full market power, to avoid unfair retail prices, we assume that both sides have reached certain agreements on the pricing policy, which is modeled through constraint (\ref{eq:App-04-RMarket-UP-2}). It includes simple lower and upper bounds as well as other bilateral contract, such as the restriction on the average price over a certain period or the price correlation among multiple goods. The inventory dynamics and other technique constraints are depicted by constraint (\ref{eq:App-04-RMarket-UP-3}).

Given the retail prices, customers solve the optimal production planning problem in the lower level 
\begin{equation}
\label{eq:App-04-RMarket-Customer}
\begin{aligned}
\min_y  ~~  &   x^T  D_C y   \\
\mbox{s.t.}~~ & F y \ge f
\end{aligned}
\end{equation}
and determine the optimal purchasing strategy. The objective function in (\ref{eq:App-04-RMarket-Customer}) is the total cost of customers, where the price vector $x$ is constant coefficient; constraints capture the demands and all other technique requirements in the production process.

Bilevel program (\ref{eq:App-04-RMarket-Retailer})-(\ref{eq:App-04-RMarket-Customer}) are not linear, although  (\ref{eq:App-04-RMarket-Customer}) is indeed an LP, because of the bilinear term $x^T  D_C y$ in (\ref{eq:App-04-RMarket-UP-1}), where both $x$ and $y$ are variables (the retailer controls $y$ indirectly through prices). The KKT condition of LP (\ref{eq:App-04-RMarket-Customer}) reads
\begin{equation}
\begin{gathered}
D^T_C x = F^T u   \\
0 \le u \bot Fy - f \ge 0
\end{gathered}  \notag
\end{equation}
where $u$ is the dual variable. The complementarity constraints can be linearized via binary variables, which has been clarified in Appendix $\ref{App-B-Sect03-05}$. Furthermore, strong duality
gives 
\begin{equation}
x^T D_C y = f^T u \notag 
\end{equation}
The right-hand side is linear in $u$. Therefore, problem (\ref{eq:App-04-RMarket-Retailer})-(\ref{eq:App-04-RMarket-Customer}) and the following MILP  
\begin{equation}
\label{eq:App-04-RMarket-MILP}
\begin{aligned}
\max_{x,y,u,v,z} ~~  &  f^T u - p^T D_M z  \\
\mbox{s.t.} ~~ &   A x \le a,~ B_1 y + B_2 z \le b \\
   & v \in \mathbb B^{N_f},~  D^T_C x = F^T u  \\
   & 0 \le u \le M(1-v) \\
   & 0 \le Fy - f \le Mv
\end{aligned}
\end{equation}
have the same optimal solution in primal variables, where $N_f$ is the dimension of $f$. 

We can learn from this example that when the problem exhibits a certain structure, the non-convexity can be eliminated without introducing additional dimensions of complexity. In problem (\ref{eq:App-04-RMarket-Retailer}), the price is a primal variable quoted by a decision maker, and is equal to the knock-down price. This scheme is called pay-as-bid. Next, we give an example of a marginal pricing market where the price is determined by the dual variables of a market clearing problem.

\vspace{12pt}
{\noindent \bf 4. A wholesale market problem}

In a wholesale market, a provider bids its offering prices to a market organizer. The organizer collects information on available resources and the bidding of the provider, and then clears the market by scheduling the production in the most economic way.  The provider is paid at the marginal cost. This problem can be modeled by a bilevel program

\begin{equation}
\label{eq:App-04-Pool-Market-Provider} 
\max_{\beta} ~ \lambda(\beta)^T p(\beta)  - f(p(\beta))  
\end{equation}
where $\beta$ is the offering price vector of the provider, $p(\beta)$ is the quantity of goods ordered by the market organizer, function $f(p) = \sum_i f_i(p_i)$, where $f_i(p_i)$ is a univariate convex function representing the production cost, and $\lambda(\beta)$ is the marginal prices of each kind of goods. Both of them depend on the value of $\beta$, and are determined from the market clearing problem in the lower level 
\begin{subequations}
\label{eq:App-04-Pool-Market-MC} 
\begin{align}
\min_{p,u} ~~ & \beta^T p + c^T u \label{eq:App-04-Pool-Market-MC-Obj}  \\
\mbox{s.t.}~~ & p_n \le p \le p_m: \eta_n, \eta_m  \label{eq:App-04-Pool-Market-Cons-1}\\
& p + F u = d: \lambda  \label{eq:App-04-Pool-Market-Cons-2}  \\ 
& A u \le a: \xi  \label{eq:App-04-Pool-Market-Cons-3}  
\end{align}
\end{subequations}
where  $u$ includes all other variables, such as the amount of each kind of goods collected from other providers or produced locally, the system operating variable, and so on;    $c$ is the coefficient including prices of goods offered by other providers, and the production cost if the organizer wishes to produce the goods by itself. Objective function (\ref{eq:App-04-Pool-Market-MC-Obj}) represents the total cost in the market to be minimized. Constraint (\ref{eq:App-04-Pool-Market-Cons-1}) defines offering limits of the upper-level provider; constraint (\ref{eq:App-04-Pool-Market-Cons-2}) is the system-wide production-demand balancing condition of each goods, the dual variable $\lambda$ at the optimal solution gives the marginal cost of each goods; (\ref{eq:App-04-Pool-Market-Cons-3}) imposes constraints which the system operation must obey, such as network flow and inventory dynamics.

In the provider's problem (\ref{eq:App-04-Pool-Market-Provider}), the offering price $\beta$ is not restricted by finite upper bounds pricing policies (but such a policy can certainly be modeled), because the competition appears in the lower level: if $\beta$ is not reasonable, the market organizer would resort to other providers or count on its own production capability.     

Compared with the situation in a retail market, problems (\ref{eq:App-04-Pool-Market-Provider})-(\ref{eq:App-04-Pool-Market-MC}) are even more complicated: the dual variable $\lambda$ appears in the objective function of the provider, and the term $\lambda^T p$ is non-convex. In the following, we reveal that it can be exactly expressed as a linear function in the primal and dual variables via (somehow tricky) algebraic transformations.

KKT conditions of the market clearing LP (\ref{eq:App-04-Pool-Market-MC}) are summarized as follows
\begin{subequations}
\label{eq:App-04-Pool-Market-MC-KKT}
\begin{gather}
\beta = \lambda + \eta_n + \eta_m  
\label{eq:App-04-Pool-Market-MC-KKT-1} \\
\eta_n^T (p-p_n) = 0
\label{eq:App-04-Pool-Market-MC-KKT-2} \\
\eta_m^T (p_m-p) = 0
\label{eq:App-04-Pool-Market-MC-KKT-3} \\
c= A^T \xi + F^T \lambda 
\label{eq:App-04-Pool-Market-MC-KKT-4} \\
\xi^T (Au - a) = 0 
\label{eq:App-04-Pool-Market-MC-KKT-5} \\
\eta_n \ge 0,~ \eta_m \le 0,~ \xi \le 0
\label{eq:App-04-Pool-Market-MC-KKT-6} \\
(\ref{eq:App-04-Pool-Market-Cons-1})-(\ref{eq:App-04-Pool-Market-Cons-3})
\label{eq:App-04-Pool-Market-MC-KKT-7} 
\end{gather}
\end{subequations}

According to (\ref{eq:App-04-Pool-Market-MC-KKT-1}), 
\begin{subequations}
\begin{equation}
\label{eq:App-04-Pool-Market-MILP-1}
\beta^T p = \lambda^T p + \eta^T_n p + \eta^T_m p
\end{equation}
From (\ref{eq:App-04-Pool-Market-MC-KKT-2}) and (\ref{eq:App-04-Pool-Market-MC-KKT-3}) we have
\begin{equation}
\label{eq:App-04-Pool-Market-MILP-2}
\eta_n^T p = \eta^T_n p_n,~ \eta_m^T p = \eta^T_m p_m,
\end{equation}
Substituting (\ref{eq:App-04-Pool-Market-MILP-2}) in (\ref{eq:App-04-Pool-Market-MILP-1}) renders
\begin{equation}
\label{eq:App-04-Pool-Market-MILP-3}
\lambda^T p = \beta^T p - \eta^T_n p_n - \eta^T_m p_m
\end{equation}
Furthermore, strong duality of LP implies the following equality
\begin{equation}
\beta^T p + c^T u = \eta^T_n p_n + \eta^T_m p_m + d^T \lambda + a^T \xi \notag
\end{equation}
or
\begin{equation}
\label{eq:App-04-Pool-Market-MILP-4}
\beta^T p - \eta^T_n p_n - \eta^T_m p_m = d^T \lambda + a^T \xi - c^T u
\end{equation}
Substituting (\ref{eq:App-04-Pool-Market-MILP-4}) in (\ref{eq:App-04-Pool-Market-MILP-3}) results in
\begin{equation}
\label{eq:App-04-Pool-Market-MILP-5}
\lambda^T p = d^T \lambda + a^T \xi - c^T u
\end{equation}
The right-hand side is a linear expression for $\lambda^T p$ in primal variable $u$ and dual variables $\lambda$ and $\xi$. 
\end{subequations}

Combining the KKT condition (\ref{eq:App-04-Pool-Market-MC-KKT}) and (\ref{eq:App-04-Pool-Market-MILP-5}) gives an MPCC which is equivalent to the bilevel wholesale market problem (\ref{eq:App-04-Pool-Market-Provider})-(\ref{eq:App-04-Pool-Market-MC})
\begin{equation}
\label{eq:App-04-Pool-Market-MPCC}
\begin{aligned}
\max ~~ & d^T \lambda + a^T \xi - c^T u - f(p(\beta)) \\
\mbox{s.t.} ~~ &  (\ref{eq:App-04-Pool-Market-MC-KKT-1})-(\ref{eq:App-04-Pool-Market-MC-KKT-7})
\end{aligned}
\end{equation}

Because complementarity conditions (\ref{eq:App-04-Pool-Market-MC-KKT-2}), (\ref{eq:App-04-Pool-Market-MC-KKT-3}), (\ref{eq:App-04-Pool-Market-MC-KKT-5}) can be linearized, and convex function $f(p)$ can be approximated by PWL functions, MPCC (\ref{eq:App-04-Pool-Market-MPCC}) can be recast as an MILP.

\subsection{Bilevel Mixed-integer Program}

Although LP can tackle many economic problems and market activities in real life, there are indeed even more decision-making problems which are beyond the reach of LP, for example, power market clearing considering unit commitment \cite{App-04-BiMIP-TEP}. KKT optimality condition or strong duality from LP theory do not apply to discrete optimization problems due to their intrinsic non-convexity. Furthermore, this is no computationally viable approach to express the optimality condition of a general discrete program in closed form, making a bilevel mixed-integer programs much more challenging to solve than a bilevel linear program. Some traditional algorithms either rely
on enumerative branch-and-bound strategies based on a weak relaxation or depends on complicated operations that are problem-specific. To our knowledge, the reformulation and decomposition algorithm proposed in \cite{App-04-BiMIP-Zeng} is the first approach that can solve general bilevel mixed-integer programs in a systematic way, and will be introduced in this section. 

The bilevel mixed-integer program has the following form
\begin{equation}
\label{eq:App-04-BiMIP-Comp-1}
\begin{aligned}
\min ~~ & f^T x + g^T y + h^T z \\
\mbox{s.t.} ~~ &  Ax \le b,~ x \in \mathbb R^{m_c} \times \mathbb B^{m_d}\\
& (y,z) \in \arg \max~~ w^T y + v^T z   \\
& \qquad \qquad \quad  \mbox{s.t.}~~  P y + Nz \le r -K x \\
& \qquad \qquad \qquad ~~~  y \in \mathbb R^{n_c},~ z \in \mathbb B^{n_d}
\end{aligned}
\end{equation}
where $x$ is the upper-level decision variable and appears in constraints of the lower-level problem; $y$ and $z$ represent lower-level continuous decision variable and discrete decision variable, respectively. We do not distinguish upper-level continuous variable and discrete variable because they have little impact on the exposition of the algorithm, unlike the ones appeared in the lower level. If the lower-level has multiple solutions,  the follower chooses the one in favor of the leader. In the current form, the upper-level constraints are independent of lower-level variables. Nevertheless, coupling constraints in the upper level can be easily incorporated \cite{App-04-BiMIP-Yue}.

In this section, we assume that the relatively complete recourse property in \cite{App-04-BiMIP-Zeng} holds, i.e., for any feasible pair $(x, z)$, the feasible set for lower-level continuous variable $y$ is non-empty. Under this premise, the optimal solution exists. This assumption is mild because we can add slack variables in the lower-level constraints and penalize constraint violation in the lower-level objective function. For instances in which the relatively complete recourse property is missing, please refer to the remedy in \cite{App-04-BiMIP-Yue}.

To eliminate $\in$ qualifier in (\ref{eq:App-04-BiMIP-Comp-1}), we duplicate decision variables and constraints of the lower-level problem and set up an equivalent formulation:
\begin{equation}
\label{eq:App-04-BiMIP-Comp-2}
\begin{aligned}
\min ~~ & f^T x + g^T y^0 + h^T z^0 \\
\mbox{s.t.} ~~ &  Ax \le b,~ x \in \mathbb R^{m_c} \times \mathbb B^{m_d}\\
& K x + P y^0 + N z^0 \le r \\
& w^T y^0 + v^T z^0 \ge \max~~ w^T y + v^T z   \\
& \qquad \qquad \qquad ~~~ \mbox{s.t.}~~  P y + Nz \le r -K x \\
& \qquad \qquad \qquad \qquad ~~  y \in \mathbb R^{n_c},~ z \in \mathbb B^{n_d}
\end{aligned}
\end{equation}
In this formulation, the leader controls its original variable $x$ as well as replicated variables $y^0$ and $z^0$. Conceptually, the leader will use $(y^0, z^0)$ to anticipate the response of follower and its impact on his objective function. Clearly, if the lower-level problem has a unique optimal solution, it must be equal to $(y^0, z^0)$. It is worth mentioning that although more variables and constraints are incorporated in (\ref{eq:App-04-BiMIP-Comp-2}), this formulation is actually an informative and convenient expression for algorithm development, as $\ge$ would be more friendly to general purpose mathematical programming solvers.

Up to now, the obstacle of solving (\ref{eq:App-04-BiMIP-Comp-2}) remains: discrete variable $z$ in the lower level, which prevents the use of optimality condition of LP. To overcome this difficulty, we treat $y$ and $z$ separately and restructure the lower-level problem as:
\begin{equation}
\label{eq:App-04-BiMIP-Obj-L}
w^T y^0 + v^T z^0 \ge \max_{z \in Z}~ v^T z + \max_y \{w^T y |Py \le r-Kx-Nz\}
\end{equation}
where $Z$ represents the set consisting of all possible values of $z$. Despite the large cardinality of $Z$, the second optimization is a pure LP, and can be replaced with its KKT  condition, resulting in: 
\begin{equation}
\label{eq:App-04-BiMIP-L-LP-KKT}
\begin{aligned}
w^T y^0 + v^T z^0 & \ge \max_{z \in Z}~ v^T z + w^T y  \\
& \qquad \mbox{s.t.}~ P^T \pi = w  \\
& \qquad \quad ~~ 0 \le \pi \bot r-Kx-Nz-Py \ge 0
\end{aligned}
\end{equation}
The complementarity constraints can be linearized via the method in Sect. \ref{App-B-Sect03-05}. Then, by enumerating $z^j$ over $Z$ with associated variables $(y^j,\pi^j)$, we arrive at an MPCC that is equivalent to problem (\ref{eq:App-04-BiMIP-Comp-2})
\begin{equation}
\label{eq:App-04-BiMIP-MPCC-Full}
\begin{aligned}
\min ~~ & f^T x + g^T y^0 + h^T z^0 \\
\mbox{s.t.} ~~ &  Ax \le b,~ x \in \mathbb R^{m_c} \times \mathbb B^{m_d}\\
& K x + P y^0 + N z^0 \le r,~ P^T \pi^j = w,~ \forall j \\
& 0 \le \pi^j \bot r - K x - N z^j - P y^j \ge 0,~ \forall j \\
& w^T y^0 + v^T z^0 \ge w^T y^j + v^T z^j,~ \forall j 
\end{aligned}
\end{equation}
Without particular mention, (\ref{eq:App-04-BiMIP-MPCC-Full}) is compatible with MILP solvers. 

Except for the KKT optimality condition, another popular approach entails applying primal-dual condition for the LP regarding the lower-level continuous variable $y$. Following this line, rewrite this LP in (\ref{eq:App-04-BiMIP-Obj-L}) by strong duality, we obtain
\begin{equation}
\label{eq:App-04-BiMIP-L-LP-PD}
\begin{aligned}
w^T y^0 + v^T z^0 & \ge \max_{z \in Z}~ v^T z + \min ~ \pi^T ( r-Kx-Nz) \\
& \qquad \qquad \qquad ~ \mbox{s.t.}~ P^T \pi = w,~ \pi \ge 0
\end{aligned}
\end{equation}
In (\ref{eq:App-04-BiMIP-L-LP-PD}), if all variables in $x$ are binary, the bilinear terms $\pi^T K x$ and $\pi^T N z$ from the leader's point of view can be linearized via the method in Sect. \ref{App-B-Sect02-02}. The min operator in the right-hand side can be omitted because the upper-level objective function is to be minimized, giving rise to
\begin{equation}
\label{eq:App-04-BiMIP-BLIP-Full}
\begin{aligned}
\min ~~ & f^T x + g^T y^0 + h^T z^0 \\
\mbox{s.t.} ~~ &  Ax \le b,~ x \in \mathbb R^{m_c} \times \mathbb B^{m_d}\\
& K x + P y^0 + N z^0 \le r  \\
& w^T y^0 + v^T z^0  \ge v^T z^j + ( r-Kx-Nz)^T \pi^j,~ \forall j  \\
& P^T \pi^j = w,~ \pi^j \ge 0,~ \forall j
\end{aligned}
\end{equation}
Clearly, (\ref{eq:App-04-BiMIP-BLIP-Full})  has fewer constraints compared to (\ref{eq:App-04-BiMIP-MPCC-Full}). Nevertheless, whenever $x$ contains continuous variables, linearizing $\pi^T K x$ would incur more binary variables.

One may think that it is hopeless to solve above enumeration forms (\ref{eq:App-04-BiMIP-MPCC-Full}) and (\ref{eq:App-04-BiMIP-BLIP-Full}) due to the large cardinality of $Z$. In a way similar to the CCG algorithm for solving robust optimization, we can start with a subset of $Z$ and solve relaxed version of problem (\ref{eq:App-04-BiMIP-MPCC-Full}), until the lower bound and upper bound of optimal value converge. The flowchart
is shown in Algorithm \ref{Ag:App-04-BiMIP-CCG}

\begin{algorithm}[!t]
\normalsize
\caption{\bf : CCG algorithm for bilevel MILP}
\begin{algorithmic}[1]
\STATE Set LB $=-\infty$, UB $= +\infty$, and $l = 0$;
\STATE Solve the following master problem
\begin{equation}
\label{eq:App-04-BiMIP-CCG-Master}
\begin{aligned}
\min ~~ & f^T x + g^T y^0 + h^T z^0 \\
\mbox{s.t.} ~~ &  Ax \le b,~ x \in \mathbb R^{m_c} \times \mathbb B^{m_d}\\
& K x + P y^0 + N z^0 \le r,~ P^T \pi^j = w,~ \forall j \le l \\
& 0 \le \pi^j \bot r - K x - N z^j - P y^j \ge 0,~ \forall j \le l \\
& w^T y^0 + v^T z^0 \ge w^T y^j + v^T z^j,~ \forall j \le l
\end{aligned}
\end{equation}
The optimal solution is $(x^*, y^{0*}, z^{0*}, y^{1*},\cdots,y^{l*}, \pi^{1*}, \cdots, \pi^{l*})$, and the optimal value is $v^*$. Update lower bound LB $=v^*$.  

\STATE Solve the following lower-level MILP with obtained $x^*$ 
\begin{equation}
\label{eq:App-04-BiMIP-CCG-SP-1}
\begin{aligned}
\theta (x^*) = \max ~~ & w^T y + v^T z  \\ 
\mbox{s.t.}~~ & P y + Nz \le r - K x^*  \\
& y \in \mathbb R^{n_c},~ z \in \mathbb B^{n_d} 
\end{aligned}
\end{equation}
The optimal value is $\theta (x^*)$. 

\STATE Solve an additional MILP to refine a solution that is favor of the leader 
\begin{equation}
\label{eq:App-04-BiMIP-CCG-SP-2}
\begin{aligned}
{\rm \Theta} (x^*) = \min ~~ & g^T y + h^T z  \\ 
\mbox{s.t.}~~ & w^T y + v^T z \ge \theta(x^*) \\
& P y + Nz \le r - K x^*  \\
& y \in \mathbb R^{n_c},~ z \in \mathbb B^{n_d} 
\end{aligned}
\end{equation}
The optimal solution is $(y^*,z^*)$, and the optimal value is ${\rm \Theta}(x^*)$. Update upper bound UB $= \min \{\mbox{UB}, f^T x^* + {\rm \Theta} (x^*)\}$.

\STATE If UB $-$ LB $=0$, terminate and report optimal solution; otherwise, set $z^{l+1}=z^*$, create new variables $(y^{l+1}, \pi^{l+1})$, adding the following cuts to master problem 
\begin{equation*}
\begin{gathered}
 w^T y^0 + v^T z^0 \ge w^T y^{l+1} + v^T z^{l+1}  \\
 0 \le \pi^{l+1} \bot r - K x - N z^{l+1} - P y^{l+1} \ge 0,~ 
 P^T \pi^{l+1} = w 
\end{gathered}
\end{equation*}
Update $l \leftarrow l+1$, and go to step 2.
\end{algorithmic}
\label{Ag:App-04-BiMIP-CCG}
\end{algorithm}

Because $Z$ has finite elements, Algorithm \ref{Ag:App-04-BiMIP-CCG} must terminate in a finite number of iterations, which is bounded by the cardinality of $Z$. When it converges, LB equals to UB without a positive gap. 

To see this, suppose that in iteration $l_1$, $(x^*, y^{0*}, z^{0*})$ is obtained in step 2 with LB $<$ UB, and $z^*$ is produced in step 4. Particularly, we assume that $z^*$ was previously derived in some iteration $l_0 < l_1$. Then, in step 5, new variables and cuts associated with $z^* = z^{l_1+1}$ will be generated and augmented with the master problem. As those variables and constraints already exist after iteration $l_0$, the augmentation is essentially redundant, and the optimal value of master problem in iteration $l_1+1$ remains the same as that in iteration $l_1$, so does LB. Consequently, in iteration $l_1+1$
\begin{equation*}
\begin{aligned}
\mbox{LB} & = f^T x^* + g^T y^{0*} + h^T z^{0*}  \\
& = f^T x^* + \min~ g^T y^0 + h^T z^0  \\ 
& \qquad \qquad ~~ \mbox{s.t.}~ P y^0 + N z^0 \le r - Kx^*,~ P^T \pi^j = w,~ \forall j \le l_1+1 \\
& \qquad \qquad \qquad 0 \le \pi^j \bot r - K x^* - N z^j - P y^j \ge 0,~ \forall j \le l_1+1 \\
& \qquad \qquad \qquad w^T y^0 + v^T z^0 \ge w^T y^j + v^T z^j,~ \forall j \le l_1+1   \\
& \ge f^T x^* + \min~ g^T y^0 + h^T z^0  \\ 
& \qquad \qquad ~~ \mbox{s.t.}~ P y^0 + N z^0 \le r - Kx^*,~ P^T \pi^j = w,~ j = l_1+1 \\
& \qquad \qquad \qquad 0 \le \pi^j \bot r - K x^* - N z^j - P y^j \ge 0,~ j = l_1+1 \\
& \qquad \qquad \qquad w^T y^0 + v^T z^0 \ge w^T y^j + v^T z^j,~ j = l_1+1\\
& \ge f^T x^* + \min~ g^T y^0 + h^T z^0  \\ 
& \qquad \qquad ~~ \mbox{s.t.}~ P y^0 + N z^0 \le r - Kx^*  \\
& \qquad \qquad \qquad w^T y^0 + v^T z^0 \ge \theta(x^*)   \\
& = f^T x^* + {\rm \Theta}(x^*) 
\end{aligned}
\end{equation*}
The second $\ge$ follows from the fact that $z^{l_1+1}$ is the optimal solution to problem (\ref{eq:App-04-BiMIP-CCG-SP-1}) and KKT condition in constraints warrants that $v^Tz^{l_1+1} + w^Ty^{l_1+1} = \theta(x^*)$. In the next iteration, the algorithm terminates since LB $\ge$ UB. 

It should be pointed out that although a large amount of variables and constraints are generated in step 5, in practice, Algorithm \ref{Ag:App-04-BiMIP-CCG} often converges to an optimal solution within a small number of iterations that could be  drastically smaller than the cardinality of $Z$, because the most critical scenarios in $Z$ can be discovered from problem (\ref{eq:App-04-BiMIP-CCG-SP-1}).

It is suggested in \cite{App-04-BiMIP-Zeng} that the master problem could be tightened by introducing variables $(\hat y, \hat \pi)$  representing the primal and dual variables of lower-level problem corresponding to $(x, z^0)$ and augmenting the following constraints 
\begin{equation*}
\begin{gathered}
w^T y^0 + v^T z^0 \ge w^T \hat y + v^T z^0  \\
 0 \le \hat \pi \bot r - K x - N z^0 - P \hat y \ge 0,~ 
 P^T \hat \pi = w
\end{gathered}
\end{equation*}
It is believed that such constraints includes some useful information that is parametric not only to $x$ but also to $z^0$, and is not available from any fixed samples $z^1,\cdots,z^l$. It is also pointed out that for instance with pure integer variables in the lower-level problem, this strategy is generally ineffective.

\section{Mathematical Programs with Equilibrium Constraints}
\label{App-D-Sect04}

A mathematical program with equilibrium constraints (MPEC) is an extension of the bilevel program by incorporating multiple followers competing with each other, resulting in a GNEP in the lower level. In this regard, an MPEC is a single-leader-multi-follower Stackelberg game. In a broader sense, MPEC is an optimization problem with variational inequalities. MPECs are difficult to solve because of the complementarity constraints.

\subsection{Mathematical Formulation}
\label{App-D-Sect04-01}

In an MPEC, the leader deploys its action $x$ prior to the followers; then each follower selects its optimal decision $y_j$ taking the decision of the leader $x$ and rivals' strategies $y_{-j}$ as given.  The MPEC can be formulated in two levels:
\begin{subequations}
\label{eq:App-04-MPEC-1}
\begin{align}
\mbox{Leader:} \qquad & \left\{
\begin{aligned}
\min_{x,\bar y, \bar \lambda, \bar \mu} ~~ & 
F(x,\bar y, \bar \lambda, \bar \mu)    \\
\mbox{s.t.} ~~ &  G(x,\bar y) \le 0    \\
               &  (\bar y, \bar \lambda, \bar \mu) \in S(x)
\end{aligned}  \right. \label{eq:App-04-MPEC-Leader}  \\
\mbox{Followers:} \qquad & \left\{
\begin{aligned}
\min_{y_j,\lambda_j, \mu_j} ~~ & f_j (x,y_j,y_{-j})  \\
\mbox{s.t.} ~~ &  g_j (x,y_j) \le 0 : \mu_j  \\
               &  h(x,y) \le 0 : \lambda_j      
\end{aligned}  \right\},~ \forall j  \label{eq:App-04-MPEC-Followers}
\end{align}
\end{subequations}

In (\ref{eq:App-04-MPEC-Leader}), the leader minimizes its payoff function $F$ which depends on the choice of its own $x$, the decisions of the followers $y$, and the dual variables $\lambda$ and $\mu$ from the lower level, because these dual variables may represent the prices of goods determined by the lower-level market clearing model. Constraints include inequalities and equalities (as a pair of opposite inequalities), as well as the optimality condition of the lower-level problem. In (\ref{eq:App-04-MPEC-Followers}), $x$ is treated as a parameter, and the competition among followers comes down to a GNEP with shared convex constraints: the payoff function $f_j(x,y_j,y_{-j})$ of follower $j$ is assumed to be convex in $y_j$; inequality $g_j(x,y_j) \le 0$ defines a local constraint of follower $j$ which is convex in $y_j$ and does not involve $y_{-j}$; inequality $h(x,y) \le 0$ is the shared constraint which is convex in $y$. Since each follower's problem is convex, the KKT condition is both necessary and sufficient for optimality. We assume that the set of GNEPs $S(x)$ is always non-empty.

The GNEP encompasses several special cases in the lower level. If the global constraint is absent, it degenerates into an NEP; moreover, if the objective functions of followers are also decoupled, the lower level reduces to independent convex optimization programs. 

By replacing the lower-level GNEP with its KKT condition (\ref{eq:App-04-GNEP-KKT}), the MPEC (\ref{eq:App-04-MPEC-1}) becomes an MPCC, which can be solved by some suitable methods explained before. As the lower level GNEP usually possesses infinitely many equilibria, the outcome found by the MPCC reformulation is the favourite one from the leader's perspective. We can also require the Lagrange multipliers for the shared constraints should be equal, so as to restrict the GNEP to VEs. If the followers' problems are linear, the primal-dual optimality condition is an alternative choice in addition to the KKT condition, as it often involves fewer constraints. Nevertheless, the strong duality may introduce products involving primal and dual variables, such as those in (\ref{eq:App-04-Linear-Bilevel-NLP}) and (\ref{eq:App-04-Linear-Bilevel-NLP-Pen}), which remain non-convex and require special treatments. 

\subsection{Linear Complementarity Problem}
\label{App-D-Sect04-02}

A linear complementarity problem (LCP) requires finding a feasible solution subject to the following constraints
\begin{equation}
\label{eq:App-04-LCP-1}
0 \le x \bot P x + q \ge 0  
\end{equation}
where $P$ is a square matrix; $q$ is a vector. Their dimensions are compatible with $x$.  

LCP is a special case of MPCC without an objective function. This type of problem  frequently arises in various disciplines including market equilibrium analysis, computational mechanics, game theory, and mathematical programming. The theory of LCPs is a well-developed field. Detailed discussions can be found in \cite{App-04-LCP-Book}. In general, an LCP is NP-hard, although it is polynomially solvable for some special cases. One situation is when the matrix $P$ is positive semidefinite. In such circumstance, problem (\ref{eq:App-04-LCP-1}) can be solved via the following convex quadratic program 
\begin{equation}
\label{eq:App-04-LCP-CQP}
\begin{aligned}
\min~~ & x^T P x + q^T x \\
\mbox{s.t.}~~ & x \ge 0,~  P x + q \ge 0    
\end{aligned}
\end{equation}
(\ref{eq:App-04-LCP-CQP}) is a CQP which is readily solvable. Its optimum must be non-negative according to the constraints.   If the optimal value of (\ref{eq:App-04-LCP-CQP}) is 0, then its optimal solution also solves LCP (\ref{eq:App-04-LCP-1}); otherwise, if  the optimal value is strictly positive, LCP (\ref{eq:App-04-LCP-1}) is infeasible. In fact, this conclusion holds no matter whether $P$ is positive semidefinite or not. However, if $P$ is indefinite, identifying the global optimum of a non-convex QP (\ref{eq:App-04-LCP-CQP}) is also NP-hard, and thus does not facilitate solving the LCP. 

There is a large body of literature discussing algorithms for solving LCPs. One of the most representative ones is the Lemke's pivoting method developed in \cite{App-04-LCP-Lemke}, and another emblematic one is the interior-point method proposed in \cite{App-04-LCP-IPO}. One drawback of the former method is its exponentially growing worst-case complexity, which makes it less efficient for large problems.  The latter approach runs in polynomial time, but it requires the positive semidefiniteness of $P$, which is a strong assumption and limits its application. In this section, we will not present comprehensive reviews on the algorithms for LCP. We will introduce MILP formulations for problem (\ref{eq:App-04-LCP-1}) devised in \cite{App-04-LCP-MILP-1,App-04-LCP-MILP-2}. They make no reference on any special structure of matrix $P$. More importantly, they offer an option to access the solutions of practical problems in a systematic way.   

Recall the MILP formulation techniques presented in Appendix \ref{App-B-Sect03-05}, it is easy to see that problem (\ref{eq:App-04-LCP-1}) can be equivalently expressed as linear constraints with additional binary variable $z$ as follows 
\begin{equation}
\label{eq:App-04-LCP-MILC}
 0 \le x \le Mz,~ 0 \le P x + q \le M(1-z)    
\end{equation}
Integrality of $z$ maintains the element-wise complementarity of $x$ and $Px+q$: at most one of $x_i$ and $(Px+q)_i$ can be strictly positive. Formulation (\ref{eq:App-04-LCP-MILC}) entails a manually specified parameter $M$, which is not instantly available at hand. On the one hand, it must be big enough to preserve all extreme points of (\ref{eq:App-04-LCP-1}). On the other hand, it is expected to be as small as possible from a computational perspective, otherwise, the continuous relaxation of (\ref{eq:App-04-LCP-MILC}) would be very loose. In this regard, (\ref{eq:App-04-LCP-MILC}) is too cursory, although it might work well. 

To circumvent above difficulty, it is proposed in \cite{App-04-LCP-MILP-1} to solve a bilinear program without a big-M parameter
\begin{equation}
\label{eq:App-04-LCP-BLP}
\begin{aligned}
\min_{x,z}~~ & z^T ( P x + q ) + ( {\bf 1} - z )^T x \\
\mbox{s.t.}~~ & x \ge 0,~ P x + q \ge 0,~ z \mbox{ binary}    
\end{aligned}
\end{equation}

If (\ref{eq:App-04-LCP-1}) has a solution $x^*$, the optimal value of (\ref{eq:App-04-LCP-BLP}) is 0: for $x^*_i>0$, we have $z^*_i=1$ and $(P x^* + q)_i=0$; for $(P x^* + q)_i>0$, we have $z^*_i=0$ and $x^*_i=0$. The optimal solution is consistent with the feasible solution of (\ref{eq:App-04-LCP-MILC}). The objective can be linearized by introducing auxiliary variables $w_{ij}=z_i x_j$, $\forall i,j$. However, applying normal integer formulation techniques in Appendix \ref{App-B-Sect02-02} on variable $w_{ij}$ again needs the upper bound of $x_i$, another interpretation of the big-M parameter. 

A parameter-free MILP formulation is suggested in \cite{App-04-LCP-MILP-1}. To understand the basic idea, recall the fact that $(1-z_i)x_i=0$; if we impose $x_i=w_{ii}=x_iz_i$, $i=1,2,\cdots$ in the constraint, $({\bf 1} - z )^T x$ in the objective can be omitted. Furthermore, multiplying both sides of $\sum_j P_{kj} x_j + q_k \ge 0$, $k=1,2,\cdots$ with $z_i$ gives $\sum_j P_{kj} w_{ij} + q_k z_i\ge 0$, $\forall i,k$. Since $z_i \in \{0,1\}$, $\sum_j P_{kj} x_j + q_k \ge \sum_j P_{kj} w_{ij} + q_k z_i$, $\forall i,k$ and $0 \le w_{ij} \le x_j$, $\forall i,j$ naturally hold. Collecting up these valid inequalities, we obtain an MILP
\begin{equation}
\label{eq:App-04-LCP-MILP-1}
\begin{aligned}
\min_{x,z,w}~~ & q^T z + \sum_i \sum_j P_{ij} w_{ij} \\
\mbox{s.t.}~~ & \sum_j P_{kj} x_j + q_k \ge \sum_j P_{kj} w_{ij} + q_k z_i \ge 0,~ \forall i,k \\
& 0 \le w_{ij} \le x_j,~ \forall i,j,~ w_{jj} = x_j,~ \forall j, ~ z \mbox{ binary}
\end{aligned}   
\end{equation}
Instead of enforcing every $\sum_j P_{kj} w_{ij} + q_k z_i$ being at 0, we relax them as inequalities and minimize their summation. More valid inequalities can be added in (\ref{eq:App-04-LCP-MILP-1}) by exploiting linear cuts of $z$. It is proved in \cite{App-04-LCP-MILP-1} that relation $w_{ij}=z_i x_j$, $\forall i,j$ is implicitly guaranteed at the optimal solution of (\ref{eq:App-04-LCP-MILP-1}). In view of this, MILP (\ref{eq:App-04-LCP-MILP-1}) is equivalent to LCP (\ref{eq:App-04-LCP-1}) in the following sense: (\ref{eq:App-04-LCP-1}) has a solution if and only if (\ref{eq:App-04-LCP-MILP-1}) has an optimal value equal to zero, and the optimal solution to (\ref{eq:App-04-LCP-MILP-1}) incurring a zero objective value is a solution of LCP (\ref{eq:App-04-LCP-1}). MILP (\ref{eq:App-04-LCP-MILP-1}) is superior compared with (\ref{eq:App-04-LCP-MILC}) and big-M linearization based MILP  formulation of MINLP (\ref{eq:App-04-LCP-BLP}) because it is parameter-free and gives tighter continuous relaxation. Nevertheless, the number of constraints in (\ref{eq:App-04-LCP-MILP-1}) is significantly larger than that in formulation (\ref{eq:App-04-LCP-MILC}). This method has been further analyzed in \cite{App-04-LCP-MILP-2} and extended to binary-constrained mixed LCPs.

Another parameter-free MILP formulation is suggested in \cite{App-04-LCP-MILP-3}, which takes the form of
\begin{equation}
\label{eq:App-04-LCP-MILP-2}
\begin{aligned}
\max_{\alpha,y,z}~~ & \alpha \\
\mbox{s.t.}~~ & 0 \le (P y)_i + q_i \alpha \le 1-z_i,~\forall i \\
& 0 \le y_i \le z_i,~ z_i \in \{0,1\},~ \forall i \\
& 0 \le \alpha \le 1
\end{aligned}   
\end{equation}
Since $\alpha =0$, $y=0$, $z=0$ is always feasible, MILP (\ref{eq:App-04-LCP-MILP-2}) is feasible and has an optimum no greater than 1. By observing the constraints, we can conclude that if MILP (\ref{eq:App-04-LCP-MILP-2}) has a feasible solution with $\bar \alpha > 0$, then $x=y/\bar \alpha$ solves problem (\ref{eq:App-04-LCP-1}). If the optimal solution $\bar \alpha= 0$, then problem (\ref{eq:App-04-LCP-1}) has no solution; otherwise, suppose $\bar x$ solves (\ref{eq:App-04-LCP-1}), and let $\bar \alpha^{-1} = \max \{\bar x_i, (P\bar x)_i+q_i,i=1,\cdots\}$, then for any $0 <\alpha \le \bar \alpha$, $\bar y = \alpha \bar x$ is feasible in (\ref{eq:App-04-LCP-MILP-2}). As a result, the optimal solution should be no less than $\bar \alpha$, rather than 0. Compared with formulation (\ref{eq:App-04-LCP-MILC}), the big-M parameter is adaptively scaled by optimizing $\alpha$. 

Because (\ref{eq:App-04-LCP-MILP-2}) works with an intermediate variable $y$, when LCP (\ref{eq:App-04-LCP-1}) should be jointly solved with other conditions on $x$, formulation (\ref{eq:App-04-LCP-MILP-2}) is not advantageous, because non-convex variable transformation $x=y/\alpha$ must be appended to link both parts. 

Robust solutions of LCPs with uncertain $P$ and $q$ are discussed in \cite{App-04-Robust-LCP}. It is found that when $P \succeq 0$, robust solutions can be extracted from an SOCP under some mild assumptions on the uncertainty set; otherwise, the more general problem with uncertainty can be reduced to a deterministic non-convex QCQP. This technique is particularly useful in uncertain traffic equilibrium problems and uncertain Nash-Cournot games. Uncertain VI problems and MPCCs can be tackled in the similar vein after some proper transformations.

It is shown in \cite{App-04-LCP-BLP-MPEC} that a linear bilevel program or its equivalent MPEC can be globally solved via a sequential LCP method. A hybrid enumerative method is suggested which substantially reduces the effort for searching a solution of the LCP or certifying that the LCP has no solution. When the LCP is easy to solve, this approach is attractive.

Several extensions of LCP, including the discretely-constrained mixed LCP, discretely-constrained Nash-Cournot game, discretely-constrained MPEC, and logic constrained equilibrium problem as well as their applications in energy markets and traffic system equilibrium have been investigated in \cite{App-04-DM-LCP-1,App-04-DM-LCP-2,App-04-DM-LCP-3,App-04-DM-LCP-4}. In a word, due to its wide applications, LCP is still an active research field, and MILP remains an appealing method for solving LCPs for practical problems.

\subsection{Linear Programs with Complementarity Constraints}

A linear program with complementarity constraints (LPCC) entails solving a linear optimization problem  with linear complementarity constraints. It is a special case of MPCC if all functions in the problem are linear, and a generalization of LCP by incorporating an objective function to be optimized. An LPCC has the following form 
\begin{equation}
\label{eq:App-04-LPCC}
\begin{aligned}
\max_{x,y}~~ & c^T x + d^T y \\
\mbox{s.t.}~~ & Ax + By \ge f \\
& 0 \le y \bot q + N x + M y \ge 0 
\end{aligned}   
\end{equation}

A standard approach for solving (\ref{eq:App-04-LPCC}) is to linearize the complementarity constraint by introducing a binary vector $z$ 
and solve the following MILP
\begin{equation}
\label{eq:App-04-LPCC-MILP-BigM}
\begin{aligned}
\max_{x,y,z}~~ & c^T x + d^T y \\
\mbox{s.t.}~~ & Ax + By \ge f \\
& 0 \le q + N x + M y \le M z \\
& 0 \le y \le M(1-z) \\
& z \in \{0,1\}^m
\end{aligned}   
\end{equation}
If both of $x$ and $y$ are bounded variables, we can readily derive the proper value of $M$ in each inequality; otherwise, finding high quality bounds  is nontrivial even if they do exist. The method in \cite{App-04-LPCC-BigM} can be used to determine proper bounds of $M$, if the NLP solver can successfully find local solutions of the bounding problems.

Using a arbitrarily large value may solve the problem correctly. Nevertheless, parameter-free method is still of great theoretical interests. A smart Benders decomposition algorithm is proposed in \cite{App-04-LPCC-Benders} to solve (\ref{eq:App-04-LPCC-MILP-BigM}) without requiring the value of $M$. The completely positive programming method developed in \cite{App-04-LPCC-CPP-Relax} can also be used to solve (\ref{eq:App-04-LPCC-MILP-BigM}). For more theory and algorithm for LPCC, please see \cite{App-04-LPCC-DCP}, \cite{App-04-LPCC-1}-\cite{App-04-LPCC-7} and references therein. Interesting connections among conic QPCCs, QCQPs, and completely positive programs are revealed in \cite{App-04-LPCC-8}.

\section{Equilibrium Programs with Equilibrium Constraints}
\label{App-D-Sect05}

An equilibrium program with equilibrium constraints (EPEC) is the most general extension of the bilevel program. It incorporates multiple leaders and multiple followers competing with each other in the upper level and the lower level, respectively, resulting in two GNEPs in both levels. In this regard, an EPEC is a multi-leader-follower Stackelberg game. 

\subsection{Mathematical model}
\label{App-D-Sect05-01}

\begin{figure}[!htp]
\centering
\includegraphics[scale=0.60]{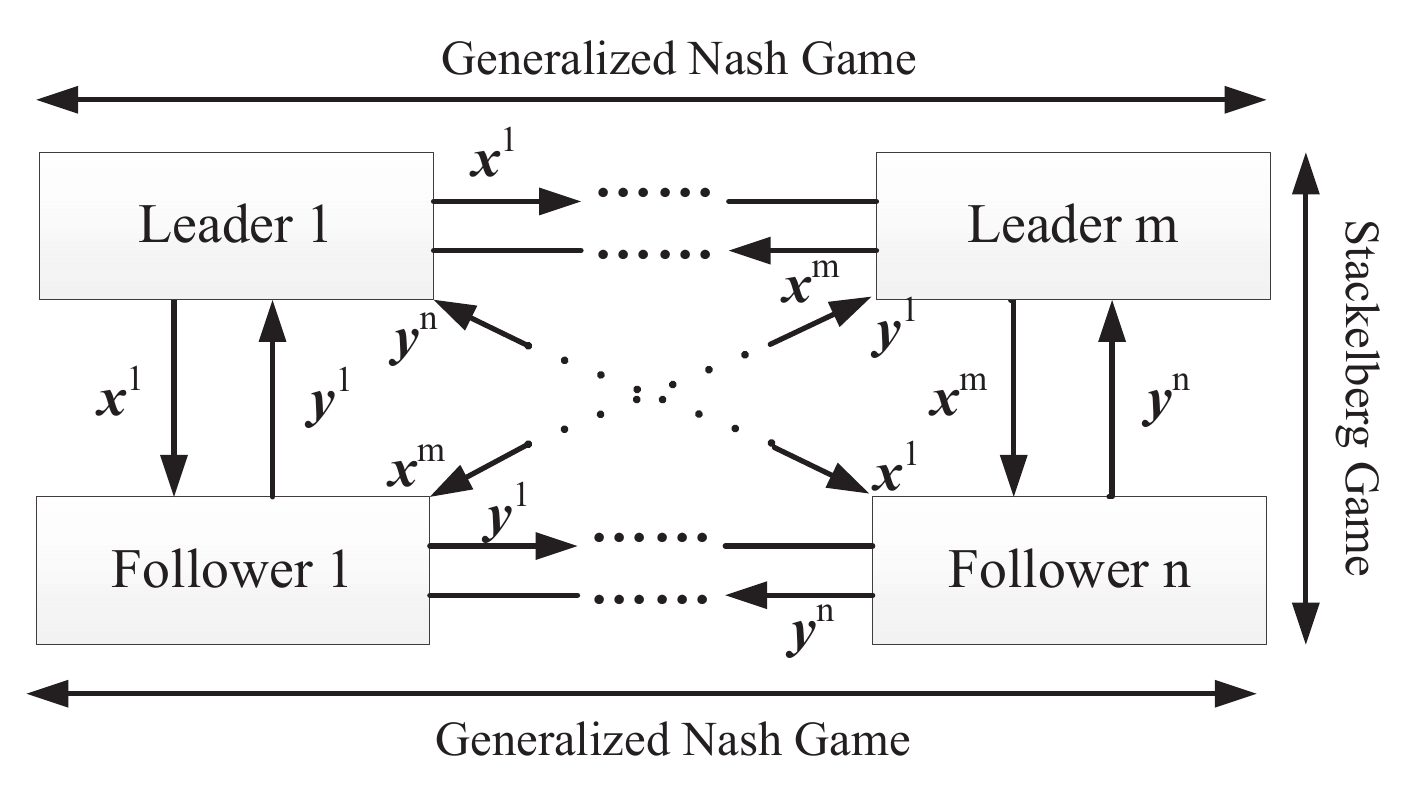}
\caption{The structure of an EPEC.}
\label{fig:App-04-04}
\end{figure}

In an EPEC, each leader $i$ deploys an action $x_i$ prior to the followers while taking movements of other leaders $x_{-i}$ into account and anticipating the best responses $y(x)$ from the followers; then each follower selects its optimal decision $y_j$ by taking the strategies of leaders $x$ and rivals' actions $y_{-j}$ as given.  The EPEC can be formulated in two levels
\begin{subequations}
\label{eq:App-04-EPEC-1}
\begin{align}
\mbox{Leaders:} \qquad & \left\{
\begin{aligned}
\min_{x_i,\bar y, \bar \lambda, \bar \mu} ~~ & 
F_i (x_i,x_{-i},\bar y, \bar \lambda, \bar \mu)    \\
\mbox{s.t.} ~~ &  G_i (x_i) \le 0    \\
               &  (\bar y, \bar \lambda, \bar \mu) \in S(x_i,x_{-i})
\end{aligned}  \right\},~\forall i 
\label{eq:App-04-EPEC-Leaders}  \\
\mbox{Followers:} \qquad & \left\{
\begin{aligned}
\min_{y_j} ~~ & f_j (x,y_j,y_{-j})  \\
\mbox{s.t.} ~~ &  g_j (x,y_j) \le 0 : \mu_j  \\
               &  h(x,y) \le 0 : \lambda_j      
\end{aligned}  \right\},~ \forall j  
\label{eq:App-04-EPEC-Followers}
\end{align}
\end{subequations}
In (\ref{eq:App-04-EPEC-Leaders}), each leader minimizes its payoff function $F_i$ which depends on its own choice $x_i$, the decisions of followers $y$, dual variables $\lambda$ and $\mu$ are parameterized in competitors' strategies $x_{-i}$. Tuple $(\bar y, \bar \lambda, \bar \mu)$ in the upper level is restricted by the optimality condition of the lower-level problem. Although the inequality constraints of leaders are decoupled, and we do not explicitly consider global constraints in the upper-level GNEP, the leaders' strategy sets as well as their payoff functions are still correlated through the best reaction map $S(x_i,x_{-i})$, and hence (\ref{eq:App-04-EPEC-Leaders}) itself is a GNEP, which is non-convex. The followers' problem (\ref{eq:App-04-EPEC-Followers}) is a GNEP with shared constraints, which is the same as the situation in an MPEC. The same convexity assumptions are made in (\ref{eq:App-04-EPEC-Followers}). The structure of EPEC (\ref{eq:App-04-EPEC-1}) is depicted in Fig. \ref{fig:App-04-04}. The equilibrium solution of EPEC (\ref{eq:App-04-EPEC-1}) is defined as the GNE among leaders' MPECs. It is common knowledge that EPECs often have no pure strategy equilibrium due to the intrinsic non-convexity of MPECs.  

\subsection{Methods for Solving an EPEC}
\label{App-D-Sect05-02}

An EPEC can be viewed as a set of coupled MPEC problems: leader $i$ is facing an MPEC composed of problem $i$ in (\ref{eq:App-04-EPEC-Leaders}) together with all followers' problems in (\ref{eq:App-04-EPEC-Followers}), which is parameterized in $x_{-i}$. By replacing lower-level GNEP with its KKT optimality conditions, it can be imaged that the GNEP among leaders have non-convex constraints which inherit the tough properties of complementarity constraints. Thus, solving an EPEC is usually extremely challenging. To our knowledge, systematic algorithms of EPEC are firstly developed in dissertations \cite{App-04-EPEC-Algorithm-1,App-04-EPEC-Algorithm-2,App-04-EPEC-Algorithm-3}. The primary application of such an equilibrium model is found in energy market problems, see \cite{App-04-EPEC-Algorithm-4} for an excellent introduction. 

Unlike the NEP and GNEP discussed in Sect. \ref{App-D-Sect01} and Sect. \ref{App-D-Sect02}, where the strategy sets are convex or jointly convex, because the lower-level problems are replaced with KKT optimality conditions and the MPEC for the leader is intrinsically non-convex, provable existence and uniqueness guarantees for the solution to EPECs are non-trivial. There is sustainable attempt on the analysis of EPEC solution properties. For example, the existence of a unique equilibrium for certain EPEC instances are discussed in \cite{App-04-EPEC-Solution-1,App-04-EPEC-Solution-2},  and in \cite{App-04-EPEC-Solution-3} for a nodal price based power market model. However, the existence and uniqueness of solution are only guaranteed under restrictive conditions. Counterexamples have been given in \cite{App-04-EPEC-Solution-4} to demonstrate that there is no general result for
the existence of solutions to EPECs due to their non-convexity. The non-uniqueness issue is studied in \cite{App-04-EPEC-Algorithm-2,App-04-EPEC-Solution-5}. It is shown that even in the simplest instances, local uniqueness of the EPEC equilibrium solution may not be guaranteed, and  a manifold of equilibria may exist. This can be understood because EPEC is a generalization of GNEP, whose solution property is illustrated in Sect. \ref{App-D-Sect02-01}. When the payoff functions possess special structures, say, a potential function exists, then the existence of a global equilibrium can be investigated using the theory of potential games \cite{App-04-EPEC-Solution-6,App-04-EPEC-Solution-7,App-04-EPEC-Potential-MPEC}. In summary, the theory of EPEC solutions are much more complicated than the single-level NEP and GNEP.

This section reviews several representative algorithms which are widely used in literature. The former two are generic and  seen in \cite{App-04-EPEC-Algorithm-1,App-04-EPEC-Algorithm-2}; the third one is motivated by the convenience brought by the property of potential games, and reported in \cite{App-04-EPEC-Potential-MPEC}; at last, a pricing game in a competitive market, which appears to be non-convex at first sight, is presented to show the hidden convexity in such a special equilibrium model.

\vspace{12pt}
{\noindent \bf 1. Best response algorithm}

Since the equilibrium of an EPEC is a GNEP among leaders' MPEC problems, the most intuitive  strategy for identifying an equilibrium  solution
is the best response algorithm. In some literature, it is also called diagonalization method or sequential MPEC method. This approach can be further categorized into Jacobian type and Gauss-Seidel type method, according to the information used when players update their strategies. 

To explain the algorithmic details, denote by MPEC($i$) the problem of leader $i$: the upper level is problem $i$ in (\ref{eq:App-04-EPEC-Leaders}), and the lower level is the GNEP described in (\ref{eq:App-04-EPEC-Followers}) given all leaders' strategies. Let $x^k_i$ be the strategy of leader $i$ in iteration $k$, and $x^k=(x^k_1,\cdots,x^k_m)$ the strategy profile of leaders. The Gauss-Seidel type algorithm proceeds as follows \cite{App-04-Complement-Book}:

\begin{algorithm}[H]
\normalsize
\caption{\bf : Best-response (Diagonalization) algorithm for EPEC}
\begin{algorithmic}[1]
\STATE Choose an initial strategy profile $x^0$ for leaders, set convergence tolerance $\varepsilon>0$, an allowed number of iterations $K$, and the iteration index $k=0$;
\STATE Let $x^{k+1}=x^k$. Loop for players $i=1,\cdots,m$:
\begin{enumerate}
\item[a.]  Solve MPEC($i$) for leader $i$ given $x^{k+1}_{-i}$.
\item[b.]  Replace $x^{k+1}_i$ with the optimal strategy of leader $i$ just obtained.  
\end{enumerate}

\STATE If $\| x^{k+1} - x^k \|_2 \le \varepsilon$, the upper level converges; solve lower-level GNEP (\ref{eq:App-04-EPEC-Followers}) with $x^*=x^{k+1}$ using the algorithms elaborated in Sect. \ref{App-D-Sect02-02}, and the equilibrium among followers is $y^*$. Report $(x^*,y^*)$ and terminate. 

\STATE If $k = K$, report failure of convergence and quit. 

\STATE Update $k \leftarrow k+1$, and go to step 2.
\end{algorithmic}
\label{Ag:App-04-EPEC-Diagonalization}
\end{algorithm} 

Without an executable criterion to judge the existence and uniqueness of solution, possible outcomes of Algorithm \ref{Ag:App-04-EPEC-Diagonalization} are discussed in three situations.

1. There is no equilibrium. Algorithm \ref{Ag:App-04-EPEC-Diagonalization} does not converge. In such circumstance, one may turn to seeking a mixed-strategy Nash equilibrium, which always exists. Examples are given in \cite{App-04-Complement-Book}: if there are two leaders, we can list possible strategy combinations and solve the lower-level GNEP among followers, then compute respective payoffs of the two leaders, and then build a bimatrix game, whose mixed-strategy Nash equilibrium can be calculated from solving an LCP, as explained in Sect. \ref{App-D-Sect01-04}.

2. There is a unique equilibrium, or there are multiple equilibria. Algorithm \ref{Ag:App-04-EPEC-Diagonalization} may converge or not, and which equilibrium will be found (if it converges) depends on the initial strategy profile offered in step 1. 

3. Algorithm \ref{Ag:App-04-EPEC-Diagonalization} may converge to a local equilibrium in the sense of \cite{App-04-EPEC-Solution-3}, if each MPEC is solved by a local NLP method which does not guarantee global optimality. The true equilibrium can be found only if each leader's MPEC can be globally solved. The MILP reformulation (if possible) offers one plausible way for this task.

\vspace{12pt}
{\noindent \bf 2. KKT system method}

To tackle the divergence issue in the best response algorithm, it is proposed to apply the KKT condition to each leader's MPEC and   solve the resulting KKT systems simultaneously \cite{App-04-EPEC-Algorithm-3,App-04-EPEC-Algorithm-5}. The solution turns out to be a strong stationary equilibrium point of EPEC (\ref{eq:App-04-EPEC-1}). There is  no convergence issue in this approach, since no iteration is deployed. However, special attention should be paid to some potential problems mentioned below.  

1. Since the EPEC is essentially a GNEP among leaders, the concentrated KKT system may have non-isolated solutions. To refine a meaningful outcome, we can manually specify a secondary objective function, which is optimized subject to the KKT system.  

2. The embedded (twice) application of KKT condition for the lower-level problems and upper-level problems inevitably introduces extensive of complementarity and slackness conditions, which greatly challenges solving the concentrated KKT system. In this regard, scalability may be a main bottleneck for this approach. If the lower-level GNEP is linear, it may be better to use primal-dual optimality condition first for followers, and then KKT condition for leaders.

3. Because each leader's MPEC is non-convex, a stationary point of the KKT condition is not necessarily an optimal solution of the leader; as a result, the solution of the  concentrated KKT system may not be an equilibrium of the EPEC. To validate the result, one can conduct the best-response method initiated at the candidate solution with a slight perturbation.   

\vspace{12pt}
{\noindent \bf 3. Potential MPEC method}

 When the upper-level problems among leaders admit a potential function satisfying (\ref{eq:App-04-PG-1}), the EPEC can be reformulated as an MPEC, and the relations of their solutions are revealed by comparing the KKT condition of the normalized Nash stationary points of the EPEC and the KKT condition of the associated MPEC \cite{App-04-EPEC-Potential-MPEC}.

For example, if the leaders' objectives are given by
\begin{equation}
F_i(x_i,x_{-i},y) = F^S_i(x_i) + H(x,y) \notag
\end{equation}
or in other words, the payoff function $F_i(x_i,x_{-i},y)$ can be decomposed as the sum of two parts: the first one $F^S_i(x_i)$ only depends on the local variable $x_i$, and the second one $H(x,y)$ is common to all leaders. In such circumstance, the potential function can be expressed as   
\begin{equation}
U(x,y) = H(x,y) + \sum_{i=1}^m F^S_i(x_i)  \notag
\end{equation}

Please see Sect. \ref{App-D-Sect01-05} for the condition under which a potential function exists and special instances in which a potential function can be easily found.

Suppose that leaders' local constraints are given by $x_i \in X_i$ which is independent of $x_{-i}$ and $y$, and the best reaction map of followers with fixed $x$ is given by $(\bar y, \bar \lambda, \bar \mu) \in S(x)$. Clearly, the solution of MPEC
\begin{equation}
\begin{aligned}
\min_{x,\bar y,\bar \lambda,\bar \mu} ~~ &  U(x,\bar y)  \\
\mbox{s.t.} ~~ & x_i \in X_i,~ i=1,\cdots,m  \\
               & (\bar y, \bar \lambda, \bar \mu) \in S(x)
\end{aligned}  \notag
\end{equation}
must be an equilibrium solution of the original EPEC.

This approach leverages the property of potential games and is superior over the previous two methods (if a potential function exists): KKT condition is applied only once to the lower level problems, and the equilibrium can be retrieved by solving MPEC only once.    

\vspace{12pt}
{\noindent \bf 4. A pricing game in a competitive market}

We consider an EPEC taken from the examples in \cite{App-04-EPEC-Hidden-Convexity-1}, which models a strategic pricing game in a competitive market. The hidden convexity in this EPEC is revealed. For ease of exposition, we study the case with two leaders and one follower. The results can be extended to the situation where more than two leaders exist. The pricing game with two leaders can be formulated by the following EPEC
\begin{subequations}
\label{eq:App-04-Two-Leader-Pricing-Game}
\begin{align}
\mbox{Leader 1:} \quad & \max_{x_1}~ \left\{  y^T(x_1,x_2) A_1 x_1 ~\middle|~ \ B_1 x_1 \le b_1 \right\}  \label{eq:App-04-TLPG-Leader-1}  \\
\mbox{Leader 2:} \quad & \max_{x_2}~ \left\{  y^T(x_1,x_2) A_2 x_2 ~\middle|~ \ B_2 x_2 \le b_2 \right\}  \label{eq:App-04-TLPG-Leader-2}  \\
\mbox{Follower:} \quad & \max_y~ \left\{ f(y) - y^T A_1 x_1 - y^T A_2 x_2
~\middle|~ Cy=d \right\}   \label{eq:App-04-TLPG-Follower} 
\end{align}
\end{subequations}
In (\ref{eq:App-04-TLPG-Leader-1}) and (\ref{eq:App-04-TLPG-Leader-2}), two leaders announce their offering prices $x_1$ and $x_2$, respectively, subject to some certain pricing policy described in their corresponding constraints. The follower then decides how many goods should be purchased from each leader, according to the optimal solution of problem (\ref{eq:App-04-TLPG-Follower}), where the profit of the follower
\begin{equation}
f(y)=-\dfrac{1}{2} y^T Q y + c^T y \notag
\end{equation}
is a strongly concave quadratic function, i.e. $Q \succ 0$, and matrix $C$ has full rank in its rows. Each player in the market wishes to maximize his own profit. The utilities of leaders are the payments from trading with the follower; the profit of follower is the revenue minus the purchasing cost.

At first sight, EPEC (\ref{eq:App-04-Two-Leader-Pricing-Game}) is non-convex, not only because the leaders' objective functions are bilinear, but also because the best response mapping is generally non-convex. In light of the strong convexity of  (\ref{eq:App-04-TLPG-Follower}), the following KKT condition:
\begin{equation}
\begin{gathered}
 c-Qy-A_1 x_1-A_2 x_2-C^T \lambda = 0 \\ Cy-d=0 
\end{gathered}  \notag
\end{equation}
is necessary and sufficient for a global optimum. Because constraints in (\ref{eq:App-04-TLPG-Follower}) are all equalities, there is no complementarity and slackness condition. Solve this set of linear equations, we can obtain the optimal solution $y$ in a closed form. To this end, substituting 
\begin{equation}
y=Q^{-1}(c-A_1 x_1-A_2 x_2-C^T \lambda) \notag 
\end{equation}
into the second equation, we have   
\begin{equation}
\lambda = M \left[N(c-A_1 x_1-A_2 x_2)-d \right] \notag
\end{equation}
where 
\begin{equation}
M = \left[C Q^{-1} C^T \right]^{-1},~ N = C Q^{-1} \notag
\end{equation}
Moreover, eliminating $\lambda$ in the expression of $y$ gives the best reaction map
\begin{equation}
\label{eq:App-04-TLPG-Follower-Best-Reaction} 
y = r + D_1 x_1 + D_2 x_2  
\end{equation}
where
\begin{equation}
\begin{gathered}
 r  =  Q^{-1} c + N^T M d - N^T M N c \\
D_1 = N^T M N A_1 -Q^{-1} A_1  \\
D_2 = N^T M N A_2 -Q^{-1} A_2  \\
\end{gathered}   \notag
\end{equation}

Substituting (\ref{eq:App-04-TLPG-Follower-Best-Reaction}) into the objective functions of leaders, EPEC (\ref{eq:App-04-Two-Leader-Pricing-Game}) reduces to a standard Nash game
\begin{equation}
\begin{aligned}
\mbox{Leader 1:} \quad & \max_{x_1}~ \left\{  \theta_1 (x_1,x_2) ~\middle|~ \ B_1 x_1 \le b_1 \right\}   \\
\mbox{Leader 2:} \quad & \max_{x_2}~ \left\{  \theta_2 (x_1,x_2) ~\middle|~ \ B_2 x_2 \le b_2 \right\}   \\
\end{aligned}    \notag
\end{equation}
where 
\begin{equation}
\begin{gathered}
\theta_1(x_1,x_2) = r^T A_1 x_1 + x_1^T D_1^T A_1 x_1 + x_2^T D_2^T A_1 x_1\\
\theta_2(x_1,x_2) = r^T A_2 x_2 + x_2^T D_2^T A_2 x_2 + x_1^T D_1^T A_2 x_2
\end{gathered}    \notag
\end{equation}

The partial Hessian matrix of $\theta_1(x_1,x_2)$ can be calculated as
\begin{equation*}
\nabla^2_{x_1} \theta_1(x_1,x_2) = 2A^T_1 (N^T M N -Q^{-1}) A_1 
\end{equation*}
As $Q \succ 0$, its inverse matrix $Q^{-1} \succ 0$; denote by $Q^{-1/2}$ the square root of $Q^{-1}$, and 
\begin{equation}
P_J = I - Q^{-1/2} C^T (CQ^{-1} C^T)^{-1} C Q^{-1/2} \notag
\end{equation}
It is easy to check that $P_J$ is a projection matrix, which is symmetric and idempotent, i.e., $P_J = P^2_J =P^3_J =\cdots$. Moreover, it can be verified that the Hessian matrix $\nabla^2_{x_1} \theta_1(x_1,x_2)$ can be expressed via 
\begin{equation}
\nabla^2_{x_1} \theta_1(x_1,x_2) = 2A^T_1 (N^T M N -Q^{-1}) A_1 
= -2 A^T_1 Q^{-1/2} P_J Q^{-1/2} A_1 \notag
\end{equation}

For any vector $z$ with a proper dimension,
\begin{equation}
\begin{aligned}
 z^T \nabla^2_{x_1} \theta_1(x_1,x_2) z ~
=~ & -2 z^T \left( A^T_1 Q^{-1/2} P_J Q^{-1/2} A_1 \right) z        \\
=~ & -2 z^T A^T_1 Q^{-1/2} P^T_J P_J Q^{-1/2} A_1 z  \\
=~ & -2 (P_J Q^{-1/2} A_1 z)^T (P_J Q^{-1/2} A_1 z) \le 0
\end{aligned}   \notag
\end{equation}
We can see that $\nabla^2_{x_1} \theta_1(x_1,x_2) \preceq 0$. The similar analysis also applies to $\nabla^2_{x_2} \theta_2(x_1,x_2)$. Therefore, the problems of leaders are actually convex programs, and a pure-strategy Nash equilibrium exists.

\section{Conclusions and Further Reading}
\label{App-D-Sect06}

Equilibrium problems entail solving interactive optimization problems simultaneously, and serve as the foundation for modeling competitive behaviors among strategic decision makers, and analyzing the stable outcome of a game. This chapter provides an overview on two kinds of equilibrium problems that frequently arise in various economic and engineering applications. 

One-level equilibrium problems, including  the NEP and GNEP, are introduced first. The existence of equilibrium can be ensured under some convexity and monotonicity assumptions. Distributed methods for solving one-level games are presented. When each player solves a strictly convex optimization problem, distributed algorithms converge with provable guarantee, and thus are preferred, whereas the KKT system renders nonlinear equations and is relatively difficult to solve. To address incomplete information and uncertainty in player's decision making, a robust optimization based game model is proposed in \cite{App-04-Robust-Game-Theory}, which is distribution-free and relaxes Harsanyi's assumptions on Bayesian games. Particularly, the robust Nash equilibrium of a bimatrix game with uncertain payoffs can be characterized  via the solution of a second-order cone complementarity problem \cite{App-04-Robust-NE-1}, and more general cases involving  $n$ players and continuous payoffs are discussed in \cite{App-04-Robust-NE-2}. Distributional uncertainty is tackled in \cite{App-04-DR-CC-Game}, in which the mixed-strategy Nash equilibrium of a distributionally robust chance-constrained game is studied. A generalized Nash game arises when the strategy sets of players are coupled. Due to practical interests from a variety of engineering disciplines, the solution method for GNEPs is still an active research area. The volume of articles is growing quickly in recent years, say, \cite{App-04-GNEP-Algorithm-1,App-04-GNEP-Algorithm-2,App-04-GNEP-Algorithm-3,App-04-GNEP-Algorithm-4,App-04-GNEP-Algorithm-5,App-04-GNEP-Algorithm-6,App-04-GNEP-Algorithm-7,App-04-GNEP-Algorithm-8,App-04-GNEP-Algorithm-9}, to name just a few. GNEPs with uncertainties are studied in \cite{App-04-GNEP-Uncertainty-1,App-04-GNEP-Uncertainty-2}.

Bilevel equilibrium problems, including  the bilevel program, MPEC, and EPEC, are investigated. These problems are intrinsically hard to solve, due to the non-convexity induced by the best reaction map of followers, and solution properties have been revealed for specific instances under restrictive assumptions. We recommend \cite{App-04-Complement-Book,App-04-BLP-Book-1,App-04-BLP-Book-2} for theoretical foundations and energy market applications of bilevel equilibrium models, and \cite{App-04-BLP-Review-Pozo} for an up-to-date survey. The theories on bilevel programs and MPEC are relatively mature. Recent research efforts have been spent on new constraint qualifications and optimality conditions, for example, the work in \cite{App-04-MPEC-CQ-1,App-04-MPEC-CQ-2,App-04-MPEC-CQ-3,App-04-MPEC-CQ-4}. The MILP reformulation is preferred by most power system applications, because the ability of MILP solvers keep improving, and a global  optimal solution can be found. Stochastic MPEC is proposed in \cite{App-04-Stochastic-MPEC-1} to model uncertainty using probability distributions. Algorithms are developed in \cite{App-04-Stochastic-MPEC-2,App-04-Stochastic-MPEC-3,App-04-Stochastic-MPEC-4,App-04-Stochastic-MPEC-5}, and a literature review can be found in \cite{App-04-Stochastic-MPEC-6}. Owing to the inherent hardness, discussions on EPEC models are limited to special cases, such as those with shared P-matrix linear complementarity constraints \cite{App-04-Muti-Leader-Follower-Game-1}, power market models \cite{App-04-Muti-Leader-Follower-Game-1,App-04-Muti-Leader-Follower-Game-2,App-04-Muti-Leader-Follower-Game-3}, those with convex quadratic objectives and linear constraints \cite{App-04-Muti-Leader-Follower-Game-1}, and Markov game models \cite{App-04-EPEC-Markov-Regularization}. Methods for solving EPEC are based on relaxing or regularizing complementarity constraints \cite{App-04-EPEC-Markov-Regularization,App-04-EPEC-Relaxation}, as well as evolutionary algorithms \cite{App-04-Muti-Leader-Follower-Game-EA}. Robust equilibria of EPEC are discussed in \cite{App-04-Robust-SNE}. An interesting connection between the bilevel program and the GNEP has been revealed in \cite{App-04-BiP-GNEP}, establishing a new look on these game models.

We believe that the equilibrium programming models will become an imperative tool for designing and analyzing interconnected energy systems and related markets, in view of the physical interdependence of heterogenous energy flows and strategic interactions among different network operators.

\backmatter

\printindex


\end{document}